\documentclass{article}

\RequirePackage{amsthm,amsmath,amsfonts,amssymb}

\RequirePackage{hyperref}
\hypersetup{
    colorlinks,
    linkcolor={red!50!black},
    citecolor={blue!50!black},
    urlcolor={blue!80!black}
}
\RequirePackage{graphicx}
\usepackage{setspace}
\usepackage{booktabs}
\usepackage{mathtools}

\usepackage[capitalise,nameinlink]{cleveref}

\crefname{supp}{Supplement}{Supplements}
\crefname{appendix}{Supplement}{Supplements}

\usepackage[shortlabels]{enumitem}

\newcounter{inappendix}

\usepackage{mathrsfs}

\usepackage{algorithm}
\usepackage{algpseudocode}
\usepackage{tikz}
\usepackage{subcaption}
\usetikzlibrary{shapes,arrows,fadings}
\usepackage{circuitikz}
\usepackage{apptools}
\usepackage{enumitem}
\usepackage{lmodern}
\usepackage{placeins}
\usepackage[savewrites,nonumberlist, seeautonumberlist,sort=def]{glossaries}
\usepackage{apptools}
\usepackage{authblk}
\usepackage{chngcntr}

\usepackage[margin=1in]{geometry}
\usepackage[shortlabels]{enumitem}
\usepackage[
natbib=true,
backend=biber,
sorting=nyt,
useprefix=true,
style=authoryear,
url=false,
doi=true,
eprint=false,
giveninits=true,
uniquename=init]{biblatex}
\AtEveryBibitem{%
  \clearfield{note}%
}
\AtEveryBibitem{\clearlist{language}}
\AtEveryBibitem{\clearfield{month}}
\AtEveryCitekey{\clearfield{month}}

\usepackage[toc, title, page,header]{appendix}
\usepackage{minitoc}
\theoremstyle{plain}

\newtheorem{theorem}{Theorem}

\newtheorem{lemma}{Lemma}
\newtheorem{proposition}{Proposition}

\theoremstyle{remark}
\newtheorem{assumption}{Assumption}
\newtheorem{definition}{Definition}
\newtheorem{example}{Example}

\newtheorem{examplerepeattwo}{Example}

\newtheorem{stheorem}{Theorem S\ignorespaces}
\newtheorem{scorollary}{Corollary S\ignorespaces}
\newtheorem{slemma}{Lemma S\ignorespaces}
\newtheorem{sproposition}{Proposition S\ignorespaces}
\newtheorem{sremark}{Remark S\ignorespaces}
\crefformat{stheorem}{#2Theorem S#1#3}
\crefformat{scorollary}{#2Corollary S#1#3}
\crefformat{slemma}{#2Lemma S#1#3}
\crefformat{sproposition}{#2Proposition S#1#3}
\crefformat{sremark}{#2Remark S#1#3}

\theoremstyle{remark}
\newtheorem{sassumption}{Assumption S\ignorespaces}

\newtheorem{sdefinition}{Definition S\ignorespaces}
\newtheorem{sexample}{Example S\ignorespaces}

\newtheorem{scounterexample}{Counterexample S\ignorespaces}
\crefformat{sassumption}{#2Assumption S#1#3}
\crefrangeformat{sassumption}{Assumptions #3S#1#4 to #5S#2#6}
\crefformat{sdefinition}{#2Definition S#1#3}
\crefformat{sexample}{#2Example S#1#3}
\crefformat{sexamplerepeat}{#2Example S#1#3}
\crefformat{sexamplerepeat2}{#2Example S#1#3}
\crefformat{scounterexample}{#2Counterexample S#1#3}

\newcommand{\alignswith}{\overset{\mathcal{C}}{\approx}}
\newcommand{\equivalign}{\overset{\mathcal{C}}{\sim}}
\newcommand{\bs}[1]{\boldsymbol{#1}}
\newcommand{\mc}[1]{\mathcal{#1}}

\DeclareMathOperator*{\argmin}{arg\,min}

\setlist[enumerate,1]{label=(\alph*)}

\renewcommand{\glossarymark}[1]{}
\makeglossaries
\newglossaryentry{integers}{
    name={$[M]$},
    description={The set \{1, \dots, M\} for any natural number $M$}
}
\newglossaryentry{rv-overline}{
    name={$\overline{X}_k$},
    description={The first $k$ elements of the random vector $X$, i.e. $\overline{X}_k = (X_1, \dots, X_k)$}
}
\newglossaryentry{support}{
    name={\ensuremath{\textnormal{$\textsf{Supp}\left[X;P\right]$}}},
    description={The support of the random variable $X$ under the distribution $P$}
}
\newglossaryentry{l2-space}{
    name={$L^2(P)$},
    description={The Hilbert space of real-valued measurable functions of a random vector $W$ distributed according to $P$ with finite second moments, equipped with the covariance inner product}
}
\newglossaryentry{l20-space}{
    name={$L^2_0(P)$},
    description={The subspace of $L^2(P)$ of mean zero functions.}
}

\newglossaryentry{l2x-space}{
    name={$L^2(X;P)$},
    description={The subspace of $L^2(P)$ of $L^2(P)$ comprised of functions that depend on $X$ only. }
}

\newglossaryentry{l20x-space}{
    name={$L^2_0(X;P)$},
    description={The subspace of $L^2(X;P)$ comprised of mean zero functions that depend on $X$ only.}
}
\newglossaryentry{projection}{
    name={$\Pi[\cdot|\mc{H}]$},
    description={The orthogonal projection onto the closed linear subspace $\mc{H}$ of a Hilbert space}
}
\newglossaryentry{adjoint}{
    name={$A^*$},
    description={The hermitian adjoint of an operator $A$}
}
\newglossaryentry{closure}{
    name={$\overline{\mc{H}}$},
    description={The closure of the set $\mc{H}$}
}
\newglossaryentry{summation}{
    name = {$+, \oplus$},
    description={For two inner product spaces that are subspaces of a common inner product space, $\mc{A}, \mc{B}$, $\mc{A} + \mc{B} = \{a + b:a\in\mc{A}, b\in\mc{B}\}$. When $\mc{A} \perp \mc{B}$, we often use $\oplus$ instead of $+$. When summing over many subspaces we instead use $\sum$ and $\bigoplus$}
}
\newglossaryentry{tangent-space}{
    name={$\mc{T}(P, \mc{M})$},
    description={The tangent space of the statistical model $\mc{M}$ at the distribution $P$}
}
\newglossaryentry{influence_function}{
    name={$\gamma^1_P$},
    description={An influence function (equivalently gradient) of the parameter $\gamma$ at $P$ (in some model $\mc{M}$)}
}
\newglossaryentry{efficient_influence_function}{
    name={$\gamma^1_{P, eff}$},
    description={The efficient influence function (equivalently canonical gradient) of the parameter $\gamma$ at $P$ (in some model $\mc{M}$)}
}
\newglossaryentry{ideal-model}{
    name={$\mc{Q}$},
    description={The statistical model of distributions for the ideal data random variable $W$}
}
\newglossaryentry{observed-model}{
    name={$\mc{P}$},
    description={The statistical model of distributions for the observed data random variable $O$}
}
\newglossaryentry{observed-model-ext}{
    name={$\mc{P}^{ext}$},
    description={The statistical model of distributions for the observed data random variable $O$ when there are no restrictions on the ideal data model, defined in \cref{def:extended-obs-data-model}}
}
\newglossaryentry{ideal-parameter}{
    name={$\psi$},
    description={The ideal data parameter of interest}
}
\newglossaryentry{observed-parameter}{
    name={$\varphi$},
    description={The observed data parameter of interest}
}
\newglossaryentry{alignments}{
    name={$\mc{C}$},
    description={The collection of alignments relating the ideal data and observed data distributions, defined in \cref{def:alignswith}}
}
\newglossaryentry{alignswith}{
    name={$P \alignswith Q$},
    description={$P$ and $Q$ are aligned relative to $\mc{C}$. See  \cref{def:alignswith}}
}

\newglossaryentry{stronglyaligns}{
     name={%
        \begin{tabular}{@{}l}    
            $(Q, U, P)$\\ 
            \ensuremath{\textnormal{$\textit{Strongly Align}$}}
        \end{tabular} 
    },
    description={See \cref{def:strongly-align}}
}

\newglossaryentry{equivalign}{
    name={$Q \equivalign \widetilde{Q}$},
    description={$Q$ and $\widetilde{Q}$ are equivalent according to the equivalence relation defined in \cref{subsec:alignment-assumptions}}
}
\newglossaryentry{equivclass}{
    name={$\xi(Q;\mc{C})$},
    description={The class of distributions equivalent to $Q$ under the relation $\equivalign$. See \cref{subsec:alignment-assumptions}}
}
\newglossaryentry{phiclass}{
    name={$\Phi(P;\mc{C})$},
    description={The equivalence class $\xi(Q;\mc{C})$ such that $P\alignswith Q$ for all $Q\in \xi(Q;\mc{C})$. See \cref{subsec:alignment-assumptions}}
}

\newglossaryentry{fused-data-model}{
    name={$(\mc{Q}, \mc{P}, \mc{C})$},
    description={A fused-data model. See  \cref{def:fdm}}
}
\newglossaryentry{fused-data-framework}{
    name={$(\mc{Q}, \mc{P}, \mc{C}, \psi, \varphi)$},
    description={A fused-data framework. See  \cref{def:fdf}}
}

\title{Towards a Unified Theory for Semiparametric Data Fusion with Individual-Level Data}
\author[1]{Ellen Graham}
\author[1]{Marco Carone}
\author[1]{Andrea Rotnitzky}
\affil[1]{Department of Biostatistics, University of Washington}
\date{}

\begin{document}

\maketitle

\begin{abstract}
We consider inference about a finite-dimensional parameter integrating samples from independent sources. A recently developed theory considers scenarios where sources align with subsets of the conditional distributions of a single factorization of the joint target distribution. While this theory applies in many settings, it falls short in important data fusion problems, such as two-sample instrumental variable analysis, settings that integrate data from epidemiological studies with diverse designs, and studies with mismeasured variables supplemented by external validation studies. In this paper, we derive a comprehensive theory that, in particular, covers these settings by allowing the integration of sources aligned with conditional distributions that do not correspond to a single factorization of the target distribution. We provide a universal characterization of the influence functions of regular and asymptotically linear estimators and the efficient influence function of a target parameter, irrespective of the parameter of interest or the statistical model for the target distribution, thus paving the way for a unified theory for machine-learning debiased, semiparametric efficient estimation.
\end{abstract}

\doparttoc 
\faketableofcontents
\part{}

\section{Introduction}

With the growing amount of data available to researchers, there has been
increasing attention to developing statistical methodology for data fusion,
which aims to effectively combine diverse data sources to estimate summaries
of interest. Studies involving data fusion are abundant. For instance, in
economics and social sciences, practitioners frequently combine large
administrative data sources routinely collected by governments (\cite{angrist_effect_1992, ridder_chapter_2007}). In vaccine research,
immunobridging research combines historical clinical trials with trials
performed in new populations or with new versions of treatments to establish
efficacy when phase three randomized trials are infeasible or unethical (%
\cite{fleming_biomarkers_2012, luedtke_partial_2017}). In settings with
measurement error, data from a primary study is often combined with data
from an external validation study to correct for the effect of measurement
error (\cite{chen_semiparametric_2008, cole_illustration_2023}).

Data fusion methods have surged in recent
years, but are tailored to specific problems. Examples include methodologies designed for transporting treatment effects
learned in one population to a new population (\cite%
{pearl_transportability_2011, rudolph_robust_2017, dahabreh_extending_2019,
dahabreh_towards_2020, luedtke_partial_2017}) and methods to estimate treatment effects with instrumental variables where the instrument and treatment are
measured in one source, and the instrument and outcome are measured
in another source (\cite%
{klevmarken_missing_1982, angrist_effect_1992, shu_improved_2020,
sun_semiparametric_2022, zhao_two-sample_2019, pierce_efficient_2013}).
A recent thread of research aims to develop general
methods that can be applied to a wide range of data structures and
summaries of interest. This includes \cite{graham_efficient_2016} who
consider parameters defined as solutions of additively
separable population moment equations using data from two sources; \cite{hu_paradoxes_2023}
who supplement individual-level data from a target
population with estimated summaries from external sources; the seminal work of \cite{li_efficient_2023} who offers a comprehensive framework for semiparametric inference when the sources align with subsets of the conditional distributions in a single factorization of the joint target distribution; and the landmark work of \cite{qiu_efficient_2024} who provide multiply robust estimation for estimands that are the means of fixed functions under the fused-data framework of \cite{li_efficient_2023}.

The work in this article aims to advance the development of a unified theory for inference with individual-level fused-data. To this end, we develop a general
theory for semiparametric efficient inference that permits the aligned conditional
distributions of different sources to correspond to the components of
different factorizations of the joint target distribution. Our formulation
includes the work of \cite{li_efficient_2023} as a special case but
it additionally addresses many common scenarios not covered by the work of these authors, including the two-sample
instrumental variables problem, measurement error problems with external
validation studies, and scenarios integrating data from diverse
epidemiological study designs, such as prospective cohorts and retrospective
case-control studies. We provide a universal characterization of the influence functions of regular and asymptotically linear estimators and the efficient one of a target parameter, regardless of the number of data sources, the parameter of interest, or the statistical model of the target distribution. This characterization paves the way for a unified theory for machine learning debiased, semiparametric efficient estimation. 

Our formulation involves a model for coarsened data that differs in a fundamental way from standard missing data models. In missing data, the sampled units are viewed as originating from the target population, and the target parameter is a functional of the distribution of the full data that would be observed for each sampled unit. The target population can be regarded as a mixture of respondent-type subpopulations, each representing a source of information, with mixing probabilities given by response probabilities.

In the data fusion setting considered here, the combined sample across sources does not represent a random sample from the target population. Instead, the target population is an abstract population that is related to each source population only through the alignments. This target population need not coincide with any individual source population, see e.g. \cref{example:transporting}. However in special cases, such as the forthcoming \cref{example:disease-prevalence}, the target population may coincide with a particular source. The distinction between the target populations of missing data and fused-data analysis parallels the difference between the target populations for the average treatment effect (ATE) and the average treatment effect among the treated (ATT) in causal inference.

There are special instances where the two formulations are the same. \cite{li_efficient_2023} noted that their fused-data models overlap with the monotone missing completely at random model when in each source, the entire distribution of the observed data aligns. The fused-data model in the present paper further includes non-monotone missing completely at random data.

The developments in the present work can be viewed as complementary to the theory characterizing the set of all influence functions and the efficient one in coarsened at random (CAR) data \citep{robins_recovery_1992, robins_estimation_1994, van_der_laan_unified_2003}. Indeed, the present work relies on deriving the score operator for the fused-data model, just as the analogous theory for CAR relies on deriving the score operator for the CAR model. 
\section{Review of semiparametric theory}

We first set the notation used throughout. A glossary is given in \cref{app:notation}. We let $[K]\coloneqq\left\{ 1,...,K\right\} $. We
use capital letters to denote random variables or vectors and lowercase letters for the corresponding
realizations. We use $P$ and $Q$ to denote
probability laws and $p$ and $q$ to denote their corresponding densities. For $x=\left( x_{1},...,x_{K}\right)$, and $k \leq K$, $\overline{x}_{k}$ denotes $(x_{1},...,x_{k})$. If $X=\left( X_{1},...,X_{K}\right)$, $X\leq x$ stands for $X_{k}\leq x_{k}$ for all $1\leq k\leq K.$ 

If $\left( X,Y\right)\sim P$, $\mathsf{Supp}\left[
X;P\right]$ denotes the support of $X$. We write $d\left( X\right) \in \mathcal{D}$ a.e.- $P$ on $\mathcal{B}$ iff $\mathcal{B}\subseteq\mathsf{Supp}\left[
X;P\right]$ and $P\left[ d\left( X\right)
\in \mathcal{D}|X\in \mathcal{B}\right] =1$. If $\mathcal{P}$ is a collection of laws
for $\left( X,Y\right)$ with common support
then $\mathsf{Supp}\left[ X;\mathcal{P}\right]$ stands for the common
support of $X$.

For $P$, $Q$ two probability laws for $(X, Y)$, with marginals for $X$ denoted by $P_X$ and $Q_X$, $\frac{dQ}{dP}(x)$ stands for $\frac{dQ_X}{dP_X}(x)$. For $x$ not in $\textsf{Supp}\left[X;P\right]$ we define $\frac{dQ}{dP}(x)$ arbitrarily.

All conditional probabilities in this article are assumed to be regular (\cite{chang_conditioning_1997}). Conditioning on lowercase letters in expectations and probabilities indicates
conditioning on the corresponding random variable taking the specific value.

If $(X, Y)\sim P$, $L^2(P)$ denotes the Hilbert space of real-valued functions of $(X, Y)$ with finite variance equipped with the covariance inner product.  $L^2(X;P)$ is the subspace of $L^2(P)$ of functions of $X$ alone. $L^2_0(P)$ and $L^2_0(X;P)$ are the subspaces of $L^2(P)$ and $L^2_0(X;P)$ respectively of mean-zero functions. For $\Lambda \left( P\right)$ a closed linear subspace of $%
L^{2}(P)$ and $f\in L^{2}(P),$ $\Pi \left[ \left. f\right\vert \Lambda
\left( P\right) \right] $ denotes the $L^{2}\left( P\right) $-orthogonal projection
of $f$ onto $\Lambda \left( P\right)$. We sometimes write $\Pi[f(X)|\Lambda(P)]$ instead of $\Pi[f|\Lambda(P)]$.

We now review key
elements of semiparametric theory, highlighting the pivotal role 
influence functions play in constructing debiased machine
learning estimators. See \cref{app:semiparametric-theory} for an expanded review. Given a model $\mathcal{M}$ for the distribution of a random vector $X$, the collection of all regular submodels of $\mathcal{M}$ at $P\in \mc{M}$ indexed by $t \in (-\varepsilon, \varepsilon)$, each with $P_{t=0} = P$, induces the tangent set $\mathcal{T}^{\circ }(P,\mc{M})$, i.e. the collection of scores at $t=0$ of all such submodels. The closed linear span of $\mathcal{T}^{\circ }(P,\mc{M})$ in $L^{2}(P)$ is called the tangent space  $\mathcal{T}(P,\mc{M})$. Models $\mathcal{M}$ with $\mc{T}(P, \mc{M}) = L^2_0(P)$ are (locally at $P$) nonparametric. When $\mc{T}(P, \mc{M})$ is finite dimensional, $\mc{M}$ is parametric. Otherwise $\mc{M}$ is semiparametric.  

A mapping $\gamma :\mathcal{M}\rightarrow \mathbb{R}$ is pathwise
differentiable at $P$ in $\mathcal{M}$ if there exists $\gamma _{P}^{1}\in
L_{0}^{2}(P)$ such that $\left. \frac{d}{dt}%
\gamma (P_{t})\right\vert _{t=0}=\left\langle \gamma _{P}^{1},h\right\rangle
_{L^{2}\left( P\right) }$ for any regular submodel of $\mathcal{M}$ indexed by $t \in (-\varepsilon, \varepsilon)$ with score $h$ at $t=0$ and with $P_{t=0}=P$. Such $\gamma^1_P$ is called a gradient. The canonical gradient $\gamma _{P,eff}^{1}$ of $\gamma $ at $P$ is the unique gradient of $\gamma $
that belongs to $\mathcal{T}(P,\mathcal{M})$.

We now turn to estimating $\gamma(P)$ based on $n$ i.i.d. draws $X_i$ from $P$. An estimator $\hat{\gamma}_n$ of $\gamma(P)$ is asymptotically linear with influence function $\Gamma_P \in L^2_0(P)$ if $n^{1/2}(\hat{\gamma}_n - \gamma(P)) = n^{-1/2}\sum_{i=1}^n \Gamma_P(X_i) + o_p(1)$. Asymptotically linear estimators are consistent and asymptotically normal with variance $var_{P}\left( \Gamma
_{P}\right)$. The estimator $\gamma _{n}$\ is regular with respect $\mc{M}$ at $P$ if its convergence to its limiting
distribution is locally uniform over laws contiguous to $P$ (\cite{bickel_efficient_1998}). An asymptotically linear estimator $\hat{\gamma} _{n}$ of $\gamma(P)$
at $P$ with influence function $\Gamma _{P}$ is regular at $P$ in $\mc{M}$ if and only if $\gamma$ is pathwise differentiable at $P$ in $\mc{M}$ and $\Gamma _{P}$ is a gradient of $\gamma $ (\cite%
{van_der_vaart_asymptotic_2000}). This is why the terms "gradients" and "influence functions" are often used interchangeably and we follow this convention. The variance of $\gamma^1_{P, eff}$ is the smallest asymptotic variance of any regular and asymptotically linear (RAL) estimator of $\gamma$. RAL estimators that achieve this bound are called semiparametric efficient. Henceforth we call the canonical gradient the efficient influence function.

RAL estimators can readily be constructed using influence functions, and an efficient estimator can be constructed using the efficient influence function. For example, given an influence function $\gamma^1_P$ and a consistent estimator $\hat{P}_n \in \mc{M}$ of $P$, the one-step estimator is $\hat{\gamma}_n(P) \coloneqq \gamma(\hat{P}_n) + n^{-1}\sum_{i=1}^n\gamma^1_{\hat{P}_n}(X_i)$ is, under regularity conditions, RAL with influence function $\gamma^1_P$ and efficient if $\gamma^1_P = \gamma^1_{P,eff}$. This is an example of a debiased estimator, so called because the correction by the influence function "debiases" the plug-in estimator $\gamma(\hat{P}_n)$, which is generally biased if $\hat{P}_n$ is constructed using flexible machine learning techniques. Other examples of debiased machine learning estimators include TMLE (\cite{van_der_laan_targeted_2006}) and estimating equation procedures based on influence functions (\cite{van_der_laan_unified_2003, tsiatis_semiparametric_2006}). 

\section{The inferential problem and fused-data framework}
\label{sec:inferential-problem}

\label{subsec:overview}

Suppose we aim to estimate the value of the functional $\psi:\mathcal{%
Q\rightarrow }\mathbb{R}$ at the distribution $Q_{0}$ of a random vector $%
W$, which takes values in $\mathcal{W}\subseteq \mathbb{R}^{K},$ under the
assumption that $Q_{0}$ belongs to $\mathcal{Q}$, a collection of mutually
absolutely continuous probabilities on $\mathcal{W}$. Suppose we lack access to
a random sample from $Q_{0}$. Instead, we have access to 
samples drawn from $J$ distinct data sources. In each unit drawn from source $j$, we observe the subvector $W^{\left( j\right) }$ of $W$. Thus, the available data are $n$
i.i.d. copies of the coarsened data vector $O\coloneqq(c(W,S),S)$
where $S$ takes values in $\left\{ 1,...,J\right\} $, with $S=j$ if the
observation originates from the $j^{th}$-source and $c(W,j)=W^{\left(
j\right) }$. For each source $j$, we assume that only specific
conditional or marginal distributions of the vector $W^{\left( j\right) }$
align, i.e. agree with, with the corresponding conditionals or marginals of $Q_{0}$. These alignments are justified either on the basis of substantive
considerations or due to the sampling design employed in each source. The
precise alignments allowed by our formulation are detailed in \cref%
{subsec:alignment-assumptions}. Throughout we refer to $O$ as the observed data vector. We let $P_{0}$ denote the law of $O$.

Fused data is distinct from missing data. In a missing data setting $Q_{0}$ stands for the law of a full
data vector $W$ drawn from a target population. We draw a random sample of
units from the target population but the entire vector $W$ is not always
observed. In contrast, in a fused-data setting, we may measure the entire data vector $W$ in each unit drawn from each source. Yet even in this instance the available data $W_1, \dots, W_n$ would not constitute a random sample from $Q_0$ unless we make the strong assumption that the full joint distribution of $W$ in each source aligns with that of $Q_0$.

Model $\mathcal{Q}$ and the alignment assumptions give rise to a model $%
\mathcal{P}$ for the observed data, where a distribution $P$ for $O$ is
included in $\mathcal{P}$ if and only if it is mutually absolutely
continuous with $P_{0}$ and there exists a distribution $Q$ in $\mathcal{Q}$
such that the assumed alignments hold between $P$ and $Q$. Multiple
distributions $Q$ in $\mathcal{Q}$ may satisfy the assumed alignments with a
single observed data law $P,$ i.e. the assumed alignments may not suffice to
identify $Q$ from $P$. However, throughout we shall assume that $\psi: \mathcal{Q}\rightarrow \mathbb{R}$ depends on $Q$ solely through the aligned
conditionals or marginals in the sense that for a given distribution $P$, $\psi$ takes a single value over all aligned distributions $Q$ in the model $\mc{Q}$. In \cref{subsec:alignment-assumptions} we will
show this implies that $\psi (Q)$ is identified by $P$ in the sense
that there exists a well-defined mapping $\varphi :\mathcal{P}\rightarrow  
\mathbb{R}$ such that $\psi (Q)=\varphi \left( P\right) $ for any $Q$ in $%
\mathcal{Q}$ that satisfies the assumed alignments with $P$.

A restricted model $\mathcal{Q},$ such as a strictly semiparametric or fully
parametric model, will sometimes place equality constraints on the aligned
conditionals or marginals of $P$ and thus will lead to a semi or fully
parametric model $\mathcal{P}$. Importantly, a non-parametric model $%
\mathcal{Q}$ may or may not result in a non-parametric model $\mathcal{P}$.
This distinction stems from the fact that, as illustrated in \cref%
{example:transporting} scenario (iii) below, the very existence of a
common distribution $Q$ aligning certain 
conditionals or marginals with those from different data sources may impose
equality restrictions on the observed data law.

Letting $\mathcal{C}$ denote the collection of assumed alignments, we refer
to the triplet $\left( \mathcal{Q},\mathcal{P},\mathcal{C}\right) $ as a
fused-data model and to the quintuplet $\left( \mathcal{Q},\mathcal{P},%
\mathcal{C}\text{\thinspace },\psi ,\varphi \right) $ as a fused-data framework. Fused-data frameworks are broad enough to encompass numerous fused-data applications. These include causal analysis combining data from multiple
sources as in examples \ref{example:tsiv-lsm} and %
\ref{example:transporting} below. In particular, it includes integrating data from prospective cohorts and retrospective case-control or case-only studies. \cref{example:transporting} below illustrates its applicability
in these contexts. It additionally encompasses studies with mismeasured variables supplemented with external validation studies as in \cref%
{example:disease-prevalence} below. 

In any data fusion application, each assumed alignment must be
well-justified. This must be done on a case-by-case basis, taking into
account substantive underpinnings and the sampling designs operating in each
source. Our paper does not aim to discuss best practices for conducting such
analyses. Instead, we assume that the analyst has already defined the
problem as a fused-data framework; we provide a general theory of inference
about the scalar summary of interest $\varphi(P)$ from that starting point onwards. All our results can be applied to a multivariate $\varphi(P)$ by applying them component wise.

In the following examples, we assume each law $Q\in\mc{Q}$ is dominated by a product measure $\mu$ and for each $j\in[J]$, $P\in \mc{P}$, $P(\cdot|S=j)$ is dominated by a product measure $\mu^{(j)}$. We write $q = \frac{dQ}{d\mu}$ and $p(\cdot|S=j) = \frac{dP(\cdot|S=j)}{d\mu^{(j)}}$.
\begin{example}[Estimating disease prevalence from misclassified disease and an external validation study]
\label{example:disease-prevalence}

Suppose $V$ is a binary indicator of disease, relatively inexpensive to
measure but prone to misclassification, $Y$ denotes the true but
costly-to-obtain binary disease indicator, and $X$ represents a $p-$vector
of baseline covariates. Let $Q_{0}$ denote the distribution of $W\coloneqq%
\left( X,V,Y\right) $ in the target population. We are interested in
estimating $\psi \left( Q_{0}\right) \coloneqq E_{Q_{0}}\left( Y\right)$, the prevalence of disease in this population, under a non-parametric model $\mathcal{Q}$. In the target
population, referred to as source 1, we obtain a random
sample of $W^{\left( 1\right) }\coloneqq\left( V,X\right) $. Thus, for $v\in
\left\{ 0,1\right\} ,$ $x\in \mathbb{R}^{p},$ $Q\coloneqq Q_{0}$ and $P%
\coloneqq P_{0},$ we have 
\begin{equation}
Q\left( V=v,X\leq x\right) =P\left( V=v,X\leq x|S=1\right).
\label{eq:prevalence-source1}
\end{equation}

Suppose we have access to an external validation study, i.e. a random sample of $W^{\left( 2\right) }\coloneqq%
\left( V,X,Y\right) $ from an external population, referred to as
source 2. Assume the sensitivity and specificity of $V$ within levels of the covariates are the same in sources 1 and 2 and the support of $X$ in source 1 and 2 are equal. Then we have that for $v,y\in \left\{ 0,1\right\}
,Q=Q_{0}$ and $P=P_{0}$ it holds that,%
\begin{equation}
Q\left( V=v|Y=y,X\right) =P\left( V=v|Y=y,X;S=2\right) \text{ a.e.- }Q.
\label{eq:prevalence-source2}
\end{equation}
Writing $\psi(Q) = E_Q[m_Q(X, V)]$ where 
\begin{equation}
m_{Q}\left( x,v\right) \coloneqq \frac{v-E_{Q}\left( V|Y=0,X=x\right) }{E_{Q}\left(
V|Y=1,X=x\right) -E_{Q}\left( V|Y=0,X=x\right) }  \label{eq:mq}
\end{equation}
we see that, if $V$ and $Y$ are dependent given $X$ under $Q$, then $\psi(Q)$ equals 
\begin{align*}
    \varphi(P) \coloneqq E_{P}\left[\frac{V-E_{P}\left[ V|Y=0,X, S=2\right] }{E_{P}\left[
V|Y=1,X, S=2\right] -E_{P}\left[ V|Y=0,X, S=2\right] }\Big|S=1\right]
\end{align*}
under alignments \eqref{eq:prevalence-source1} and \eqref{eq:prevalence-source2}. In fact, the entire distribution $Q$ is determined by $P$ and the model $\mc{P}$ is nonparametric under these assumptions as shown in \cref{app:examples-extras}. 

We conclude that this setting gives rise to a fused-data framework $\left( 
\mathcal{Q},\mathcal{P},\mathcal{C},\psi ,\varphi \right) $ in which models $%
\mathcal{Q}$ and $\mathcal{P}$ are non-parametric. The collection $%
\mathcal{C}$ of alignments is given by $\left( \ref{eq:prevalence-source1}%
\right) $ and $\left( \ref{eq:prevalence-source2}\right)$. \cite{cole_illustration_2023} consider this fused-data framework without baseline covariates. 

By viewing $S$ as a missingness indicator for $Y$ instead of a source indicator, this fused-data framework is closely related to assuming that $Y$ is missing not at random, but that we have access to a so-called shadow variable $V$ for $Y$ \citep{dhaultfoeuille_new_2010, wang_instrumental_2014, miao_identification_2018, li_non-parametric_2023, wang_identification_2024, park_single_2024}. The distinction between the two frameworks lies in the target distribution; in data fusion, this distribution is the distribution of $(X, Y)$ among individuals where $Y$ is missing. In shadow variables, the target distribution is the distribution of $(X, Y)$ had $Y$ always been observed. Despite the differences, the theory in this paper applies with slight modification to shadow variables analysis. In \cref{app-subsec:shadow-variables} we discuss the relationship in detail and expand this example to accommodate the subtleties arising when $V$ and $Y$ are non-binary.
\end{example}

\begin{example}[Two-Sample instrumental variables under a linear structural
equation model]
\label{example:tsiv-lsm}

Let $W=\left( L,X,Y\right) $ where $L,X$ and $Y$ are scalar random
variables. Suppose the law $Q_{0}$ of $W$ belongs to the collection $%
\mathcal{Q}$ of distributions $Q$ satisfying: (i) $X$ and $L$ and correlated
under $Q$ and, (ii) neither $Y$ nor $X$ is a deterministic function of the other (iii) there exist unique scalars $\alpha (Q)$ and $\psi (Q)$
solving 
\begin{equation}
E_{Q}\left[ Y-\alpha -\psi X|L\right] =0\text{ a.e.- }Q.  \label{a1}
\end{equation}%
We are interested in estimating $\psi(Q_0)$. In causal inference, model \eqref{a1} arises from
assuming a linear structural equation model for an outcome $Y$ on an
endogenous treatment $X$ with an instrument $L$ (\cite{anderson_estimation_1949}). In the case where $L$ is
binary, $\psi \left( Q\right) $ admits a different causal interpretation
under a different set of structural assumptions (\cite%
{imbens_identification_1994}). For a review see Chapter 12 of \cite%
{hansen_econometrics_2022}.

For $L$ polytomous or continuous, model $\mathcal{Q}$ is semiparametric. When $L$ is binary, $\mc{Q}$ is only restricted by the already assumed non-zero correlation between $X$ and $%
L$.

Importantly, $\psi \left( Q\right) $ depends on $Q$ solely through the
conditional distributions of $X|L$ and of $Y|L$. Therefore, it is possible
to identify $\psi \left( Q\right) $ from two separate samples, each
providing information about one of these conditional expectations (\cite%
{klevmarken_missing_1982}). Specifically, suppose we obtain a random sample
of $W^{\left( 1\right) }\coloneqq\left( L,Y\right) $ from source 1 and a random sample of $W^{\left( 2\right) }%
\coloneqq\left( L,X\right) $ from source 2  and we assume for all $y,x\in \mathbb{R}$
\begin{equation}
Q\left( Y\leq y|L\right) =P\left( Y\leq y|L;S=1\right) \text{ and }Q\left(
X\leq x|L\right) =P\left( X\leq x|L;S=2\right) \text{ a.e.- }Q.  \label{a2}
\end{equation}
This scenario fits the fused-data
framework $\left( \mathcal{Q},\mathcal{P},\mathcal{C}\text{\thinspace },\psi
,\varphi \right) $ if we assume that the supports of the instrument $L$ in
sources 1 and 2 include the support of $L$ under $Q_{0}$. For a non-dichotomous $L,$ the model $\mathcal{P}$ for the law $P$ of the observed data is semiparametric because the model $\mathcal{Q}$ restricts the distribution of $Y$ given $L$. The
alignments in $\mathcal{C}$ are as in $\left( \ref{a2}\right) $
and $\varphi \left( P\right)$ is the solution to $E_{P}\left(
Y|L,S=1\right) -\tau -\varphi E_{P}\left( X|L,S=2\right) =0$. This
fused-data framework is known as the Two-Sample Instrumental Variables
(TSIV) Model and has been extensively studied (\cite{zhao_two-sample_2019,
pierce_efficient_2013, pacini_two-sample_2019, graham_efficient_2016,
ridder_chapter_2007, shi_data_2023, shu_improved_2020}).
\end{example}

\begin{example}[Transporting average treatment effects]
\label{example:transporting}

Suppose that we have access to data drawn from two populations, the second one being the target population of interest. The data
available from the first population comes from a prospective cohort study in
which we measure $W =  (A, L, Y)$ where $A$ is a binary treatment, $Y$ is a binary disease outcome, and $L$ is a $p-$vector of baseline covariates. For the data available from the target
population we shall consider the following three scenarios. In scenario (i), we draw a random sample and only measure the covariates. In scenario (ii) we draw a random sample from a disease registry and measure covariates and treatments. In scenario (iii) we perform a case-control study and measure outcomes, covariates, and treatments.

Suppose we are interested in estimating the average treatment effect (ATE)
in the target population. The
data at our disposal under any of the three scenarios is insufficient to
identify it, even when $L$ suffices to control for confounding. So, we wish
to incorporate data from the cohort study in the first population with the hope of identifying the ATE in the target population. In \cref{app-subsec:causal-identification} we discuss the causal identification of the ATE in the target population under each of these scenarios. In all scenarios, under appropriate causal assumptions, we can recast the causal effect of interest as the evaluation at a law $Q_{0}$ for $W$ of the
functional $\psi :\mathcal{Q\rightarrow }\mathbb{R}$ defined as 
\begin{equation}
\psi \left( Q\right) \coloneqq E_{Q}\left[ Q\left( Y=1|L,A=1\right) -Q\left(
Y=1|L,A=0\right) \right]  \label{eq:psi-example3}
\end{equation}%
where $\mathcal{Q}$ is a nonparametric model of laws on $W$ and the definition of $Q_{0}$ varies depending on the scenario. The first scenario has been extensively studied (\cite%
{pearl_transportability_2011, rudolph_robust_2017, dahabreh_extending_2020,
dahabreh_towards_2020, shi_data_2023, li_efficient_2023}). \cite{jia_identifiability_2006, chatterjee_constrained_2016} discuss estimation of parameters with fused data assuming alignments as in the third scenario. Their approaches differ from ours in that we allow for a nonparametric ideal data model as opposed to the parametric setting of those authors. 

In the three aforementioned scenarios, the available data comprises a random sample
of observed data $O=\left( c\left( W,S\right) ,S\right) ,$ where $S=1$ if
the observation originates from the cohort study and $S=2$ if it stems from
the study conducted on the target population. On the other hand, $%
c\left( W,j\right) $ indicates the subset of the ideal vector $W,$ denoted
with $W^{\left( j\right) },$ that is available when $S=j,j=1,2.$ Thus, $%
c\left( W,1\right) \coloneqq W^{\left( 1\right) }=\left( L,A,Y\right) $ and $%
c\left( W,2\right) \coloneqq W^{\left( 2\right) }$ varies based on the data
available from each study design.

\textbf{Scenario (i):} We observe only $L$ in source $2$ and hence $c(W, 2) = L$. Under the assumptions in \cref{app-subsec:causal-identification}, the supports of $\left( L,A\right) $ under $P_0(\cdot|S=1)$ and under $Q_{0}$ are equal and the following alignments hold for $Q=Q_0$ and $P=P_0$ 
\begin{align}
Q\left( Y=1|L,A\right) =&P\left( Y=1|L,A,S=1\right)
 \text{ a.e.-}Q\label{eq:scenario-1-align-1} \\
Q\left( L\leq l\right) =&P\left( L\leq l|S=2\right) \text{ for all $l\in \mathbb{R}^p$}\label{eq:scenario-1-align-2}.
\end{align}
 Model $\mathcal{P}$ is non-parametric because the
existence of a law $Q$ satisfying $\eqref{eq:scenario-1-align-1}$ and $\eqref{eq:scenario-1-align-2} $ does not place equality constraints on $P$. Moreover, $\psi
\left( Q_{0}\right) =\varphi \left( P_{0}\right) $ where for any $P$ 
\begin{equation}
\varphi \left( P\right) \coloneqq E_{P}\left[ P\left( Y=1|L,A=1,S=1\right)
-P\left( Y=1|L,A=0,S=1\right) |S=2\right].  \label{phi-def-example3}
\end{equation}%
We conclude that the present scenario fits the fused-data framework $%
\left( \mathcal{Q},\mathcal{P},\mathcal{C}\text{\thinspace },\psi ,\varphi
\right) $ with $\mathcal{C}$ the collection of alignments defined by \eqref{eq:scenario-1-align-1} and \eqref{eq:scenario-1-align-2}. Note that we have assumed
neither that $P_{0}\left( Y=1|L,A,S=1\right) =P_{0}\left( Y=1|L,A,S=2\right) 
$ nor that $P_{0}\left( L|S=2\right) =P_{0}\left( L|S=1\right) $ hold,
implying that $Q_{0}$ is not necessarily the distribution of the complete
data $\left( L,A,Y\right) $ in either source. In fact, there exist
infinitely many distributions $Q_{0}$ that satisfy \eqref{eq:scenario-1-align-1}
and \eqref{eq:scenario-1-align-2} for $Q=Q_{0}$ and $P=P_{0}$ because these
equations do not restrict the propensity score $Q_{0}\left( A=1|L\right) .$
This is an example of a scenario where the alignments ensure that $\psi
\left( Q_{0}\right) $ is identified by $P_{0}$, but do not ensure that $Q_0$
is identified by $P_{0}$.\medskip

\textbf{Scenario (ii):} We observe $W=(L, A, Y)$ in source $2$ and hence $c(W, 2) = W$. Under the assumptions in \cref{app-subsec:causal-identification}, the following alignments hold for $Q=Q_0, P=P_0$ 
\begin{align}
\label{eq:scenario-2-align-1}Q\left(Y=1|L, A\right) =&P\left(Y=1|L, A,S=1\right) \text{ a.e.-$Q$}\\
Q\left( L\leq l,A=a|Y=1\right) =&P\left( L\leq l,A=a|Y=1,S=2\right) \text{ for all $l \in \mathbb{R}^p$, $a\in \{0, 1\}$}.
\label{eq:scenario-2-align-2}
\end{align}
Assuming $Q\left( Y=1|L,A\right) >0$ a.e.- $Q$, in \cref{app:examples-extras} we show that under the alignments \eqref{eq:scenario-2-align-1} and \eqref{eq:scenario-2-align-2}, $q(l) = \sum_{a=0}^1 \beta(l, a, 1;P)$ where
\begin{align*}
\beta\left(a,l, y;P\right)\coloneqq&\left. \frac{p\left( l,a|Y=y,S=2\right) }{P\left(
Y=y|l,a,S=1\right) }\right/ \sum_{a=0}^{1}\int \frac{p\left( l^{\prime
},a|Y=y,S=2\right) }{P\left( Y=y|l^{\prime },a,S=1\right) }dl^{\prime }.
\end{align*}%
Furthermore, with $\alpha(l;P) \coloneqq  \sum_{a=0}^1\beta(l, a, 1;P)$, $\psi \left( Q_{0}\right) $ is equal to 
\begin{equation}
\varphi \left( P_{0}\right) \coloneqq\int \left\{ P_{0}\left(
Y=1|l,A=1,S=1\right) -P_{0}\left( Y=1|l,A=0,S=1\right) \right\} \alpha\left(
l;P_{0}\right) dl  \label{idenphi}
\end{equation}

Since we make no assumptions on $Q_{0}$ other than assumptions on its support, model $\mathcal{Q}$ is non-parametric. Furthermore,
since the existence of $Q_{0}$ satisfying \eqref{eq:scenario-2-align-1} and \eqref{eq:scenario-2-align-2} places only support restrictions on $P_{0}$ (see 
\cite{arnold_specification_1996}, \cref{app:examples-extras}), we conclude that the only equality restriction model $\mathcal{P}$
imposes on the observed data law $O$ is that the law of $Y|S=2$ assigns
probability 1 to $Y=1$. Letting $\mathcal{C}$ be the collection of
alignments \eqref{eq:scenario-2-align-1} and \eqref{eq:scenario-2-align-2} and $\varphi \left( P\right) $ defined as in $\left( \ref{idenphi}\right) $ with $%
P$ instead of $P_{0},$ we thus conclude that this scenario also fits the
fused-data framework $\left( \mathcal{Q},\mathcal{P},\mathcal{C}\text{%
\thinspace },\psi ,\varphi \right) $.\medskip

\textbf{Scenario (iii):} In this scenario, we observe $W=(L, A, Y)$ in source $2$ and hence $c(W, 2) = W$. It is well known that the causal risk difference cannot be identified from case-control data, even if the covariates $L$ are sufficient to control for confounding. One strategy for identification involves supplementing the case-control study with a consistent estimator of the prevalence of disease in the target population (\cite{van_der_laan_estimation_2008}). In cases where such an estimator is unavailable, we are still able to identify the causal risk difference in the target population by integrating the case-control study with data from the prospective cohort study if we are prepared to make additional assumptions that ultimately, allow us to identify the prevalence in the target population. Specifically, under the assumptions discussed in \cref{app-subsec:causal-identification}, the following alignments hold for $Q=Q_0$ and $P=P_0$ 
\begin{align}
\label{eq:scenario-3a-align-1}Q\left(Y=1|L, A\right) =&P\left(Y=1|L, A,S=1\right) \text{ a.e.-$Q$}\\
Q\left( L\leq l,A=a|Y\right) =&P\left( L\leq l,A=a|Y,S=2\right) \text{ a.e.-$Q$}
\label{eq:scenario-3a-align-2}
\end{align}
for all $l \in \mathbb{R}^p$, $a\in \{0, 1\}$. Similarly to scenario (ii), under the additional assumption that $0 <Q\left(Y=1|L, A\right) < 1$ a.e.-$Q$, $\psi(Q_0)$ is equal to $\varphi(P_0)$ as defined in \eqref{idenphi}. In contrast to scenario (ii), the additional assumption that the conditional distribution $L, A|Y=0$ aligns in source 2 implies the equality constraint $\beta(l, a, 1;P) = \beta(l, a, 0;P)$ for all $l\in \mathbb{R}^p, a\in \{0, 1\}$ on the observed data law $P$. Thus, the model $\mc{P}$ is strictly semiparametric. To recap, this scenario fits the fused-data framework $\left( \mathcal{Q},%
\mathcal{P},\mathcal{C}\text{\thinspace },\psi ,\varphi \right) $, now with
a non-parametric model $\mathcal{Q},$ $\mathcal{C}$ being comprised of the
alignments \eqref{eq:scenario-3a-align-1} and \eqref{eq:scenario-3a-align-2} and $\psi (Q_{0})$ and $\varphi (P_{0})$ as defined in
scenario (ii). However, unlike the preceding two scenarios, although model $%
\mathcal{Q}$ remains non-parametric, the assumption of the very existence of
a common distribution $Q$ that satisfies \eqref{eq:scenario-3a-align-1} and \eqref{eq:scenario-3a-align-2} gives rise to a strict
semiparametric model $\mathcal{P}$. 
\end{example}

\cite{li_efficient_2023} developed a comprehensive theory for semiparametric estimation in a special type of fused-data framework. They consider scenarios where the aligned conditionals in
each source $j$ corresponds to a source-specific subset of the factors in
\begin{equation*}
q_{0}\left( W\right) =q_{0}\left( W_{1}\right) \times q_{0}\left(
W_{2}|W_{1}\right) \times q_{0}\left( W_{3}|\overline{W}_{2}\right) \times
\cdots \times q_{0}\left( W_{K}|\overline{W}_{K-1}\right)
\end{equation*}%
and the observed data vector $W^{\left( j\right) }$ in each source $j$
suffices to identify the aligned components. While their framework applies in significant contexts, such as in \cref{example:transporting}, scenario (i), it does not apply to any other examples discussed above. We aim to develop a versatile theory applicable to a broader class of fused-data frameworks. As in \cite%
{li_efficient_2023}, we shall establish a template for
calculating influence functions and in particular, the efficient influence
function, of $\varphi :\mathcal{P}\mapsto \mathbb{R}$ from influence
functions of $\psi :\mathcal{Q}\mapsto \mathbb{R}$. In \cref{app:semiparametric-theory} we review the central role played by influence functions for constructing
semiparametric efficient, debiased machine learning estimators.

\subsection{Alignment assumptions and the fused-data model definition}

\label{subsec:alignment-assumptions}

In this subsection, we give the alignment
assumptions made by our theory. We first provide an informal description to facilitate understanding.

We assume that within each data source $j$, one or more conditional
distributions of the subvector $W^{\left( j\right) }$ align with the
corresponding conditionals under $Q_{0}$. While we stipulate that within
each source $j$, these aligned conditionals correspond to some or all
factors of a particular factorization of the joint law of $W^{\left(
j\right) }$, we do not mandate that across all sources, the aligning
conditionals correspond to a single factorization of the distribution of the
full data vector $W$. Specifically, we assume that for each $j\in \left[ J\right]
$, there exists a partition of the set $\omega^{(j)}\coloneqq \{W^{(j)}_k\}_{1\leq k\leq \text{dim}(W^{(j)})}$ such that with $Z_k^{(j)}$, $k=1, \dots, K^{(j)}$, denoting the distinct subsets of the partition of $\omega^{(j)}$, and with $Z^{(j)} \coloneqq (Z_1^{(j)}, \dots, Z_{K^{(j)}}^{(j)})$, the aligned distributions in source $j$ correspond to a subset of the factors in
\begin{equation*}
q_{0}( Z^{( j) }) =q_{0}( Z_{1}^{( j)
}) \times q_{0}( Z_{2}^{( j) }|Z_{1}^{( j)
}) \times q_{0}( Z_{3}^{( j) }|\overline{Z}%
_{2}^{( j) }) \times \cdots \times q_{0}( Z_{K^{(
j) }}^{( j) }|\overline{Z}_{K^{( j) }-1}^{(
j) }).
\end{equation*}
Here and throughout, in a slight abuse of notation each $Z_k^{(j)}$ is regarded as a random vector of the elements of the corresponding subset of $\omega^{(j)}$, with the entries of $Z_k^{(j)}$ sorted in an arbitrary order. For purposes of exposition, we have assumed a product measure $\mu$ that dominates $Q_0$ exists and let $q_0 \coloneqq \frac{dQ_0}{d\mu}$, but our theory does not require this. 

In addition, to accommodate settings such as \cref{example:transporting}
scenario (ii), where alignment occurs only for a subset of the elements in
the support of the conditioning variables, for each $k$ in $\left\{
2,...,K^{\left( j\right) }\right\} ,$ we let $\overline{\mathcal{Z}}%
_{k-1}^{\left( j\right) }$ be a set such that the distributions $Z_k^{(j)}|\overline{Z}_{k-1}^{(j)} = \overline{z}_{k-1}^{(j)}$ under $Q_0$ and in source $j$ align if and only if  $
\overline{z}_{k-1}^{\left( j\right) }$ is in $\overline{\mathcal{Z}}%
_{k-1}^{\left( j\right) }$. We additionally assume that if $\overline{\mathcal{Z}}%
_{k-1}^{\left( j\right) }\not=\emptyset $, $\overline{\mathcal{Z}}%
_{k-1}^{\left( j\right) }$ has positive probability under $Q_0$ and $P_0(\cdot|S=j)$. With this definition $\overline{\mathcal{Z}}_{k-1}^{\left( j\right)
}=\emptyset$, for some $k\geq 2$, indicates the distribution $Z_k^{(j)}|\overline{Z}_{k-1}^{(j)}$ under $Q_{0}$ does not align with the corresponding conditional  in the $j^{th}$ source. Additionally, since $\overline{Z}_{0}^{(j)}$ does not exist, we require a notational convention to distinguish between alignment
and no alignment of the marginal distributions of $Z_{1}^{(j)}$ under $%
P_{0}(\cdot |S=j)$ and $Q_{0}$. We define $\overline{z}%
_{0}^{\left( j\right) }\coloneqq\ast $ and let $\overline{\mathcal{Z}}%
_{0}^{(j)}$ be either $\left\{ \ast \right\}$ or $\emptyset$. Thus, $\overline{Z}_0^{(j)} = \{\ast\}$ if and only if the marginal distribution of $Z_0^{(j)}$ under $Q_0$ aligns with the corresponding marginal in the $j^{th}$ source.

The collection of assumed alignments is thus fully determined by the
collection 
\begin{equation}
\mathcal{C}\coloneqq\{ ( W^{( j) },\{ {Z}_{k}^{( j)
}\} _{k=1,...K^{( j) }},\{ \overline{\mathcal{Z}}_{k-1}^{( j)
}\} _{k=1,...K^{( j) }}) :j=1,...,J\}.
\label{colect}
\end{equation}

\begin{figure}[!ht]
\caption{Illustration of a fused-data model}
\centering
\resizebox{.8\textwidth}{!}{%
\fbox{\begin{circuitikz}
\tikzstyle{every node}=[font=\large]
\draw  (2.75,13) rectangle (4.25,11.5);
\draw  (4.25,13) rectangle (5.75,11.5);
\draw  (5.75,13) rectangle (7.25,11.5);
\draw  (7.25,13) rectangle (8.75,11.5);
\draw  (8.75,13) rectangle (10.25,11.5);
\draw  (10.25,13) rectangle (11.75,11.5);
\draw  (11.75,13) rectangle (13.25,11.5);
\draw  (13.25,13) rectangle (14.75,11.5);
\node [font=\large, align=center] at (0, 12.25) {Ideal Data\\Vector};
\node [font=\large, align=center] at (0, 8.75) {Aligned\\Components};
\node [font=\large, align=center] at (0, 6.5) {Source};
\node [font=\large] at (3.5,12.25) {$W_1$};
\node [font=\large] at (5,12.25) {$W_2$};
\node [font=\large] at (6.5,12.25) {$W_3$};
\node [font=\large] at (8,12.25) {$W_4$};
\node [font=\large] at (9.5,12.25) {$W_5$};
\node [font=\large] at (11,12.25) {$W_6$};
\node [font=\large] at (12.5,12.25) {$\dots$};
\node [font=\large] at (14,12.25) {$W_{K}$};
\draw [short] (3.5,11.5) -- (3.5,10.5);
\draw [short] (9.5,11.5) -- (3.5,10.5);
\node [font=\large] at (4.5,8.5) {$Q(\underbrace{W_1, W_5}_{Z_2^{(1)}}|\underbrace{W_2=w_2, W_3}_{Z_1^{(1)}})$};
\draw [->, >=Stealth] (3.5,10.5) -- (3.25,9.5);
\draw [short] (5,11.5) -- (5.5,10.25);
\draw [short] (6.5,11.5) -- (5.5,10.25);
\draw [->, >=Stealth] (5.5,10.25) -- (5.5,9.5);
\node [font=\large] at (11,8.5) {$Q(\underbrace{W_2}_{Z_3^{(2)}}|\underbrace{W_1}_{Z_2^{(2)}}, \underbrace{W_6}_{Z_1^{(2)}})$};
\node [font=\large] at (14.75,8.5) {$Q(\underbrace{W_6}_{Z_1^{(2)}})$};
\draw [->, >=Stealth, dashed] (5,11.5) -- (9.75,9.5);
\draw [->, >=Stealth, dashed] (3.5,11.5) -- (10.75,9.5);
\draw [->, >=Stealth, dashed] (11,11.5) -- (12,9.5);
\draw [->, >=Stealth, dashed] (11,11.5) -- (14.75,9.5);

\node [font=\large] at (4.5,6.5) {$S=1$};
\node [font=\large] at (12.5,6.5) {$S=2$};
\end{circuitikz}}
}%
\label{fig:example-figure}
\end{figure}

To illustrate the complex alignments allowed by a fused-data model consider \cref{fig:example-figure} which corresponds to a model with two sources. In the first source, we observe the subvector $W^{(1)} = (W_1, W_2, W_3, W_5)$ which we decompose as $Z^{(1)}=(Z^{(1)}_1, Z^{(1)}_2)$ where $Z_1^{(1)} = (W_2, W_3)$ and $Z_2^{(2)} = (W_1, W_5)$. In the second source, we observe the subvector $W^{(2)} = (W_1, W_2, W_6)$ which we decompose as $Z^{(2)} =(Z^{(2)}_1, Z^{(2)}_2, Z^{(2)}_3)$ where $Z_1^{(2)} = W_6$, $Z_2^{(2)} = W_1$, and $Z_3^{(2)} = W_2$. We then assume the conditional distribution of $Z_2^{(1)}$ given $Z_1^{(1)} = z_1^{(1)}$ for $z_1^{(1)}$ in $\{w_2\}\times \textsf{Supp}\left[W_3;Q\right]$ in source 1 aligns with the corresponding conditional distribution of $Q$. Likewise, the conditional distribution of $Z^{(2)}_3$ given $(Z^{(2)}_2, Z^{(2)}_1)$ and the marginal distribution of $Z^{(2)}_1$ in source 2 align with the corresponding conditional and marginal distributions of $Q$. These alignments are encoded in the class 
$\mc{C} = \{(W^{(1)}, \{{Z}_1^{(1)}, {Z}_2^{(1)}\}, \{\overline{\mc{Z}}_0^{(1)}, \overline{\mc{Z}}_1^{(1)}\}), (W^{(2)}, \{{{Z}}_1^{(2)}, {{Z}}_2^{(2)}, {{Z}}_3^{(3)}\}, \{\overline{\mc{Z}}_0^{(2)}, \overline{\mc{Z}}_1^{(2)}, \overline{\mc{Z}}_1^{(3)}\})\}$
where $\overline{\mc{Z}}_0^{(1)} = \emptyset$, $\overline{\mc{Z}}_1^{(1)} = \left\{w_2\right\}\times\textsf{Supp}\left[W_6;Q\right]$, $\overline{\mc{Z}}_0^{(2)} = \{*\}$, $\overline{\mc{Z}}_1^{(2)} = \emptyset$, and $\mc{Z}_3^{(2)} = \textsf{Supp}\left[W_6;Q\right]\times \textsf{Supp}\left[W_1;Q\right]$.

We are now ready to give the formal definition of alignment in our theory. 
\begin{definition}
\label{def:alignswith}
We say that a law $P$ on the sample space of $O\coloneqq( c(
W,S) ,S) $ and a law $Q$ on the sample space of $W$ are \emph{aligned}
relative to a collection $\mathcal{C}$ defined as in $( \ref{colect}%
) $ if and only if $S\in [ J] ,c( W,j)
=W^{( j) }$ and
\begin{enumerate}
\item $\overline{\mathcal{Z}}_{k-1}^{( j) }\subseteq \mathsf{Supp}%
[ \overline{Z}_{k-1}^{( j) };P( \cdot |S=j) %
] \cap \mathsf{Supp}[ \overline{Z}_{k-1}^{( j) };%
Q] $ for $k\in \{ 2,...,K^{( j) }\} ,$ 
$j\in \lbrack J],$

\item If $\overline{\mathcal{Z}}_{k-1}^{( j) }\not=\emptyset ,$
then $P( \overline{Z}_{k-1}^{( j) }\in \overline{\mathcal{Z}}%
_{k-1}^{( j) }|S=j) >0$ and $Q( \overline{Z}%
_{k-1}^{( j) }\in \overline{\mathcal{Z}}_{k-1}^{( j)
}) >0;$ and

\item For all $z_{k}^{(j)}\in \mathbb{R}^{\dim ( Z_{k}^{(j)}) },$ 
$j\in [ J] $ and $k\in [ K^{( j) }] ,$ $P( Z_{k}^{( j) }\leq z_{k}^{(j)}\vert \overline{%
Z}_{k-1}^{( j) },S=j) =Q( Z_{k}^{(
j) }\leq z_{k}^{(j)}\vert \overline{Z}_{k-1}^{( j)
})$ a.e.-$Q$ on $\overline{\mathcal{Z}}_{k-1}^{(
j) }$.
\end{enumerate}
\end{definition}
Hereafter, the symbol $P\alignswith Q$ denotes alignment of $P$ and $Q$ relative to $%
\mathcal{C}$.

We can now give a precise definition of a fused-data model. To ease future
notation, given a law $Q\in \mc{Q}$ and a law $P_{0}$ on
the sample space of $O\coloneqq(c(W,S),S),$ we define%
\begin{align*}
\mathcal{P}_{Q}\coloneqq\{ P:P \alignswith Q, P\text{ mutually absolutely continuous with }P_{0}\} .
\end{align*}%
Note that even though $\mathcal{P}_{Q}$ depends on $P_{0}$ through the
requirement of mutual absolute
continuity with $P_{0},$ we do not make this dependence explicit in the
notation.

\begin{definition}
\label{def:fdm}
Given $Q_{0}$ a law for $W$ and a law $P_{0}$ for $O$, the triplet $( \mathcal{Q},\mathcal{P},\mathcal{C}%
) $ is a \emph{fused-data model} with respect to $( Q_{0},P_{0}) $
if and only if $P_{0}\alignswith Q_{0}$ and
\begin{enumerate}
\item $\mathcal{Q}$ is a collection of laws $Q$ on the sample space of $W$
that are mutually absolutely continuous with $Q_{0},$ and $Q_{0}\in \mathcal{%
Q}$.

\item $\mathcal{C}$ is defined as in $( \ref{colect}) $ for some $%
W^{( j) }\subseteq W$, $\bigcup_{k\in [K^{(j)}]} Z_k^{(j)} = \{W^{(j)}_l\}_{1\leq l \leq \text{dim}(W^{(j)})}$ and $Z_k^{(j)} \cap Z_{k'}^{(j)} = \emptyset$ for $k \not= k'$,  $\overline{\mathcal{Z}}_{k-1}^{(
j) }\subseteq \mathsf{Supp}[ \overline{Z}_{k-1}^{( j)
};\mathcal{Q}] $, $Q( \overline{Z}_{k-1}^{( j) }\in 
\overline{\mathcal{Z}}_{k-1}^{( j) }) >0$ whenever $\overline{\mathcal{Z}}_{k-1}^{( j) } \not= \emptyset$ for all $Q\in 
\mathcal{Q}$, $k$ in $\{ 2,...,K^{( j) }\} $ and $j\in
\lbrack J]$.

\item $\mathcal{P}=\bigcup_{Q\in \mathcal{Q}}\mathcal{P}_{Q}$ and $P_{0}\in 
\mathcal{P}.$
\end{enumerate}
\end{definition}

As discussed earlier, model $\mathcal{P}$ may impose equality constraints on 
$P$ for two reasons: first, if the law $P$ inherits equality constraints
imposed by the model $\mathcal{Q}$, and second, if the mere existence of a
single $Q$ that aligns with the specific conditionals in each source places
restrictions on $P$, as seen in \cref{example:transporting} scenario
(iii). This distinction is crucial because it affects the structure of the
set of influence functions of pathwise differentiable functionals of the
observed data law $P$. To facilitate a clear distinction in later sections
between the two sources of restrictions in model $\mathcal{P}$, we introduce
the following definition.
\begin{definition}
\label{def:extended-obs-data-model}
Given a fused-data model $\left( \mathcal{Q},\mathcal{P},\mathcal{C}\right) $
with respect to $\left( Q_{0},P_{0}\right) $, the \emph{extended observed data
model} is defined as $\mathcal{P}^{ext}=\bigcup_{Q\in \mathcal{Q}^{np}}%
\mathcal{P}_{Q}$ where $\mathcal{Q}^{np}\supseteq \mathcal{Q}$ is the
collection of all laws on the sample space of $W$ that are mutually
absolutely continuous with $Q_{0}$.
\end{definition}
Note that $\mathcal{P}^{ext}$ imposes equality constraints on $P$
if and only if the mere existence of some law $Q$ that aligns with $P$ according to $\mathcal{C}$ creates
equality constraints on $P$.

Next we will formalize the assertion that $P_{0}$ identifies the summary $\psi \left( Q_{0}\right) $ if $%
P_{0} \alignswith Q_{0}$ and $\psi $ depends on $Q$ only through
the aligned conditional distributions. To formally address this concept, we
begin by defining an equivalence relation $\overset{\mathcal{C}}{\sim }$ on $%
\mathcal{Q}$. Given a fused-data model $(\mc{Q}, \mc{P}, \mc{C})$, for a pair of laws $Q$ and $\widetilde{Q}\in \mathcal{Q}$ we
write $Q\overset{\mathcal{C}}{\sim }\widetilde{Q}$ whenever for all $%
z_{k}^{(j)}\in \mathbb{R}^{\dim \left( Z_{k}^{(j)}\right) },k\in \left[
K^{\left( j\right) }\right] $ and $j\in \left[ J\right] $ it holds that%
\begin{equation}
Q\left( \left. Z_{k}^{\left( j\right) }\leq z_{k}^{(j)}\right\vert \overline{%
Z}_{k-1}^{\left( j\right) }\right) =\widetilde{Q}\left( \left. Z_{k}^{\left(
j\right) }\leq z_{k}^{(j)}\right\vert \overline{Z}_{k-1}^{\left( j\right)
}\right) \text{ a.e.- }Q\text{ on }\overline{\mathcal{Z}}_{k-1}^{\left(
j\right) }.  \label{equiv}
\end{equation}

By the assumption of mutual absolute continuity of the laws in $\mathcal{Q}$%
, $\left( \ref{equiv}\right) $ holds if and only if the same equality holds
a.e.- $\widetilde{Q}$ on $\overline{\mathcal{Z}}_{k-1}^{\left( j\right) }.$
This, in turn, implies that the relation $\overset{\mathcal{C}}{\sim }$ is
transitive. Consequently, defining for any $Q\in \mathcal{Q}$ the
equivalence class $\xi \left( Q;\mathcal{C}\right) \coloneqq\left\{ \widetilde{Q}:Q\overset{%
\mathcal{C}}{\sim }\widetilde{Q},\widetilde{Q}\in \mathcal{Q}\right\}$
we have that given $Q\in \mathcal{Q},$ $P \alignswith Q$ if and only if $P \alignswith \widetilde{Q}$ for every $\widetilde{Q}$ in $\xi \left( Q;\mathcal{C}\right) $%
. Note that even though $\xi \left( Q;\mathcal{C}\right) $ depends on model $%
\mathcal{Q},$ we do not make this dependence explicit in the notation. The
next assumption formalizes the assertion that $\psi \left( Q\right) $
depends on $Q$ solely through its aligned components.

\begin{assumption}
\label{as:identification} $Q \overset{\mathcal{C}}{\sim} \widetilde{Q}$
implies that $\psi \left( Q\right) =\psi \left( \widetilde{Q}\right) .$
\end{assumption}

Under \cref{as:identification}, there exists a mapping from $\mathcal{P}$
to $\mathbb{R}$ such that when evaluated at $P\in \mathcal{P}$ it is equal
to $\psi \left( Q\right) $ for any $Q\in \mathcal{Q}$ such that $P \alignswith Q$. Specifically, let $\Xi $ denote the partition
of $\mathcal{Q}$ into equivalence classes $\xi \left( Q;\mathcal{C}\right)$ and define the map $\Psi :\Xi \mapsto \mathbb{R}$ such that $\Psi \left[ \xi
\left( Q;\mathcal{C}\right) \right] $ is equal to the unique value $\psi (%
\widetilde{Q})$ taken by all $\widetilde{Q}\in \xi \left( Q;\mathcal{C}%
\right) $. Next, define the map $\Phi(\cdot;\mc{C}) :\mathcal{P\rightarrow }\Xi $, which
assigns to every $P\in \mathcal{P}$ the unique equivalence class $\xi \left(
Q;\mathcal{C}\right) $ such that $P \alignswith \widetilde{Q}$ for all $\widetilde{%
Q}\in \xi (Q,\mathcal{C})$. Finally, define $\varphi :\mathcal{P}$ $\mapsto 
\mathbb{R}$ as $\varphi \coloneqq\Psi \circ \Phi(\cdot;\mc{C})$

\begin{theorem}
\label{thrm:identifiability} Given a fused-data model $\left( \mathcal{Q},%
\mathcal{P},\mathcal{C}\right) $ with respect to $\left( Q_{0},P_{0}\right)
, $ if \cref{as:identification} holds then $\psi \left( Q\right) =\varphi
\left( P\right) $ for any $Q\in \mathcal{Q}$ and $P\in \mathcal{P}$ such
that $P \alignswith Q$.
\end{theorem}

We are now ready to give the precise definition of a fused-data framework. In \cref{app:examples-extras}, we define the fused-data frameworks in examples  \ref{example:disease-prevalence}-%
\ref{example:transporting}. 

\begin{definition}
\label{def:fdf} The quintuplet $\left( \mathcal{Q},\mathcal{P},\mathcal{C}%
,\psi ,\varphi \right) $ is a \emph{fused-data framework}
with respect to $\left( Q_{0},P_{0}\right) $ if and only if $\left( \mathcal{%
Q},\mathcal{P},\mathcal{C}\right) $ is a fused-data model with respect to $%
\left( Q_{0},P_{0}\right) ,$ $\psi :\mathcal{Q\rightarrow }\mathbb{R}$
satisfies \cref{as:identification} and $\varphi :\mathcal{P\rightarrow }%
\mathbb{R}$ defined as above.
\end{definition}

The identification result in \cref{thrm:identifiability} is general, but satisfying \cref{as:identification} in a fused-data framework can delicately depend on the choice of $\psi$, the ideal data model $\mc{Q}$, and the collection of alignments $\mc{C}$. In the simplest case, the full distribution $Q$ is determined by the aligned marginal and conditional distributions. \citet{gelman_characterizing_1993} provides sufficient (though not necessary; see \citet{gelman_corrigendum_1999}) conditions for identifying $Q$ from these aligned distributions. 

In other settings, the full ideal data distribution $Q$ might not be identified, but it is possible to directly express the target functional as a function of certain conditional or marginal distributions of $Q$ that are identified. This is the case in \cref{example:transporting}(i), where the treatment propensity $q(A|L)$ under $Q$ is unspecified, but the target functional can be directly expressed as a functional of the conditional distribution of $Y|A, L$ and the marginal distribution of $L$ under $Q$, each of which is aligned with a distinct source. \cite{jia_identifiability_2006} give sufficient conditions under which certain aligned conditionals and marginals are sufficient to identify other conditionals and marginals, from which the target parameter might be naturally expressed by. 

There is a substantial literature on the deriving graphical criteria under which counterfactual distributions are identified when combining experimental and observational results from a common (target) population \citep{bareinboim_transportability_2014, bareinboim_causal_2016, hunermund_causal_2023, kivva_revisiting_2022, lee_general_2020, lee_general_2024, pearl_transportability_2011, bareinboim_recovering_2015}. These works provide necessary and sufficient conditions for identifying counterfactual distributions using graphical criteria under a latent causal directed acyclic graph (DAG) model. In such settings, the alignments considered in this work must be made on the level of the latent variables. When all latent variables are observed, the present work applies directly to identified counterfactual parameters. When some latent variables are not observed, there is no guarantee that the implied alignments on the observed variables follow the structure considered in this work. We leave an extension of the present results to this setting for future work.


Intricate identification scenarios arise when the target parameter $\psi$ is identified by the alignments $\mc{C}$ only under additional restrictions on the ideal-data model $\mc{Q}$. This occurs in \cref{example:disease-prevalence}, where \cref{as:identification} requires that $Y$ and $V$ are correlated within levels of $X$ for every $Q \in \mc{Q}$ in order to identify nontrivial target functionals. \cite{li_non-parametric_2023} derive an analogous condition when $Y$ and $V$ are non-binary. We discuss this condition and implications for efficient inference in \cref{app-subsec:shadow-variables}.

\section{Main results}
\label{sec:main-results} 
\subsection{The score operator}

\label{subsec:score-operator-theory}

In this section we provide a high-level discussion of the strategy that we employ to derive all influence functions $\varphi_{P}^{1}$ of $\varphi 
$, in terms of the influence functions $\psi _{Q}^{1}$ of $\psi $. In what
follows, we refer to $\varphi _{P}^{1}$ as observed data influence
functions and $\psi _{P}^{1}$ as ideal data influence functions. Our
strategy builds on calculations that involve the so-called score operator as
discussed in Section 25.5 of \cite{van_der_vaart_asymptotic_2000}. A key
point driving our strategy is the observation that $\mathcal{P%
}$ in a fused-data model $\left( \mathcal{Q},\mathcal{P},\mathcal{C}\right) $
can be indexed by probability measures, one of which is the
ideal data law $Q$ as \cref{lemma:alt-model-characterization} below establishes. 

Here and throughout given $\left( \mathcal{Q},\mathcal{P},%
\mathcal{C}\right) $ a fused-data model with respect to $\left(
Q_{0},P_{0}\right) $ we let $\Lambda $ denote the set of all probabilities $%
\lambda $ on the sample space $\left[ J\right] $ for $S$ such that $\lambda
(S=j)>0$ for all $j\in \left[ J\right] $, and let $\mathcal{U}$ denote the
collection of $J$-tuplets $U\coloneqq(U^{(1)},\dots ,U^{(J)})$ where each $%
U^{(j)}$ is a law on the sample space of $Z^{\left( j\right) },$ mutually
absolutely continuous with the law of $Z^{\left( j\right) }$ under $%
P_{0}\left(\cdot|S=j\right) ,$ but otherwise unrestricted. Furthermore, given $Q\in \mathcal{Q%
},U\in \mathcal{U}$, and $\lambda \in \Lambda$, we define $P_{Q,U,\lambda }$ as the law on the sample space of $O$ such that for each $j\in[J]$
\begin{enumerate}
\item 
\begin{align*}
    P_{Q,U,\lambda }\left( Z_{1}^{\left( j\right) }\leq z_{1}^{(j)}|S=j\right) %
\coloneqq Q\left( Z_{1}^{\left( j\right) }\leq z_{1}^{(j)}\right) 
\end{align*}
if $\overline{\mathcal{Z}}_{0}^{(j)}=\left\{ \ast \right\}$ and 
\begin{align*}
    {P_{Q,U,\lambda }\left( Z_{1}^{\left( j\right) }\leq z_{1}^{(j)}|S=j\right) %
\coloneqq U^{\left( j\right) }\left( Z_{1}^{\left( j\right) }\leq
z_{1}^{(j)}\right)}
\end{align*} 
if $\overline{\mathcal{Z}}_{0}^{(j)}=\emptyset$ for all $z_{1}^{(j)}\in \mathbb{R}^{\dim \left( Z_{1}^{\left(
j\right) }\right) }$.

\item For each $k=2,...,K^{\left( j\right) }$ and all $z_{k}^{(j)}\in \mathbb{R}^{\dim \left( Z_{k}^{\left( j\right)
}\right)}$
\begin{align*}
    {P_{Q,U,\lambda }\left( Z_{k}^{\left( j\right) }\leq z_{k}^{(j)}|\overline{Z}%
_{k-1}^{\left( j\right) }=\overline{z}_{k-1}^{\left( j\right) },S=j\right) %
\coloneqq Q\left( Z_{k}^{\left( j\right) }\leq z_{k}^{(j)}|\overline{Z}%
_{k-1}^{\left( j\right) }=\overline{z}_{k-1}^{\left( j\right) }\right)}    
\end{align*}
if $\overline{z}_{k-1}^{\left( j\right) }\in \overline{%
\mathcal{Z}}_{k-1}^{\left( j\right) }$ and 
\begin{align*}
    {P_{Q,U,\lambda }\left( Z_{k}^{\left( j\right) }\leq z_{k}^{(j)}|\overline{Z}%
_{k-1}^{\left( j\right) }=\overline{z}_{k-1}^{\left( j\right) },S=j\right) %
\coloneqq U^{\left( j\right) }\left( Z_{k}^{\left( j\right) }\leq z_k^{(j)}|\overline{Z}%
_{k-1}^{\left( j\right) }=\overline{z}_{k-1}^{\left( j\right) }\right)}
\end{align*}
if ${\overline{z}_{k-1}^{\left( j\right) }\in \mathsf{Supp}%
\left[ \overline{Z}_{k-1}^{\left( j\right) };P_{0}\left( \cdot|S=j\right) \right]
\backslash \overline{\mathcal{Z}}_{k-1}^{\left( j\right) }}$.
\item $P_{Q,U,\lambda }\left( S=j\right) \coloneqq\lambda \left( S=j\right)
$.
\end{enumerate}

\begin{lemma}
\label{lemma:alt-model-characterization} Let $\left( \mathcal{Q},\mathcal{P},%
\mathcal{C}\right) $ be a fused-data model. Then $\mathcal{P}=\left\{ P_{Q,U,\lambda }:Q\in \mathcal{Q},U\in \mathcal{U}%
,\lambda \in \Lambda \right\}$.
\end{lemma}

While under \cref{as:identification} alignment of $Q$ and $P$ relative to $\mc{C}$ is sufficient to identify $\psi(Q)$ with $\varphi(P)$, it does not ensure that $\varphi$ is pathwise differentiable at $P$. As discussed in \cref{app:semiparametric-theory}, pathwise differentiability is a necessary condition for the existence of regular and asymptotically linear estimators of $\varphi(P)$ based on a random sample from $P$. The following definition encodes regularity conditions we will assume throughout to derive necessary and sufficient conditions for the pathwise differentiability of $\varphi$ at $P$. 
\begin{definition}
\label{def:strongly-align}
Let $\left( \mathcal{Q},\mathcal{P},\mathcal{C}\right) $ and $\left( Q,U,P\right) $ be in $\mathcal{Q}\times \mathcal{U\times P}$.
\begin{enumerate}
    \item $\left( Q,P\right) $ is \emph{strongly aligned} with respect $%
\mathcal{C}$ iff $P \alignswith Q$ and there exists $\delta >0$ such that 
$\delta ^{-1}\leq \frac{dP(\cdot|S=j)}{dQ}(\overline{Z}%
_{k-1}^{(j)})\leq \delta \text{ \ a.e.- }Q\text{ on }\overline{\mathcal{Z}}%
_{k-1}^{(j)}$ for all 
$j\in \left[ J\right] ,k\in \left\{ 2,...,K^{\left( j\right) }\right\} $.

\item $\left( U,P\right) $ is \emph{strongly aligned} with respect $%
\mathcal{C}$ iff $P=P_{Q^{\prime },U,\lambda }$ for some  $Q^{\prime }$ in $\mathcal{Q}$ and $\lambda \in
\Lambda $ and there exists $%
\epsilon >0$  such that $\epsilon ^{-1}\leq \frac{dP(\cdot|S=j)}{dU^{\left(
j\right) }}(\overline{Z}_{k-1}^{(j)})\leq \epsilon$ a.e.-$
U^{\left( j\right) }$ on $\mathsf{Supp}\left[ \overline{Z}%
_{k-1}^{(j)},P\left( \cdot|S=j\right) \right] \setminus\overline{\mathcal{Z}}%
_{k-1}^{(j)}$ for all $j\in \left[ J\right]$. 

\item $\left(Q,U, P\right) $ is a \emph{strongly aligned triplet} with
respect to $\mathcal{C}$ iff $\left( Q,P\right) $ and $\left( U,P\right) $
are strongly aligned with respect to $\mathcal{C}$.
\end{enumerate}
\end{definition}

Suppose that $\left( Q,U, P\right) $ is a strongly aligned triplet. In the proof of \cref{lemma:score-operator} we show 
that regular submodels of $\mathcal{Q},\mathcal{U},$ and $\Lambda 
$ respectively indexed by $t \in (-\varepsilon, \varepsilon)$ such that $Q_{t=0}=Q,U_{t=0}=U$ and $%
\lambda _{t=0}=\lambda $, with scores at $t=0$ denoted by $h^{\left(
Q\right) },h^{\left( U\right) }\coloneqq\left( h^{\left( U^{\left( 1\right)
}\right) },...,h^{\left( U^{\left( J\right) }\right) }\right) $ and $%
h^{\left( \lambda \right) },$ induce a  submodel $%
\left\{P_{t}\coloneqq P_{Q_{t},U_{t},\lambda _{t}}:t \in (-\varepsilon, \varepsilon)\right\}$ of model $\mathcal{%
P}$ differentiable in quadratic mean with score $g$ at $t=0$. Additionally, two 
submodels for $\left( Q,U,\lambda \right) $ with the same scores $h^{\left(
Q\right) },h^{\left( U\right) }$ and $h^{\left( \lambda \right) }$ induce
the same score $g.$ We can thus define the map $A_{Q,U,\lambda }:%
\mathcal{H}^{\circ }$ $\mathcal{\rightarrow }$ $L_{0}^{2}\left( P\right) $
such that 
\begin{equation*}
A_{Q,U,\lambda }h=g
\end{equation*}%
where $h\coloneqq$ $\left( h^{\left( Q\right) },h^{\left( U\right)
},h^{\left( \lambda \right) }\right) $ and $\mathcal{H}^{\circ }\coloneqq%
\mathcal{T}^{\circ }(Q,\mathcal{Q})\times \prod_{j\in \lbrack
J]}L_{0}^{2}\left( U^{\left( j\right) }\right) \times L_{0}^{2}\left(
\lambda \right) $ is the Cartesian product of the maximal tangent sets of
model $\mathcal{Q}$ and of the unrestricted models $\mathcal{U}$ and $%
\Lambda .$ The range of this map is the maximal tangent set $\mathcal{T}%
^{\circ }\left( P,\mathcal{P}\right) $ for model $\mathcal{P}$ at $P$ and
the $L^2(P)$-closed linear span of $\mathcal{T}^{\circ }\left( P,\mathcal{P}%
\right) $ is the maximal tangent space $\mathcal{T}(P,\mathcal{P}).$

The Cartesian product $\mathcal{H}\coloneqq\mathcal{T}(Q,\mathcal{Q})\times
\prod_{j\in \lbrack J]}L_{0}^{2}\left( U^{\left( j\right) }\right) \times
L_{0}^{2}\left( \lambda \right) ,$ where $\mathcal{T}(Q,\mathcal{Q})$ is the 
$L^{2}\left( Q\right) -$closed linear span of $\mathcal{T}^{\circ }(Q,%
\mathcal{Q}),$ endowed with the inner product \\
${\left\langle \left( h_{1}^{\left( Q\right) },h_{1}^{\left( U\right) },h_{1}^{\left( \lambda \right) }\right) ,\left( h_{2}^{\left( Q\right)
},h_{2}^{\left(U\right) },h_{2}^{\lambda }\right) \right\rangle _{\mathcal{H}}} 
\coloneqq\left\langle h_{1}^{\left( Q\right) },h_{2}^{\left( Q\right)
}\right\rangle _{L_{0}^{2}\left( Q\right) }+\sum_{j=1}^{J}\left\langle
h_{1}^{\left( U^{\left( j\right) }\right) },h_{2}^{\left( U^{\left( j\right)
}\right) }\right\rangle _{L_{0}^{2}\left( U^{\left( j\right) }\right)
}+\left\langle h_{1}^{\left( Q\right) },h_{2}^{\left( Q\right)
}\right\rangle _{L_{0}^{2}\left( \lambda \right) }$%
is a Hilbert space. Defining $\nu :\mathcal{Q\times U}\times \Lambda
\rightarrow \mathbb{R}$ as $\nu \left( Q,U,\lambda \right) \coloneqq\psi \left( Q\right)$ 
we then conclude that when $\psi $ is pathwise differentiable at $Q$ and $%
\psi _{Q}^{1}$ denotes any one of its influence functions, then for any regular
parametric submodels of $\mathcal{Q},\mathcal{U},$ and $\Lambda $
respectively indexed by a scalar $t$ such that $Q_{t=0}=Q,U_{t=0}=U$ and $%
\lambda _{t=0}=\lambda $, with scores at $t=0$ denoted by $h^{\left(
Q\right) },h^{\left( U\right) }\coloneqq\left( h^{\left( U^{\left( 1\right)
}\right) },...,h^{\left( U^{\left( J\right) }\right) }\right) $ and $%
h^{\left( \lambda \right) }$ it holds that 
\begin{align}
\left. \frac{d}{dt}\psi \left( Q_{t}\right) \right\vert _{t=0} =\left. 
\frac{d}{dt}\nu \left( Q_{t},U_{t},\lambda _{t}\right) \right\vert _{t=0}
\label{bp1} 
=\left\langle \left( \psi _{Q}^{1},\mathbf{0}_{J},0\right) ,h\right\rangle
_{\mathcal{H}} 
\end{align}%
where $h\coloneqq\left( h^{\left( Q\right) },h^{\left( U\right) },h^{\left(
\lambda \right) }\right) $ and $\mathbf{0}_{J}$ is the vector of dimension $%
J $ with all its entries equal to $0$. On the other hand, $\varphi $ is
pathwise differentiable at $P$ with respect to the tangent space $\mc{T}(P, \mc{P})$ if and only if there exists $\varphi
_{P}^{1} $ in $L_{0}^{2}\left( P\right) $ such that for all such submodels, 
\begin{eqnarray}
\left. \frac{d}{dt}\varphi \left( P_{Q_{t},U_{t},\lambda _{t}}\right)
\right\vert _{t=0} =\left. \frac{d}{dt}\nu \left( Q_{t},U_{t},\lambda
_{t}\right) \right\vert _{t=0}  \label{bp2} =\left\langle \varphi _{P}^{1},A_{Q,U,\lambda }h\right\rangle
_{L^{2}\left( P\right) }.
\end{eqnarray}

In \cref{lemma:score-operator} below we establish that when $\left(
Q,U, P\right) $ is a strongly aligned triplet with respect to $\mathcal{C}$%
, the map $A_{Q,U,\lambda }$ is bounded, meaning that sup$_{h:\left\Vert
h\right\Vert _{\mathcal{H}}=1}\left\Vert A_{Q,U,\lambda }h\right\Vert
_{L^{2}\left( P\right) }<\infty $. Because $A_{Q,U,\lambda }$ is both linear
and bounded, it can be extended to a linear bounded operator on $\mathcal{H}$%
, the closed linear span of $\mathcal{H}^{0}$. For this extension, which we will
continue to denote as $A_{Q,U,\lambda }$, there exists an adjoint operator $%
A_{Q,U,\lambda }^{\ast }:L_{0}^{2}\left( P\right) \rightarrow \mathcal{H}$
that satisfies 
$\left\langle \varphi _{P}^{1},A_{Q,U,\lambda }h\right\rangle _{L^{2}\left(
P\right) }=\left\langle A_{Q,U,\lambda }^{\ast }\varphi
_{P}^{1},h\right\rangle _{_{\mathcal{H}}}$ for all $h\in \mathcal{H}$. 
Therefore, equating the rightmost hand sides of $\left( \ref{bp1}\right) $
and $\left( \ref{bp2}\right) ,$ we conclude that if $\psi $ is pathwise
differentiable at $Q$, then $\varphi $ is pathwise differentiable at $P$ if
and only if there exists $\varphi _{P}^{1}\in L_{0}^{2}\left( P\right) $
solving 
$
\left\langle A_{Q,U,\lambda }^{\ast }\varphi _{P}^{1}-\left( \psi _{Q}^{1},%
\mathbf{0}_{J},0\right) ,h\right\rangle _{\mathcal{H}}=0$ for all $h\in 
\mathcal{H}$ in which case, $\varphi _{P}^{1}$ is a influence function of $\varphi $ at $P.$ This equality holds if and only if 
\begin{equation}
A_{Q,U,\lambda }^{\ast }\varphi _{P}^{1}=\left( \psi _{Q,eff}^{1},\mathbf{0}%
_{J},0\right).  \label{ap4}
\end{equation}

We conclude that when $\left(Q,U, P\right) $ is a strongly aligned triplet
with respect to $\mathcal{C}$, the existence of a solution to equation $%
\left( \ref{ap4}\right) ,$ equivalently the condition that $\left( \psi
_{Q,eff}^{1},\mathbf{0}_{J},0\right) $ lies in the range of $A_{Q,U,\lambda
}^{\ast },$ is the necessary and sufficient condition for pathwise
differentiability of $\varphi $ with respect to $\mathcal{T}^{\circ }(P,%
\mathcal{P}).$ Furthermore, the collection of all solutions to $\left( \ref%
{ap4}\right) $ is precisely the class of observed data influence functions.
The efficient influence function $\varphi _{P,eff}^{1}$ of $\varphi $ is the
unique element of $\mathcal{T}(P,\mathcal{P})$ that satisfies the preceding
equation. Equation $\left( \ref{ap4}\right) $ then gives the fundamental
equation that determines all observed data influence functions $\varphi
_{P}^{1}$.

Some important subtle points warrant distinction at this juncture. Suppose
that $\left( \mathcal{Q},\mathcal{P},\mathcal{C},\psi ,\varphi \right) $ is
a fused-data framework and we are interested in inference about $\psi \left(
Q_{0}\right) $ at a particular $Q_{0}$ based on $n$ i.i.d. draws from $%
P_{0}. $ There always exist $\left( Q,U\right) $ and $\left(
Q_{0},U_{0}\right) \in \mathcal{Q\times U}$ such that $\left(
Q_{0},U_{0}\right) \not=\left( Q,U\right) $ and with $P_{Q_{0},U_{0},\lambda
_{0}}=P_{Q,U,,\lambda _{0}}=P_{0}$ for some $\lambda _{0}\in \Lambda .$ In
particular, $Q$ may be chosen to be any element of the equivalence class $\xi
\left( Q_{0};\mathcal{C}\right) $ which will be comprised of more than one
element when the alignments in $\mathcal{C}$ do not determine $Q_{0}.$ For
such $\left( Q,U\right) $ and $\left( Q_{0},U_{0}\right) $ it may happen
that both $\left( Q,U,P_{0}\right) $ and $\left( Q_{0},U_{0},P_{0}\right) $
are strongly aligned with respect to $\mathcal{C}$. Therefore, the score
operators $A_{Q_{0},U_{0},\lambda _{0}}$ and $A_{Q,U,\lambda _{0}}$ are
bounded. Thus, following the preceding argument with $\left(
Q_{0},U_{0},\lambda _{0}\right) $ or $\left( Q,U,\lambda _{0}\right) $
instead of $\left( Q,U,\lambda \right) $ we arrive at the conclusion that
the set of observed data influence functions $\varphi
_{P_{0}}^{1}$ is the set of solutions to $\left( \ref{ap4}\right) $ with $%
\left( Q,U,\lambda \right) $ replaced by either $\left( Q_{0},U_{0},\lambda
_{0}\right) $ or $\left( Q,U,\lambda _{0}\right) .$ This is true even though
the spaces $\mathcal{H}_{0}=\mathcal{T}(Q_{0},\mathcal{Q})\times \prod_{j\in
\lbrack J]}L_{0}^{2}\left( U_{0}^{\left( j\right) }\right) \times
L_{0}^{2}\left( \lambda_0 \right) $ and $\mathcal{H}=\mathcal{T}(Q,\mathcal{Q}%
)\times \prod_{j\in \lbrack J]}L_{0}^{2}\left( U^{\left( j\right) }\right)
\times L_{0}^{2}\left( \lambda_0 \right) $ are different and the score
operators $A_{Q_{0},U_{0},\lambda _{0}}$ and $A_{Q,U,\lambda _{0}}$ are also different.

If $\left( Q_{0},U_{0},P_{0}\right) $ is not strongly aligned with respect
to $\mathcal{C},$ we cannot ensure that the score operator $%
A_{Q_{0},U_{0},\lambda _{0}},$ assuming it exists, is bounded. Thus, the
adjoint of $A_{Q_{0},U_{0},\lambda _{0}}$ may not be everywhere defined and therefore the
equation $\left( \ref{ap4}\right) $ is not available to characterize the set
of observed data influence functions. However, we can always take $U_{0}$
such that $U_{0}^{\left( j\right) }=P_0\left( \cdot|S=j\right) $ and for such $%
U_{0},$ it holds that the pair $\left( U_{0},P_{0}\right) $ is strongly
aligned. For such choice, failure of the strong alignment of the triplet $%
\left( Q_{0},U_{0},P_{0}\right) $ can only be due to failure of strong
alignment of the pair $\left( Q_{0},P_{0}\right) .$ Nevertheless, if $\xi
\left( Q_{0};\mathcal{C}\right) $ has more than one element, it may be
possible to find a $Q\in $ $\xi \left( Q_{0};\mathcal{C}\right) $ such that $%
\left( Q,U_{0},P_{0}\right) $ is strongly aligned and we can therefore use
that $Q\,\ $to derive the set of all observed data influence functions as
the set of solutions of equation $\eqref{ap4}$. See the continuation of \cref{example:transporting2} scenario (i)
in \cref{subsec:examples IF} for an illustration of this point.

For a strongly aligned triplet $\left( Q,U,P\right) ,$ the closure of the
range of the score operator $A_{Q,U,\lambda }:\mathcal{H}\rightarrow
L_{0}^{2}\left( P\right) $ is the tangent space $\mathcal{T}\left( P,%
\mathcal{P}\right) $. Thus, if the range of the score operator is not dense
in $L_{0}^{2}(P)$, then $\mathcal{T}\left( P,\mathcal{P}\right) $ is
strictly included in $L_{0}^{2}\left( P\right) $, implying the existence of
infinitely many observed data influence functions of which, the unique
influence function $\varphi _{P,eff}^{1}$ lying in $\mathcal{T}\left( P,%
\mathcal{P}\right) $ has the smallest variance. In the context of a
non-parametric model $\mathcal{Q}$, the score operator's range can only fail
to be dense in $L_{0}^{2}(P)$ due to alignment assumptions that impose
constraints on the distribution $P$. This implies that relaxing some
alignment assumptions might still allow for parameter identification, and
illustrates the usual trade-off between bias and variance: reducing the
number of alignment assumptions may decrease the risk of bias if some
assumptions prove to be invalid, yet it also reduces the efficiency of
parameter estimation. \cref{example:transporting2} scenario (iii) in %
\cref{subsec:examples IF} illustrates these concepts.

We will now establish that when the triplet $\left( Q,U,P\right) $ is
strongly aligned, the score operator $A_{Q,U,\lambda }$ exists, is
bounded, and is linear. We define the subspace of $
L^{2}(W;Q)$ 
\begin{equation*}
\mathcal{D}_{k}^{(j)}\left( Q\right) \coloneqq\left\{ I(\overline{z}%
_{k-1}^{(j)}\in \overline{\mc{Z}}_{k-1}^{(j)})\{d(\overline{z}_{k}^{(j)})-E_{Q}%
[ d(\overline{Z}_{k}^{(j)})|\overline{z}%
_{k-1}^{(j)}] \}:d\in L^{2}(\overline{Z}_{k}^{(j)};Q)\right\} ,
\end{equation*}%
and the subspaces of $L^{2}(Z^{(j)};P_{Q,U,\lambda }\left(\cdot |S=j\right) )$
for $k\in [ K^{\left( j\right) }] ,j\in [ J] ,$ 
\begin{align*}
\mathcal{R}_{k}^{(j)}&\left( P_{Q,U,\lambda }\right) \coloneqq\\
\big\{& I(%
\overline{z}_{k-1}^{(j)}\not\in \overline{\mc{Z}}_{k-1}^{(j)})\{r(\overline{z}%
_{k}^{(j)})-E_{U^{(j)}}[r(\overline{Z}_{k}^{(j)})|\overline{z}%
_{k-1}^{(j)},S=j]\}:r\in L^{2}(\overline{Z}_{k}^{(j)};U^{(j)} )\big\}
\end{align*}%
where hereafter $\overline{z}_{k-1}^{(j)}\not\in \overline{\mathcal{Z}}%
_{k-1}^{(j)}$ is a shortcut for $\overline{z}_{k-1}^{\left(
j\right) }\in{\textsf{Supp}\left[ \overline{Z}_{k-1}^{\left( j\right)
};P_{Q,U,\lambda }\left( \cdot |S=j\right) \right] \setminus \overline{\mathcal{Z}}_{k-1}^{\left( j\right) }}$.

We express the score operator with orthogonal projections of
ideal data scores into the spaces $\mathcal{D}_{k}^{(j)}( Q) $
and $\mathcal{R}_{k}^{(j)}( P_{Q,U,\lambda })$. As such, note that for any $h\in L^{2}(W;Q),$ its projection into $\mathcal{D}_{k}^{(j)}( Q)$ is given by 
\begin{align*}
    \Pi [ h\vert \mathcal{D}_{k}^{(j)}( Q) ]
( \overline{z}_{k}^{(j)}) =I(\overline{z}_{k-1}^{(j)}\in 
\overline{\mc{Z}}_{k-1}^{(j)})\{E_{Q}[ h(W)|\overline{%
z}_{k}^{(j)}] -E_{Q}[ h(W)|\overline{z}%
_{k-1}^{(j)}]\}
\end{align*}
and for any $u^{( j) }\in L^{2}(Z^{(j)};P_{Q,U,\lambda }(
|S=j) )$, its projection into $\mathcal{R}%
_{k}^{(j)}( P_{Q,U,\lambda })$ is given by 
\begin{align*}
    \Pi [ u^{( j) }\vert \mathcal{R}%
_{k}^{(j)}( P_{Q,U,\lambda }) ] ( \overline{z}%
_{k}^{(j)}) =I(\overline{z}_{k-1}^{(j)}\not\in \overline{\mc{Z}}%
_{k-1}^{(j)})\{E_{U^{(j)}}[ u^{( j) }(Z^{( j) })|\overline{z}_{k}^{(j)}] -E_{U^{(j)}}[ u^{(
j) }(Z^{( j) })|\overline{z}_{k-1}^{(j)}]\}.
\end{align*}

The following lemma provides the precise expression for $A_{Q,U,\lambda }$
and its adjoint $A_{Q,U,\lambda }^{\ast }$. These expressions invoke the
decomposition of $g\in L_{0}^{2}(P_{Q,U,\lambda })$ as 
\begin{equation}
g(o)=\sum_{j\in \lbrack J]}I(s=j)\sum_{k\in \left[ K^{\left( j\right) }%
\right] }\left\{ m_{k}^{(j)}(\overline{z}_{k}^{(j)})+n_{k}^{(j)}(\overline{z}%
_{k}^{(j)})\right\} +\gamma \left( s\right)  \label{decompose}
\end{equation}%
for $m_{k}^{(j)}\in \mathcal{D}_{k}^{(j)}\left( Q\right) ,$ $n_{k}^{(j)}\in 
\mathcal{R}_{k}^{(j)}\left( P_{Q,U,\lambda }\right) ,k\in \left[ K^{\left(
j\right) }\right] ,$ $j\in \lbrack J]$ and $\gamma \in L_{0}^{2}(\lambda )$. Whenever $(Q, U, P)$ are strongly aligned 
with respect to $\mathcal{C}$ this decomposition exists and is unique as the spaces $\{I(S=j)m_k^{(j)}(\overline{Z}_k^{(j)}):  m_k^{(j)} \in \mc{D}_k^{(j)}(Q)\}$, $\{I(S=j)n_k^{(j)}(\overline{Z}_k^{(j)}):  n_k^{(j)} \in \mc{R}_k^{(j)}(P_{Q, U, \lambda})\}$, $L^2_0(\lambda)$ are mutually orthogonal subspaces of $L^2(P_{Q, U, \lambda})$ for $j \in [J], k\in[K^{(j)}]$. Their direct sum is equal to the space $L^2_0(P_{Q, U, \lambda})$. The decomposition can be computed by taking the $L^2(P_{Q, U, \lambda})$ orthogonal projection of $g$ onto each subspace.

\begin{lemma}
\label{lemma:score-operator} Let $\left( \mathcal{Q},\mathcal{P},\mathcal{C}%
\right) $ be a fused-data model.
Let $\left( Q,U,P\right) $ be strongly aligned with respect to $\mathcal{C}.$
Let $\lambda(S=j) \coloneqq P(S=j)$. Then, the score operator $%
A_{Q,U,\lambda }:\mathcal{H}\rightarrow L_{0}^{2}(P)$ exists, is bounded and
linear, and for any $h\coloneqq\left(
h^{(Q)},h^{(U^{(1)})},...,h^{(U^{(J)})},h^{(\lambda )}\right) \in \mathcal{H}$, $A_{Q,U,\lambda }h\coloneqq A_{Q}h^{(Q)}+\sum_{j\in \lbrack
J]}A_{U^{(j)}}h^{(U^{(j)})}+A_{\lambda }h^{(\lambda )}$
where $(A_{\lambda }h^{(\lambda )})(o)\coloneqq h^{(\lambda )}(s)$ and
\begin{align*}
(A_{Q}h^{(Q)})(o)\coloneqq& \sum_{j\in \lbrack J]}I(s=j)\sum_{k\in \left[
K^{\left( j\right) }\right] }\Pi \left[ \left. h^{\left( Q\right)
}\right\vert \mathcal{D}_{k}^{(j)}\left( Q\right) \right] \left( \overline{z}%
_{k}^{(j)}\right) \\
(A_{U^{(j)}}h^{(U^{(j)})})(o)\coloneqq& I(s=j)\sum_{k\in \left[ K^{\left(
j\right) }\right] }\Pi \left[ \left. h^{(U^{(j)})}\right\vert \mathcal{R}%
_{k}^{(j)}\left( P_{Q,U,\lambda }\right) \right] \left( \overline{z}%
_{k}^{(j)}\right).
\end{align*}
The adjoint $A_{Q,U,\lambda }^{\ast}\coloneqq(A_{Q}^{\ast
},A_{U^{(1)}}^{\ast },\dots ,A_{U^{(J)}}^{\ast },A_{\lambda }^{\ast
}):L_{0}^{2}(P)\rightarrow \mathcal{H}$ satisfies for any $g\in L_{0}^{2}(P)$
\begin{align*}
(A_{Q}^{\ast }g)(w) =&\sum_{j\in \lbrack J]}\sum_{k\in \left[ K^{\left(
j\right) }\right] }\Pi \left[ \left. \frac{dP(\cdot|S=j)}{dQ}(\overline{Z}%
_{k-1}^{(j)})\lambda (S=j)m_k^{(j)}(\overline{Z}%
_k^{(j)})\right\vert \mathcal{T}(Q,%
\mathcal{Q})\right] \left( w\right)  \\
(A_{U^{(j)}}^{\ast }g)(z^{(j)}) =&\sum_{k\in \left[ K^{\left( j\right) }%
\right] }\frac{dP(\cdot|S=j)}{dU^{\left( j\right)}}(\overline{z}_{k-1}^{(j)})\lambda
(S=j)n_k^{(j)}(\overline{z}_k^{(j)})  \notag
\end{align*}%
and $(A_{\lambda }^{\ast }g)(s) =\gamma \left( s\right)$ with  $m_{k}^{(j)},n_{k}^{(j)}$ and $\gamma $ being the components of the
decomposition $\left( \ref{decompose}\right) $ of $g$ and  $\Pi \left[ \left.
\frac{dP(\cdot|S=j)}{dQ}(\overline{Z}%
_{k-1}^{(j)})\lambda (S=j)m_k^{(j)}(\overline{Z}%
_k^{(j)})\right\vert \mathcal{T}(Q,\mathcal{Q})\right]$ the $%
L^{2}\left( Q\right) $-projection of the function $\overline{z}%
_k^{(j)} \mapsto \frac{dP(\cdot|S=j)}{dQ}(\overline{z}%
_{k-1}^{(j)})\lambda (S=j)m_k^{(j)}(\overline{z}%
_k^{(j)})$ into $\mathcal{T}(Q,\mathcal{Q}).$
\end{lemma}

Fused-data models where the mere fact that the alignments $\mc{C}$ restrict the observed data model, i.e. $\mc{T}(P, \mc{P}^{ext}) \subsetneq L^2_0(P)$, are of particular interest as the alignment assumptions could potentially be weakened without impacting identification. The following lemma provides a necessary and sufficient condition for $\mc{T}(P, \mc{P}^{ext}) = L^2_0(P)$.
\begin{lemma}
    \label{lemma:extended-model-nonparametric}
    Let $(\mc{Q}, \mc{P}, \mc{C})$ be a fused-data model. Suppose that $(Q, P)$ is strongly aligned with respect to $\mc{C}$. Then, $\mc{T}(P, \mc{P}^{ext}) =  L^2_0(P)$ if and only if the spaces $\mc{D}_k^{(j)}$, $k\in[K^{(j)}], j\in[J]$ are linearly independent in the sense that  $0 = \sum_{j\in[J]}\sum_{k\in[K^{(j)}]}m_k^{(j)}$ for $m_k^{(j)} \in \mc{D}_k^{(j)}(Q)$ if and only if $m_k^{(j)} = 0$ a.e.-$Q$ for all $k\in[K^{(j)}], j\in[J]$. 
\end{lemma}
When $J=2$, the condition reduces to checking if $\bigcap_{j=1}^2\bigoplus_{k\in[K^{(j)}]}\mathcal{D}_{k}^{(j)}\left( Q\right) = \{0\}$.
\subsection{Characterizing observed data influence functions}

\label{subsec:characterize-ifs}

In the next lemma, we
invoke \cref{lemma:score-operator} and equation 
\eqref{ap4}
to derive two equivalent necessary and sufficient conditions for the
pathwise differentiability of $\varphi $.

\begin{lemma}
\label{lemma:pathwise-differentiabilty} Let $\left( \mathcal{Q},\mathcal{P},%
\mathcal{C}\text{\thinspace },\psi ,\varphi \right) $ be a fused-data
framework. Let $P\in \mathcal{P}
$. Suppose there exists $Q$ in $\mc{Q}$ such that $\left(
Q,P\right) $ is strongly aligned with respect to $\mathcal{C}$ and $\psi $
is pathwise differentiable at $Q$ in model $\mathcal{Q}$. Then each of the following assertions is equivalent to $\varphi $ being pathwise differentiable at $P$ in model $\mathcal{P}$:

\begin{enumerate}
\item \label{item:pathwise-diff-2} There exists $m_{k}^{(j)}\in \mathcal{D}%
_k^{(j)}(Q),k\in[K^{(j)}], j\in[J]$ such that 
\begin{equation*}
\psi _{Q,eff}^{1}(w)=\sum_{j\in \lbrack J]}\sum_{k\in \lbrack K^{(j)}]}\Pi 
\left[ \frac{dP(\cdot|S=j)}{dQ}(\overline{Z}%
_{k-1}^{(j)})\lambda (S=j)m_k^{(j)}(\overline{Z}%
_k^{(j)})\Big|\mathcal{T}(Q,\mathcal{Q})\right] (w)
\end{equation*}

where $\psi _{Q,eff}^{1}$ is the efficient influence function of $\psi $ at $%
Q$ in model $\mathcal{Q}.$

\item \label{item:pathwise-diff-3} There exists an influence function $\psi
_{Q}^{1}$ of $\psi $ at $Q$ in model $\mathcal{Q}$ and
$m_{k}^{(j)}\in \mathcal{D}%
_k^{(j)}(Q),k\in[K^{(j)}], j\in[J]$ such that 
\begin{equation}
\psi _{Q}^{1}=\sum_{j\in \lbrack J]}\sum_{k\in \lbrack K^{(j)}]}m_{k}^{(j)}.
\label{eq:if-decomposition}
\end{equation}
\end{enumerate}
\end{lemma}

Part \ref{item:pathwise-diff-2} of the preceding Lemma is equivalent to the assertion that $\psi
_{Q,eff}^{1}$ is in the range of the operator $A_{Q}^{\ast }$ defined in \cref{lemma:score-operator} and thus is equivalent to the assertion that $\left( \psi _{Q,eff}^{1},\mathbf{0}%
_{J},0\right) $ is in the range of $A_{Q,U,\lambda }^{\ast },$ the adjoint
of the score operator, i.e. the necessary and sufficient condition for
pathwise differentiability of $\varphi $ with respect to $\mathcal{T}^{\circ
}(P,\mathcal{P})$ discussed in \cref{subsec:score-operator-theory}.

According to the preceding Lemma, pathwise differentiability of $\varphi \,$%
\ at $P$ can be confirmed by exhibiting the decomposition $\left( \ref%
{eq:if-decomposition}\right) $ for some ideal data influence function $\psi _{Q}^{1}$,
for any $Q$ such that the pair $\left( Q,P\right) $ strongly aligns. Note
that it may be the case that the target ideal data law is a $\widetilde{Q}$
that aligns with $P$ but does not strongly align with $P.$ The theorem
establishes that strong alignment of the target $\widetilde{Q}$ is not
needed to derive the observed data influence functions. It suffices to find
a $Q$ in the equivalence class $\xi \left( \widetilde{Q};\mathcal{C}\right) $
such that $\left( Q,P\right) $ strongly aligns and $\psi $ is pathwise
differentiable at that $Q.$

A natural question is whether the mere fact that the pathwise differentiable
functional $\psi :\mathcal{Q\rightarrow }\mathbb{R}$ depends on $Q$ solely
through the aligned conditionals, implies that $\varphi $ is pathwise
differentiable at a $P$ that strongly aligns with $Q$. In \cref{sec:appendix-score-operator}, we show that the answer is in general negative.
However, for particular fused-data frameworks it can be shown that the
decomposition $\left( \ref{eq:if-decomposition}\right) $ holds whenever $%
\psi $ is pathwise differentiable. This is the case in fused-data frameworks 
$\left( \mathcal{Q},\mathcal{P},\mathcal{C}\text{\thinspace },\psi ,\varphi
\right) $ where the fused-data models $\left( \mathcal{Q},\mathcal{P},%
\mathcal{C}\text{\thinspace }\right) $ are those in \cref%
{example:disease-prevalence} and \cref{example:transporting} scenarios
(ii) and (iii), see \cref{app:examples-extras}. It is also the case in the fused-data frameworks
considered by \cite{li_efficient_2023}. These frameworks include the framework in \cref{example:transporting} scenario (i). 

The following theorem characterizes all observed data influence
functions in terms of ideal data influence functions when the observed data
parameter is pathwise differentiable. 

\begin{theorem}
\label{Theorem:influence functions} Let $\left( \mathcal{Q},\mathcal{P},%
\mathcal{C}\text{\thinspace },\psi ,\varphi \right) $ be a fused-data
framework. Let $P\in \mathcal{P}
$ satisfy that $\varphi $ is pathwise differentiable at $P$ in model $%
\mathcal{P}$. Suppose there exists $Q$ in $\mc{Q}$ such that $%
\left( Q,P\right) $ is strongly aligned with respect to $\mathcal{C}$ and $%
\psi $ is pathwise differentiable at $Q$ in model $\mathcal{Q}$. Then,

\begin{enumerate}
\item \label{item:influence functions1} $\varphi _{P}^{1}$ is an influence
function of $\varphi $ at $P$ in model $\mathcal{P}$ iff $\varphi
_{P}^{1}$ can be expressed as 
\begin{equation}
\varphi _{P}^{1}\left( o\right) =\sum_{j=1}^{J}\frac{I\left( s=j\right) }{%
P(S=j)}\sum_{k\in \left[ K^{\left( j\right) }\right] }\frac{dQ}{dP(\cdot|S=j)}(\overline{z}%
_{k-1}^{(j)})m_{k}^{(j)}(\overline{z}%
_{k}^{(j)})  \label{newIF}
\end{equation}%
for some $m_{k}^{(j)}\in \mathcal{D}_{k}^{(j)}\left( Q\right),k\in[K^{(j)}], j\in[J]$ such that there exists
an influence function $\psi _{Q}^{1}$ of $\psi $ at $Q$ in model $\mathcal{Q}$ that satisfies the decomposition \eqref{eq:if-decomposition}.

\item\label{item:influence functions2} If $\mathcal{T}\left( P,\mathcal{P}%
^{ext}\right) =L_{0}^{2}\left( P\right) $ then for every ideal data
influence function $\psi _{Q}^{1}$ there exists at most one collection of
functions $\left\{ m_{k}^{(j)}\in \mathcal{D}_{k}^{(j)}\left( Q\right) :k\in %
\left[ K^{\left( j\right) }\right] ,j\in \left[ J\right] \right\} $
satisfying the decomposition \eqref{eq:if-decomposition}. If $\mathcal{T}\left( P,\mathcal{P}%
^{ext}\right) \subsetneq L_{0}^{2}\left( P\right)$ then there exists 0 or infinitely many such collections for each $\psi^1_Q$. 
\end{enumerate}
\end{theorem}
Throughout we say that an
observed data influence function $\varphi _{P}^{1}$ corresponds to
the ideal data influence function $\psi _{Q}^{1}$ if there exists a
collection $\big\{ m_{k}^{(j)}\in \mathcal{D}_{k}^{(j)}\left( Q\right)
:k\in \left[ K^{\left( j\right) }\right] ,j\in \left[ J\right] \big\} $
such that $\varphi _{P}^{1}$ decomposes as $\left( \ref{newIF}\right) $ and $%
\psi _{Q}^{1}$ decomposes as $\left( \ref{eq:if-decomposition}\right) .$ Part \ref{item:influence functions1} of the preceding theorem establishes that every observed data
influence function $\varphi _{P}^{1}$ corresponds to an ideal data influence
function $\psi^1_Q$. However, it is not true that every $\psi _{Q}^{1}$ necessarily corresponds to a $\varphi^1_P$. There are fused-data frameworks where only
a strict subset of all ideal data influence functions $\psi _{Q}^{1}$ can
be decomposed as $\left( \ref{eq:if-decomposition}\right) $ for functions $%
m_{k}^{(j)}$ in $\mathcal{D}_{k}^{(j)}\left( Q\right) .$ Such subset is
included in the set of ideal data influence functions that are orthogonal to
the null space of the component $A_{Q}$ of the score operator $%
A_{Q,U,\lambda }$. For instance, suppose that in \cref%
{example:transporting} scenario (i), model $\mathcal{Q}$ for the ideal data
law $Q$ restricts the propensity score to a fixed and known $Q_{0}\left(
A=1|L\right) $. Then, it is well known (\cite{robins_estimation_1994}) that $\widetilde{\psi }_{Q}^{1}\left( l,a,y\right) =\psi _{Q}^{1}\left(
l,a,y\right) +d\left( l\right) \left\{ a-Q_{0}\left( A=1|l\right) \right\}$
is an ideal data influence function for any $d\in L^{2}\left( L;Q\right) ,$ where $\psi _{Q}^{1}$ is given in
equation $\left( \ref{AIPW}\right) $ in \cref{subsec:examples IF}. While $%
\psi _{Q}^{1}$ can be decomposed as in $\left( \ref{eq:if-decomposition}%
\right) $, $\widetilde{\psi }_{Q}^{1}\left( l,a,y\right) $ cannot when $%
d\not=0\,\ $because $d\left( l\right) \left\{ a-Q_{0}\left( A=1|l\right)
\right\} $ is orthogonal to the spaces $\mathcal{D}_{k}^{(j)}\left( Q\right) 
$ for all $k\in[K^{(j)}]$, $j\in\{1, 2\}$.

Part \ref{item:influence functions1} of \cref{Theorem:influence functions}
implies that if we can express an ideal data influence function as in %
\eqref{eq:if-decomposition}, we can find an observed data
influence function without needing to calculate new pathwise derivatives. The following proposition shows that finding the decomposition reduces to solving a single integral equation when there are only two sources.

\begin{proposition}
\label{prop:two-source-solution} Let $\left( \mathcal{Q},\mathcal{P},%
\mathcal{C}\text{\thinspace },\psi ,\varphi \right) $ be a fused-data
framework. Suppose there exists $Q \in \mc{Q}$ such that $\left(
Q,P\right) $ is strongly aligned with respect to $\mathcal{C}$ and $\psi $
is pathwise differentiable at $Q$ in model $\mathcal{Q}$. Suppose $J=2$%
. Then

\begin{enumerate}
\item \label{item:two-source-solution1}$\varphi $ is pathwise differentiable at $P$ in $\mc{P}$
iff there exists an influence function $\psi _{Q}^{1}$ for $\psi $ at $Q$ in $\mathcal{Q%
}$ such that the following equation has a solution $m^{(2)}$ in $\bigoplus_{k=1}^{K^{\left( 2\right) }}\mathcal{D}_{k}^{(2)}(Q)$,
\begin{align}
&m^{(2)}(Z^{(2)})-\sum_{k=1}^{K^{\left( 1\right) }}I(\overline{z}%
_{k-1}^{\left( 1\right) } \in \overline{\mathcal{Z}}_{k-1}^{\left( 1\right)
})\{E_{Q}[m^{(2)}(Z^{(2)})|\overline{z}_{k}^{\left( 1\right)
}]-E_{Q}[m^{(2)}(Z^{(2)})|\overline{z}_{k-1}^{\left( 1\right) }]\}
\label{eq:algorithm-op-eq-simple} \\
&=\psi _{Q}^{1}(w)-\sum_{k=1}^{K^{\left( 1\right) }}I(\overline{z}%
_{k-1}^{\left( 1\right) }\in \overline{\mathcal{Z}}_{k-1}^{\left( 1\right)
})\{E_{Q}[\psi _{Q}^{1}(W)|\overline{z}_{k}^{\left( 1\right)
}]-E_{Q}[\psi_{Q}^{1}(W)|\overline{z}_{k-1}^{\left( 1\right) }]\}.  \notag
\end{align}%

\item \label{item:two-source-solution2} Suppose %
\eqref{eq:algorithm-op-eq-simple} has a solution $m^{(2)}\in
\bigoplus_{k=1}^{K^{\left( 2\right) }}\mathcal{D}_{k}^{(2)}(Q)$ for some
influence function $\psi _{Q}^{1}$ for $\psi $ at $Q$ in $\mathcal{Q}$. Then, the following is an influence function for $\varphi $ at $P$ in $\mathcal{P}$,
\begin{align}
\varphi _{P}^{1}(o) =&\frac{I(s=1)}{P(S=1)}\sum_{k=1}^{K^{\left( 1\right) }}%
\frac{dQ}{dP(\cdot|S=1)}(\overline{z}_{k-1}^{(1)}) \Pi\left[\psi _{Q}^{1}-m^{(2)}|\mc{D}_k^{(1)}(Q)\right](\overline{z}_k^{(1)})
\label{phi-solution-IF} \\
&+\frac{I(s=2)}{P(S=2)}\sum_{k=1}^{K^{\left( 2\right) }}\frac{dQ}{dP(\cdot|S=2)}(\overline{z}%
_{k-1}^{(2)})\Pi\left[m^{(2)}|\mc{D}_k^{(2)}(Q)\right](\overline{z}_k^{(2)}). \notag
\end{align}%

\item \label{item:two-source-solution3} Suppose $\varphi $ is pathwise
differentiable at $P$ in $\mathcal{P}$, $\psi _{Q}^{1}$ is an ideal data
influence function such that \cref{eq:algorithm-op-eq-simple} has a
solution, and $\varphi _{P}^{1}$ is defined as in $\left( \ref%
{phi-solution-IF}\right)$. Then 
\begin{equation*}
\left\{ \varphi _{P}^{1}(o)\right\} +\left\{ \sum_{j=1}^{2}\left( -1\right)
^{j+1}\frac{I(s=j)}{P(S=j)}\sum_{k=1}^{K^{\left( j\right) }}\frac{dQ}{dP(\cdot|S=1)}(\overline{z}_{k-1}^{(j)})\Pi \left[ f|%
\mathcal{D}_{k}^{(j)}\left( Q\right) \right] \left( \overline{z}%
_{k}^{(j)}\right) :f\in \mc{F} \right\}
\end{equation*}%
is the set of all influence
functions for $\varphi $ at $P$ in $\mathcal{P}$ that correspond to $\psi _{Q}^{1}$ where $\mc{F} = \bigcap_{j=1}^2\bigoplus_{k\in[K^{(j)}]}\mathcal{D}_{k}^{(j)}\left( Q\right)$.
\end{enumerate}
\end{proposition}

In \cref{app:examples-extras} we show that a closed-form solution to \cref{eq:algorithm-op-eq-simple} exists in the fused-data frameworks of examples \ref{example:disease-prevalence}, \ref{example:tsiv-lsm} or \ref{example:transporting}. In fact, under the alignments of those examples, for any ideal data parameter satisfying \cref{as:identification} whose corresponding observed data parameter is pathwise differentiable, a closed-form solution of \cref{eq:algorithm-op-eq-simple} will exist. \cref{prop:example-ifs} in \cref{app:examples-extras} establishes the general expression for
the observed data influence functions in such settings. We note that the latter set in part \ref{item:two-source-solution3} of the above proposition is the orthogonal complement of the extended observed data tangent space $\mc{T}(P, \mc{P}^{ext})$. As shown in \cref{lemma:extended-model-nonparametric}, this space is non-empty whenever $\mc{F}$ is non-empty. In Section \cref{app:decomposing-ideal} we discuss an extension of \cref{prop:two-source-solution} to fused-data frameworks with more than two sources. In such cases, we must solve $J-1$ linear operator equations analogous to \eqref{eq:algorithm-op-eq-simple} sequentially. 

\subsection{The observed data efficient influence
function}

\label{subsec:compute-EIF}

\begin{theorem}
\label{theorem:eif} Let $\left( \mathcal{Q},\mathcal{P},\mathcal{C}\text{%
\thinspace },\psi ,\varphi \right) $ be a fused-data framework. Let $P\in \mathcal{P} $ and suppose that $%
\varphi $ is pathwise differentiable at $P$ in model $\mathcal{P}$. Suppose
there exists $Q$ in $\mc{Q}$ such that $\left( Q,P\right) $
is strongly aligned with respect to $\mathcal{C}$ and $\psi $ is pathwise
differentiable at $Q$ in model $\mathcal{Q}$. The following statements are equivalent, where all limits are in $L^2(Q)$ norm
\begin{enumerate}
\item \label{item:eif1} $\varphi _{P,eff}^{1}$ is the efficient influence
function of $\varphi $ at $P$ in model $\mathcal{P}$.

\item \label{item:eif2} 
\begin{align*}
    \varphi _{P,eff}^{1}=\sum_{j\in[J]}\frac{I\left(
s=j\right) }{P(S=j)}\sum_{k\in[K^{\left( j\right)}]}\frac{dQ}{dP(\cdot|S=j)}(\overline{z}%
_{k-1}^{(j)})m_{k}^{(j)}(\overline{z}%
_{k}^{(j)})
\end{align*}
where $m_{k}^{(j)}\in \mathcal{D}_{k}^{(j)}\left(
Q\right),k\in \left[ K^{\left( j\right) }\right] ,j\in \left[ J\right]$ satisfy \eqref{eq:if-decomposition} for some ideal data influence function $\psi^1_Q$ and there exists $h_{n}^{\left( Q\right) }\in 
\mathcal{T}\left( Q;\mathcal{Q}\right) ,n=1,2,...,$ satisfying for all $k\in \left[ K^{\left( j\right) }\right] ,j\in \left[ J\right]$
\begin{equation*}
m_{k}^{(j)}=\lim_{n\rightarrow \infty }\frac{dP(\cdot|S=j)}{dQ}(\overline{z}_{k-1}^{(j)})P(S=j)\Pi \left[ \left. h_{n}^{\left( Q\right)
} \right\vert \mathcal{D}_{k}^{(j)}\left( Q\right) \right](\overline{z}_{k}^{(j)}).
\end{equation*}%

\item \label{item:eif3} 
\begin{align*}
    \varphi _{P,eff}^{1}(o)=\sum_{j\in[J]}I\left(
s=j\right) \sum_{k\in \left[ K^{\left( j\right) }\right] }\lim_{n\rightarrow
\infty }\Pi \left[ \left. h_{n}^{\left( Q\right) }
\right\vert \mathcal{D}_{k}^{(j)}\left( Q\right) \right](\overline{z}_k^{(j)})
\end{align*} 
where $h_{n}^{\left( Q\right) }\in \mathcal{T}\left( Q;\mathcal{Q}\right)
,n=1,2,...$ satisfy
\begin{equation}
\psi _{Q,eff}^{1}=\lim_{n\rightarrow \infty }\sum_{j\in \lbrack
J]}\sum_{k\in \lbrack K^{(j)}]}\Pi \left\{ \left. \frac{dP(\cdot|S=j)}{dQ}(\overline{Z}_{k-1}^{(j)})P(S=j)\Pi \left[
h_{n}^{\left( Q\right) }|\mathcal{D}_{k}^{(j)}\left( Q\right) \right](\overline{Z}_k^{(j)})
\right\vert \mathcal{T}\left( Q;\mathcal{Q}\right) \right\}
\label{information-equation}
\end{equation}
and  $\psi _{Q,eff}^{1}$ the efficient ideal data influence function.
\end{enumerate}
\end{theorem}

Of course part \ref{item:eif3} of the preceding theorem implies that if we find $h^{\left( Q\right) }$ in $\mathcal{T}\left( Q;%
\mathcal{Q}\right) $ that satisfies simultaneously for all $k\in \left[ K^{\left( j\right) }\right] ,j\in \left[
J\right]$,
\begin{equation}
\psi _{Q,eff}^{1}=\sum_{j\in \lbrack J]}\sum_{k\in \lbrack K^{(j)}]}\Pi
\left\{ \left. \frac{dP(\cdot|S=j)}{dQ}(\overline{Z}_{k-1}^{(j)})P(S=j)\Pi \left[
h^{\left( Q\right) }|\mathcal{D}_{k}^{(j)}\left( Q\right) \right](\overline{Z}_k^{(j)}) \right\vert \mathcal{T}\left( Q;\mathcal{Q}\right) \right\}\label{new-big-eq}
\end{equation}
then $\varphi _{P,eff}^{1}(o)=\sum_{j=1}^{J}I\left( s=j\right)
\sum_{k\in \left[ K^{\left( j\right) }\right] }\Pi \left[ \left. h^{\left(
Q\right) }\left( W\right) \right\vert \mathcal{D}_{k}^{(j)}\left( Q\right) %
\right](\overline{z}_k^{(j)}) $ is the efficient observed data influence function. If $\varphi$ is pathwise differentiable, a sufficient condition for such an $h^{(Q)}$ to exist is that the range of $A_Q^*$ is closed.  

A necessary and
sufficient condition for such an $h^{\left( Q\right) }$ to exist is that $%
\left( \psi _{Q,eff}^{1},\boldsymbol{0}_{J},0\right) $ is in the range of
the information operator $A_{Q,U,\lambda }^{\ast }A_{Q,U,\lambda }:%
\mathcal{H\rightarrow H}$ when $\left( Q,U,P\right) $ is strongly aligned. As discussed in \cite%
{van_der_vaart_asymptotic_2000} Chapter 25.5, when this condition holds $\varphi _{P,eff}^{1}=A_{Q,U,\lambda }\left( A_{Q,U,\lambda }^{\ast
}A_{Q,U,\lambda }\right) ^{-}\left( \psi _{Q,eff}^{1},\boldsymbol{0}%
_{J},0\right)$ is the efficient influence function. Here $\left( A_{Q,U,\lambda }^{\ast
}A_{Q,U,\lambda }\right) ^{-}$ is a generalized inverse.

We have been unable to derive a simple sufficient condition under which  $%
\left( \psi _{Q,eff}^{1},\boldsymbol{0}_{J},0\right) $ is in the range of
the information operator for an arbitrary fused-data framework. In contrast, in coarsening at random models the information operator is invertible under a strong positivity assumption on the coarsening
mechanism. In fact, for coarsening at random models, the inverse of the information operator can be computed with the
method of successive approximations because the identity minus the
information operator is a contraction (\cite{robins_estimation_1994, van_der_laan_unified_2003}). Unfortunately, this technique
cannot be used in general in fused-data models with information
operators with bounded inverses. In \cref{sec:appendix-score-operator}, we exhibit a fused-data model where $A_{Q,U,\lambda }^{\ast
}A_{Q,U,\lambda }$ has a bounded inverse on the appropriate domain but $I-$ $A_{Q,U,\lambda }^{\ast }A_{Q,U,\lambda }$
is not a contraction.

The variance of the observed data efficient influence function is an
efficiency bound that quantifies the information about $\varphi $ in the
observed data. This information comes from two
distinct sets of restrictions imposed on $P$ by the model in fused-data frameworks. The first is the
set of restrictions inherited by $P$ from constraints on $Q$ imposed by the
ideal data model $\mathcal{Q}$. The second is the set of equality
constraints imposed on $P$ by the mere existence of an ideal data
distribution that aligns on the marginals and conditionals dictated by $%
\mathcal{C}$. In \cref{app:model-types-eif} we classify fused-data models according to these two types of restrictions and discuss the inherent difficulties in computing the efficient influence function under each fused-data model type.

\subsection{Examples revisited}

\label{subsec:examples IF}

We now compute the observed data influence functions and the efficient one 
for examples \ref{example:disease-prevalence}-\ref{example:transporting}, providing the derivations for \cref{example:disease-prevalence} and leaving all other derivations to \cref{app:examples-extras} as well as extensions to generalizations of the examples.

\begin{examplerepeattwo}
\label{example:disease-prevalence2} Suppose $%
P\in \mathcal{P}$ is such that the unique $Q \in \mc{Q}$ that aligns with $P$ is strongly aligned. For $\psi \left( Q\right) \coloneqq E_{Q}\left(
Y\right)$, $\psi _{Q}^{1}\left( X,V,Y\right) =Y-\psi \left( Q\right) $
is its unique influence function.

By \cref{prop:two-source-solution} to derive $\varphi^1_P$ we must solve the integral equation 
\begin{align}
E_{Q}[m^{(1)}(X,V)|X=x,Y=y]& =E_{Q}[\psi^1_Q(X, Y, V)|X=x,Y=y] = y - \psi(Q)
\label{eq:example-1-prop-1-equation}\end{align}%
for $m^{(1)}\in L_{0}^{2}(X,V;Q)$. For $m_Q(X, V)$ defined as in \eqref{eq:mq} it holds that $E_Q[m_Q(X, V)|x, y] = y$ and $E_Q[m_Q(X, V)] = \psi(Q)$. Then $m^{(1)}(x, v) \coloneqq m_Q(x, v) - \psi(Q) \in L^2_0(X, V;Q)$ solves \eqref{eq:example-1-prop-1-equation}. It follows by \cref{prop:two-source-solution} that $\varphi\left( P\right)$
is pathwise differentiable at $P$ and 
\begin{equation*}
\varphi _{P}^{1}(o)=\frac{I(s=1)}{P(S=1)}\left\{ m_{Q}(x,v)-\psi (Q)\right\}
+\frac{I(s=2)}{P(S=2)}\frac{q(x,y)}{p(x,y|S=2)}\left\{ y-m_{Q}(x,v)\right\}
\end{equation*}%
is its unique influence function which is a function of $P$ because $Q$ is identified by $P$. This influence function is related to but distinct from the ones derived in  \cite{li_nonparametric_2022} and \cite{park_single_2024} for the shadow variables target parameter. As we discussed above, this distinction arises because the target functionals are different. See \cref{app-subsec:shadow-variables} for more details. 
\end{examplerepeattwo}

\begin{examplerepeattwo}
\label{example:tsiv-lsm2} Suppose that $P\in \mathcal{P}$ and $\delta ^{-1}\leq \frac{p(L|S=2)}{p(L|S=1)}\leq \delta $ a.e. - $P.$ for some $\delta >0$. Then, any $Q \in \Phi\left( P;\mc{C}\right)$ such that $Q(L)=P\left( L|S=1\right)$ is strongly aligned with $P$ and
\begin{align*}
\nu _{P}^{1}\left( o\right) \coloneqq& B_{Q}\left( g\right) ^{-1}g(l)\frac{q\left(
l\right) }{p(l|S=2)}\Big[\frac{I(s=1)}{P(S=1)}\frac{p(l|S=2)}{p(l|S=1)}%
\{y-E_{Q}\left( Y|L=l\right) \} \\
& +\frac{I(s=2)}{P(S=2)}\{E_{Q}\left( Y|L=l\right) -\alpha \left( Q\right)
-\psi \left( Q\right) x\}\Big]
\end{align*}%
is an influence function for $\nu(P) \coloneqq (\tau(P), \varphi(P))'$ for every $g:\mathbb{R}\rightarrow\mathbb{R}^2$ a function of $L$ such that $B_Q(g)\coloneqq E_Q[g(L)(1, X)]$ is non-singular. All influence functions take this form. The efficient influence function of $\nu$ is 
\begin{align*}
\nu _{P,eff}^{1}\left( o\right) =&B_{P\left( \cdot|S=2\right) }\left(
t_{P,eff}\right) ^{-1}t_{P,eff}\left( l\right) \Big[\frac{I(s=1)}{P(S=1)}\frac{p(l|S=2)}{p(l|S=1)}\{y-E_{P\left( \cdot|S=1\right)
}\left( Y|L=l\right) \}\\
&+{\frac{I(s=2)}{P(S=2)}\{E_{P\left( \cdot|S=1\right)
}\left( Y|L=l\right) -\tau \left( P\right) -\varphi \left( P\right) x\}}\Big]
\end{align*}
where $t_{P,eff}\left( L\right) \coloneqq\sigma ^{-2}\left( L\right) E_{P}\left(
U|L\right) $, $U \coloneqq \frac{I(S=2)}{P(S=2)}(1, X)'$. 
The second entries of $\nu _{P}^{1}$ and $\nu _{P,eff}^{1}$
are the influence function and the efficient
one of $\varphi$. This agrees with \cite%
{zhao_two-sample_2019} who derived a class of estimating equations whose
solutions are, up to asymptotic equivalence, all RAL estimators of $\nu$. Their influence functions are $\nu^1_P$ for some $g$.
\end{examplerepeattwo}

\begin{examplerepeattwo}
In all scenarios $\psi \left( Q\right) \coloneqq E_{Q}\left[ E_{Q}\left(
Y|A=1,L\right) -E_{Q}\left( Y|A=0,L\right) \right] $ and $\mathcal{T}\left(
Q,\mathcal{Q}\right) =L_{0}^{2}\left( Q\right) $. The unique influence function of $\psi$ is (\cite{robins_estimation_1994, hahn_role_1998}) 
\begin{equation}
\psi _{Q}^{1}\left( l,a,y\right) =E_{Q}\left( Y|l,A=1\right) -E_{Q}\left(
Y|l,A=0\right) -\psi \left( Q\right) +\frac{2a-1}{q\left( a|l\right) }%
\left\{ y-E_{Q}\left( Y|l,a\right) \right\}. \label{AIPW}
\end{equation}
\label{example:transporting2}
\textbf{Scenario (i):} Suppose $P\in \mathcal{P}$ is such that $\delta ^{-1}<\frac{p(L|S=2)}{p(L|S=1)}<\delta $ a.e.-$P(\cdot|S=2)$  for some $\delta >0$. Then, $Q\in \Phi \left( P;\mc{C}\right) $ with $Q\left( A=1|L\right) =P\left( A=1|L,S=1\right)$ a.e.-$Q$ is strongly aligned with $P$. As discussed in \cref{subsec:score-operator-theory}, even if the target distribution $Q_0$ is aligned but not strongly aligned with $P$, we may use $Q$ instead of $Q_0$ to derive the observed data influence function. The unique influence function of $\varphi$ at $P$ is
\begin{align*}
\varphi _{P}^{1}(o)=&\frac{I(s=1)}{P(S=1)}\frac{p(l|S=2)}{p(l|S=1)}\frac{2a-1}{p\left(
a|l,S=1\right) }\left\{ y-E_{P}\left[Y|a,l, S=1\right]
\right\} \\
&+\frac{I(s=2)}{P(S=2)}\left\{ E_{P}\left[
Y|A=1,l, S=1\right] -E_{P}\left[ Y|A=0,l, S=1\right] -\varphi
\left( P\right) \right\}
\end{align*}
which coincides with the influence function was derived in \cite{rudolph_robust_2017}. See also
\cite{dahabreh_extending_2020} and \cite{li_efficient_2023}.

\textbf{Scenario (ii)} Suppose $P\in \mathcal{P}$ is such that the
unique aligned $Q$ satisfies $\delta ^{-1}\leq \frac{p(L,A|S=1)}{q(L,A)}%
\leq \delta $ a.e.-$Q$ and $Q(Y=1|L, A) \geq \delta^{-1}$ a.e.-$Q$ for some $\delta >0$. Then, in \cref{app:examples-extras} we show $\varphi$ has unique influence function
\begin{align*}
\varphi _{P}^{1}( o)  =&\frac{I( s=1) }{P(S=1)}\frac{%
q( l,a) }{p( l,a|S=1) }\{ \psi _{Q}^{1}(l,a,y)-%
\frac{y}{Q(Y=1|l,a)}E_{Q}[ \psi _{Q}^{1}( L,A,Y) |l,a]
\} \\
& +\frac{I( s=2) }{P(S=2)}\frac{Q(Y=1)}{P(Y=1|S=2)}\frac{y}{%
Q(Y=1|l,a)}E_{Q}[ \psi _{Q}^{1}( L,A,Y) |l,a].
\end{align*}%
The righthand side of the above display is a function of $P$ because $Q$ is identified by $P$. 

\textbf{Scenario (iii):} Suppose $P\in \mathcal{P}$ is such that the
unique aligned $Q$ satisfies $\delta ^{-1}\leq \frac{p(L,A|S=1)}{q(L,A)}%
\leq \delta $ a.e.- $Q$ and $\frac{q(L, A)q(Y)}{q(L,A, Y)} <\delta$ a.e.-$Q$ for some $\delta >0$. Suppose $\mathsf{Supp}\left[ \left( L,A,Y\right) ;Q\right] =\mathsf{Supp}\left[
\left( L,A\right) ;Q\right] \times \mathsf{Supp}\left[ Y;Q\right]$. Let
\begin{align*}
    f_t(l, a, y) \coloneqq \frac{q\left( l,a\right) q\left( y\right) }{q\left(
l,a,y\right) }\{ t\left( l,a,y\right) -E_{Q^{\ast }}\left[ \left.
t\left( L,A,Y\right) \right\vert l,a\right] -E_{Q^{\ast }}\left[ \left.
t\left( L,A,Y\right) \right\vert y\right] +E_{Q^{\ast }}\left[ t\left(
L,A,Y\right) \right] \}
\end{align*} where $E_{Q^*}$ denotes expectation under the law $Q^*$ with density $q^*(l, a, y) \coloneqq q(l, a)q(y)$. In \cref{app:examples-extras} we show that for any $t \in L^2(Q)$,
\begin{align*}
\varphi^1_P\left( o\right) =&\frac{I\left( s=1\right) }{P(S=1)}\frac{%
q\left( l,a\right) }{p\left( l,a|S=1\right) }\left\{ \psi _{Q}^{1}(l,a,y)-%
\frac{q\left( l,a\right) q\left( y\right) }{q\left( l,a,y\right) }E_{Q}\left[
\psi _{Q}^{1}\left( L,A,Y\right) |l,a\right] + f_t(l, a, y)\right\}\\
&+\frac{I\left( s=2\right) }{P(S=2)}\frac{q(y)}{p(y|S=2)}\left\{ \frac{%
q\left( l,a\right) q\left( y\right) }{q\left( l,a,y\right) }E_{Q}\left[ \psi
_{Q}^{1}\left( L,A,Y\right) |l,a\right] -f_t(l, a, y)\right\}
\end{align*}%
is an influence functions of $\varphi$. All influence functions for $\varphi$ take this form for some $t \in L^2(Q)$. The efficient influence function of $\varphi$ is 
\begin{align*}
\varphi _{P,eff}^{1}( o) =& I(s=1)\{ h^{( Q)
}(l,a,y)-E_{Q}[h^{( Q) }(l,a,Y)|l,a]\}+I(s=2)\{ h^{( Q) }(l,a,y)-E_{Q}[h^{( Q)
}(L,A,y)|y]\}
\end{align*}
where $h^{\left( Q\right) }$ solves the integral equation 
\begin{align*}
    \frac{q(l,a,y)}{p(l,a,y)}&\psi _{Q}^{1}(l,a,y)=h(l,a,y)-p(S=1|l,a,y)E_{Q}[h(l,a,Y)|l,a] 
 -p(S=2|l,a,y)E_{Q}[h(L,A,y)|y].
\end{align*}
In fact, $h^{\left( Q\right)}$ admits a closed form expression (see \cref{prop:example-eifs} of the \cref{app:examples-extras}).

In this scenario, the distribution $Q$ and consequently the target parameter $\psi(Q)$ is identified under the weaker alignments \eqref{eq:scenario-2-align-1} and \eqref{eq:scenario-2-align-2} assumed in scenario (ii). Alternatively, $Q$ is also identified when weakening \eqref{eq:scenario-3a-align-1} to $Q(Y=1|l_0, a=0) = P(Y=1|l_0, a=0, S=1)$ a.e.-$Q$, as we study in scenario (iv) in \cref{app:examples-extras}.
These relaxations decrease the efficiency with which $\psi(Q)$ can be estimated. To investigate this phenomenon, we computed the asymptotic relative efficiency of semiparametric efficient estimators for the average treatment effect under a data-generating process that agrees with scenarios (ii), (iii), and (iv) simultaneously for several different values of $P(S=1)$, the probability of observing data from the prospective cohort study. This data-generating process is described in \cref{app-subsubsec:dgp}. \cref{fig:asym-var-main} summarizes the results of this investigation. The degree of variance reduction under scenario (iii) illustrates that the alignment assumptions for this scenario impose strong restrictions on the observed data model. Recall that the observed data models in scenarios (ii) and (iv) do not impose equality constraints. As usual, relaxing assumptions broadens the data-generating processes under which efficient estimators of $\varphi \left( P\right)$ are asymptotically unbiased. 
\begin{figure}[ht]
\caption{Asymptotic relative efficiency of efficient estimators of the ATE under the scenarios (ii), (iii), and (iv) of \cref{example:transporting}}
\label{fig:asym-var-main}\centering
\includegraphics[width=.9\textwidth]{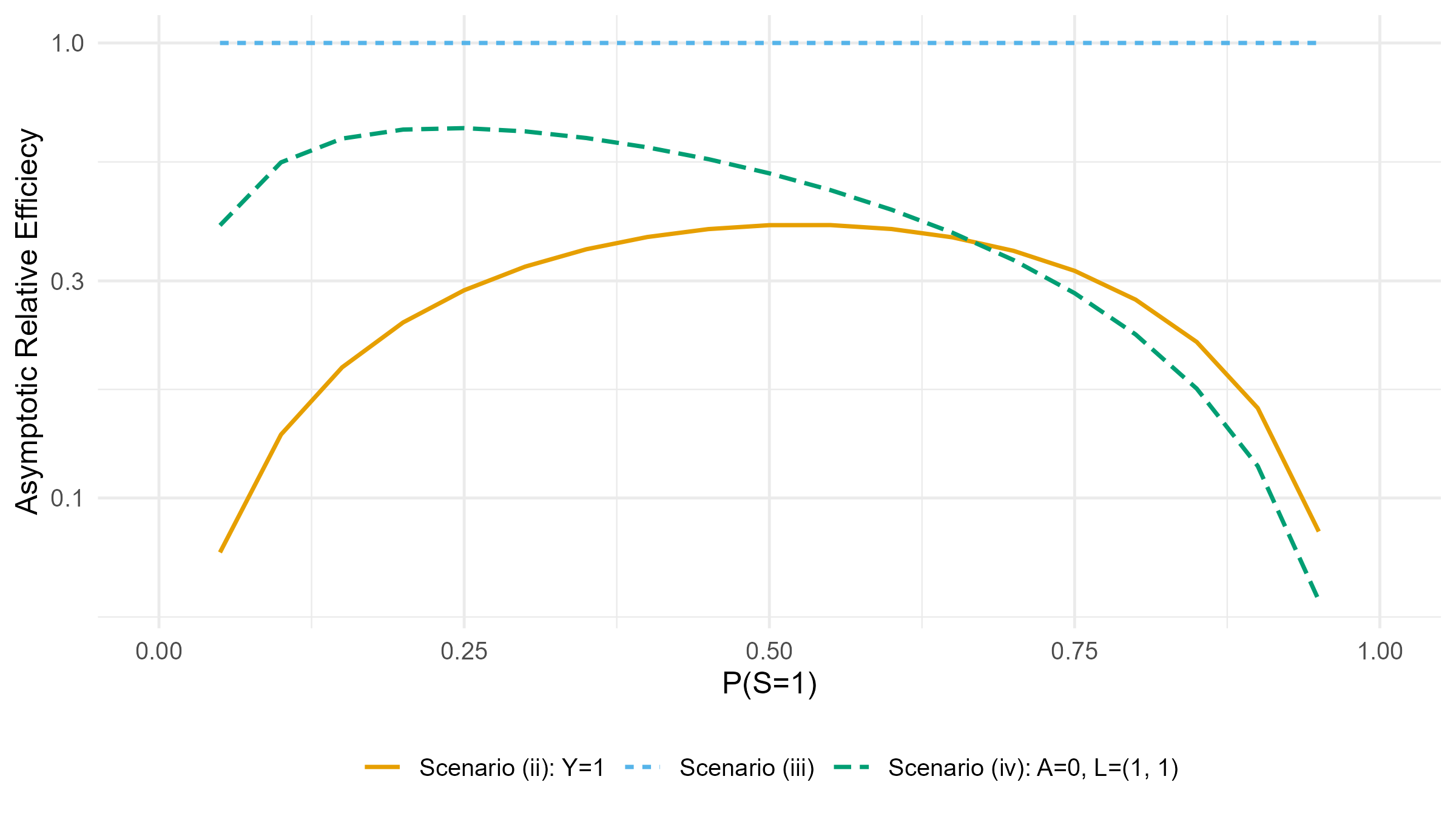}
\end{figure}
\end{examplerepeattwo}

\section{Discussion}

\label{subsec:discussion}

We have introduced a comprehensive framework for integrating
individual-level data from multiple sources. Our framework assumes that
certain conditional or marginal distributions from each source align with
those of the target distribution and that the finite-dimensional parameter
of interest depends on the target distribution only through the aligned
components and is therefore identified by the fused data. Our first main
contribution is a characterization that allows one to directly find the class of all observed data
influence functions from ideal data influence functions without needing to calculate new pathwise derivatives. Our second main contribution is a universal characterization of the structure of the efficient influence functions under our general class of fused-data frameworks. These characterizations pave
the way for conducting machine learning debiased estimation. They also highlight
the challenges of conducting efficient inference in settings
where the alignments themselves impose restrictions on the observed data
distribution.

Our third main contribution is framing inference with fused data in a manner that allows applying the powerful theory outlined in Chapter 25.5 of \cite%
{van_der_vaart_asymptotic_2000}. This framing opens the path to machine learning debiased, semiparametric efficient inference under fused data under alignments other than conditional and marginal distributions, for instance alignments of conditional means, or alignments defined via copulas as in \cite{evans_parameterizing_2024}, \cite{manela_testing_2024}. To compute the influence functions and the semiparametric efficient one in such settings, one would compute the score operator and its adjoint as in \cref{subsec:score-operator-theory} but under the new fused-data models. 

\cite{qiu_efficient_2024} noted that within the framework of \cite{li_efficient_2023}, the influence functions for linear target functionals have a particular multiply robust structure. Our preliminary results indicate that this multiply robust structure is preserved in our more general framework. A full investigation of this topic is beyond the scope of this paper. 

Another direction for future work is to develop sensitivity analysis methods for data fusion. One setting where such sensitivity analysis is already covered by the present theory is as follows. Consider the case in which the target parameter depends on a conditional distribution that is not identified from the observed sources, in addition to aligned conditionals and marginals. One approach to sensitivity analysis is to restrict the model $\mc{Q}$ by fixing the unidentified conditional at a specified value and then estimating the target parameter under this restriction. By varying the fixed conditional distribution over a set of scientifically plausible values and observing the resulting variation in the estimated parameter, one obtains a sensitivity analysis. The theoretical results developed here apply directly in this setting, since within each restricted ideal data model $\mc{Q}$ the target functional is identified.

Our framework addresses data fusion when analysts have access to
individual-level data. However, in many studies, practitioners may have
access to individual-level data from some sources and summary statistics
from others. For example, certain relevant summary statistics may be
available from published material. Additionally, institutions may release
only summary statistics from earlier studies to external researchers to
protect the privacy of study participants. Recent work by \cite%
{hu_paradoxes_2023} develops semiparametric efficient estimation assuming a
random sample from the target population and only summary statistics from
external sources. Developing a unified theory that accommodates the
possibility of using individual-level data from some sources and summary
statistics from others would be of great practical importance and deserves
further study.

\section*{Acknowledgements}

This material is based upon work supported by the National Science Foundation Graduate Research Fellowship Program under Grant No. GE-2140004, by National Heart, Lung, and Blood Institute grant R01-
HL137808, and by National Institute of Allergy and Infectious Diseases grants UM1-AI068635 and R37-
AI029168. The content is solely the responsibility of the authors. Any opinions, findings, and conclusions or
recommendations expressed in this material are those of the authors and do not necessarily reflect the views
of the National Science Foundation or the National Institutes of Health.

\printbibliography
\pagebreak
\appendix

\counterwithin{figure}{inappendix}
\counterwithin{algorithm}{inappendix}
\counterwithin{table}{inappendix}

\renewcommand{\thesection}{S\arabic{section}}
\renewcommand{\thefigure}{S\arabic{figure}}
\renewcommand{\thealgorithm}{S\arabic{algorithm}}
\renewcommand{\thetable}{S\arabic{table}}
\stepcounter{inappendix}
\setcounter{equation}{0}
\renewcommand{\theequation}{S\arabic{equation}}

\begin{refsection}
\addcontentsline{toc}{section}{Supplementary Material} 
\part{Supplementary Material} 
To distinguish the results in the supplement from the main text, we preface section numbers, lemma numbers, equation numbers, etc. that appear for the first time in the supplement by S. 
\parttoc

\section{Glossary of notation}
\label{app:notation}
\glsaddall
\setlength\glsdescwidth{0.8\hsize}
\printglossary[style = long]

\section{Proofs of main text results}

\label{sec:proofs} We first state a lemma that will be useful in the
subsequent proofs.

\begin{slemma}
\label{lemma:f-in-Q} Let $(\mathcal{Q},\mathcal{P},\mathcal{C})$ be a
fused-data model with respect to $\left( Q_{0},P_{0}\right) $. Suppose that $%
\left(Q, U, P\right) \in \mc{Q}\times\mc{U}\times{P}$ is strongly aligned with respect
to $\mathcal{C}$. Let $m_{k}^{(j)}\in \mathcal{D}_{k}^{(j)}(Q)$ and $n_{k}^{(j)}\in \mathcal{R}_{k}^{(j)}(Q)$ for some $k\in [{K}^{(j)}], j\in[J]$. Then the functions
\begin{align*}
\overline{z}_k^{(j)} \mapsto& \frac{dP(\cdot|S=j)}{dQ}(\overline{z}_{k-1}^{(j)})P(S=j)m_{k}^{(j)}(\overline{z}_{k}^{(j)})\\
\overline{z}_k^{(j)} \mapsto& \frac{dQ}{dP(\cdot|S=j)}(\overline{z}_{k-1}^{(j)})P(S=j)^{-1}m_{k}^{(j)}(\overline{z}_{k}^{(j)})
\end{align*}%
are in $\mathcal{D}_{k}^{(j)}(Q)$ and the functions
\begin{align*}
\overline{z}_k^{(j)} \mapsto& \frac{dP(\cdot|S=j)}{dU^{(j)}}(\overline{z}_{k-1}^{(j)})P(S=j)n_{k}^{(j)}(\overline{z}_{k}^{(j)})\\
\overline{z}_k^{(j)} \mapsto& \frac{dU^{(j)}}{dP(\cdot|S=j)}(\overline{z}_{k-1}^{(j)})P(S=j)^{-1}n_{k}^{(j)}(\overline{z}_{k}^{(j)})
\end{align*}%
are in $\mathcal{R}_{k}^{(j)}(P)$ for all $k\in \lbrack K^{(j)}],j\in \lbrack
J]$.
\end{slemma}

\begin{proof}[Proof of \cref{lemma:f-in-Q}.]
Let $f_k^{(j)}(\overline{z}_{k-1}^{(j)})\coloneqq \frac{dP(\cdot|S=j)}{dQ}(\overline{z}_{k-1}^{(j)})P(S=j)m_{k}^{(j)}(\overline{z}_{k}^{(j)})$. We aim to show that $f_k^{(j)} \in \mc{D}_k^{(j)}(Q)$. First note that $f_{k}^{(j)}(\overline{z}_{k}^{(j)})=0$ if $\overline{z}%
_{k}^{(j)}\not\in \overline{\mathcal{Z}}_{k}^{\left( j\right) }$. Next, 
\begin{align*}
& E_{Q}[f_{k}^{(j)}(\overline{Z}_{k}^{(j)})|\overline{Z}_{k-1}^{(j)}] \\
=& E_{Q}\left[\frac{dP(\cdot|S=j)}{dQ}(\overline{Z}_{k-1}^{(j)})P(S=j)m_k^{(j)}(\overline{Z}^{(j)}_k)\}\big|\overline{Z}%
_{k-1}^{(j)}\right] \\
=& \frac{dP(\cdot|S=j)}{dQ}(\overline{z}_{k-1}^{(j)})P(S=j)E_{Q}\left[ m_k^{(j)}(\overline{Z}^{(j)}_k)|\overline{Z}_{k-1}^{(j)}%
\right] \\
=& 0.
\end{align*}%
Additionally, 
\begin{align*}
& E_{Q}\left[ f_{k}^{(j)}(\overline{Z}_{k}^{(j)})^{2}\right] \\
=& E_{Q}\left[ \left\{ \frac{dP(\cdot|S=j)}{dQ}(\overline{Z}_{k-1}^{(j)})P(S=j)m_k^{(j)}(\overline{Z}_{k-1}^{(j)})\right\} ^{2}\right]
\\
\leq& P(S=j)^{2}\delta^2E_{Q}\left[m_k^{(j)}(\overline{Z}_{k-1}^{(j)})^{2}%
\right] \\
<& \infty
\end{align*}%
where the first inequality follows because $(Q, P)$ is strongly aligned. The proofs for the remaining three assertions follow the same lines and we omit them. This
concludes the proof of \cref{lemma:f-in-Q}.
\end{proof}

\begin{proof}[Proof of \cref{lemma:alt-model-characterization}]
Let $\mc{P}_{op} \coloneqq \{P_{Q, U, \lambda}:Q \in \mc{Q}, U \in \mc{U},
\lambda \in \Lambda\}$.

\textbf{Proof that $\mc{P} \subseteq \mc{P}_{op}$:}  Let $P \in \mc{P}$. Let 
$Q\in\mc{Q}$ such that $P\alignswith Q$. Then $Q\in \mc{Q}$ and $P \in \mc{P}%
_{Q}$. For each $j\in[J]$ let $U^{(j)} \coloneqq P(\cdot|S=j)$ and $%
\lambda(S=j) \coloneqq P(S=j)$. From the definition of $%
P_{Q, U, \lambda}$ and the fact that $P\alignswith Q$ we know that $P = P_{Q, U,
\lambda}$. It remains to show $U \in \mc{U}$ and $\lambda \in \Lambda$.
First note that $P(\cdot|S=j) = U^{(j)}$ is mutually
absolutely continuous with $P_0(\cdot|S=j)$ because $P \in \mc{P}_Q$. Hence $U^{(j)} \in \mc{U}^{(j)}$
and so $U \in \mc{U}$. $\lambda \in \Lambda$ because $P(S=j)>0$ for each $j%
\in[J]$. This proves $\mc{P} \subseteq \mc{P}_{op}$.

\textbf{Proof that $\mc{P}_{op} \subseteq \mc{P}$:}  Let $Q \in \mc{Q}, U
\in \mc{U}, \lambda \in \Lambda$. We will show $P_{Q, U, \lambda} \in \mc{P}%
_{Q}\subseteq \mc{P}$. First, it follows from $U^{(j)}\in \mc{U}^{(j)}$ that $U^{(j)}$
is mutually absolutely continuous with $P_0(\cdot|S=j)$. Next, we note that $\lambda \in \Lambda
$ means that $\lambda(S=j) > 0$ for $j\in[J]$. Finally, $P_{Q, U, \lambda}$
and $P_0$ are mutually absolutely continuous because $Q \in \mc{Q}$ is
mutually absolutely continuous with $Q_0$.

It remains to show $P_{Q, U,\lambda}\alignswith Q$. Conditions (a) and (b) in \cref{def:alignswith} hold because $P_0%
\alignswith Q_0$ and $P_0$ and $P_{Q, U, \lambda}$ are mutually absolutely
continuous. Condition (c) of that definition holds by the construction of $%
P_{Q, U, \lambda}$. Hence, $P_{Q, U, \lambda} \in \mc{P}_{Q} \subseteq \mc{P}
$ and so $\mc{P}_{op} \subseteq \mc{P}$. This concludes the proof.
\end{proof}

\begin{proof}[Proof of \cref{thrm:identifiability}]
As in the main text, let $\Xi$ be the partition of $\mc{Q}$ into equivalence
classes $\xi(Q;\mc{C})$. The map $\Psi:\Xi\rightarrow\mathbb{R}$ such that $%
\Psi(\xi(Q;\mc{C})) \coloneqq \psi(\widetilde{Q})$ for all $Q\in \mc{Q}$, $%
\widetilde{Q}\in \xi(Q;\mc{C})$ is well defined by \cref{as:identification}.
The map $\Phi(\cdot;\mc{C}):\mc{P}\rightarrow \Xi$ such that $\Phi(P;\mc{C})\coloneqq \xi(Q;\mc{C}%
)$ for any $Q\in \mc{Q}$ where $P\alignswith Q$ is well defined because for
all $Q, \widetilde{Q}\in \mc{Q}$ such that $P \alignswith Q$ and $P%
\alignswith \widetilde{Q}$, $\xi(Q;\mc{C}) = \xi(\widetilde{Q};\mc{C})$. Let 
$Q \in \mc{Q}$, $P\in \mc{P}$ such that $P\alignswith Q$. Then, $\varphi(P)
=\Psi(\Phi(P;\mc{C})) = \Psi(\xi(Q;\mc{C})) = \psi(Q)$. This concludes the proof of %
\cref{thrm:identifiability}.
\end{proof}

\begin{proof}[Proof of \cref{lemma:score-operator}]

In what follows we define $I_k^{(j)}(\overline{z}_{k-1}^{(j)}) \coloneqq I(\overline{z}_{k-1}^{(j)} \in \mc{Z}_k^{(j)})$. By assumption, $\left( Q,U,P\right) ~$are strongly aligned with respect to $%
\mc{C}$. Therefore, letting $\lambda (S=j)\coloneqq P(S=j),$ it follows from %
\cref{lemma:alt-model-characterization} that $P_{Q,U,\lambda }=P$. Let $%
t:(-\varepsilon ,\varepsilon )\rightarrow Q_{t},t:(-{\varepsilon },{%
\varepsilon })\rightarrow U_{t}$ and $t:(-{\varepsilon },{\varepsilon }%
)\rightarrow \lambda _{t},$ be regular parametric submodels of $\mathcal{Q},%
\mathcal{U}$, and $\Lambda ,$ with $Q_{t=0}=Q,U_{t=0}^{\left( j\right)
}=U^{\left( j\right) }$ and $\lambda _{t=0}=\lambda $ and scores at $t=0$
denoted by $h^{\left( Q\right) },h^{\left( U\right) }$, and $h^{\left(
\lambda \right) }.$ The induced submodel $t:(-\varepsilon ,\varepsilon
)\rightarrow P_{Q_{t},U_{t},\lambda _{t}}$ of $\mathcal{P}$, is
differentiable in quadratic mean if there exists $g\in L^{2}(P_{Q,U,\lambda
})$ that satisfies 
\begin{equation}
\lim_{t\rightarrow 0}\left\Vert \frac{\sqrt{\frac{dP_{Q_{t},U_{t},\lambda
_{t}}}{dP_{Q,U,\lambda }}}-1}{t}-\frac{1}{2}g\right\Vert
_{L^{2}(P_{Q,U,\lambda })}=0  \label{eq:score-operator-condition}
\end{equation}%
i.e. $\frac{1}{2}g$ is the Frechet derivative at $t=0$ of $t\mapsto \sqrt{%
\frac{dP_{Q_{t},U_{t},\lambda _{t}}}{dP_{Q,U,\lambda }}}$ mapping $%
(-\varepsilon ,\varepsilon )$ into $L^{2}(P_{Q,U,\lambda }).$

We will now show that such $g$ exists and it is equal to $A_{Q,U,\lambda }h,$
where $h\coloneqq\left( h^{\left( Q\right) },h^{\left( U\right) },h^{\left(
\lambda \right) }\right) .$ Write 
\begin{equation*}
\frac{dP_{Q_{t},U_{t},\lambda _{t}}}{dP_{Q,U,\lambda }}\left( o\right) =%
\frac{d\lambda _{t}}{d\lambda }\left( s\right) \prod_{j\in \lbrack
J]}\prod_{k\in \left[ K^{\left( j\right) }\right] }r_{k}^{\left( j\right)
}\left( \overline{z}_{k-1}^{(j)},s\right) \widetilde{r}_{k}^{\left( j\right)
}\left( \overline{z}_{k-1}^{(j)},s\right) 
\end{equation*}%
where for $k\in \lbrack K^{(j)}],j\in \lbrack J]$ 
\begin{equation*}
r_{k,t}^{\left( j\right) }\left( \overline{z}_{k}^{(j)},s\right) \coloneqq\left[ 
\frac{\frac{dQ_{t}}{dQ}(\overline{z}_{k}^{(j)})}{\frac{dQ_{t}}{dQ}(\overline{z}%
_{k-1}^{(j)})}\right] ^{I_k^{(j)}(\overline{z}_{k-1}^{(j)})I\left( s=j\right) }
\end{equation*}%
and 
\begin{equation*}
\widetilde{r}_{k,t}^{\left( j\right) }\left( \overline{z}_{k}^{(j)},s\right) %
\coloneqq\left[ \frac{\frac{dU_{t}^{(j)}}{dU^{(j)}}(\overline{z}_{k}^{(j)})}{%
\frac{dU_{t}^{(j)}}{dU^{(j)}}(\overline{z}_{k-1}^{(j)})}\right] ^{I(\overline{z}%
_{k-1}^{(j)}\notin \overline{\mathcal{Z}}_{k-1}^{(j)})I\left( s=j\right) }.
\end{equation*}%
Note that $r_{k,t}^{\left( j\right) }$ and $\widetilde{r}_{k,t}^{\left(
j\right) }$ are well defined for all $k\in \lbrack K^{(j)}],j\in \lbrack J]$
because, by assumption, the laws in $\mathcal{Q}$ are mutually absolutely
continuous and likewise the laws in $\mathcal{U}^{\left( j\right) }$ are
also mutually absolutely continuous. The condition \ref%
{eq:score-operator-condition} is the same as the condition that $\frac{1}{2}g
$ is the Frechet derivative at $t=0$ of the map%
\begin{equation}
t\mapsto \sqrt{\frac{d\lambda _{t}}{d\lambda }}\prod_{j\in \lbrack
J]}\prod_{k\in \left[ K^{\left( j\right) }\right] }\sqrt{r_{k,t}^{\left(
j\right) }}\sqrt{\widetilde{r}_{k,t}^{\left( j\right) }}. \label{eq:bigmap}
\end{equation}%
Now, $\sqrt{\frac{d\lambda _{t}}{d\lambda }},\sqrt{r_{k,t}^{\left( j\right) }%
}$ and $\sqrt{\widetilde{r}_{k,t}^{\left( j\right) }}$ are in $L^{2}\left(
P_{Q,U,\lambda }\right) .$ Then a sufficient condition for the map \ref%
{eq:bigmap} to be Frechet differentiable at $t=0$ is that each of the maps 
\begin{equation}
t\mapsto \sqrt{\frac{d\lambda _{t}}{d\lambda }},\text{ }t\mapsto \sqrt{%
r_{k,t}^{\left( j\right) }}\text{ and }t\mapsto \sqrt{\widetilde{r}%
_{k,t}^{\left( j\right) }}  \label{eq:mapslist}
\end{equation}%
viewed as a map from $(-\varepsilon ,\varepsilon )$ to $L^{2}(P_{Q,U,\lambda
}),$ is Frechet differentiable at $t=0$. Furthermore, in that case, by the product
rule, the Frechet derivative of the map \eqref{eq:bigmap} is equal to 
\begin{equation*}
\left. \frac{d}{dt}\sqrt{\frac{d\lambda _{t}}{d\lambda }}\right\vert
_{t=0}+\sum\limits_{j\in \left[ J\right] }\sum\limits_{k\in \left[ K^{\left(
j\right) }\right] }\left. \frac{d}{dt}\sqrt{r_{k,t}^{\left( j\right) }}%
\right\vert _{t=0}+\sum\limits_{j\in \left[ J\right] }\sum\limits_{k\in %
\left[ K^{\left( j\right) }\right] }\left. \frac{d}{dt}\sqrt{\widetilde{r}%
_{k,t}^{\left( j\right) }}\right\vert _{t=0}
\end{equation*}%
where in a slight abuse of notation $\frac{d}{dt}\left( \cdot \right) $ in
the last display denotes Frechet differentiation. Thus, to show that $g$
exists and is equal to $A_{Q,U,\lambda }h$ it suffices to show that the maps in %
\eqref{eq:mapslist} are Frechet differentiable at $t=0$ and satisfy%
\begin{equation}
\left. \frac{d}{dt}\sqrt{\frac{d\lambda _{t}}{d\lambda }}\right\vert _{t=0}=%
\frac{1}{2}h^{\left( \lambda \right) }  \label{eq:lambdaFrechet}
\end{equation}%
\begin{equation}
\left. \frac{d}{dt}\sqrt{r_{k,t}^{\left( j\right) }}\right\vert _{t=0}\left( 
\overline{z}_{k}^{(j)},s\right) =\frac{1}{2}I\left( s=j\right) \Pi \left[
h^{\left( Q\right) }|\mathcal{D}_{k}^{\left( j\right) }\left( Q\right) %
\right] \left( \overline{z}_{k}^{(j)}\right) \text{ }  \label{eq:rjkFrechet}
\end{equation}%
and 
\begin{equation}
\left. \frac{d}{dt}\sqrt{\widetilde{r}_{k,t}^{\left( j\right) }}\right\vert
_{t=0}=\frac{1}{2}I\left( s=j\right) \Pi \left[ h^{\left( U^{\left( j\right)
}\right) }|\mathcal{R}_{k}^{\left( j\right) }\left( P_{Q,U,\lambda }\right) %
\right] \left( \overline{z}_{k}^{(j)}\right)   \label{eq:rjkTildeFrechet}.
\end{equation}

Now, the equality $\eqref{eq:lambdaFrechet}$ holds because $h^{(\lambda )}$ is
the score at $t=0$ of the regular parametric submodel $t\rightarrow $ $%
\lambda _{t}$ with $\lambda _{t=0}=\lambda $ and  
\begin{equation*}
\left\Vert \frac{\sqrt{\frac{d\lambda _{t}}{d\lambda }}-1}{t}-h^{(\lambda
)}/2\right\Vert _{L^{2}(\lambda )}=\left\Vert \frac{\sqrt{\frac{d\lambda _{t}%
}{d\lambda }}-1}{t}-h^{(\lambda )}/2\right\Vert _{L^{2}(P_{Q,U,\lambda })}
\end{equation*}%
since $\frac{d\lambda _{t}}{d\lambda }$ depends on $o$ only through the
source indicator $s.$ On the other hand, to show the equality \eqref{eq:rjkFrechet} it suffices to show that the maps 
\begin{equation*}
t\mapsto \sqrt{\left[ \frac{dQ_{t}}{dQ}(\overline{z}_{k}^{(j)})\right] ^{I_k^{(j)}(\overline{z}_{k-1}^{(j)})I\left( s=j\right) }}
\end{equation*}%
and%
\begin{equation*}
t\mapsto \sqrt{\left[ \frac{dQ_{t}}{dQ}(\overline{z}_{k-1}^{(j)})\right] ^{I(\overline{%
z}_{k-1}^{(j)}\in \overline{\mathcal{Z}}_{k-1}^{(j)})I\left( s=j\right) }}
\end{equation*}%
viewed as maps from $(-\varepsilon ,\varepsilon )$ into $L^{2}(P_{Q,U,%
\lambda })$ are Frechet differentiable with derivatives at $t=0$ equal to  
\begin{equation*}
\left. \frac{d}{dt}\sqrt{\left[ \frac{dQ_{t}}{dQ}(\overline{z}_{k}^{(j)})\right]
^{I_k^{(j)}(\overline{z}_{k-1}^{(j)})I\left( s=j\right) }%
}\right\vert _{t=0}=\frac{1}{2}I\left( s=j\right) I_k^{(j)}(\overline{z}_{k-1}^{(j)})E\left[ h^{\left( Q\right) }\left( W\right) |\overline{z}_{k}^{(j)}\right] 
\end{equation*}%
and 
\begin{equation}
\label{eq:lemma2-proof-intermediate-1}
\left. \frac{d}{dt}\sqrt{\left[ \frac{dQ_{t}}{dQ}(\overline{z}_{k-1}^{(j)})\right]
^{I_k^{(j)}(\overline{z}_{k-1}^{(j)})I\left( s=j\right) }%
}\right\vert _{t=0}=\frac{1}{2}I\left( s=j\right) I_k^{(j)}(\overline{z}_{k-1}^{(j)})E\left[ h^{\left( Q\right) }\left( W\right) |\overline{z}_{k-1}^{(j)}\right] 
\end{equation}%
since in such case, \eqref{eq:rjkFrechet} follows by an application of the
rule for the derivatives of ratios, to the map 
\begin{equation*}
t\mapsto \left. \sqrt{\left[ \frac{dQ_{t}}{dQ}(\overline{z}_{k}^{(j)})\right] ^{I_k^{(j)}(\overline{z}_{k-1}^{(j)})I\left( s=j\right) }}%
\right/ \sqrt{\left[ \frac{dQ_{t}}{dQ}(\overline{z}_{k}^{(j)})\right] ^{I_k^{(j)}(\overline{z}_{k-1}^{(j)})\left( s=j\right) }}.
\end{equation*}%
Now,       
\begin{align*}
&\lim_{t\rightarrow 0}\int \left\{ \frac{\sqrt{\left[ \frac{dQ_{t}}{dQ}(\overline{z}_{k}^{(j)})\right]
^{I_k^{(j)}(\overline{z}_{k-1}^{(j)})I\left( s=j\right) }%
}-1}{t}-\frac{1}{2}I\left( s=j\right) I_k^{(j)}(\overline{z}_{k-1}^{(j)})E\left[ h^{\left( Q\right) }\left( W\right) |\overline{z}_{k}^{(j)}\right] \right\} ^{2}dP_{Q,U,\lambda }\left(
o\right)  \\
=&\lim_{t\rightarrow 0}\int I\left( s=j\right) I_k^{(j)}(\overline{z}_{k-1}^{(j)})\left\{ \frac{\sqrt{\left[ \frac{dQ_{t}}{dQ}(\overline{z}_{k}^{(j)})%
\right] ^{I_k^{(j)}(\overline{z}_{k-1}^{(j)})I\left(
s=j\right) }}-1}{t}-\frac{1}{2}E_{Q}\left[ h^{\left( Q\right) }\left(
W\right) |\overline{z}_{k}^{(j)}\right] \right\}
^{2}dP_{Q,U,\lambda }\left( o\right)  \\
=&\lim_{t\rightarrow 0}\lambda\left( S=j\right) \int I_k^{(j)}(\overline{z}_{k-1}^{(j)})\left\{ \frac{\sqrt{\frac{dQ_{t}}{dQ}(\overline{z}%
_{k}^{(j)})}-1}{t}-\frac{1}{2}E_{Q}\left[ h^{\left( Q\right) }\left(
W\right) |\overline{z}_{k}^{(j)}\right] \right\}
^{2}dP\left( \overline{z}_{k}^{(j)}|S=j\right)  \\
\leq& \lim_{t\rightarrow 0}\delta \lambda \left( S=j\right) \int \left\{ \frac{\sqrt{\frac{dQ_{t}%
}{dQ}(\overline{z}_{k}^{(j)})}-1}{t}-\frac{1}{2}E_{Q}\left[ h^{\left( Q\right)
}\left( W\right) |\overline{z}_{k}^{(j)}\right] \right\}
^{2}dQ\left( \overline{z}_{k}^{(j)}\right)  \\
=&0
\end{align*}%
where the inequality follows because $dP\left( \cdot |S=j\right) /dQ\left( \overline{z}_{k}^{(j)}\right)I_k^{(j)}(\overline{z}_{k-1}^{(j)}) \leq \delta $ by the strong alignment assumption, and the limit to 0 follows from the known expression for the scores in
information loss models (\cite{van_der_vaart_asymptotic_2000} Section 25.5.2). \eqref{eq:lemma2-proof-intermediate-1} follows analogously which then implies \eqref{eq:rjkFrechet}. The proof of \eqref{eq:rjkTildeFrechet} follows along the same lines with $U^{(j)}$ replacing $Q$ and $\epsilon$ replacing $\delta$ due to the strong alignment of $(U, P)$ and we omit it.

The boundedness of the operator $A_{Q,U,\lambda }$ follows once again from the strong alignment of $(Q,U,P)$ with respect to $\mc{C}$, since
\begin{align*}
&\sum_{j\in \lbrack J]}\sum_{k\in \left[ K^{\left( j\right) }\right]
}E_{P_{Q,U,\lambda }}\left[ I(S=j)I_k^{(j)}(\overline{z}_{k-1}^{(j)})\left\{ E_{Q}[h^{(Q)}(W)|\overline{Z}%
_{k}^{(j)}]-E_{Q}[h^{(Q)}(W)|\overline{Z}_{k-1}^{(j)}]\right\} ^{2}\right] 
\\
=&\sum_{j\in \lbrack J]}\sum_{k\in \left[ K^{\left( j\right) }\right] }E_{Q}%
\Bigg[ \frac{dP_{Q,U,\lambda }(\cdot |S=j)}{dQ}\left( \overline{Z}%
_{k-1}^{(j)}\right) \lambda \left( S=j\right) \\
&\times I_k^{(j)}(\overline{z}_{k-1}^{(j)})\left\{ E_{Q}[h^{(Q)}(W)|\overline{Z}%
_{k}^{(j)}]-E_{Q}[h^{(Q)}(W)|\overline{Z}_{k-1}^{(j)}]\right\} ^{2}\Bigg] 
\\
\leq &\delta \sum_{j\in \lbrack J]}\sum_{k\in \left[ K^{\left( j\right) }%
\right] }E_{Q}\left[ \left\{ I_k^{(j)}(\overline{z}_{k-1}^{(j)})\left[ E_{Q}[h^{(Q)}(W)|\overline{Z}%
_{k}^{(j)}]-E_{Q}[h^{(Q)}(W)|\overline{Z}_{k-1}^{(j)}]\right] \right\} ^{2}%
\right]  \\
\leq &J\delta K\left\Vert h^{(Q)}\right\Vert _{L^{2}\left( Q\right) }
\end{align*}%
with the first inequality in last display holding by the strong alignment of $(Q, P)$ and
the last inequality holding because $I_k^{(j)}(\overline{z}_{k-1}^{(j)})\left\{ E_{Q}[h^{(Q)}(W)|\overline{Z}%
_{k}^{(j)}]-E_{Q}[h^{(Q)}(W)|\overline{Z}_{k-1}^{(j)}]\right\} $ is the
projection of $h^{(Q)}$ into  $\mathcal{D}_{k}^{\left( j\right) }\left(
Q\right) $. Similarly, 
\begin{equation*}
\sum_{j\in \lbrack J]}\left\Vert A_{U^{\left( j\right) }}h^{\left( U^{\left(
j\right) }\right) }\right\Vert _{L^{2}\left( P_{Q,U,\lambda }\right)
}^{2}\leq \sum_{j\in \lbrack J]}\epsilon \left\Vert h^{\left( U^{\left( j\right) }\right) }\right\Vert _{L^{2}\left( U^{\left(
j\right) }\right) }^{2}
\end{equation*}%
by the strong alignment of $(U, P)$. Consequently, 
\begin{eqnarray*}
&&\left. \left\Vert A_{Q,U,\lambda }h\right\Vert _{L^{2}\left(
P_{Q,U,\lambda }\right) }^{2}=\right.  \\
&=&\left\Vert A_{Q}h^{\left( Q\right) }\right\Vert _{L^{2}\left(
P_{Q,U,\lambda }\right) }^{2}+\sum_{j\in \lbrack J]}\left\Vert A_{U^{\left(
j\right) }}h^{\left( U^{\left( j\right) }\right) }\right\Vert _{L^{2}\left(
P_{Q,U,\lambda }\right) }^{2}+\left\Vert A_{\lambda }h^{\left( \lambda
\right) }\right\Vert _{L^{2}\left( P_{Q,U,\lambda }\right) }^{2} \\
&\leq &\delta JK \left\Vert h^{(Q)}\right\Vert _{L^{2}\left( Q\right)
}^{2}+\sum_{j\in \lbrack J]}\epsilon \left\Vert h^{(U^{\left( j\right)
})}\right\Vert _{L^{2}\left( U^{\left( j\right) }\right) }^{2}+\left\Vert
h^{\left( \lambda \right) }\right\Vert _{L^{2}\left( \lambda \right) }^{2} \\
&\leq &\max \left\{ J\delta K ,\epsilon ,1\right\} \left\{ \left\Vert
h^{(Q)}\right\Vert _{L^{2}\left( Q\right) }^{2}+\sum_{j\in \lbrack
J]}\left\Vert h^{(U^{\left( j\right) })}\right\Vert _{L^{2}\left( U^{\left(
j\right) }\right) }^{2}+\left\Vert h^{\left( \lambda \right) }\right\Vert
_{L^{2}\left( \lambda \right) }^{2}\right\}  \\
&=&\max \left\{ J\delta K ,\epsilon ,1\right\} \left\Vert h\right\Vert _{\mathcal{H}%
}^{2}.
\end{eqnarray*}

Next, we show that $A_{Q,U,\lambda }^{\ast }:L_{0}^{2}\left( P_{Q,U,\lambda
}\right) \rightarrow \mathcal{H},$ as defined in the Lemma, is the adjoint
of $A_{Q,U,\lambda }$. Let $g(o) \coloneqq \sum_{j\in[J]}I(s=j)\sum_{k\in[K^{(j)}]}\left\{m_k^{(j)}(\overline{z}_k^{(j)}) + n_k^{(j)}(\overline{z}_k^{(j)})\right\} + \gamma(s)\in L^2_0(P)$ for $m_k^{(j)} \in \mc{D}_k^{(j)}(Q)$, $n_k^{(j)} \in \mc{R}_k^{(j)}(P_{Q, U, \lambda})$, for $k\in[K^{(j)}]$, $j\in[J]$ and $\gamma \in L^2_0(\lambda)$.

We first argue the codomain of $A^*_{Q, U, \lambda}$ is $\mc{H}$. Recall $A^*_{Q, U, \lambda}g = (A^*_Qg, A^*_{U^{(1)}}g, \dots, A^*_{U^{(J)}}g, A^*_\lambda g)$.  
\begin{align*}
    A_{Q}^{\ast
}g\coloneqq\sum_{j\in \lbrack J]}\sum_{k\in \left[ K^{\left( j\right) }%
\right] }\Pi \left[ \left. \frac{dP(\cdot|S=j)}{dQ}(\overline{Z}_{k-1}^{(j)})P(S=j)m_{k}^{(j)}(\overline{Z}_{k}^{(j)})\right\vert \mathcal{T}(Q,\mathcal{Q})%
\right] \in \mathcal{T}(Q,\mathcal{Q})
\end{align*}
by construction where the projection is well defined by \cref{lemma:f-in-Q}. Again by \cref{lemma:f-in-Q}, 
\begin{align*}
    A_{U^{(j)}}^{\ast
}g\coloneqq\sum_{k\in \left[ K^{\left( j\right) }%
\right] }\frac{dP(\cdot|S=j)}{dU^{(j)}}(\overline{z}_{k-1}^{(j)})P(S=j)n_{k}^{(j)}(\overline{z}_{k}^{(j)})\in L^2_0(U^{(j)}) = \mathcal{T}(U^{(j)},\mc{U}^{(j)}) 
\end{align*}
for all $j\in[J]$. Additionally, $A_{\lambda }^{\ast }g\coloneqq\gamma \in
L_{0}^{2}\left( \lambda \right) =\mathcal{T}(\lambda ,\Lambda )$. We then conclude $A^*_{Q, U, \lambda}$ maps $L^2_0(P_{Q, U, \lambda})$ into $\mc{H}$. 

Now let $h = (h^{(Q)}, h^{(U^{(1)})}, \dots, h^{(U^{(J)})}, h^{\lambda}) \in \mc{H}$. The proof that $A_{Q,U,\lambda
}^{\ast }$ is the
adjoint of $A_{Q,U,\lambda }$ is then completed if we show that 
\begin{equation}
\langle g,A_{Q}h^{\left( Q\right) }\rangle _{L^{2}\left( P_{Q,U,\lambda
}\right) }=\left\langle A_{Q}^{\ast }g,h^{(Q)}\right\rangle _{L^{2}\left(
Q\right) }  \label{p1}
\end{equation}%
\begin{equation}
\langle g,A_{U^{\left( j\right) }}h^{\left( U^{(j)}\right) }\rangle
_{L^{2}\left( P_{Q,U,\lambda }\right) }=\left\langle A_{U^{\left( j\right)
}}^{\ast }g,h^{(U^{(j)})}\right\rangle _{L^{2}\left( U^{\left( j\right)
}\right) }\text{ for }j\in \left[ J\right]   \label{p2}
\end{equation}%
and 
\begin{equation*}
\langle g,A_{\lambda }h^{\left( \lambda \right) }\rangle _{L^{2}\left(
P_{Q,U,\lambda }\right) }=\left\langle A_{\lambda }^{\ast }g,h^{(\lambda
)}\right\rangle _{L^{2}\left( \lambda \right) }
\end{equation*}%
since then 
\begin{eqnarray*}
\langle g,A_{Q,U,\lambda }h\rangle _{L^{2}\left( P_{Q,U,\lambda }\right) }
&=&\langle g,A_{Q}h^{\left( Q\right) }\rangle _{L^{2}\left( P_{Q,U,\lambda
}\right) }+\sum_{j\in \left[ J\right] }\langle g,A_{U^{\left( j\right)
}}h^{\left( U^{\left( j\right) }\right) }\rangle _{L^{2}\left(
P_{Q,U,\lambda }\right) }+\langle g,A_{\lambda }h^{\left( \lambda \right)
}\rangle _{L^{2}\left( P_{Q,U,\lambda }\right) } \\
&=&\left\langle A_{Q}^{\ast }g,h^{(Q)}\right\rangle _{L^{2}\left( Q\right)
}+\sum_{j\in \left[ J\right] }\left\langle A_{u^{\left( j\right) }}^{\ast
}g,h^{\left( U^{\left( j\right) }\right) }\right\rangle _{L^{2}\left(
U^{\left( j\right) }\right) }+\left\langle A_{\lambda }^{\ast }g,h^{(\lambda
)}\right\rangle _{L^{2}\left( \lambda \right) } \\
&=&\langle A_{Q,U,\lambda }^{\ast }g,h\rangle _{\mathcal{H}}.
\end{eqnarray*}%
We show $\left( \ref{p1}\right) $ next. Let $f_k^{(j)}(\overline{z}_k^{(j)}) \coloneqq \frac{dP_{Q,U,\lambda }(\cdot |S=j)}{dQ}(\overline{z}%
_{k-1}^{(j)})\lambda (j)m_{k}^{(j)}(\overline{z}_{k}^{(j)})$. Then 
\begin{align*}
& \left. \langle g,A_{Q}h^{\left( Q\right) }\rangle _{P_{Q,U,\lambda
}}\right.  \\
& =\sum_{j\in \lbrack J]}\sum_{k\in \left[ K^{\left( j\right) }\right]
}E_{P_{Q,U,\lambda }}\left[ I(S=j)m_{k}^{(j)}(\overline{Z}_{k}^{(j)})\Pi %
\left[ h^{(Q)}|\mc{D}_{k}^{(j)}(Q)\right] (\overline{Z}_{k}^{(j)})\right]  \\
& =\sum_{j\in \lbrack J]}\sum_{k\in \left[ K^{\left( j\right) }\right]
}E_{P_{Q,U,\lambda }}\left[ \left. \lambda (j)m_{k}^{(j)}(\overline{Z}%
_{k}^{(j)})\Pi \left[ h^{(Q)}|\mc{D}_{k}^{(j)}(Q)\right] (\overline{Z}%
_{k}^{(j)})\right\vert S=j\right]  \\
& =\sum_{j\in \lbrack J]}\sum_{k\in \left[ K^{\left( j\right) }\right] }E_{Q}%
\left[ \frac{dP_{Q,U,\lambda }(\cdot |S=j)}{dQ}(\overline{Z}%
_{k-1}^{(j)})\lambda (j)m_{k}^{(j)}(\overline{Z}_{k}^{(j)})\Pi \left[
h^{(Q)}|\mc{D}_{k}^{(j)}(Q)\right] (\overline{Z}_{k}^{(j)})\right]  \\
& =\sum_{j\in \lbrack J]}\sum_{k\in \left[ K^{\left( j\right) }\right] }E_{Q}%
\left[ f_{k}^{(j)}(\overline{Z}_{k}^{(j)})\Pi \left[ h^{(Q)}|\mc{D}%
_{k}^{(j)}(Q)\right] (\overline{Z}_{k}^{(j)})\right]  \\
& =\sum_{j\in \lbrack J]}\sum_{k\in \left[ K^{\left( j\right) }\right] }E_{Q}%
\left[ f_{k}^{(j)}(\overline{Z}_{k}^{(j)})h^{(Q)}(W)\right]  \\
& =\sum_{j\in \lbrack J]}\sum_{k\in \left[ K^{\left( j\right) }\right] }E_{Q}%
\left[ \Pi \left[ \left. f_{k}^{(j)}\right\vert \mathcal{T}(Q,\mathcal{Q})%
\right] h^{(Q)}(W)\right]  \\
& =\left\langle \sum_{j\in \lbrack J]}\sum_{k\in \left[ K^{\left( j\right) }%
\right] }\Pi \left[ \left. f_{k}^{(j)}\right\vert \mathcal{T}(Q,\mathcal{Q})%
\right] ,h^{(Q)}\right\rangle _{L^{2}\left( Q\right) } \\
& =\left\langle A_{Q}^{\ast }g,h^{(Q)}\right\rangle _{L^{2}\left( Q\right) }.
\end{align*}%
In the preceding display the third equality follows because $P_{Q,U,\lambda }%
\alignswith Q$ and the fifth equality because $f_{k}^{(j)}\in \mc{D}%
_{k}^{(j)}(Q)$ by \cref{lemma:f-in-Q}. The proof of $\left( \ref{p2}\right) $ follows along the
same lines as that of $\left( \ref{p1}\right) $ and we omit it. Finally,%
\begin{eqnarray*}
\langle g,A_{\lambda }h^{\left( \lambda \right) }\rangle _{P_{Q,U,\lambda }}
&=&E_{P_{Q,U,\lambda }}\left[ \gamma \left( S\right) h^{\left( \lambda
\right) }\left( S\right) \right]  \\
&=&E_{\lambda }\left[ \gamma \left( S\right) h^{\left( \lambda \right)
}\left( S\right) \right]  \\
&=&\left\langle \gamma ,h^{\left( \lambda \right) }\right\rangle
_{L^{2}\left( \lambda \right) }
\end{eqnarray*}%
This concludes the proof of \cref{lemma:score-operator}.
\end{proof}

\begin{proof}[Proof of \cref{lemma:extended-model-nonparametric}]
    Let $U$ be such that $(Q, U, P)$ is strongly aligned. Define the extended score
    operator $A_{Q,U,\lambda }^{ext}:\mc{H}^{ext}\rightarrow L_{0}^{2}(P)$ where $\mc{H}^{ext}=L_{0}^{2}(Q)\times \prod_{j\in
    \lbrack J]}L_{0}^{2}(U^{(j)})\times L_{0}^{2}(\lambda )$ and $A_{P,U,\lambda
    }^{ext}$ is the natural extension of $A_{Q,U,\lambda }$ to $\mc{H}^{ext}$ defined identically as $A_{Q,U,\lambda }$ but with domain $\mc{H}^{ext}$.
    
    Let 
    \begin{align*}
        g(o) \coloneqq \sum_{j\in[J]}I(s=j)\sum_{k\in[K^{(j)}]}\left\{m_k^{(j)}(\overline{z}_k^{(j)}) + n_k^{(j)}(\overline{z}_k^{(j)})\right\} + \gamma(s)
    \end{align*}
    for some $m_k^{(j)} \in \mc{D}_k^{(j)}(Q)$, $n_k^{(j)} \in \mc{R}_k^{(j)}(P_{Q, U, \lambda})$, for $k\in[K^{(j)}]$, $j\in[J]$ and $\gamma \in L^2_0(\lambda)$. Then $g\in L^2_0(P)$. Similar arguments as in the proof of \cref{lemma:score-operator} show the adjoint of $A_{Q, U, \lambda}^{ext}$ is given by 
    \begin{equation*}
    A_{Q,U,\lambda }^{ext,\ast }g\coloneqq(A_{Q}^{ext,\ast }g,A_{U^{(1)}}^{\ast
    }g,\dots ,A_{U^{(J)}}^{\ast }g,A_{\lambda }^{\ast }g)
    \end{equation*}%
    where 
    \begin{equation*}
    (A_{Q}^{ext,\ast }g)(w)\coloneqq\sum_{j\in \lbrack J]}\sum_{k\in \lbrack
    K^{(j)}]} \frac{dP(\cdot|S=j)}{dQ}(\overline{z}_{k-1}^{(j)})P(S=j)m_k^{(j)}(\overline{z}_k^{(j)}).
    \end{equation*}%
    Note that for all $g\in L^2_0(P)$, $\Pi[A_{Q}^{ext,\ast }g|\mc{T}(Q, \mc{Q})] = A_Q^*g$. 

    It follows from standard results in linear operator theory that $\mc{T}(P, \mc{P}^{ext}) \coloneqq \overline{\text{Range}(A^{ext}_{Q, U, \lambda})} = \text{Null}(A^{ext, *}_{Q, U, \lambda})^\perp$. Hence, $\mc{T}(P, \mc{P}^{ext}) = L^2_0(P) $ if and only if $\text{Null}(A^{ext, *}_{Q, U, \lambda}) = 0$. 
    
    We first argue that $\text{Null}(A^{ext, *}_{Q, U, \lambda}) = \{0\}$ if the spaces $\mc{D}_k^{(j)}$, $k\in[K^{(j)}], j\in[J]$ are linearly independent in the sense that  $0 = \sum_{j\in[J]}\sum_{k\in[K^{(j)}]}m_k^{(j)}$ for $m_k^{(j)} \in \mc{D}_k^{(j)}(Q)$ if and only if $m_k^{(j)} = 0$ a.e.-$Q$ for all $k\in[K^{(j)}], j\in[J]$. 
    
    Suppose these spaces are linearly independent. Let $g \in \text{Null}(A^{ext, *}_{Q, U, \lambda})$, which we write as 
    \begin{equation*}
        g(o)=\sum_{j\in \lbrack J]}I(s=j)\sum_{k\in \left[ K^{\left( j\right) }%
    \right] }\left\{ m_{k}^{(j)}(\overline{z}_{k}^{(j)})+n_{k}^{(j)}(\overline{z}%
    _{k}^{(j)})\right\} +\gamma \left( s\right).
    \end{equation*}
    Then, 
    \begin{align*}
        (A_{Q}^{ext,\ast }g)(w)\coloneqq\sum_{j\in \lbrack J]}\sum_{k\in \lbrack
    K^{(j)}]} \frac{dP(\cdot|S=j)}{dQ}(\overline{z}_{k-1}^{(j)})P(S=j)m_k^{(j)}(\overline{z}_k^{(j)}) = 0.
    \end{align*}
    But, $\frac{dP(\cdot|S=j)}{dQ}(\overline{z}_{k-1}^{(j)})P(S=j)m_k^{(j)}(\overline{z}_k^{(j)}) \in \mc{D}_k^{(j)}(Q)$ for each $k\in[K^{(j)}], j\in[J]$, and so by assumption $\frac{dP(\cdot|S=j)}{dQ}(\overline{z}_{k-1}^{(j)})P(S=j)m_k^{(j)}(\overline{z}_k^{(j)}) = 0$ a.e.-$Q$. But the strong alignment of $(Q, P)$ then implies that $m_k^{(j)} = 0$ a.e.-$Q$. 

    Furthermore, for each $j\in[J]$,
    \begin{align*}
        (A_{U^{(j)}}^{\ast }g)(w) =& \sum_{k\in[K^{(j)}]}\frac{dP(\cdot|S=j)}{dU^{(j)}}(\overline{z}_{k-1}^{(j)})P(S=j)n_k^{(j)}(\overline{z}_k^{(j)}) = 0
    \end{align*}
    because $g \in \text{Null}(A^{ext, *}_{Q, U, \lambda})$. By strong alignment of $(Q, U, P)$ and the orthogonality of $n_k^{(j)}(\overline{z}_k^{(j)})$ and $n_{k'}^{(j)}(\overline{z}_{k'}^{(j)})$ for $k \not=k'$, it follows that $n_k^{(j)} = 0$ a.e.-$P$ for all $k\in[K^{(j)}], j\in[J]$. 

    Finally, 
    \begin{align*}
        (A_{\lambda}^*g)(s) = \gamma(s) = 0
    \end{align*}
    because $g \in \text{Null}(A^{ext, *}_{Q, U, \lambda})$. Then, $g = 0$, and so $\text{Null}(A^{ext, *}_{Q, U, \lambda}) = 0$. Hence $\mc{T}(P, \mc{P}^{ext}) = L^2_0(P)$. 

    We now prove the converse implication. Suppose that $\text{Null}(A^{ext, *}_{Q, U, \lambda}) = 0$, which is the case if and only if $\mc{T}(P, \mc{P}^{ext}) = L^2_0(P)$. We aim to show that $0 = \sum_{j\in[J]}\sum_{k\in[K^{(j)}]}m_k^{(j)}$ for $m_k^{(j)} \in \mc{D}_k^{(j)}(Q)$ if and only if $m_k^{(j)} = 0$ a.e.-$Q$ for all $k\in[K^{(j)}], j\in[J]$. The if part of this last assertion is trivially true. To prove the only if part, suppose that $0 = \sum_{j\in[J]}\sum_{k\in[K^{(j)}]}m_k^{(j)}$ for $m_k^{(j)} \in \mc{D}_k^{(j)}(Q)$. Then, 
    \begin{align*}
        0 =& \sum_{j\in[J]}\sum_{k\in[K^{(j)}]}\left\{\frac{dP(\cdot|S=j)}{dU^{(j)}}(\overline{z}_{k-1}^{(j)})P(S=j)\right\}\left\{\frac{dP(\cdot|S=j)}{dU^{(j)}}(\overline{z}_{k-1}^{(j)})P(S=j)\right\}^{-1}m_k^{(j)}(\overline{z}_k^{(j)})\\
        =& (A_{Q}^{ext,\ast }g)(w)
    \end{align*}
    where 
    \begin{align*}
        g(o)=\sum_{j\in \lbrack J]}I(s=j)\sum_{k\in \left[ K^{\left( j\right) }%
    \right] }\left\{\frac{dP(\cdot|S=j)}{dU^{(j)}}(\overline{z}_{k-1}^{(j)})P(S=j)\right\}^{-1}m_k^{(j)}(\overline{z}_{k}^{(j)})
    \end{align*}
    because $\left\{\frac{dP(\cdot|S=j)}{dU^{(j)}}(\overline{z}_{k-1}^{(j)})P(S=j)\right\}^{-1}m_k^{(j)}(\overline{z}_k^{(j)}) \in \mc{D}_k^{(j)}(Q)$ for each $j\in[J], k\in[K^{(j)}]$. Since this shows that $g \in \text{Null}(A^{ext, *}_{Q, U, \lambda})$, then $g = 0$ by assumption. This then implies 
    \begin{align*}
        \left\{\frac{dP(\cdot|S=j)}{dU^{(j)}}(\overline{z}_{k-1}^{(j)})P(S=j)\right\}^{-1}m_k^{(j)}(\overline{z}_{k}^{(j)}) = 0
    \end{align*}
    because the spaces $\{I(S=j)\widetilde{m}_k^{(j)}(\overline{Z}_k^{(j)}):  \widetilde{m}_k^{(j)} \in \mc{D}_k^{(j)}(Q)\}$ are mutually orthogonal in $L^2_0(P)$ for $j\in[J]$, $k\in[K^{(j)}]$ . Then, by the strong alignment of $(Q, P)$, $m_k^{(j)} = 0$ a.e.-$Q$. This shows that $0 = \sum_{j\in[J]}\sum_{k\in[K^{(j)}]}m_k^{(j)}$ if and only if $m_k^{(j)} = 0$, hence the spaces $\mc{D}_k^{(j)}(Q)$, $k\in[K^{(j)}], j\in[J]$ are linearly independent. This completes the proof. 
\end{proof}

\begin{proof}[Proof of \cref{lemma:pathwise-differentiabilty}]
Define $%
U^{(j)}\coloneqq P(\cdot |S=j)$. Then, by construction, $(Q,U,P)$ is
strongly aligned with respect to $\mc{C}$. Define and $\lambda (S=j)\coloneqq
P(S=j).$ 

\textbf{Proof that $\varphi$ pathwise differentiable 
$\Leftrightarrow \text{statement \ref{item:pathwise-diff-2}. }$}
This result
follows almost immediately from \cref{lemma:score-operator}. As argued in \cref{subsec:score-operator-theory} and in Theorem
25.31 of \cite{van_der_vaart_asymptotic_2000}, $\varphi $ will be pathwise
differentiable if and only if $(\psi _{Q,eff}^{1},\bs{0}_{J},0)$ is in the
range of $A_{Q,U,\lambda }^{\ast }$. $\bs{0}_{J}$ and $0$ are in the range
of $A_{U}^{\ast }$ and $A_{\lambda }^{\ast }$ respectively because $%
A_{U}^{\ast }$ and $A_{\lambda }^{\ast }$ are linear operators. Next, by %
\cref{lemma:score-operator}, $\psi _{Q,eff}^{1}$ is in the range of $%
A_{Q}^{\ast }$ if and only if there exists $\left\{m_{k}^{(j)}\in \mc{D}%
_{k}^{(j)}(Q):k\in[K^{(j)}], j\in[J]\right\}$ such that 
\begin{equation*}
\psi _{Q;eff}^{1}(w)=\sum_{j\in \lbrack J]}\sum_{k\in \lbrack K^{(j)}]}\Pi 
\left[ \frac{dP(\cdot|S=j)}{dQ}(\overline{Z}_{k-1}^{(j)})P\left( S=j\right) m_{k}^{(j)}(\overline{Z}%
_{k}^{(j)})\Big|\mc{T}(Q,\mc{Q})\right] (w)
\end{equation*}%
where $\psi _{Q;eff}^{1}$ is the efficient influence function of $\psi $ at $%
Q$ in $\mc{Q}$ proving the desired result.

\textbf{Proof that $\text{statement \ref{item:pathwise-diff-2}}\Rightarrow 
\text{statement \ref{item:pathwise-diff-3}}$} Now, suppose statement \ref%
{item:pathwise-diff-2} holds. Then $\varphi$ is pathwise differentiable at $P$. Let $f_k^{(j)}(\overline{z}%
_{k}^{(j)})\coloneqq \frac{dP(\cdot|S=j)}{dQ}(\overline{z}_{k-1}^{(j)})P\left( S=j\right) m_{k}^{(j)}(\overline{z}%
_{k}^{(j)})$. Therefore, it follows from \ref{item:pathwise-diff-2} that%
\begin{equation*}
\psi _{Q}^{1}(w)\coloneqq\sum_{j\in \lbrack J]}\sum_{k\in \lbrack
K^{(j)}]}f_k^{(j)}(\overline{z}%
_{k}^{(j)}).
\end{equation*}%
is an influence function for $\psi $ at $Q$ in $\mc{Q}$ because $\Pi (\psi
_{Q}^{1}|\mc{T}(Q,\mc{Q}))$ is the efficient influence function of $\psi $
at $Q$ in $\mc{Q}$. This proves the desired result because by \cref%
{lemma:f-in-Q}, $f_k^{(j)}\in \mc{D}_{k}^{(j)}(Q).$

\textbf{Proof that $\text{statement \ref{item:pathwise-diff-3}}\Rightarrow 
\text{statement \ref{item:pathwise-diff-2}}$:} Let 
\begin{equation*}
\widetilde{m}_{k}^{(j)}(\overline{z}_{k}^{(j)})\coloneqq\frac{dQ}{dP(\cdot|S=j)}(\overline{z}%
_{k-1}^{(j)})P(S=j)^{-1}m_{k}^{(j)}(\overline{z%
}_{k}^{(j)}).
\end{equation*}%
Thus $
\widetilde{m}_{k}^{(j)}\in \mc{D}_{k}^{(j)}(Q)$ by \cref{lemma:f-in-Q}. Then, $\left\{\widetilde{m}_k^{(j)}:k\in[K^{(j)}], j\in[J]\right\}$ satisfies 
\begin{align*}
    \psi _{Q,eff}^{1}(w)=\sum_{j\in \lbrack J]}\sum_{k\in \lbrack K^{(j)}]}\Pi 
\left[ \frac{dP(\cdot|S=j)}{dQ}(\overline{Z}%
_{k-1}^{(j)})\lambda (S=j)
\widetilde{m}_k^{(j)}(\overline{Z}%
_k^{(j)})\Big|\mathcal{T}(Q,\mathcal{Q})\right] (w)
\end{align*}
because 
\begin{align*}
    \psi^1_{Q, eff} =& \Pi[\psi^1_Q|\mc{T}(Q, \mc{Q})]\\
    =& \sum_{j\in[J]}\sum_{k\in[K^{(j)}]}\Pi[m_k^{(j)}|\mc{T}(Q, \mc{Q})]\\
    =& \sum_{j\in \lbrack J]}\sum_{k\in \lbrack K^{(j)}]}\Pi 
\left[ \frac{dP(\cdot|S=j)}{dQ}(\overline{Z}%
_{k-1}^{(j)})\lambda (S=j)\widetilde{m}_k^{(j)}(\overline{Z}%
_k^{(j)})\Big|\mathcal{T}(Q,\mathcal{Q})\right].
\end{align*}

This concludes the proof of \cref{lemma:pathwise-differentiabilty}.
\end{proof}

\begin{proof}[Proof of \cref{Theorem:influence functions}]
Define $U^{(j)}\coloneqq P(\cdot |S=j)$ and $\lambda (S=j)\coloneqq P(S=j)$.
By construction, $(Q,U,P)$ is strongly aligned with respect to $\mc{C}$.

\textbf{Proof of part \ref{item:influence functions1} $(\Rightarrow)$}

Suppose that $\varphi _{P}^{1}$ is an influence function of $\varphi $ at $%
P\in \mc{P}$. By \cref{lemma:score-operator} and Theorem 25.31 of \cite%
{van_der_vaart_asymptotic_2000}, this implies that $A_{Q,U,\lambda }^{\ast
}\varphi _{P}^{1}=(\psi _{Q,eff}^{1},\bs{0}_{J},0)$. $\varphi
_{P}^{1}$ admits a decomposition  
\begin{equation*}
\varphi _{P}^{1}=\sum_{j\in \lbrack J]}I(s=j)\sum_{k\in \lbrack
K^{(j)}]}\left\{ \widetilde{m}_{k}^{(j)}(\overline{z}_{k}^{(j)})+n_{k}^{(j)}(%
\overline{z}_{k}^{(j)})\right\} +\gamma (s)
\end{equation*}%
for some $\widetilde{m}_{k}^{(j)}\in \mc{D}_{k}^{(j)}(Q)$, $n_{k}^{(j)}\in %
\mc{R}_{k}^{(j)}(P_{Q,U,\lambda })$, $\gamma \in L_{0}^{2}(\lambda )$ because $\varphi^1_{P} \in L^2_0(P)$.
However, $A_{U^{\left( j\right) }}^{\ast }\varphi _{P_{0}}^{1}=0$ implies $%
n_{k}^{(j)}=0$ (a.e.-$U^{(j)}$) for each $k\in \lbrack K^{(j)}]$
and $j\in \lbrack J]$ and $A_{\lambda }^{\ast }\varphi
_{P_{0}}^{1}=0$ implies that $\gamma =0$. Therefore,%
\begin{equation*}
\varphi _{P}^{1}=\sum_{j\in \lbrack J]}I(s=j)\sum_{k\in \lbrack K^{(j)}]}%
\widetilde{m}_{k}^{(j)}(\overline{z}_{k}^{(j)}).
\end{equation*}%
Now, the functions $\widetilde{m}_{k}^{(j)}\in %
\mc{D}_{k}^{(j)}(Q)$ must satisfy the equation 
\begin{equation*}
\psi _{Q,eff}^{1}=\sum_{j\in \lbrack J]}\sum_{k\in \left[ K^{\left( j\right)
}\right] }\Pi \left[ \left. f_{k}^{(j)}\right\vert \mathcal{T}(Q,\mathcal{Q})%
\right] 
\end{equation*}%
with 
\begin{equation*}
f_{k}^{(j)}\left( \overline{z}_{k}^{(j)}\right) \coloneqq\frac{dP(\cdot|S=j)}{dQ}(\overline{z}_{k-1}^{(j)})P(S=j)\widetilde{m}_{k}^{(j)}(%
\overline{z}_{k}^{(j)})
\end{equation*}
because $A_{Q}^{\ast }\varphi _{P_{Q,U,\lambda }}^{1}=\psi
_{Q,eff}^{1}$. 

But then $\psi _{Q}^{1}\coloneqq\sum_{j\in \lbrack J]}\sum_{k\in \left[
K^{\left( j\right) }\right] }f_{k}^{(j)}\left( \overline{z}_{k}^{(j)}\right) $ is
an influence function of $\psi $ because its projection onto the tangent
space $\mc{T}(Q,\mc{Q})$ is the efficient influence function. The proof is
completed by taking  $m_{k}^{(j)}\coloneqq f_{k}^{(j)}$ since then, by %
\cref{lemma:f-in-Q}, $m_{k}^{(j)}\in \mc{D}_{k}^{(j)}(Q)$ and 
\begin{align*}
\varphi _{P}^{1}(o)=& \sum_{j\in \lbrack J]}I(s=j)\sum_{k\in \lbrack
K^{(j)}]}\widetilde{m}_{k}^{(j)}(\overline{z}_{k}^{(j)}) \\
 =&\sum_{j\in \lbrack J]}\frac{I(s=j)}{P(S=j)}\sum_{k\in \lbrack K^{(j)}]}%
\frac{dQ}{dP(\cdot|S=j)}(\overline{z}_{k-1}^{(j)})%
f_{k}^{(j)}\left( \overline{z}_{k}^{(j)}\right)  \\
=& \sum_{j\in \lbrack J]}\frac{I(s=j)}{P(S=j)}\sum_{k\in \lbrack K^{(j)}]}%
\frac{dQ}{dP(\cdot|S=j)}(\overline{z}_{k-1}^{(j)})%
m_{k}^{(j)}(\overline{z}_{k}^{(j)}).
\end{align*}

\textbf{Proof of part \ref{item:influence functions1} $(\Leftarrow)$}

Now, suppose that for some collection  $\left\{ m_{k}^{(j)}\in \mc{D}%
_{k}^{(j)}(Q),k\in \lbrack K^{(j)}],j\in \lbrack J]\right\} $, there exists
an influence function $\psi^1_Q$  for $\psi $ at $Q$ in $\mc{Q}$ such that 
\begin{equation}
\psi _{Q}^{1}=\sum_{j\in \lbrack J]}\sum_{k\in \lbrack K^{(j)}]}m_{k}^{(j)}.
\label{eq:ifPSImjk}
\end{equation}
Let 
\begin{equation*}
\varphi _{P}^{1}\left( o\right) \coloneqq\sum_{j\in \lbrack J]}\frac{I(S=j)}{%
P(S=j)}\sum_{k\in \lbrack K^{(j)}]}\frac{dQ}{dP(\cdot |S=j)}(\overline{z}%
_{k-1}^{(j)})m_{k}^{(j)}(\overline{z}_{k}^{(j)}).
\end{equation*}
We will show $\varphi _{P}^{1}$ satisfies $%
A_{Q,U,\lambda }^{\ast }\varphi _{P}^{1}=\left( \psi _{Q,eff}^{1},\mathbf{0}%
_{J},0\right)$ which will then prove $\varphi^1_P$ is an influence function of $\varphi$ by Theorem 25.31 of \cite{van_der_vaart_asymptotic_2000}. 
Let 
\begin{equation*}
\widetilde{m}_{k}^{(j)}(\overline{z}_{k}^{(j)})\coloneqq\frac{dQ}{dP(\cdot|S=j)}(\overline{z}_{k-1}^{(j)})P(S=j)^{-1}m_{k}^{(j)}(\overline{z}_{k}^{(j)}).
\end{equation*}%
We re-express $\varphi _{P}^{1}$ as $\varphi _{P}^{1}=\sum_{j\in \lbrack
J]}I(S=j)\sum_{k\in \lbrack K^{(j)}]}\widetilde{m}_{k}^{(j)}$. By \cref{lemma:f-in-Q} $\widetilde{m}%
_{k}^{\left( j\right) }\in \mc{D}_{k}^{(j)}(Q)$. Then, $\varphi _{P}^{1}$
satisfies the decomposition \eqref{decompose} with $n_{k}^{(j)}=\gamma
^{(j)}=0$. On the other hand, the expression \ref{eq:ifPSImjk} implies that
the ideal data efficient influence function satisfies 
\begin{equation*}
\psi _{Q,eff}^{1}=\sum_{j\in \lbrack J]}\sum_{k\in \lbrack K^{(j)}]}\Pi 
\left[ m_{k}^{(j)}|\mc{T}(Q,\mc{Q})\right] .
\end{equation*}%
Then, it follows from the expression $A_{Q,U,\lambda }^{\ast }$ established
in \cref{lemma:score-operator}, that $\varphi _{P}^{1}$ satisfies $%
A_{Q,U,\lambda }^{\ast }\varphi _{P}^{1}=\left( \psi _{Q,eff}^{1},\mathbf{0}%
_{J},0\right) $, and as such it is an influence function for $\varphi $ at $P
$ in $\mc{P}$.

\textbf{Proof of part \ref{item:influence functions2}}

Recall the extended score
operator $A_{Q,U,\lambda }^{ext}:\mc{H}^{ext}\rightarrow L_{0}^{2}(P)$ where $\mc{H}^{ext}=L_{0}^{2}(Q)\times \prod_{j\in
\lbrack J]}L_{0}^{2}(U^{(j)})\times L_{0}^{2}(\lambda )$ and $A_{P,U,\lambda
}^{ext}$ is the natural extension of $A_{Q,U,\lambda }$ to $\mc{H}^{ext}$ defined identically as $A_{Q,U,\lambda }$ but with domain $\mc{H}^{ext}$.

Let 
\begin{align*}
    g(o) \coloneqq \sum_{j\in[J]}I(s=j)\sum_{k\in[K^{(j)}]}\left\{m_k^{(j)}(\overline{z}_k^{(j)}) + n_k^{(j)}(\overline{z}_k^{(j)})\right\} + \gamma(s)
\end{align*}
for some $m_k^{(j)} \in \mc{D}_k^{(j)}(Q)$, $n_k^{(j)} \in \mc{R}_k^{(j)}(P_{Q, U, \lambda})$, for $k\in[K^{(j)}]$, $j\in[J]$ and $\gamma \in L^2_0(\lambda)$. Then $g\in L^2_0(P)$. Similar arguments as in the proof of \cref{lemma:score-operator} show the adjoint of $A_{Q, U, \lambda}^{ext}$ is given by 
\begin{equation*}
A_{Q,U,\lambda }^{ext,\ast }g\coloneqq(A_{Q}^{ext,\ast }g,A_{U^{(1)}}^{\ast
}g,\dots ,A_{U^{(J)}}^{\ast }g,A_{\lambda }^{\ast }g)
\end{equation*}%
where 
\begin{equation*}
(A_{Q}^{ext,\ast }g)(w)\coloneqq\sum_{j\in \lbrack J]}\sum_{k\in \lbrack
K^{(j)}]} \frac{dP(\cdot|S=j)}{dQ}(\overline{z}_{k-1}^{(j)})P(S=j)m_k^{(j)}(\overline{z}_k^{(j)}).
\end{equation*}%
Note that for all $g\in L^2_0(P)$, $\Pi[A_{Q}^{ext,\ast }g|\mc{T}(Q, \mc{Q})] = A_Q^*g$. 

We first show that $\varphi^1_P$ is an influence function of $\varphi$ at $P$ in model $\mc{P}$ if and only if there exists an influence function $\psi^1_Q$ at $Q$ in model $\mc{Q}$ such that $A^{*, ext}_Q\varphi^1_P = \psi^1_Q$ and $A^*_{U^{(j)}}\varphi^1_P = A^*_{\lambda}\varphi^1_P = 0$ for all $j\in[J]$. Recall that $\varphi^1_P$ is an influence function if and only if $A^*_{Q, U, \lambda}\varphi^1_P = (\psi^1_{Q, eff}, \bs{0}_J, 0)$ which holds if and only if $A^*_Q\varphi^1_P = \psi^1_{Q, eff}$ and $A^*_{U^{(j)}}\varphi^1_P = A^*_{\lambda}\varphi^1_P = 0$. But  $A^*_Q\varphi^1_P = \Pi\left[A^{ext, *}_Q\varphi^1_P|\mc{T}(Q, \mc{Q})\right]$. Rearranging terms this implies that $\Pi\left[A^{ext, *}_Q\varphi^1_P-\psi^1_{Q, eff}|\mc{T}(Q, \mc{Q})\right] = 0$ or equivalently that $A^{ext, *}_Q\varphi^1_P = \psi^1_{Q, eff} + \widetilde{h}^{(Q)}$ with $\widetilde{h}^{(Q)} \in \mc{T}(Q, \mc{Q})^\perp$. But this means that $A^{ext, *}_Q\varphi^1_P \eqqcolon \psi^1_Q$ is an influence function of $\psi$. 

Suppose the extended model $\mc{P}^{ext}$ is nonparametric at $P$. Then, $\mc{T}(P,\mc{P}^{ext})=L_{0}^{2}(P)$, and so the range of the extended score
operator $A_{Q,U,\lambda }^{ext}$ is
dense in $L_{0}^{2}(P)$. Now, $A^{ext, *}_{Q, U, \lambda}$ is injective since $A^{ext}_{Q, U, \lambda}$ has a dense range in $L^2_0(P)$. This in turn implies that $A^{*, ext}_{Q}$ is injective.  Then, given $\psi^1_Q$ an influence function for $\psi$, the influence function $\varphi^1_P$ for $\varphi$ solving $A^{ext, *}_{Q}\varphi^1_P = \psi^1_Q$ will be unique, if it exists. But this means the collection $\left\{m_k^{(j)}\in \mc{D}_k^{(j)}(Q):k\in[K^{(j)}], j\in[J]\right\}$ satisfying \eqref{newIF} and \eqref{eq:if-decomposition} for $\varphi^1_P$ and $\psi^1_Q$ respectively must be unique because two distinct collections would lead to two distinct influence functions for $\varphi$ solving $A^{ext, *}_{Q}\varphi^1_P = \psi^1_Q$. 

Now, suppose instead the extended model is strictly semiparametric at $P$. Then, $\mc{T}(P, \mc{P}^{ext}) \subsetneq L_0^2(P)$. It follows that $\{0\}\subsetneq \mc{T}(P, \mc{P}^{ext})^\perp = \text{Null}(A^{ext, *}_{Q, U, \lambda})$. Let $\psi^1_Q$ be an influence function for $\psi$ at $Q$ in model $\mc{Q}$. First suppose there does not exist $\varphi^1_P \in L^2_0(P)$ such that $A^{ext, *}_{Q, U, \lambda}\varphi^1_{P} = \psi^1_Q$. Then, there exists no collection $\{m_k^{(j)} \in \mc{D}_k^{(j)}(Q):k\in [K^{(j)}], j\in[J]\}$ satisfying the decomposition \eqref{eq:if-decomposition} because otherwise there would be an observed data influence function $\varphi^1_{P} \in L^2_0(P)$ solving $A^{ext, *}_{Q, U, \lambda}\varphi^1_{P} = \psi^1_Q$ by part \ref{item:influence functions1} of this lemma. Now suppose there does exist $\varphi^1_P \in L^2_0(P)$ such that $A^{ext, *}_{Q, U, \lambda}\varphi^1_{P} = \psi^1_Q$. By part \ref{item:influence functions1} of this lemma, there exists a collection $\{m_k^{(j)} \in \mc{D}_k^{(j)}(Q):k\in [K^{(j)}], j\in[J]\}$ satisfying \eqref{newIF} and \eqref{eq:if-decomposition} for $\varphi^1_P$ and $\psi^1_Q$ respectively. Let $f \in \text{Null}(A^{ext, *}_{Q, U, \lambda})\setminus \{0\}$ which is non-empty by the above argument. Then, $A^{ext, *}_{Q, U, \lambda}(\varphi^1_{P}+ f) = \psi^1_Q$, and so $\varphi^1_{P}+ f$ is an observed data influence function. But then, there must exist a distinct collection $\{\widetilde{m}_k^{(j)} \in \mc{D}_k^{(j)}(Q):k\in [K^{(j)}], j\in[J]\}$ satisfying \eqref{newIF} and \eqref{eq:if-decomposition} for $\varphi^1_P + f$ and $\psi^1_Q$ respectively. That $\{\widetilde{m}_k^{(j)} \in \mc{D}_k^{(j)}(Q):k\in [K^{(j)}], j\in[J]\} \not= \{{m}_k^{(j)} \in \mc{D}_k^{(j)}(Q):k\in [K^{(j)}], j\in[J]\}$ follows because $\varphi^1_P \not= \varphi^1_P + f$. There are infinitely collections satisfying \eqref{eq:if-decomposition} because each choice $f\in \text{Null}(A^{ext, *}_{Q, U, \lambda})\setminus \{0\}$ will lead to a distinct collection of functions satisfying \eqref{newIF} and \eqref{eq:if-decomposition} for $\varphi^1_P + f$ and $\psi^1_Q$ respectively. 
\end{proof}

\begin{proof}[Proof of \cref{prop:two-source-solution}]
Parts \ref{item:two-source-solution1} and \ref{item:two-source-solution2}
are direct corollaries of \cref{lemma:compute-ifs-algorithm} in \cref{app:decomposing-ideal} with $J=2$.
As such we only provide a proof for part \ref{item:two-source-solution3}.
Let 
\begin{align*} 
\mc{G} =& \left\{ \varphi _{P}^{1}(o)\right\} \\
+&\Bigg\{ \sum_{j=1}^{2}\left(
-1\right) ^{j+1}\frac{I(s=j)}{P(S=j)}\sum_{k=1}^{K^{\left( j\right) }}\frac{%
dQ}{dP(\cdot|S=1)}(\overline{z}_{k-1}^{(j)})\Pi \left[ f|\mathcal{D}%
_{k}^{(j)}\left( Q\right) \right] \left( \overline{z}_{k}^{(j)}\right) :f\in 
\mathcal{D}^{(1)}\left( Q\right) \cap \mathcal{D}^{(2)}\left( Q\right)
\Bigg\}.
\end{align*}
where $\mc{D}^{(j)} \coloneqq \bigoplus_{k=1}^{K^{(j)}}\mc{D}_k^{(j)}(Q)$. Let $\left\{m_k^{(j)}\in \mc{D}_k^{(j)}(Q):k\in[K^{(j)}], j\in[J]\right\}$ be
such that 
\begin{align*}
\varphi^1_P(o) = \sum_{j\in[J]}\frac{I(s=j)}{P(S=j)}\sum_{k\in[K^{(j)}]}%
\frac{dQ}{dP(\cdot|S=j)}(\overline{z}_{k-1}^{(j)})m_k^{(j)}(\overline{z}_{k})
\end{align*}
which exists by part \ref{item:influence functions1} of \cref%
{Theorem:influence functions} and let $m^{(j)} \coloneqq \sum_{k\in[K^{(j)}]%
}m_k^{(j)}$. Then $m^{(2)}$ solves \eqref{eq:algorithm-op-eq-simple} and $%
m^{(1)} = m^{(2)} - \psi^1_Q$.

We first show that every element in $\mc{G}$ is an influence function of $%
\varphi $ at $P$ in model $\mc{P}$. Let $f\in \mathcal{D}^{(1)}\left(
Q\right) \cap \mathcal{D}^{(2)}\left( Q\right) $ and 
\[
g_{f}(o)\coloneqq\varphi _{P}^{1}(o)+\sum_{j=1}^{2}\left( -1\right) ^{j+1}%
\frac{I(s=j)}{P(S=j)}\sum_{k=1}^{K^{\left( j\right) }}\frac{dQ}{dP(\cdot
|S=1)}(\overline{z}_{k-1}^{(j)})\Pi \left[ f|\mathcal{D}_{k}^{(j)}\left(
Q\right) \right] \left( \overline{z}_{k}^{(j)}\right) .
\]%
Clearly $g_{f}\in \mc{G}$.

We have that $m^{(2)} + f \in \mc{D}^{(2)}(Q)$ because $f \in \mc{D}^{(2)}(Q)
$. We now show $m^{(2)} + f$ solves \eqref{eq:algorithm-op-eq-simple}. Note
we may rewrite \eqref{eq:algorithm-op-eq-simple} as 
\begin{align*}
\Pi\left[m^{(2)}\left\vert \left(\bigoplus_{k\in[K^{(2)}]} \mc{D}%
_k^{(1)}(Q)\right)^\perp\right.\right] = \Pi\left[\psi^1_Q\left\vert
\left(\bigoplus_{k\in[K^{(2)}]} \mc{D}_k^{(1)}(Q)\right)^\perp\right.\right].
\end{align*}
But $\Pi\left[m^{(2)} + f\left\vert \left(\bigoplus_{k\in[K^{(2)}]} \mc{D}%
_k^{(1)}(Q)\right)^\perp\right.\right] = \Pi\left[m^{(2)}\left\vert
\left(\bigoplus_{k\in[K^{(2)}]} \mc{D}_k^{(1)}(Q)\right)^\perp\right.\right]$
because $f \in \mc{D}^{(1)}(Q) = \bigoplus_{k\in[K^{(2)}]} \mc{D}_k^{(1)}(Q)$%
. But this means $m^{(2)}+f$ also solves \eqref{eq:algorithm-op-eq-simple}.
Hence, 
\begin{align*}
\varphi^1_P(o) + \sum_{j=1}^{2}\left( -1\right) ^{j+1}\frac{I(s=j)}{P(S=j)}%
\sum_{k=1}^{K^{\left( j\right) }}\frac{dQ}{dP(\cdot|S=1)}(\overline{z}%
_{k-1}^{(j)})\Pi \left[ f|\mathcal{D}_{k}^{(j)}\left( Q\right) \right]
\left( \overline{z}_{k}^{(j)}\right)
\end{align*}
is an influence function for $\varphi$ at $P$ in model $\mc{P}$ by part \ref%
{item:two-source-solution2} of this lemma.

We now show every observed data influence function is an element of $\mc{G}$%
. Let $\widetilde{\varphi }_{P}^{1}$ be an influence function for $\varphi $
at $P$ in model $\mc{P}$ that corresponds to $\psi _{Q}^{1}$. Let $\{%
\widetilde{m}_{k}^{(j)}\in \mc{D}_{k}^{(j)}(Q):k\in \lbrack K^{(j)}],j\in
\lbrack J]\}$ be the collection such that 
\[
\widetilde{\varphi }_{P}^{1}(o)=\sum_{j\in \lbrack J]}\frac{I(s=j)}{P(S=j)}%
\sum_{k\in \lbrack K^{(j)}]}\frac{dQ}{dP(\cdot |S=j)}(\overline{z}_{k}^{(j)})%
\widetilde{m}_{k}^{(j)}(\overline{z}_{k}^{(j)}).
\]%
Let $\widetilde{m}^{(j)}\coloneqq\sum_{k\in \lbrack K^{(j)}]}\widetilde{m}%
_{k}^{(j)}\in \mc{D}^{(j)}(Q)$ for $j\in \{1,2\}$. Let $f\coloneqq\widetilde{m}%
^{(2)}-m^{(2)}\in \mc{D}^{(2)}(Q)$. $\widetilde{\varphi }_{P}^{1}$ will belong $%
\mc{G}$ if $f\in \mc{D}^{(1)}(Q)$ because 
\[
\widetilde{\varphi }_{P}^{1}(o)=\varphi _{P}^{1}(o)+\sum_{j=1}^{2}\left(
-1\right) ^{j+1}\frac{I(s=j)}{P(S=j)}\sum_{k=1}^{K^{\left( j\right) }}\frac{%
dQ}{dP(\cdot |S=1)}(\overline{z}_{k-1}^{(j)})\Pi \left[ f|\mathcal{D}%
_{k}^{(j)}\left( Q\right) \right] \left( \overline{z}_{k}^{(j)}\right) .
\]%
But we have that $\widetilde{m}^{(2)}$ is also solution to %
\eqref{eq:algorithm-op-eq-simple} by part \ref{item:compute-ifs-algorithm3}
of \cref{lemma:compute-ifs-algorithm}. Additionally, 
\[
\Pi \left[ m^{(2)}\left\vert \left\{ \bigoplus_{k\in \lbrack K^{(1)}]}\mc{D}%
_{k}^{(1)}(Q)\right\} ^{\perp }\right. \right] =\Pi \left[ \widetilde{m}%
^{(2)}\left\vert \left\{ \bigoplus_{k\in \lbrack K^{(1)}]}\mc{D}%
_{k}^{(1)}(Q)\right\} ^{\perp }\right. \right] 
\]%
because $m^{(2)}$ and $\widetilde{m}^{(2)}$ both solve 
\[
\Pi \left[ f^{(1)}\left\vert \left\{ \bigoplus_{k\in \lbrack K^{(1)}]}\mc{D}%
_{k}^{(1)}(Q)\right\} ^{\perp }\right. \right] =\Pi \left[ \psi
_{Q}^{1}\left\vert \left\{ \bigoplus_{k\in \lbrack K^{(1)}]}\mc{D}%
_{k}^{(1)}(Q)\right\} ^{\perp }\right. \right].
\]%
This implies that, $m^{(2)}-\widetilde{m}^{(2)}\in \mc{D}^{(1)}(Q)$ giving the
desired result. This concludes the proof of this lemma.
\end{proof}

Before proving \cref{theorem:eif}, we provide a lemma characterizing the tangent space $\mc{T}(P, \mc{P})$ of the observed data model $\mc{P}$. In the following for $\mc{E}$
a Hilbert space we will use the notation $\lim_{n\rightarrow\infty}^\mc{E}$
to denote the limit in the space $\mc{E}$.
\begin{slemma}
    \label{lemma:tangent-space}
     Let $\left( \mathcal{Q},\mathcal{P},\mathcal{C}%
\right) $ be a fused-data model with respect to $\left( Q_{0},P_{0}\right) .$
Let $\left( Q,U,P\right) $ be strongly aligned with respect to $\mathcal{C}.$
Let $\lambda(S=j) = P(S=j)$. Then the tangent space $\mc{T}(P, \mc{P})$ of model $\mc{P}$ at $P$ is
\begin{align*}  
    \Bigg\{&o \mapsto \gamma(s) + \sum_{j\in[J]}I(s=j)\sum_{k\in[K^{(j)}]}\left(n_k^{(j)}(\overline{z}_k^{(j)}) + \lim_{n\rightarrow\infty}^{L^2(Q)}\Pi\left[h_n^{(Q)}|\mc{D}_k^{(j)}(Q)\right](\overline{z}_k^{(j)})\right): \\
    &h_n^{(Q)} \in \mc{T}(Q, \mc{Q}) \text{$n\in \{1, 2, \dots\}$ such that the limits exist, }\\
    &n_k^{(j)} \in \mc{R}_k^{(j)}(P) \text{ for all } k\in[K^{(j)}], j\in[J], \gamma \in L^2_0(\lambda)\Bigg\}
\end{align*}
or equivalently 
\begin{align*}
    \mc{T}(P, \mc{P}) = \overline{A_{Q}\mc{T(Q, \mc{Q})}} \oplus \bigoplus_{j\in[J]}A_{U^{(j)}}L^2_0(U^{(j)}) \oplus A_{\lambda}L^2_0(\lambda).
\end{align*}
\end{slemma}

\begin{proof}[Proof of \cref{lemma:tangent-space}]
    Recall from the discussion in \cref{subsec:score-operator-theory} that $\mc{T}(P, \mc{P}) = \overline{A_{Q, U, \lambda}\mc{H}}$. Then
    \begin{align*}
        \mc{T}(P, \mc{P}) =& \overline{A_{Q, U, \lambda}\mc{H}}\\
        =&\overline{A_{Q}\mc{T(Q, \mc{Q})} \oplus \bigoplus_{j\in[J]}A_{U^{(j)}}L^2_0(U^{(j)}) \oplus A_{\lambda}L^2_0(\lambda)}\\
        =& \overline{A_{Q}\mc{T(Q, \mc{Q})}} \oplus \bigoplus_{j\in[J]}A_{U^{(j)}}L^2_0(U^{(j)}) \oplus A_{\lambda}L^2_0(\lambda)
    \end{align*}
    where the second equality follows from the expression of the score operator in \cref{lemma:score-operator} and the third because $A_{Q}\mc{T(Q, \mc{Q})}$, $A_{U^{(1)}}L^2_0(U^{(1)}), \dots, A_{U^{(1)}}L^2_0(U^{(1)})$,  and $A_{\lambda}L^2_0(\lambda)$ are mutually orthogonal and $A_{U^{(j)}}L^2_0(U^{(j)})$ $j\in[J]$, $A_{\lambda}L^2_0(\lambda)$ are closed spaces. By the expression of $A_{U^{(j)}}$ it follows that
    \begin{align*}
        \bigoplus_{j\in[J]}A_{U^{(j)}}L^2_0(U^{(j)}) = \left\{\sum_{j\in[J]}\sum_{k\in[K^{(j)}]}n_k^{(j)}:n_k^{(j)} \in \mc{R}_k^{(j)}(P) \text{ for all } k\in[K^{(j)}], j\in[J]\right\}.
    \end{align*}
    Furthermore, $A_\lambda L^2_0(\lambda) = L^2_0(\lambda)$. 

    The lemma will then be proved if we show that
    \begin{align*}
        \overline{A_{Q}\mc{T(Q, \mc{Q})}} = \left\{o \mapsto \sum_{j\in[J]}I(s=j)\sum_{k\in[K^{(j)}]}\lim_{n\rightarrow \infty}^{L^2(Q)}\Pi\left[h_n^{(Q)}|\mc{D}_k^{(j)}(Q)\right](\overline{z}_k^{(j)})\right\}.
    \end{align*}
    Let $\mc{A}$ be the set in the right-hand side of the above display. Let $f \in \overline{A_{Q}\mc{T(Q, \mc{Q})}}$. Then there exists $h_n^{(Q)} \in \mc{T}(Q, \mc{Q})$, $n\in\{1, 2, \dots\}$ such that  
    \begin{align}
        \label{eq:closure-equivalence}
        \overline{A_{Q}\mc{T(Q, \mc{Q})}} \ni f(o) =& \lim_{n\rightarrow \infty}^{L^2(P)}\sum_{j\in[J]}I(s=j)\sum_{k\in[K^{(j)}]}\Pi\left[h_n^{(Q)}|\mc{D}_k^{(j)}(Q)\right](\overline{z}_k^{(j)})\\
        =&\sum_{j\in[J]}I(s=j)\lim_{n\rightarrow \infty}^{L^2(P(\cdot|S=j))}\sum_{k\in[K^{(j)}]}\Pi\left[h_n^{(Q)}|\mc{D}_k^{(j)}(Q)\right](\overline{z}_k^{(j)})\nonumber\\
        =& \sum_{j\in[J]}I(s=j)\lim_{n\rightarrow \infty}^{L^2(Q)}\sum_{k\in[K^{(j)}]}\Pi\left[h_n^{(Q)}|\mc{D}_k^{(j)}(Q)\right](\overline{z}_k^{(j)})\nonumber\\
        =& \sum_{j\in[J]}I(s=j)\sum_{k\in[K^{(j)}]}\lim_{n\rightarrow \infty}^{L^2(Q)}\Pi\left[h_n^{(Q)}|\mc{D}_k^{(j)}(Q)\right](\overline{z}_k^{(j)}) \in \mc{A}\nonumber
    \end{align}
    where the second equality follows because $I(s=j)f_1(\overline{z}^{(j)})$
    and $I(s=j^{\prime })f_2(\overline{z}^{(j^{\prime })})$ are orthogonal in $%
    L^2(P)$ for $j, j^{\prime }\in[J]$ such that $j\not=j^{\prime }$ and $f_1
    \in L^2(P(\cdot|S=j))$, $f_2 \in L^2(P(\cdot|S=j^{\prime }))$ and the fourth equality follows because $\mc{D}_k^{(j)}(Q), \mc{D}_{k'}^{(j)}(Q)$ are orthogonal in $L^2(Q)$ for $k\not=k'$. The third
    equality follows because $(Q, P)$ strongly aligned implies that a sequence in $\sum_{k\in[K^{(j)}]}\mc{D}_k^{(j)}(Q)$ converges with respect to the 
    $L^2(P(\cdot|S=j))$ norm if and only if it converges with respect to the $L^2(Q)$ norm. To see
    this, note for any element in $m^{(j)} \in \bigoplus_{k\in[K^{(j)}]}\mc{D}_k^{(j)}(Q)
    $,  
    \begin{align*}
    \|m^{(j)}\|_{L^2(P(\cdot|S=j))} =& E_P\left[\left.\left\{\sum_{k\in[K^{(j)}]%
    }m_k^{(j)}(\overline{Z}_k^{(j)})\right\}^2\right|S=j\right] \\
    =&\sum_{k\in[K^{(j)}]}E_P\left[\left.m_k^{(j)}(\overline{Z}%
    _k^{(j)})^2\right|S=j\right] \\
    =&\sum_{k\in[K^{(j)}]}E_Q\left[\frac{dP(\cdot|S=j)}{dQ}(\overline{Z}_{k-1}^{(j)})m_k^{(j)}(\overline{Z}_k^{(j)})^2\right] \\
    \leq&\delta\sum_{k\in[K^{(j)}]}E_Q\left[m_k^{(j)}(\overline{Z}%
    _k^{(j)})^2\right] \\
    =&\delta E_Q\left[\left\{\sum_{k\in[K^{(j)}]}m_k^{(j)}(\overline{Z}%
    _k^{(j)})\right\}^2\right] \\
    =&\delta \|m^{(j)}\|_{L^2(Q)}.
    \end{align*}
    A similar argument shows that $\|m^{(j)}\|_{L^2(P(\cdot|S=j))} \geq \delta^{-1}\|m^{(j)}\|_{L^2(Q)}$. Hence $m^{(j)}_n\in \bigoplus_{k\in[K^{(j)}]}\mc{D}_k^{(j)}(Q)$ converges to $m^{(j)}$ with respect to the $L^2(P(\cdot|S=j))$ norm if and only if it converges to $m^{(j)}$ with respect to the $L^2(Q)$ norm. 

    The above arguments show that $\overline{A_Q\mc{T}(Q, \mc{Q})} \subseteq \mc{A}$. Repeating the steps in \eqref{eq:closure-equivalence} but in reverse starting with $h_n^{(Q)} \in \mc{T}(Q, \mc{Q})$, $n\in\{1, 2, \dots\}$, such that the limits $\lim_{n\rightarrow \infty}^{L^2(Q)}\Pi\left[h_n^{(Q)}|\mc{D}_k^{(j)}(Q)\right](\overline{z}_k^{(j)})$ for $k\in[K^{(j)}]$, $j\in[J]$ exist shows $\mc{A} \subseteq \overline{A_Q\mc{T}(Q, \mc{Q})}$. This completes the proof. 
\end{proof}

\begin{proof}[Proof of \cref{theorem:eif}]
Let $U^{(j)} \coloneqq P(\cdot|S=j)$ and $\lambda(S=j) \coloneqq P(S=j)$ for $j\in[J]$. 

\textbf{Proof of that part \ref{item:eif1} $\Rightarrow$ part \ref%
{item:eif3}}

Suppose $\varphi^1_{P,eff}$ is the efficient influence function of $\varphi$ at $P$ in model $\mc{P}$. We first show that $\varphi _{P,eff}^{1}(o)=\sum_{j\in[J]}I\left(
s=j\right) \sum_{k\in \left[ K^{\left( j\right) }\right] }\lim_{n\rightarrow
\infty }^{L^2(Q)}\Pi \left[ \left. h_{n}^{\left( Q\right) }
\right\vert \mathcal{D}_{k}^{(j)}\left( Q\right) \right](\overline{z}_k^{(j)})$ for some $%
h_{n}^{\left( Q\right) }$ $\in \mathcal{T}\left( Q;\mathcal{Q}\right)
,n=1,2,...$. We have that $A^*_{U^{(j)}}\varphi^1_{P, eff} = 0$ for $j\in[J]$ and $A^*_{\lambda}\varphi^1_{P, eff}=0$ because $\varphi^1_{P, eff}$ is an influence function of $\varphi$. Then $\varphi^1_{P, eff} \in \left\{\bigoplus_{j\in[J]}A_{U^{(j)}}L^2_0(U^{(j)}) \oplus A_{\lambda}L_0^2(\lambda)\right\}^\perp$. Hence $\varphi^1_{P, eff} \in \overline{A_Q\mc{T}(Q, \mc{Q})}$ by \cref{lemma:tangent-space} because $\varphi^1_{P, eff} \in \mc{T}(P, \mc{P})$. But all elements in $\overline{A_Q\mc{T}(Q, \mc{Q})}$ may be written as 
\begin{align*}
    \sum_{j\in[J]}I\left(
s=j\right) \sum_{k\in \left[ K^{\left( j\right) }\right] }\lim_{n\rightarrow
\infty }^{L^2(Q)}\Pi \left[ \left. h_{n}^{\left( Q\right) }
\right\vert \mathcal{D}_{k}^{(j)}\left( Q\right) \right](\overline{z}_k^{(j)})
\end{align*}
for some $%
h_{n}^{\left( Q\right) }$ $\in \mathcal{T}\left( Q;\mathcal{Q}\right)
,n=1,2,...$, proving the desired representation.  

Now we demonstrate that any $h_n^{(Q)}$ that corresponds with $\varphi^1_{P, eff}$ must solve \eqref{information-equation}. Again, $A^*_{Q, U, \lambda}\varphi^1_P = (\psi^1_Q, \bs{0}_J, 0)$ because $\varphi^1_{P, eff}$ is an influence function. This implies that $A^*_Q\varphi^1_{P,eff} = A^*_Q\lim_{n\rightarrow\infty}^{L^2(P)}A_Qh_n^{(Q)} = \psi^1_{Q, eff}$. This in turn implies that $\lim_{n\rightarrow\infty}^{L^2(Q)}A^*_QA_Qh_n^{(Q)} = \psi^1_{Q, eff}$ because $A^*_Q$ is continuous. But this expression may be rewritten as 
\begin{align*}
\psi _{Q,eff}^{1}=\lim_{n\rightarrow\infty}^{L^2(Q)}\sum_{j\in \lbrack
J]}\sum_{k\in \lbrack K^{(j)}]}\Pi \left\{ \left. \frac{dP(\cdot|S=j)}{dQ}(\overline{Z}_{k-1}^{(j)})P(S=j)\Pi \left[ h_n^{\left(
Q\right) }|\mathcal{D}_{k}^{(j)}\left( Q\right) \right](\overline{Z}_k^{(j)}) \right\vert \mathcal{%
T}\left( Q;\mathcal{Q}\right) \right\}
\end{align*}
proving the desired result. 

\textbf{Proof of that part \ref{item:eif3} $\Rightarrow$ part \ref%
{item:eif1}}

Suppose $\varphi _{P,eff}^{1}(o)=\sum_{j\in[J]}I\left(
s=j\right) \sum_{k\in \left[ K^{\left( j\right) }\right] }\lim_{n\rightarrow
\infty }\Pi \left[ \left. h_{n}^{\left( Q\right) }
\right\vert \mathcal{D}_{k}^{(j)}\left( Q\right) \right](\overline{z}_k^{(j)})$ where $%
h_{n}^{\left( Q\right) }$ $\in \mathcal{T}\left( Q;\mathcal{Q}\right)
,n=1,2,...,$ satisfies 
\begin{equation}
\label{eq:eif-intermediary2}
\psi _{Q,eff}^{1}=\lim_{n\rightarrow \infty }\sum_{j\in \lbrack
J]}\sum_{k\in \lbrack K^{(j)}]}\Pi \left\{ \left. \frac{dP(\cdot|S=j)}{dQ}(\overline{Z}_{k-1}^{(j)})P(S=j)\Pi \left[
h_{n}^{\left( Q\right) }|\mathcal{D}_{k}^{(j)}\left( Q\right) \right](\overline{Z}_k^{(j)})
\right\vert \mathcal{T}\left( Q;\mathcal{Q}\right) \right\}.
\end{equation}
Clearly $\varphi^1_{P, eff} \in \mc{T}(P, \mc{P})$ by \cref{lemma:tangent-space}. Additionally, $A^*_{U^{(j)}}\varphi^1_{P, eff} = 0$ for $j\in[J]$ and $A^*_{\lambda}\varphi^1_{P, eff} = 0$. Also $\varphi^1_{P, eff} = \lim_{n\rightarrow \infty}^{L^2(Q)} A_Qh_n^{(Q)}$ from the expression for $A_Q$ in \cref{lemma:score-operator}. It remains to show $A^*_Q\varphi^1_{P, eff} = A^*_Q\lim_{n\rightarrow \infty}^{L^2(Q)}A_Qh_n^{(Q)} = \psi^1_{Q, eff}$. But this follows because the continuity of $A^*_Q$ allows us to move it inside the limit and that $\lim_{n\rightarrow \infty}^{L^2(Q)}A^*_QA_Qh_n^{(Q)}$ is equal to the right-hand side of \eqref{eq:eif-intermediary2}. 

\textbf{Proof of that part \ref{item:eif1} $\Rightarrow$ part \ref{item:eif2}}

Suppose $\varphi^1_{P, eff}$ is the efficient influence function of $\varphi$ at $P$ in model $\mc{P}$. Because part $\ref{item:eif1}\Leftrightarrow $ part \ref{item:eif3}, we know that $\varphi _{P,eff}^{1}(o)=\sum_{j\in[J]}I\left(
s=j\right) \sum_{k\in \left[ K^{\left( j\right) }\right] }\lim_{n\rightarrow
\infty }\Pi \left[ \left. h_{n}^{\left( Q\right) }
\right\vert \mathcal{D}_{k}^{(j)}\left( Q\right) \right](\overline{z}_k^{(j)})$ where $%
h_{n}^{\left( Q\right) }$ $\in \mathcal{T}\left( Q;\mathcal{Q}\right)
,n=1,2,...,$ satisfies 
\begin{equation*}
\psi _{Q,eff}^{1}=\lim_{n\rightarrow \infty }^{L^2(Q)}\sum_{j\in \lbrack
J]}\sum_{k\in \lbrack K^{(j)}]}\Pi \left\{ \left. \frac{dP(\cdot|S=j)}{dQ}(\overline{Z}_{k-1}^{(j)})P(S=j)\Pi \left[
h_{n}^{\left( Q\right) }|\mathcal{D}_{k}^{(j)}\left( Q\right) \right](\overline{Z}_k^{(j)})
\right\vert \mathcal{T}\left( Q;\mathcal{Q}\right) \right\}.
\end{equation*}
We may equivalently write 
\begin{align*}
    &\varphi _{P,eff}^{1}(o)\\
    &=\sum_{j\in[J]}\frac{I\left(
s=j\right)}{P(S=j)} \sum_{k\in \left[ K^{\left( j\right) }\right] }\frac{dQ}{dP(\cdot|S=j)}(\overline{z}_{k-1}^{(j)})\frac{dP(\cdot|S=j)}{dQ}(\overline{z}_{k-1}^{(j)})P(S=j)\lim_{n\rightarrow
\infty }^{L^2(Q)}\Pi \left[ \left. h_{n}^{\left( Q\right) }
\right\vert \mathcal{D}_{k}^{(j)}\left( Q\right) \right](\overline{z}_k^{(j)})\\
&=\sum_{j\in[J]}\frac{I\left(
s=j\right)}{P(S=j)} \sum_{k\in \left[ K^{\left( j\right) }\right] }\frac{dQ}{dP(\cdot|S=j)}(\overline{z}_{k-1}^{(j)})\lim_{n\rightarrow
\infty }^{L^2(Q)}\frac{dP(\cdot|S=j)}{dQ}(\overline{z}_{k-1}^{(j)})P(S=j)\Pi \left[ \left. h_{n}^{\left( Q\right) }
\right\vert \mathcal{D}_{k}^{(j)}\left( Q\right) \right](\overline{z}_k^{(j)})\\
&=\sum_{j\in[J]}\frac{I\left(
s=j\right)}{P(S=j)} \sum_{k\in \left[ K^{\left( j\right) }\right] }\frac{dQ}{dP(\cdot|S=j)}(\overline{z}_{k-1}^{(j)})m_k^{(j)}(\overline{z}_{k}^{(j)})
\end{align*}
with 
\begin{align*}
    m_k^{(j)}(\overline{z}_{k}^{(j)}) \coloneqq \lim_{n\rightarrow \infty}^{L^2(Q)}\frac{dP(\cdot|S=j)}{dQ}(\overline{z}_{k-1}^{(j)})P(S=j)\Pi \left[ \left. h_{n}^{\left( Q\right) }
\right\vert \mathcal{D}_{k}^{(j)}\left( Q\right) \right](\overline{z}_k^{(j)})
\end{align*}
where the second equality follows because $\frac{dP(\cdot|S=j)}{dQ}(\overline{z}_{k-1}^{(j)})P(S=j)$ is bounded by the strong alignment of $(Q, P)$ and as such can be brought inside the limit. 

It remains to show that $\psi^1_Q = \sum_{j\in[J]}\sum_{k\in[K^{(j)}]}m_k^{(j)}$ for some influence function $\psi^1_Q$ of $\psi$. $\varphi^1_{P, eff}$ satisfies $A^*_{Q, U, \lambda}\varphi^1_{P, eff} = (\psi^1_{Q, eff}, \bs{0}_J, 0)$ because it is the efficient influence function. It particular, it satisfies $A^*_Q \varphi^1_{P, eff} = \psi^1_{Q, eff}$. Recall that $A^*_Qg = \Pi\left[A^{ext,*}_{Q}g|\mc{T}(Q, \mc{Q})\right]$ for all $g\in L^2_0(P)$ with $A^{ext,*}_Q$ defined as in part \ref{item:influence functions2} of \cref{Theorem:influence functions}. Then, $\Pi\left[A^{ext, *}_Q\varphi^1_{P, eff} - \psi^1_{Q, eff}|\mc{T}(Q, \mc{Q})\right] = 0$. Equivalently, $A^{ext, *}_Q\varphi^1_{P, eff} =  \psi^1_{Q, eff} + \widetilde{h}^{(Q)}$ with $\widetilde{h}^{(Q)} \in \mc{T}(Q, \mc{Q})^\perp$. This means $A^{ext, *}_Q\varphi^1_{P, eff}$ is an influence function of $\psi$. But, $A^{ext, *}_Q\varphi^1_{P, eff} = \sum_{j\in[J]}\sum_{k\in[K^{(j)}]}m_k^{(j)}$ proving the desired result.

\textbf{Proof of that part \ref{item:eif2} $\Rightarrow$ part \ref{item:eif1}}

We first show $\varphi^1_{P, eff} \in \mc{T}(P, \mc{P})$. We may write
\begin{align*}
    &\varphi^1_{P, eff}(o) \\
    =& \sum_{j\in[J]}\frac{I\left(
s=j\right)}{P(S=j)} \sum_{k\in \left[ K^{\left( j\right) }\right] }\frac{dQ}{dP(\cdot|S=j)}(\overline{z}_{k-1}^{(j)})m_k^{(j)}(\overline{z}_{k}^{(j)})\\
=&\sum_{j\in[J]}\frac{I\left(
s=j\right)}{P(S=j)} \sum_{k\in \left[ K^{\left( j\right) }\right] }\frac{dQ}{dP(\cdot|S=j)}(\overline{z}_{k-1}^{(j)})\lim_{n\rightarrow
\infty }\frac{dP(\cdot|S=j)}{dQ}(\overline{z}_{k-1}^{(j)})P(S=j)\Pi \left[ \left. h_{n}^{\left( Q\right) }
\right\vert \mathcal{D}_{k}^{(j)}\left( Q\right) \right](\overline{z}_k^{(j)})\\
=& \sum_{j\in[J]}I\left(
s=j\right) \sum_{k\in \left[ K^{\left( j\right) }\right] }\lim_{n\rightarrow
\infty }\Pi \left[ \left. h_{n}^{\left( Q\right) }
\right\vert \mathcal{D}_{k}^{(j)}\left( Q\right) \right](\overline{z}_k^{(j)})
\end{align*}
where the third equality follows because $\frac{dQ}{dP(\cdot|S=j)}(\overline{z}_{k-1}^{(j)})P(S=j)^{-1}$ is bounded by the strong alignment of $(Q, P)$ and as such can be brought inside the limit. Hence, $\varphi^1_{P, eff} \in \mc{T}(P, \mc{P})$ by \cref{lemma:tangent-space}. 

Now we show $\varphi^1_{P,eff}$ is an influence function of $\varphi$. This is the case if $A^*_{Q, U, \lambda}\varphi^1_{P, eff} = (\psi^1_{Q, eff}, \bs{0}_J, 0)$. Clearly $A^*_{U^{(j)}}\varphi^1_{P, eff} = 0$ for $j\in[J]$ and $A^*_\lambda \varphi^1_{P, eff} = 0$. It remains to show $A^*_Q\varphi^1_{P, eff} = \psi^1_{Q, eff}$. Let $A^{ext,*}_Q$ be defined as in part \ref{item:influence functions2} of \cref{Theorem:influence functions}. It follows that $A^{ext, *}_Q\varphi^1_{P, eff} = \sum_{j\in[J]}\sum_{k\in[K^{(j)}]}m_k^{(j)}$. But we have that $A^*_Qg = \Pi\left[A^{ext, *}_Qg|\mc{T}(Q, \mc{Q})\right]$ for all $g\in L^2_0(P)$. Then, 
\begin{align*}
    A^*_Q\varphi^1_{P, eff} =& \Pi\left[A^{ext, *}_Q\varphi^1_{P, eff}|\mc{T}(Q, \mc{Q})\right]\\
    =& \Pi\left[\sum_{j\in[J]}\sum_{k\in[K^{(j)}]}m_k^{(j)}|\mc{T}(Q, \mc{Q})\right]\\
    =& \Pi\left[\psi^1_Q|\mc{T}(Q, \mc{Q})\right]\\
    =& \psi^1_{Q, eff}
\end{align*}
completing the proof of this Theorem.
\end{proof}

\section{Decomposing ideal data influence functions}
\label{app:decomposing-ideal}

\cref{prop:two-source-solution} in \cref{subsec:characterize-ifs} provides a method to compute the observed data influence functions corresponding to ideal data influence functions when there are two data sources. This method relies on solving a particular linear operator equation, which then provides the necessary decomposition of the ideal data influence function into functions belonging to the spaces $\mc{D}_k^{(j)}$ for $k\in[K^{(j)}]$, $j\in[J]$. We now provide an algorithm generalizing the procedure outline in \cref{prop:two-source-solution}. This algorithm requires solving $J-1$ non-trivial linear operator equations sequentially. It receives as input an ideal data influence function $\psi _{Q}^{1}$ and the spaces $\left\{ \mathcal{D}%
_{k}^{(j)}\left( Q\right) :k\in \left[ K^{\left( j\right) }\right] ,j\in %
\left[ J\right] \right\} $ and returns a set of functions $\left\{
m_{k}^{(j)}\in \mathcal{D}_{k}^{(j)}\left( Q\right) :k\in \left[ K^{\left(
j\right) }\right] ,j\in \left[ J\right] \right\} \,$\ satisfying $\left( \ref%
{eq:if-decomposition}\right) $ if one such class exists. The algorithm is
complete in the sense that any such class can be found as the output of the
algorithm for any given $\psi _{Q}^{1}.$

In what follows we define $\sum_{l=c}^{c-1}\cdot \coloneqq0$ for any
non-negative integer $c.$

\begin{algorithm}[ht]
    \caption{\texttt{DECOMPOSE}}
    \label{algorithm:decompose}
    \hspace*{\algorithmicindent} \textbf{Input: $\left\{\mc{D}_k^{(j)}(Q):k\in[K^{(j)}], j\in[J]\right\}$}, $\psi^1_Q$.\\
 \hspace*{\algorithmicindent} \textbf{Output:} A class {$\left\{m_k^{(j)}\in \mc{D}_k^{(j)}(Q):k\in[K^{(j)}], j\in[J]\right\}$}; Or FAIL.
    \begin{algorithmic}[1]
        \State $j\gets 1$
        \While{$j \leq J$}
            \If{the operator equation
            \begin{equation}
\label{eq:algorithm-op-eq}
\Pi \left[ f^{(j)}\left\vert \left\{\bigoplus_{k=1}^{K^{(j)}}\mathcal{D}_k^{(j)}(Q)\right\}^{\perp }\right.\right]
=\Pi \left[ \left. \psi _{Q}^{1}-\sum_{l=1}^{j-1}\sum_{k=1}^{K^{\left(
l\right) }}m_{k}^{(l)}\right\vert \left\{\bigoplus_{k=1}^{K^{(j)}}\mathcal{D}_k^{(j)}(Q)\right\}^{\perp }\right] 
\end{equation}
\hspace*{\algorithmicindent} does not have a solution on $\{0\}\cup\left\{\sum_{l={j+1}}^J\bigoplus_{k=1}^{K^{(l)}}\mathcal{D}_k^{(l)}(Q)\right\}$}
            \Return{FAIL} 
            \Else \State Let $\widetilde{f}^{(j)}\in \{0\}\cup\left\{\sum_{l={j+1}}^J\bigoplus_{k=1}^{K^{(l)}}\mathcal{D}_k^{(l)}(Q)\right\}$ be a solution to \eqref{eq:algorithm-op-eq}
            \State $k\gets 1$
            \While{$k\leq K^{(j)}$}
                \State Let 
                \begin{align*}
                    m_k^{(j)} \coloneqq \Pi\left[\left.\psi^1_Q - \widetilde{f}^{(j)} - \sum_{l=1}^{j-1}\sum_{k=1}^{K^{(l)}}m_k^{(l)} \right|\mc{D}_k^{(j)}(Q)\right]
                \end{align*}
                \State $k\gets k+1$
            \EndWhile
            \EndIf
        \State $j\gets j+1$
        \EndWhile
        \State \Return{$\left\{m_k^{(j)}\in \mc{D}_k^{(j)}(Q):k\in[K^{(j)}], j\in[J]\right\}$}
    \end{algorithmic}
\end{algorithm}
\FloatBarrier
\begin{slemma}[Computing observed data influence functions]
\label{lemma:compute-ifs-algorithm} Let $\left( \mathcal{Q},\mathcal{P},%
\mathcal{C}\text{\thinspace },\psi ,\varphi \right) $ be a fused-data
framework with respect to $\left( Q_{0},P_{0}\right) $. Let $P\in \mathcal{P}
$. Suppose there exists $Q$ in $\Phi \left( P;\mc{C}\right) $ such that $\left(
Q,P\right) $ is strongly aligned with respect to $\mathcal{C}$ and, $\psi $
is pathwise differentiable at $Q$ in model $\mathcal{Q}$. Then

\begin{enumerate}
\item \label{item:compute-ifs-algorithm1} $\varphi $ is pathwise
differentiable if and only if there exists an ideal data influence function $%
\psi _{Q}^{1}$ at $Q$ in model $\mathcal{Q}$ such that Algorithm \texttt{%
DECOMPOSE} does not return FAIL

\item \label{item:compute-ifs-algorithm2} If the Algorithm \texttt{DECOMPOSE} returns a collection $\left\{
m_{k}^{(j)}\in \mathcal{D}_{k}^{(j)}\left( Q\right) :k\in \left[ K^{\left(
j\right) }\right] ,j\in \left[ J\right] \right\} $ then 
\begin{equation}
\varphi _{P}^{1}\left( o\right) =\sum_{j=1}^{J}\frac{I\left( s=j\right) }{%
P(S=j)}\sum_{k=1}^{K^{\left( j\right) }}\frac{dQ}{dP(\cdot|S=j)}(\overline{z}_{k-1}^{(j)})m_{k}^{(j)}(\overline{z}_{k}^{(j)})
\label{observedIF-new}
\end{equation}%
is an influence function of $\varphi $ at $P$ in model $\mathcal{P}$

\item \label{item:compute-ifs-algorithm3} Any collection $\left\{ m_{k}^{(j)}\in \mathcal{D}_{k}^{(j)}\left(
Q\right) :k\in \left[ K^{\left( j\right) }\right] ,j\in \left[ J\right]
\right\} $ such that the right hand side of $\left( \ref{observedIF-new}%
\right) $ is an observed data influence function of $\varphi $ at $P$ in
model $\mathcal{P\,\ }$is the output of the Algorithm \texttt{DECOMPOSE} for
some $\psi _{Q}^{1}$ and some choice of solutions $\widetilde{f}%
^{(j)},j=1,...,J,$ to equation \eqref{eq:algorithm-op-eq}.
\end{enumerate}
\end{slemma}

Recall that \cite{li_efficient_2023} considered a fused-data framework in
which the aligned distributions from different sources correspond to factors
of a single factorization of the joint law of $W$. In this special case, the
algorithm simplifies significantly and no longer requires solving integral equations. This simplification arises because, in the fused-data frameworks they consider, for all $k\in [dim(W)]$, $j, j' \in \{j: K^{(j)} \geq k\}$, $m_k^{(j)}\in \mc{D}_k^{(j)}(Q)$, 
\begin{align*}
    \Pi\left[m_k^{(j)}\left\vert\mc{D}_{k}^{(j')}(Q)\right.\right](\overline{z}^{(j')}_k) = I(\overline{z}^{(j')}_{k-1} \in \mc{Z}_{k-1}^{(j')})m_k^{(j)}(\overline{z}^{(j)}_k)
\end{align*}
and, for all $k, k'\in [dim(W)]$ such that $k\not=k'$, $j, j'\in \{j: K^{(j)} \geq k\}\cap \{j: K^{(j)} \geq k'\}$, $m_k^{(j)}\in \mc{D}_k^{(j)}(Q)$,  $\Pi\left[m_k^{(j)}\left\vert\mc{D}_{k'}^{(j')}(Q)\right.\right] = 0$. As such, there are no integral equations to solve when applying \cref{algorithm:decompose} \texttt{DECOMPOSE}. This property results from the assumption that all aligned components are factors in a single factorization of the joint distribution $Q$ of $W$. Applying \cref{algorithm:decompose} \texttt{DECOMPOSE} like in part \ref{item:compute-ifs-algorithm3} of the above lemma with this fact in hand gives the influence function in Theorem 2 of \cite{li_efficient_2023}, as well as all other influence functions in their fused-data frameworks. 

\subsection{Proofs for \texorpdfstring{\cref{app:decomposing-ideal}}{Supplement C}}

\begin{proof}[Proof of \cref{lemma:compute-ifs-algorithm}]
$ $\newline
\textbf{Proof of part \ref{item:compute-ifs-algorithm1} ($\Rightarrow$):}
Suppose $\varphi$ is pathwise differentiable at $P$ in $\mc{P}$. We fix $J$
and will use induction on $j \in \{1, \dots, J\}$ to prove this result.

Suppose $\varphi $ is pathwise differentiable at $P$ in model $\mathcal{P}$.
Let 
\[
\psi _{Q}^{1}=\sum_{l\in \lbrack J]}\sum_{k\in \lbrack K^{(l)}]}\widetilde{m}%
_{k}^{(l)}
\]
be an arbitrary ideal data influence function at $Q$ in model $\mathcal{Q}$
such that there exists $\{\widetilde{m}_{k}^{(j)}\in \mc{D}_{k}^{(j)}(Q):k\in \lbrack
K^{(l)}],l\in \lbrack J]$ satisfying the above display. That such a $\psi_{Q}^{1}$ exists
follows from \cref{lemma:pathwise-differentiabilty}.  First, set $j=1$. We will argue that
a solution $\widetilde{f}^{(1)}$ to \eqref{eq:algorithm-op-eq} exists that
satisfies $\widetilde{f}^{(1)}\in \left\{ 0\right\} \cup
\sum_{l=2}^{J}\bigoplus_{k=1}^{K^{(l)}}\mc{D}_{k}^{(l)}(Q)$. Then, we will
demonstrate that for any solution to \eqref{eq:algorithm-op-eq} in $\left\{
0\right\} \cup \sum_{l=2}^{J}\bigoplus_{k=1}^{K^{(l)}}\mc{D}_{k}^{(l)}(Q)$,
there exists a collection $\left\{ m_{k}^{\dag (j)}\in \mc{D}_{k}^{(j)}(Q):k\in
\lbrack K^{(j)}],j=2,...,J\right\} $ such that 
\[
\psi _{Q}^{1}=\sum_{k\in \lbrack
K^{(1)}]}m_{k}^{(1)}+\sum_{j=2}^{J}\sum_{k\in \lbrack K^{(j)}]}m_{k}^{\dag
(j)}
\]
where $m_{k}^{(1)}\in \mc{D}_{k}^{(1)}(Q)$, $k\in \lbrack K^{(1)}],$ are
defined in the algorithm. 

Define $\widetilde{f}^{(1)}\coloneqq\sum_{l=2}^{J}\sum_{k\in \lbrack
K^{(l)}]}\widetilde{m}_{k}^{(l)}$. Then, 
\begin{align*}
\Pi \left[ \widetilde{f}^{(1)}\left\vert \left\{ \bigoplus_{k=1}^{K^{(1)}}%
\mc{D}_{k}^{(1)}(Q)\right\} ^{\perp }\right. \right] =& \Pi \left[ \psi
_{Q}^{1}-\sum_{k\in \lbrack K^{(1)}]}\widetilde{m}_{k}^{(1)}\left\vert
\left\{ \bigoplus_{k=1}^{K^{(1)}}\mc{D}_{k}^{(1)}(Q)\right\} ^{\perp
}\right. \right]  \\
=& \Pi \left[ \psi _{Q}^{1}\left\vert \left\{ \bigoplus_{k=1}^{K^{(1)}}\mc{D}%
_{k}^{(1)}(Q)\right\} ^{\perp }\right. \right] .
\end{align*}%
Hence, $\widetilde{f}^{(1)}$ is a solution to \eqref{eq:algorithm-op-eq}
which satisfies that $\widetilde{f}^{(1)}\in \left\{ 0\right\} \cup
\sum_{l=2}^{J}\bigoplus_{k=1}^{K^{(l)}}\mc{D}_{k}^{(l)}(Q)$.

Next, let $\widetilde{\widetilde{f}}^{(1)}$ be an arbitrary solution to %
\eqref{eq:algorithm-op-eq} in $\left\{ 0\right\} \cup
\sum_{l=2}^{J}\sum_{k\in \lbrack K^{(l)}]}\mc{D}_{k}^{(l)}(Q)$ for $j=1$.
Let $m_{k}^{(1)}\coloneqq\Pi \left[ \psi _{Q}^{1}-\widetilde{\widetilde{f}}%
^{(1)}|\mc{D}_{k}^{(1)}(Q)\right] $ for each $k\in \lbrack K^{(1)}],$ be the
functions defined in the algorithm for step $j=1$. Then $\psi _{Q}^{1}-%
\widetilde{\widetilde{f}}^{(1)}\in \bigoplus_{k=1}^{K^{(1)}}\mc{D}%
_{k}^{(1)}(Q)$ because $\Pi \left[ \psi _{Q}^{1}-\widetilde{\widetilde{f}}%
^{(1)}\left\vert \left\{ \bigoplus_{k=1}^{K^{(1)}}\mc{D}_{k}^{(1)}(Q)\right%
\} ^{\perp }\right. \right] =0$ by $\widetilde{\widetilde{f}}^{(1)}$ being a
solution of \eqref{eq:algorithm-op-eq}. Additionally, $%
\sum_{k=1}^{K^{(1)}}m_{k}^{(1)}=\psi _{Q}^{1}-\widetilde{\widetilde{f}}^{(1)}
$ because $\mc{D}_{k}^{(1)}(Q)$ and $\mc{D}_{k^{\prime }}^{(1)}(Q)$ are
orthogonal for $k\not=k^{\prime }$. Finally, we note that we can write $%
\widetilde{\widetilde{f}}^{(1)}=\sum_{l=2}^{J}\sum_{k\in \lbrack
K^{(l)}]}m_{k}^{\dag (l)}$ for some $m_{k}^{\dag (l)}\in \mc{D}_{k}^{(l)}(Q)$%
, $k\in \lbrack K^{(l)}]$, $l\in \{2,\dots ,J\}$ because $\widetilde{%
\widetilde{f}}^{(1)}$ is in $\left\{ 0\right\} \cup \sum_{l=2}^{J}\sum_{k\in
\lbrack K^{(l)}]}\mc{D}_{k}^{(l)}(Q)$. Thus, 
\begin{align*}
\psi _{Q}^{1}=& \widetilde{\widetilde{f}}^{(1)}+\sum_{k\in \lbrack
K^{(j)}]}m_{k}^{(j)} \\
=& \sum_{l=2}^{J}\sum_{k\in \lbrack K^{(l)}]}m_{k}^{\dag (l)}+\sum_{k\in
\lbrack K^{(1)}]}m_{k}^{(1)}.
\end{align*}%
This concludes the proof for the step $j=1$.

We now move to the inductive step. Suppose that a solution $\widetilde{f}%
^{(r)}$ in $\left\{ 0\right\} \cup \sum_{r=j+1}^{J}\sum_{k\in \lbrack
K^{(r)}]}\mc{D}_{k}^{(r)}(Q)$ exists for all $r<j$ and let $m_{k}^{(l)}$'s $%
k\in \lbrack K^{(l)}]$, $l\in \{1,\dots ,j-1\}$ be defined as in the
algorithm. Suppose also there exists a collection $\{m_{k}^{\dag (l)}\in %
\mc{D}_{k}^{(l)}(Q):k\in \lbrack K^{(l)}],l=j,...,J\}$ such that $\psi
_{Q}^{1}=\sum_{l=1}^{j-1}\sum_{k\in \lbrack
K^{(l)}]}m_{k}^{(l)}+\sum_{l=j}^{J}\sum_{k\in \lbrack K^{(l)}]}m_{k}^{\dag
(l)}$. We will again show that a solution $\widetilde{f}^{(j)}$exists to %
\eqref{eq:algorithm-op-eq} such that $\widetilde{f}^{(j)}\in \left\{
0\right\} \cup \sum_{l=j+1}^{J}\sum_{k\in \lbrack K^{(l)}]}\mc{D}%
_{k}^{(l)}(Q)$. Next, we will demonstrate that for any solution $\widetilde{%
\widetilde{f}}^{(j)}$ to \eqref{eq:algorithm-op-eq} such that $\widetilde{%
\widetilde{f}}^{(j)}\in \left\{ 0\right\} \cup \sum_{l=j+1}^{J}\sum_{k\in
\lbrack K^{(l)}]}\mc{D}_{k}^{(l)}(Q)$, there exists $\left\{ m_{k}^{\dag
\dag (l)}\in \mc{D}_{k}^{(j)}(Q):k\in \lbrack K^{(l)}],l=j+1,...,J\right\} $
such that for the $m_{k}^{(l)}$, $k\in \lbrack K^{(l)}]$, $l\in \left[ j%
\right] $ defined by the algorithm using $\widetilde{\widetilde{f}}^{(j)}$
it holds that 
\[
\psi _{Q}^{1}=\sum_{l=j+1}^{J}\sum_{k\in \lbrack K^{(l)}]}m_{k}^{\dag \dag
(l)}+\sum_{l=1}^{j}\sum_{k\in \lbrack K^{(l)}]}m_{k}^{(l)}.
\]

To show that a solution $\widetilde{f}^{(j)}$ to \eqref{eq:algorithm-op-eq}
in $\left\{ 0\right\} \cup \sum_{l=j+1}^{J}\sum_{k\in \lbrack K^{(l)}]}\mc{D}%
_{k}^{(l)}(Q)\,$\ exists, let $\widetilde{f}^{(j)}\coloneqq%
\sum_{l=j+1}^{J}\sum_{k\in \lbrack K^{(l)}]}m_{k}^{\dag (l)}$. Clearly $%
\widetilde{f}^{(j)}\in \sum_{l=j+1}^{J}\bigoplus_{k\in \lbrack K^{(l)}]}%
\mc{D}_{k}^{(j)}(Q)$. Then, by the inductive assumption $\widetilde{f}%
^{(j)}=\psi _{Q}^{1}-\sum_{k=1}^{K^{(j)}}m_{k}^{\dag
(j)}-\sum_{l=1}^{j-1}\sum_{k=1}^{K^{(l)}}{m}_{k}^{(l)}$ and consequently,   
\begin{align*}
\Pi \left[ \left. \widetilde{f}^{(j)}\right\vert \left\{
\bigoplus_{k=1}^{K^{(j)}}\mc{D}_{k}^{(j)}(Q)\right\} ^{\perp }\right] =& \Pi %
\left[ \psi _{Q}^{1}-\sum_{k=1}^{K^{(j)}}m_{k}^{\dag
(j)}-\sum_{l=1}^{j-1}\sum_{k=1}^{K^{(l)}}{m}_{k}^{(l)}\left\vert \left\{
\bigoplus_{k=1}^{K^{(j)}}\mc{D}_{k}^{(j)}(Q)\right\} ^{\perp }\right. \right]
\\
=& \Pi \left[ \psi _{Q}^{1}-\sum_{l=1}^{j-1}\sum_{k=1}^{K^{(l)}}{m}%
_{k}^{(l)}\left\vert \left\{ \bigoplus_{k=1}^{K^{(j)}}\mc{D}%
_{k}^{(j)}(Q)\right\} ^{\perp }\right. \right] .
\end{align*}%
This shows that there exists $\widetilde{f}^{(j)}$ solving %
\eqref{eq:algorithm-op-eq} such that $\widetilde{f}^{(j)}$ $\in \left\{
0\right\} \cup \sum_{l=j+1}^{J}\sum_{k\in \lbrack K^{(l)}]}\mc{D}%
_{k}^{(l)}(Q)$.

Now, take $\widetilde{\widetilde{f}}^{(j)}\in \left\{ 0\right\} \cup
\sum_{l=j+1}^{J}\sum_{k\in \lbrack K^{(l)}]}\mc{D}_{k}^{(l)}(Q)$ to be an
arbitrary solution to \eqref{eq:algorithm-op-eq}. Let $m_{k}^{(j)}\coloneqq%
\Pi \left[ \left. \psi _{Q}^{1}-\widetilde{\widetilde{f}}^{(j)}-%
\sum_{l=1}^{j-1}\sum_{k=1}^{K^{(l)}}{m}_{k}^{(l)}\right\vert \mc{D}%
_{k}^{(j)}(Q)\right] $ for each $k\in \lbrack K^{(j)}],$ be the functions
defined in the algorithm for step $j$. Then $\psi _{Q}^{1}-\widetilde{%
\widetilde{f}}^{(j)}-\sum_{l=1}^{j-1}\sum_{k=1}^{K^{(l)}}{m}_{k}^{(l)}\in
\bigoplus_{k=1}^{K^{(j)}}\mc{D}_{k}^{(j)}(Q)$ because 
\begin{align*}
\Pi \left[ \psi
_{Q}^{1}-\widetilde{\widetilde{f}}^{(j)}-\sum_{l=1}^{j-1}\sum_{k=1}^{K^{(l)}}%
{m}_{k}^{(l)}\left\vert \left\{ \bigoplus_{k=1}^{K^{(j)}}\mc{D}%
_{k}^{(j)}(Q)\right\} ^{\perp }\right. \right] =0    
\end{align*}
by $\widetilde{\widetilde{%
f}}^{(j)}$ being a solution of \eqref{eq:algorithm-op-eq}. Additionally, $%
\sum_{k=1}^{K^{(j)}}m_{k}^{(j)}=\psi _{Q}^{1}-\widetilde{\widetilde{f}}%
^{(j)}-\sum_{l=1}^{j-1}\sum_{k=1}^{K^{(l)}}{m}_{k}^{(l)}$ because $\mc{D}%
_{k}^{(j)}(Q)$ and $\mc{D}_{k^{\prime }}^{(j)}(Q)$ are orthogonal for $%
k\not=k^{\prime }$. Finally, we note that we can write $\widetilde{%
\widetilde{f}}^{(j)}=\sum_{l=j+1}^{J}\sum_{k\in \lbrack K^{(l)}]}m_{k}^{\dag
\dag (l)}$ for some $m_{k}^{\dag \dag (l)}\in \mc{D}_{k}^{(l)}(Q)$, $k\in
\lbrack K^{(l)}]$, $l\in \{j+1,\dots ,J\}$ because $\widetilde{\widetilde{f}}%
^{(j)}$ is in $\left\{ 0\right\} \cup \sum_{l=j+1}^{J}\sum_{k\in \lbrack
K^{(l)}]}\mc{D}_{k}^{(l)}(Q)$. Thus, 
\begin{align*}
\psi _{Q}^{1}=& \widetilde{\widetilde{f}}^{(j)}+\sum_{l=1}^{j}\sum_{k\in
\lbrack K^{(l)}]}m_{k}^{(l)} \\
=& \sum_{l=j+1}^{J}\sum_{k\in \lbrack K^{(l)}]}m_{k}^{\dag \dag
(l)}+\sum_{l=1}^{j}\sum_{k\in \lbrack K^{(l)}]}m_{k}^{(l)}
\end{align*}%
This concludes the proof for the inductive step $j\Rightarrow j+1$.

\textbf{Proof of part \ref{item:compute-ifs-algorithm1} ($\Leftarrow$):}

Suppose \cref{algorithm:decompose} \texttt{DECOMPOSE} returns FAIL for
all influence functions $\psi _{Q}^{1}$ of $\psi $ at $Q$ in model $\mc{Q}$.
The arguments in the proof of the forward direction of part \ref%
{item:compute-ifs-algorithm1} of this lemma show that the algorithm can only
return fail at $j=1$. Suppose the algorithm returned FAIL at $j=1$. Then,
there does not exists $\widetilde{f}^{(1)}\in
\sum_{l=2}^{J}\bigoplus_{k=1}^{K^{(j)}}\mc{D}_{k}^{(j)}(Q)$ such that 
\[
\Pi \left[ \widetilde{f}^{(j)}-\psi _{Q}^{1}\left\vert \left\{
\bigoplus_{k=1}^{K^{(j)}}\mc{D}_{k}^{(j)}(Q)\right\} ^{\perp }\right. \right]
=0
\]%
for any influence function $\psi _{Q}^{1}$ of $\psi $ at $Q$ in model $\mc{Q}
$. Hence, there cannot exist $\{m_{k}^{(l)}\in \mc{D}_{k}^{(j)}(Q):k\in
\lbrack K^{(l)}],l\in \lbrack J]\}$ such that $\psi _{Q}^{1}=\sum_{l\in
\lbrack J]}\sum_{{k\in \lbrack K^{(l)}]}}m_{k}^{(l)}$ for otherwise $%
\widetilde{f}^{(1)}=\sum_{l=2}^{J}\sum_{k\in \lbrack K^{(l)}]}m_{k}^{(l)}$
would solve the above equation. As this does not exist for any ideal data
influence function $\psi _{Q}^{1}$, it follows from \cref%
{lemma:pathwise-differentiabilty}, $\varphi $ cannot pathwise differentiable
at $P$ in model $\mc{P}$.

\textbf{Proof of part \ref{item:compute-ifs-algorithm2}} Suppose that %
\cref{algorithm:decompose} \texttt{DECOMPOSE} returned a collection $%
\{m_{k}^{(l)}\in \mc{D}_{k}^{(j)}(Q):k\in \lbrack K^{(l)}],l\in \lbrack J]\}$.
Examining the algorithm at the $j=J$ step, we see that $\widetilde{f}^{(J)}=0
$ solves \eqref{eq:algorithm-op-eq}. Hence, the functions $m_{k}^{(J)}$
defined in the algorithm for step $J\,\ $satisfy $m_{k}^{(J)}=\Pi \left[
\psi _{Q}^{1}-\sum_{l=1}^{J-1}\sum_{k=1}^{K^{(J)}}m_{k}^{(l)}\left\vert %
\mc{D}_{k}^{(J)}(Q)\right. \right] $. Identical arguments as in the proof of
part \ref{item:compute-ifs-algorithm1} of this lemma show that $\sum_{k\in
\lbrack K^{(J)}]}m_{k}^{(J)}=\psi
_{Q}^{1}-\sum_{l=1}^{j-1}\sum_{k=1}^{(K^{(j)}}m_{k}^{(l)}$. Hence, 
\[
\psi _{Q}^{1}=\sum_{l\in \lbrack J]}\sum_{k\in \lbrack K^{(j)}]}m_{k}^{(l)}.
\]%
Then, by part \ref{item:influence functions1} of \cref{Theorem:influence
functions}, 
\[
\varphi _{P}^{1}(o)=\sum_{j=1}^{J}\frac{I(s=j)}{P(S=j)}\sum_{k=1}^{K^{(j)}}%
\frac{dQ}{dP(\cdot |S=j)}(\overline{z}_{k-1}^{(j)})m_{k}^{(j)}
\]%
is an influence function of $\varphi $ at $P$ in model $\mc{P}$.

\textbf{Proof of part \ref{item:compute-ifs-algorithm3}:}

Suppose the right-hand side of \eqref{observedIF-new} is an observed data
influence function of $\varphi$ at $P$ in model $\mc{P}$ for a collection $%
\{m_k^{(l)}\in \mc{D}_k^{(j)}(Q):k\in[K^{(l)}],l\in[J]\}$. By part \ref%
{item:influence functions1} of \cref{Theorem:influence functions}, $%
\psi^1_Q = \sum_{l\in[J]}\sum_{k\in[K^{(l)}]}m_k^{(l)}$ is an influence
function of $\psi$ at $Q$ in model $\mc{Q}$. Run \cref%
{algorithm:decompose} \texttt{DECOMPOSE} with input $\psi^1_Q$ and for each $%
j\in[J]$ choose the solution $\widetilde{f}^{(j)} = \sum_{l=j+1}\sum_{k\in[%
K^{(l)}]}m_k^{(j)}$ to equation \eqref{eq:algorithm-op-eq}. Then, the
algorithm will output $\{m_k^{(l)}\in \mc{D}_k^{(j)}(Q):k\in[K^{(l)}],l\in[J]\}$%
.
\end{proof}

\section{Relationship between this work and \texorpdfstring{\protect\cite{li_efficient_2023}}{Li and Luedtke 2023}}\label{sec:li-luedtke}

In this section, we examine the connection between the theory presented in \cite{li_efficient_2023} and the framework proposed in this work. As noted earlier, our work generalizes their model and extends their results. Specifically, 1) their models fall within the class of fused-data frameworks defined in our paper; 2) \cref{Theorem:influence functions}, applied to their models, characterizes the full set of observed data influence functions, thus extending their Theorem 2 which provides a subset of all observed data influence functions; 3) \cref{theorem:eif}, applied to their models, characterizes the observed data efficient influence function (EIF) thus extending their Corollary 1, as corrected in an upcoming Corrigendum (which we learned about through personal communication with the authors) and in the updated preprint version of that work (\cite{li_efficient_2025}), which characterizes the observed data EIF in a subclass of their models. 4) we provide a closed-form expression for the observed data EIF in an extension of the models assumed in their Corollary 1 that allows scenarios where alignment of conditional distributions may occur only on a subset of the supports of the conditioned variables.

We first demonstrate the models of \cite{li_efficient_2023} are a type of fused-data framework, which, adapting the terminology of \cite{qiu_efficient_2024}, we refer to as sequential conditional fused-data frameworks. In an abuse of notation, given $\widetilde{Q}, Q\in \mc{Q}$, we define $\frac{d\widetilde{Q}}{dQ}(z^{(j)}_k|\overline{z}^{(j)}_{k-1}) \coloneqq \frac{d\widetilde{Q}}{dQ}(\overline{z}^{(j)}_k)\big/\frac{d\widetilde{Q}}{dQ}(\overline{z}^{(j)}_{k-1})$.
\begin{sdefinition}
    \label{def:vi-fdf}
    We say a $(\mc{Q}, \mc{P}, \mc{C}, \psi, \varphi)$ is a sequential conditional fused-data framework with respect to $(Q_0, P_0)$ if it is a fused-data framework with respect to $(Q_0, P_0)$ and the following conditions hold:
    \begin{enumerate}
        \item \label{item:def:vi-fdf1} For all $j, j' \in [J]$, $k, k' \in \min(K^{(j)}, K^{(j)})$, it holds that $Z_k^{(j)} = Z_k^{(j')}$.
        \item \label{item:def:vi-fdf3} Given $Q, \widetilde{Q}\in \mc{Q}$, $k\in [K^{(j)}], j\in[J]$, define the law $\overline{Q} \in \mc{Q}$ whose Radon-Nikodym derivative with respect to $Q$ is
        \begin{align*}
            \frac{d\overline{Q}}{dQ}(w) = \frac{d\widetilde{Q}}{dQ}(z^{(j)}_k|\overline{z}^{(j)}_{k-1})^{I(\overline{z}^{(j)}_{k-1} \in \overline{\mc{Z}}^{(j)}_{k-1})}.
        \end{align*}
        Then $\overline{Q} \in \mc{Q}$. 
    \end{enumerate}
\end{sdefinition}

Implicit in condition \ref{item:def:vi-fdf1} of \cref{def:vi-fdf} is the assumption that it is possible to sort the sources in such a way that the observed data vectors are nested: $Z^{(1)} \subseteq Z^{(2)} \subseteq \dots \subseteq Z^{(J)}$. Letting $Z \coloneqq Z^{(J)}$, condition \ref{item:def:vi-fdf1} additionally requires that the aligned conditional and marginal distributions in each source correspond to a subset of the factors in the single factorization of the joint distribution of $Z$
\begin{align*}
    q(z) = q(z_{K^{(J)}}|\overline{z}_{K^{(J)}-1})\times \dots  \times q(z_2|z_1)\times q(z_1).  
\end{align*}
This factorization motivates our choice of terminology of "sequential conditional" fused-data frameworks. 

Condition \ref{item:def:vi-fdf3} of \cref{def:vi-fdf} is the assumption that aligned components are variation independent in the model $\mc{Q}$ in the sense given a law $Q \in \mc{Q}$, we may replace the $k, j$-th aligned conditional with the corresponding conditional of any other law in $\mc{Q}$ and remain in the model $\mc{Q}$. This condition is analogous to Condition 2 in \cite{li_efficient_2023} additionally accounting for the fact that alignment may occur only on a subset of the supports of the conditioned variables. 

We will now argue that the models assumed by \cite{li_efficient_2023} are particular sequential conditional fused-data frameworks. Specifically, \cite{li_efficient_2023} assume that the entire vector $Z$ is available in each source and for a subset of the sources, whose indices define the set $\mc{S}_k$, their conditional distribution of $Z_k|\overline{Z}_{k-1}$ aligns with the corresponding conditional distribution of $Q$.  Furthermore, they assume the aligned conditionals are variation independent. These assumptions then correspond to a sequential conditional fused-data framework in which for any $k \in[K^{(j)}], j\in[J]$, $\overline{\mc{Z}}^{(j)}_{k-1}$ is either empty or equal to $\textsf{Supp}(\overline{Z}_k;Q)$ and $K^{(1)} = K^{(2)} = \dots =K^{(J)}$. Note that with our definitions, the aforementioned set $\mc{S}_k$ is equal to $\{j \in [J]: \mc{Z}_{k-1}^{(j)} =  \textsf{Supp}(\overline{Z}_k;Q)\}$. Sequential conditional fused-data frameworks generalize the models of \cite{li_efficient_2023} because they allow for the possibility that only a subset of $Z$ is observed in every source and, more importantly, that the conditional distribution of $Z_k|\overline{Z}_{k-1}$ aligns with the corresponding conditional distribution of $Q$ potentially only for a subset of the support of $\overline{Z}_{k-1}$ under $Q$. We note that the first extension has no impact as inference is unchanged whether or not one considers observing nested subsets of $Z$. 

\cref{example:transporting} scenario (i) is an example of a sequential conditional fused-data framework. However, the fused-data frameworks for all other examples in this work do not qualify as sequential conditional frameworks.

With \cref{def:vi-fdf} established, we can now characterize the set of all observed data influence functions for a sequential conditional fused-data framework. To simplify notation, we first define $ \mc{J}_k \coloneqq \{ j \in [J] : k \leq K^{(j)} \} $ and  $K^{[J]} \coloneqq \max(K^{(1)}, \dots, K^{(J)})$. To better align with the notation used in \cite{li_efficient_2023}, many of the subsequent results are expressed in terms of the random vector $Z$ and its realization $z$, rather than the random vectors $ Z^{(j)}$ and realizations $z^{(j)}$.  
\begin{slemma}
    \label{lemma:li-ifs}
    Let $(\mc{Q}, \mc{P}, \mc{C}, \psi, \varphi)$ be a sequential conditional fused-data framework with respect to $(Q_0, P_0)$. Let $P\in \mc{P}$. Suppose there exists $Q\in \Phi(P;\mc{C})$ such that $(Q, P)$ is strongly aligned with respect to $\mc{C}$ and $\psi$ is pathwise differentiable at $Q$ in model $\mc{Q}$. Then,
    \begin{enumerate}
        \item\label{item:lemma:li-ifs1} $\varphi$ is pathwise differentiable at $P$ in model $\mc{P}$
        \item\label{item:lemma:li-ifs2} $\varphi^1_P$ is an influence function of $\varphi$ at $P$ in model $\mc{P}$ if and only if there exists an influence function $\psi^1_Q$ of $\psi$ at $Q$ in model $\mc{Q}$ such that 
        \begin{align*}
            \varphi^1_P(o) = \sum_{j\in[J]}\frac{I(s=j)}{P(S=j)}\sum_{k\in[K^{(j)}]}\frac{dQ}{dP(\cdot|S=j)}(\overline{z}_{k-1})m_k^{(j)}(\overline{z}_k)
        \end{align*}
        where $m_k^{(j)} \in \mc{D}_k^{(j)}(Q)$ for all $k\in[K^{(j)}]$, $j\in[J]$, and, or all $k\in [K^{[J]}]$,
        \begin{align}
            \label{eq:li-ifs-m-condition}
            \sum_{j\in\mc{J}_k}m_k^{(j)}(\overline{z}_k) = I\left(\overline{z}_{k-1}^{(j)} \in \bigcup_{j\in\mc{J}_k}\mc{\overline{Z}}^{(j)}_{k-1}\right) \left\{E_Q\left[\psi^1_Q(W)|\overline{z}_{k}\right] - E_Q\left[\psi^1_Q(W)|\overline{z}_{k-1}\right]\right\}
        \end{align}
    \end{enumerate}
\end{slemma}

In the special case in which the sequential conditional fused-data framework corresponds to the model in \cite{li_efficient_2023}, the assumption that there exists $Q \in \Phi(P; \mc{C}) $ such that $(Q, P)$ is strongly aligned with respect to $ \mc{C}$ is equivalent to Conditions 1 and 3 in \cite{li_efficient_2023}. Part \ref{item:lemma:li-ifs1} of the above lemma corresponds to Lemma 1 in \cite{li_efficient_2023}. Part \ref{item:lemma:li-ifs2} of the above lemma corresponds to Theorem 2 of \cite{li_efficient_2023}. In that theorem, the characterize set of observed data influence functions which coincides with the set of all observed data influence functions in the special case in which different sources align with disjoint factors of the ideal data likelihood factorization, but not otherwise.

The observed data influence functions provided in Theorem 2 of \cite{li_efficient_2023} are of the form
\begin{align}
    \label{eq:li-luedtke-if-form}
    \sum_{k=1}^{K^{([J])}}\frac{I(s \in \mc{S}_k)}{P(S \in \mc{S}_k)}\frac{dQ}{dP(\cdot|S \in \mc{S}_k)}(\overline{z}_{k-1})\left\{E_Q\left[\psi^1_Q(W)|\overline{z}_{k}\right] - E_Q\left[\psi^1_Q(W)|\overline{z}_{k-1}\right]\right\}.
\end{align}
In particular, if the ideal data model is non-parametric, that theorem characterizes only a single observed data influence function. However, when there exists any pair of sources whose alignments agree on at least one factor of the ideal data likelihood, the observed data model $\mc{P}$ is strictly semiparametric. Hence, there are infinitely many observed data influence functions. The preceding lemma, \cref{lemma:li-ifs}, characterizes them all. Note that regardless of whether or not the ideal data model is non-parametric, \eqref{eq:li-luedtke-if-form} is equal to 
\begin{align*}
    \sum_{j\in[J]}\frac{I(s=j)}{P(S=j)}\sum_{k\in[K^{(j)}]}\frac{dQ}{dP(\cdot|S=j)}(\overline{z}_{k-1})m_k^{(j)}(\overline{z}_k)
\end{align*}
where 
\begin{align*}
    m_k^{(j)}(\overline{z}_k) \coloneqq & \frac{dP\left(\cdot|S=j\right)}{dP\left(\cdot|S\in \mc{S}_k\right)}(\overline{z}_{k-1})\frac{P(S=j)}{P\left(S\in \mc{S}_k\right)}\\
    &\times I(\overline{\mc{Z}}_{k-1}^{(j)} = \textsf{Supp}(\overline{Z}_{k-1}^{(j)};Q))\left\{E_Q\left[\psi^1_Q(W)|\overline{z}_{k}\right] - E_Q\left[\psi^1_Q(W)|\overline{z}_{k-1}\right]\right\}
\end{align*}
for each $k\in K^{(j)}, j\in[J]$. It is easily checked that these $m_k^{(j)}$'s satisfy \eqref{eq:li-ifs-m-condition} of part \ref{item:lemma:li-ifs2} of \cref{lemma:li-ifs} when each $\overline{\mc{Z}}_{k-1}^{(j)}$ is either empty or the support of $\overline{Z}_{k-1}^{(j)}$ under $Q$. This proves that the class of influence functions in Theorem 2 of \cite{li_efficient_2023} is a subset of the class characterized in part \ref{item:lemma:li-ifs2} of \cref{lemma:li-ifs}. We note that when $\mc{Q}$ is non-parametric, the observed data influence function provided by Theorem 2 of \cite{li_efficient_2023} is efficient even when there exists any pair of sources whose alignments agree on at least one factor of the ideal data likelihood.

Part \ref{item:lemma:li-ifs2} of \cref{lemma:li-ifs} follows by applying \cref{lemma:compute-ifs-algorithm} to the special case of sequential conditional fused-data frameworks. For these frameworks, it happens that the operator equations that must be solved in \cref{algorithm:decompose} (\texttt{DECOMPOSE}) have closed-form solutions which yields the direct computation of observed data influence functions from ideal data influence functions established in part \ref{item:lemma:li-ifs2} of \cref{lemma:li-ifs}.

That the operator equations in \cref{algorithm:decompose} (\texttt{DECOMPOSE}) have closed form solutions follows because, for sequential conditional fused-data frameworks, it holds that for all \( k \in [K^{[J]}] \), \( j, j' \in \mc{J}_k \), and \( m_k^{(j)} \in \mc{D}_k^{(j)}(Q) \), 
\[
\Pi\left[m_k^{(j)} \middle\vert \mc{D}_k^{(j')}(Q)\right](\overline{z}_k) = I(\overline{z}_{k-1} \in \mc{Z}_{k-1}^{(j')}) m_k^{(j)}(\overline{z}_k).
\]
and additionally, for all \( k, k' \in [K^{[J]}] \) such that \( k \neq k' \), \( j \in \mc{J}_k \), \( j' \in \mc{J}_{k'} \), and \( m_k^{(j)} \in \mc{D}_k^{(j)}(Q) \),
\[
\Pi\left[m_k^{(j)} \middle\vert \mc{D}_{k'}^{(j')}(Q)\right] = 0.
\]  
These properties follow from part \ref{item:def:vi-fdf1} of \cref{def:vi-fdf}.

In Corollary 1 of \cite{li_efficient_2023}, as corrected in an upcoming Corrigendum (which we learned about through personal communication with the authors) and in the updated preprint version of that work (\cite{li_efficient_2025}), the authors establish that the observed data 
efficient influence function under sequential conditional fused-data frameworks has a closed-form expression for certain ideal data models \( \mc{Q} \). Here, we extend these results to allow for cases where the aligned conditional distributions may align only on a subset of the support of \( \overline{Z}_{k-1}^{(j)} \) under \( Q \), for each \( k \in [K^{(j)}] \) and \( j \in [J] \), rather than on the entire support of \( \overline{Z}_{k-1}^{(j)} \) under \( Q \). In what follows we define $\mc{S}_k(\overline{z}_{k-1}) \coloneqq \{j \in \mc{J}_k: \overline{z}_{k-1} \in \overline{\mc{Z}}_{k-1}^{(j)}\}$. 
\begin{scorollary}
    \label{cor:li-eif}
    Let $(\mc{Q}, \mc{P}, \mc{C}, \psi, \varphi)$ be a sequential conditional fused-data framework with respect to $(Q_0, P_0)$. Let $P\in \mc{P}$. Suppose there exists $Q\in \Phi(P;\mc{C})$ such that $(Q, P)$ is strongly aligned with respect to $\mc{C}$ and $\psi$ is pathwise differentiable at $Q$ in model $\mc{Q}$. Let $\psi^1_{Q, eff}$ be the efficient influence function of $\psi$ at $Q$ in model $\mc{Q}$.
    Suppose
    \begin{align}
        \label{eq:cor:li-eif-corrigendum}&\frac{1}{P(S \in \mc{S}_k(\overline{z}_{k-1}))}\frac{dQ}{dP(\cdot|S \in \mc{S}_k(\overline{z}_{k-1}))}(\overline{z}_{k-1})\\
        \nonumber &\times I\left(\overline{z}_{k-1}\in \bigcup_{j\in \mc{J}_k}\overline{\mc{Z}}_{k-1}^{(j)}\right)\left\{E_Q\left[\psi^1_{Q, eff}(W)|\overline{z}_{k}\right] - E_Q\left[\psi^1_{Q, eff}|\overline{z}_{k-1}\right]\right\} \in \mc{T}(Q;\mc{Q})
    \end{align}
    for all $k\in[K^{[J]}]$. Then, 
    \begin{align}
        \label{eq:li-eif}
        \varphi^1_{P, eff}(o) = \sum_{k\in[K^{[J]}]}&\frac{I(s \in \mc{S}_k(\overline{z}_{k-1}))}{P(S \in \mc{S}_k(\overline{z}_{k-1}))}\frac{dQ}{dP(\cdot|S \in \mc{S}_k(\overline{z}_{k-1}))}(\overline{z}_{k-1})\\
        \nonumber&\times \{E_Q\left[\psi^1_{Q, eff}(W)|\overline{z}_{k}\right] - E_Q\left[\psi^1_{Q, eff}(W)|\overline{z}_{k-1}\right]\}
    \end{align}
    is the efficient influence function of $\varphi$ at $P$ in $\mc{P}$.
\end{scorollary}
It is unclear whether \eqref{eq:cor:li-eif-corrigendum} is necessary to guarantee in general that the efficient observed data influence function has a closed form. In the special case where $\overline{\mc{Z}}_{k-1}^{(j)}$ is either empty or equal to the support of $\overline{Z}_{k-1}^{(j)}$ under $Q$, a discussion of when this condition holds can be found in the aforementioned forthcoming Corrigendum to \cite{li_efficient_2023}.

\subsection{Proofs of \texorpdfstring{\cref{sec:li-luedtke}}{this section}}

In the following proofs we let $\mc{I}_{k-1} \coloneqq \bigcup_{j\in \mc{J}_k} \overline{\mc{Z}}_{k-1}^{(j)}$. For a subset $\mc{E}$ of $L^2(Q)$ and a closed linear subspace $\mc{D}$ of $L^2_0(Q)$ we let $\Pi[\mc{E}|\mc{D}] \coloneqq \{\Pi[e|\mc{D}]:e\in \mc{E}\}$. We first state a useful lemma. 
\begin{slemma}
    \label{lemma:li-luedtke-tangent}
    Let $(\mc{Q}, \mc{P}, \mc{C}, \psi, \varphi)$ be a sequential conditional fused-data framework with respect to $(Q_0, P_0)$. Let $P\in \mc{P}$. Suppose there exists $Q\in \Phi(P;\mc{C})$ such that $(Q, P)$ is strongly aligned with respect to $\mc{C}$.  Then,
    \begin{enumerate}
        \item \label{lemma:li-luedtke-tangent:item:1} 
        \begin{align*}
            \Pi\left[\mc{T}(Q, \mc{Q})\left|\left(\sum_{j\in[J]}\bigoplus_{k\in[K^{(j)}]}\mc{D}_k^{(j)}(Q)\right)^\perp\right.\right] &\subseteq \mc{T}(Q, \mc{Q})
        \end{align*}
        \item \label{lemma:li-luedtke-tangent:item:2} For all $k\in[K^{(j)}], j\in[J]$,
        \begin{align*}
            \Pi\left[\mc{T}(Q, \mc{Q})\left|\mc{D}_k^{(j)}(Q)\right.\right] &\subseteq \mc{T}(Q, \mc{Q})
        \end{align*}
        \item \label{lemma:li-luedtke-tangent:item:3} For all $k \in [K^{[J]}]$,
        \begin{align*}
            \Pi\left[\mc{T}(Q, \mc{Q})\left|\sum_{j\in\mc{J}_k}\mc{D}_k^{(j)}(Q)\right.\right] &\subseteq \mc{T}(Q, \mc{Q})
        \end{align*}
    \end{enumerate}
\end{slemma}

\begin{proof}[Proof of \cref{lemma:li-luedtke-tangent}]
    We begin with part \ref{lemma:li-luedtke-tangent:item:1}. Let $\{\widetilde{Q}_t:t\in(-\varepsilon,\varepsilon)\}$ be a regular parametric submodel in $\mc{Q}$ with $\widetilde{Q}_t\vert_{t=0} = Q$ and score $\widetilde{h}^{(Q)}$ at $t=0$. For every $t \in (-\varepsilon, \varepsilon)$, let $Q_t$ be a law on $W$ such that
    \begin{align*}
        \frac{dQ_t}{dQ}(z) = \prod_{k\in[K^{[J]}]}\frac{d\widetilde{Q}_t}{dQ}(z_k\vert\overline{z}_{k-1})^{1-I(\overline{z}_{k-1} \in \mc{I}_{k-1})}.
    \end{align*}
    That is, the conditional distribution of $Z_k|\overline{Z}_{k-1} = \overline{z}_{k-1}$ under $Q_t$ equals that conditional distribution under $Q$ whenever $\overline{z}_{k-1} \in \mc{I}_{k-1}$, i.e. whenever there is at least one source such that the conditional distribution of $Z_k|\overline{Z}_{k-1} = \overline{z}_{k-1}$ aligns with that of $Q$ at $\overline{z}_{k-1}$. Furthermore, the conditional distribution of $Z_k|\overline{Z}_{k-1} = \overline{z}_{k-1}$ under $Q_t$ equals the conditional distribution under $\widetilde{Q}_t$ otherwise.  
    
    By part \ref{item:def:vi-fdf3} of \cref{def:vi-fdf}, $\{Q_{t}:t \in (-\varepsilon, \varepsilon)\}$ is also a regular parametric submodel in $\mc{Q}$ with $Q_t\vert_{t=0} = Q$. Recall from the proof of \cref{lemma:score-operator} that the score of $Q_t$ at $t=0$ will be the element $h^{(Q)}$ in $ L^2_0(Q)$ such that $h^{(Q)}/2$ is the Frechet derivative at $t=0$ of the map $t \mapsto \sqrt{\frac{dQ_t}{dQ}}$ as a mapping from $(-\varepsilon, \varepsilon)$ to $L^2(Q)$. Let $r_{k,t}(\overline{z}_k) \coloneqq \sqrt{\frac{d\widetilde{Q}_t}{dQ}(z_k\vert\overline{z}_{k-1})^{1-I(\overline{z}_{k-1} \in \mc{I}_{k-1})}}$ for $k\in [K^{[J]}]$. Then, $r_{k,t} \in L^2(Q)$ for each $k \in [K^{[J]}]$. By noticing that $t \mapsto \sqrt{\frac{dQ_t}{dQ}} = \prod_{k\in[\widetilde{K}]}r_{k,t}$ and applying the product rule similarly to the proof of \cref{lemma:score-operator}, we conclude that $Q_t$ is quadratic mean differentiable with score at $t=0$ given by 
    \begin{align*}
        h^{(Q)}(z) =& \sum_{k\in[K^{[J]}]}(1-I(\overline{z}_{k-1}\in \mc{I}_{k-1}))\{E_Q[\widetilde{h}^{(Q)}(Z)|\overline{z}_k] - E_Q[\widetilde{h}^{(Q)}(Z)|\overline{z}_{k-1}]\}\\
        =&\Pi\left[\widetilde{h}^{(Q)}\left|\left(\bigoplus_{k\in[K^{[J]}]}\sum_{j\in\mc{J}_k}\mc{D}_k^{(j)}(Q)\right)^\perp\right.\right](z)\\
        =&\Pi\left[\widetilde{h}^{(Q)}\left|\left(\sum_{j\in[J]}\bigoplus_{k\in[K^{(j)}]}\mc{D}_k^{(j)}(Q)\right)^\perp\right.\right](z).
    \end{align*}
    As $\{\widetilde{Q}_t:t\in (-\varepsilon, \varepsilon)\}$ was arbitrary this proves part \ref{lemma:li-luedtke-tangent:item:1}. 

    Turn now to the proof of part \ref{lemma:li-luedtke-tangent:item:2}. Fix $j\in[J]$ and $k\in[K^{(j)}]$ and let $\{\widetilde{Q}_t:t\in (-\varepsilon, \varepsilon)\}$ be a regular parametric submodel in $\mc{Q}$ with $\widetilde{Q}_t\vert_{t=0} = Q$ and score $\widetilde{h}^{(Q)}$ at $t=0$. For every $t \in (-\varepsilon, \varepsilon)$, let $Q_t$ be a law on $W$ such that
    \begin{align*}
        \frac{dQ_t}{dQ}(z) = \frac{d\widetilde{Q}_t}{dQ}(z_k\vert\overline{z}_{k-1})^{I(\overline{z}_{k-1} \in \mc{Z}_k^{(j)})}.
    \end{align*}
    That is, the conditional distribution of $Z_k|\overline{Z}_{k-1} = \overline{z}_{k-1}$ under $Q_t$ equals that conditional distribution under $\widetilde{Q}_t$ whenever $\overline{z}_{k-1} \in \mc{Z}^{(j)}_{k}$ , i.e. whenever the conditional distribution of $Z_k|\overline{Z}_{k-1} = \overline{z}_{k-1}$ in source $j$ aligns with that of $Q$ at $\overline{z}_{k-1}$. Furthermore, the conditional distribution of $Z_k|\overline{Z}_{k-1} = \overline{z}_{k-1}$ under $Q_t$ equals the conditional distribution under $Q$ otherwise.  

    By part \ref{item:def:vi-fdf3} of \cref{def:vi-fdf}, $\{Q_t:t\in (-\varepsilon, \varepsilon)\}$ is also a regular parametric submodel in $\mc{Q}$ with $Q_t\vert_{t=0} = Q$. Similar arguments as above show that $Q_t$ is quadratic mean differentiable with score at $t=0$ given by 
    \begin{align*}
        h^{(Q)}(\overline{z}_k) = I(\overline{z}_{k-1} \in \overline{\mc{Z}}_k^{(j)})\{E_Q[\widetilde{h}^{(Q)}(Z)|\overline{z}_k] - E_Q[\widetilde{h}^{(Q)}(Z)|\overline{z}_{k-1}]\}.
    \end{align*}
    As $k$, $j$, $\{\widetilde{Q}_t:t\in (-\varepsilon, \varepsilon)\}$ were arbitrary this proves part \ref{lemma:li-luedtke-tangent:item:2}.

    Turn now to the proof of part \ref{lemma:li-luedtke-tangent:item:3}. Fix $k\in [K^{[J]}]$ and note that 
    \begin{align*}
        \sum_{j \in \mc{J}_k}\mc{D}_k^{(j)}(Q) =& \bigoplus_{i =1}^{|\mc{J}_k|}\left\{I\left(\overline{z}_k \in \overline{\mc{Z}}_k^{(j_i)}\setminus \left(\bigcup_{l < i} \overline{\mc{Z}}_k^{(j_l)}\right)\right)\{d(\overline{z}_k) - E_Q[d(\overline{Z}_{k})|\overline{z}_{k}]\}:d\in L^2(Q)\right\}
    \end{align*}
    where $\{j_1, \dots, j_{\mc{J}_k}\} = \mc{J}_k$. It then follows that 
    \begin{align*}
        &\Pi\left[\mc{T}(Q, \mc{Q})\left|\sum_{j\in\mc{J}_k}\mc{D}_k^{(j)}(Q)\right.\right]\\
        =&\Pi\left[\mc{T}(Q, \mc{Q})\left|\bigoplus_{i =1}^{|\mc{J}_k|}\left\{I\left(\overline{z}_k \in \overline{\mc{Z}}_k^{(j_i)}\setminus \left(\bigcup_{l < i} \overline{\mc{Z}}_k^{(j_l)}\right)\right)\{d(\overline{z}_k) - E_Q[d(\overline{Z}_{k})|\overline{z}_{k}]\}:d\in L^2(Q)\right\}\right.\right]\\
        =&\bigoplus_{i =1}^{|\mc{J}_k|}\Pi\left[\mc{T}(Q, \mc{Q})\left|\left\{I\left(\overline{z}_k \in \overline{\mc{Z}}_k^{(j_i)}\setminus \left(\bigcup_{l < i} \overline{\mc{Z}}_k^{(j_l)}\right)\right)\{d(\overline{z}_k) - E_Q[d(\overline{Z}_{k})|\overline{z}_{k}]\}:d\in L^2(Q)\right\}\right.\right]\\
        \subseteq& \sum_{j\in \mc{J}} \Pi\left[\mc{T}(Q, \mc{Q})\left|\mc{D}_k^{(j)}(Q)\right.\right]\\
         \subseteq& \sum_{j\in \mc{J}} \mc{T}(Q; \mc{Q}) \\
        \subseteq& \mc{T}(Q;\mc{Q}).
    \end{align*}
    As $k$ was arbitrary this completes the proof. 
\end{proof}

\begin{proof}[Proof of \cref{lemma:li-ifs}]
    $ $\newline
    \textbf{Proof of part \ref{item:lemma:li-ifs1}:}
    To prove this result, we show the efficient ideal data influence function $\psi^1_{Q, eff}$ is an element of $\sum_{j\in[J]}\bigoplus_{k\in[K^{(j)}]}\mc{D}_k^{(j)}(Q)$.  Then, the pathwise differentiability of $\varphi$ follows by \cref{lemma:pathwise-differentiabilty}. 
    
    Let $\{\widetilde{Q}_t:t\in(-\varepsilon,\varepsilon)\}$ be a regular parametric submodel in $\mc{Q}$ with $\widetilde{Q}_t\vert_{t=0} = Q$ and score $\widetilde{h}^{(Q)}$ at $Q$. For every $t \in (-\varepsilon, \varepsilon)$, let $Q_t$ be a law on $W$ such that
    \begin{align*}
        \frac{dQ_t}{dQ}(z) = \prod_{k\in[K^{[J]}]}\frac{d\widetilde{Q}_t}{dQ}(z_k\vert\overline{z}_{k-1})^{1-I(\overline{z}_{k-1} \in \mc{I}_{k-1})}.
    \end{align*}
    That is, the conditional distribution of $Z_k|\overline{Z}_{k-1} = \overline{z}_{k-1}$ under $Q_t$ equals that conditional distribution under $Q$ whenever $\overline{z}_{k-1} \in \mc{I}_{k-1}$, i.e. whenever there is at least one source such that the conditional distribution of $Z_k|\overline{Z}_{k-1} = \overline{z}_{k-1}$ aligns with that of $Q$ at $\overline{z}_{k-1}$. Furthermore, the conditional distribution of $Z_k|\overline{Z}_{k-1} = \overline{z}_{k-1}$ under $Q_t$ equals the conditional distribution under $\widetilde{Q}_t$ otherwise.  
    
    By part \ref{item:def:vi-fdf3} of \cref{def:vi-fdf}, $\{Q_t:t\in (-\varepsilon, \varepsilon)\}$ is also a regular parametric submodel in $\mc{Q}$ with $Q_t\vert_{t=0} = Q$. From the proof of \cref{lemma:li-luedtke-tangent} we have that $\{Q_t:t\in(-\varepsilon, \varepsilon)\}$ is quadratic mean differentiable with score at $t=0$ given by 
    \begin{align*}
        h^{(Q)}(z) =&\Pi\left[\widetilde{h}^{(Q)}\left|\left(\sum_{j\in[J]}\bigoplus_{k\in[K^{(j)}]}\mc{D}_k^{(j)}(Q)\right)^\perp\right.\right](z)
    \end{align*}
    Now, we have that 
    \begin{align*}
        \frac{d}{dt}\psi(Q_t)\big\vert_{t=0} = \frac{d}{dt}\psi(Q)\big\vert_{t=0} = 0.
    \end{align*}
    because the conditional distribution $Z_k^{(j)}\vert \overline{Z}_{k-1} = \overline{z}_{k-1}^{(j)}$ for $z_{k-1}^{(j)} \in \overline{\mc{Z}}_k^{(j)}$ under $Q_t$ is the same for all $t\in(-\varepsilon,\varepsilon)$, $k\in[K^{(j)}]$, $j\in[J]$ and, because by \cref{as:identification}, $\psi$ depends only on these conditional distributions. On the other hand,
    \begin{align}
        \label{eq:li-luedtke-pd-proof-orthogonal}
        \frac{d}{dt}\psi(Q_t)\big\vert_{t=0} = \langle \psi^1_Q, h^{(Q)} \rangle_Q = 0
    \end{align}
    for every influence function $\psi^1_Q$ for $\psi$ at $Q$ in $\mc{Q}$ because $\psi$ is pathwise differentiable. In particular, this holds for the efficient influence function $\psi^1_{Q, eff} \in \mc{T}(Q, \mc{Q})$. Now, it follows there exists a collection of submodels in $\mc{Q}$ through $Q$ constructed as above with scores at $t=0$ that are dense in $\Pi\left[\mc{T}(Q, \mc{Q})\left|\left(\sum_{j\in[J]}\bigoplus_{k\in[K^{(j)}]}\mc{D}_k^{(j)}(Q)\right)^\perp\right.\right]$ because $\widetilde{Q}_t$ was arbitrary.  Additionally, \eqref{eq:li-luedtke-pd-proof-orthogonal} holds for each of these submodels and their corresponding scores. Hence, we have that $\psi^1_{Q, eff}$ is orthogonal to a dense subset of $\Pi\left[\mc{T}(Q, \mc{Q})\left|\left(\sum_{j\in[J]}\bigoplus_{k\in[K^{(j)}]}\mc{D}_k^{(j)}(Q)\right)^\perp\right.\right]$, and that $\psi^1_{Q, eff} \in \mc{T}(Q, \mc{Q})$. But this means that $\psi^1_{Q, eff} \in \Pi\left[\mc{T}(Q, \mc{Q})\left|\sum_{j\in[J]}\bigoplus_{k\in[K^{(j)}]}\mc{D}_k^{(j)}(Q)\right.\right]$, and in particular that $\psi^1_{Q, eff} \in \sum_{j\in[J]}\bigoplus_{k\in[K^{(j)}]}\mc{D}_k^{(j)}(Q)$. By \cref{lemma:pathwise-differentiabilty}, this implies that $\varphi$ is pathwise differentiable at $P$ in model $\mc{P}$.

    \textbf{Proof of part \ref{item:lemma:li-ifs2} ($\Rightarrow$):}
    
    Suppose $\varphi^1_P$ is an influence function of $\varphi$ at $P$ in model $\mc{P}$. Then there exists $\psi^1_Q$ an influence function of $\psi$ at $Q$ in $\mc{Q}$ such that $\varphi^1_P$ and $\psi^1_Q$ are decomposed as in part \ref{item:influence functions1} of \cref{Theorem:influence functions}. Let $\left\{m_k^{(j)} \in \mc{D}_k^{(j)}(Q):k\in[K^{(j)}], j\in[J]\right\}$ be the collection such that $\varphi^1_P(o)$ satisfies \eqref{newIF} for this collection and
    \begin{align*}
        \psi^1_Q &= \sum_{j\in[J]}\sum_{k\in[K^{(j)}]}m_k^{(j)}\\
        &=\sum_{k\in[K^{[J]}]}\sum_{j\in\mc{J}_k}m_k^{(j)}
    \end{align*} 
     where the second equality follows by definition of $\mc{J}_k$. The spaces $\sum_{j\in\mc{J}_k}\mc{D}_k^{(j)}(Q)$ and $\sum_{j\in\mc{J}_{k'}}\mc{D}_{k'}^{(j)}$ are orthogonal for all $k, k'\in [K^{[J]}]$ such that $k\not=k'$ by part \ref{item:def:vi-fdf1} of \cref{def:vi-fdf} and 
     \begin{align*}
         \Pi\left[f\Bigg|\sum_{j\in \mc{J}_k} \mc{D}_k^{(j)}(Q)\right](w) = I(\overline{z}_{k-1} \in \mc{I}_{k-1})\{E_Q[f(W)|\overline{z}_{k}] - E_Q[f(W)|\overline{z}_{k-1}]\}.
     \end{align*}
     for all $f \in L^2(Q)$. This implies that 
    \begin{align*}
        \psi^1_Q(z) = \sum_{k\in [K^{[J]}]}I\left(\overline{z}_{k-1} \in \mc{I}_{k-1}\right)\{E_Q[\psi^1_Q(W)|\overline{z}_k] - E_Q[\psi^1_Q(W)|\overline{z}_{k-1}]\}
    \end{align*}
    But this in turn implies that 
    \begin{align*}
        I\left(\overline{z}_{k-1} \in \mc{I}_{k-1}\right)\{E_Q[\psi^1_Q(Z)|\overline{z}_k] - E_Q[\psi^1_Q(Z)|\overline{z}_{k-1}]\} = \sum_{j\in{J}_k}m_k^{(j)}
    \end{align*}
    by orthogonality of the subspaces $\sum_{j\in\mc{J}_k}\mc{D}_k^{(j)}(Q)$ and $\sum_{j\in\mc{J}_{k'}}\mc{D}_{k'}^{(j)}(Q)$ for $k\not=k'$ and the fact that $\sum_{j\in\mc{J}_k}m_k^{(j)} \in \sum_{j\in\mc{J}_k}\mc{D}_k^{(j)}(Q)$. This completes this direction of the proof. 

    \textbf{Proof of part \ref{item:lemma:li-ifs2} ($\Leftarrow$):}
    Let $\psi^1_Q$ be an influence function of $\psi$ at $Q$ in model $\mc{Q}$. Let
    \begin{align*}
        \varphi^1_P(o) \coloneqq  \sum_{j\in[J]}\frac{I(s=j)}{P(S=j)}\sum_{k\in[K^{(j)}]}\frac{dQ}{dP(\cdot|S=j)}(\overline{z}_{k-1})m_k^{(j)}(\overline{z}_k)
    \end{align*}
    for some $\left\{m_k^{(j)}\in \mc{D}_k^{(j)}(Q):k\in[K^{(j)}], j\in[J]\right\}$ satisfying, for all $k\in [K^{[J]}]$,
    \begin{align}
        \label{eq:m_k-decomp-li-luedtke}
        \sum_{j\in\mc{J}_k}m_k^{(j)}(\overline{z}_k) = I\left(\overline{z}_{k-1} \in \mc{I}_{k-1}\right) \left\{E_Q\left[\psi^1_Q(W)|\overline{z}_{k}\right] - E_Q\left[\psi^1_Q(W)|\overline{z}_{k-1}\right]\right\}
    \end{align}
    
    Define 
    \begin{align*}
        \widetilde{\psi}^1_Q \coloneqq \sum_{k\in [K^{[J]}]}I\left(\overline{z}_{k-1} \in \mc{I}_{k-1}\right) \left\{E_Q\left[\psi^1_Q(W)|\overline{z}_{k}\right] - E_Q\left[\psi^1_Q(W)|\overline{z}_{k-1}\right]\right\}.
    \end{align*}
    We will next show that $\widetilde{\psi}^1_Q$ is an influence function of $\psi$ at $Q$ in model $\mc{Q}$ and that $\widetilde{\psi}^1_Q = \sum_{j\in[J]}\sum_{k\in[K^{(j)}]}m_k^{(j)}$. We will then conclude by part \ref{item:influence functions1} of \cref{Theorem:influence functions} that $\varphi^1_P$ is an influence function of $\varphi$ at $P$ in model $\mc{P}$. 
    
    Turn first to the proof that $\widetilde{\psi}^1_Q$ is an influence function of $\psi$ at $Q$ in model $\mc{Q}$. First note that $\widetilde{\psi}^1_Q = \sum_{k\in[K^{[J]}]}\Pi\left[\psi^1_Q\left|\sum_{j\in\mc{J}_k}\mc{D}_k^{(j)}(Q)\right.\right]$. Then, letting $\psi^1_{Q, eff}$ be the efficient influence function of $\psi$ at $Q$ in model $\mc{Q}$ we have that 
    \begin{align*}
        \widetilde{\psi}^1_Q =& \sum_{k\in[K^{[J]}]}\Pi\left[\psi^1_Q\left|\sum_{j\in\mc{J}_k}\mc{D}_k^{(j)}(Q)\right.\right]\\
        =&\sum_{k\in[\widetilde{K}]}\Pi\left[\psi^1_{Q, eff} + f\left|\sum_{j\in\mc{J}_k}\mc{D}_k^{(j)}(Q)\right.\right]\\
        =&\psi^1_{Q, eff} + \sum_{k\in[K^{[J]}]}\Pi\left[ f\left|\sum_{j\in\mc{J}_k}\mc{D}_k^{(j)}(Q)\right.\right]
    \end{align*}
    where the third equality follows from the proof of part \ref{item:lemma:li-ifs1} in which we showed that $\psi^1_{Q, eff} \in \bigoplus_{k\in[K^{[J]}]}\sum_{j\in\mc{J}_k}\mc{D}_k^{(j)}(Q)$. Then, $\widetilde{\psi}^1_Q$ will be an influence function of $\psi$ at $Q$ in model $\mc{Q}$ if $\Pi\left[ f\left|\sum_{j\in\mc{J}_k}\mc{D}_k^{(j)}(Q)\right.\right] \in \mc{T}(Q, \mc{Q})^\perp$ for each $k\in [K^{[J]}]$. 
    Let $h^{(Q)} \in \mc{T}(Q, \mc{Q})$. For each $k \in [K^{[J]}]$,
    \begin{align*}
        \left\langle \Pi\left[ f\left|\sum_{j\in\mc{J}_k}\mc{D}_k^{(j)}(Q)\right.\right], h^{(Q)}\right\rangle_Q =& \left\langle \Pi\left[ f\left|\sum_{j\in\mc{J}_k}\mc{D}_k^{(j)}(Q)\right.\right], \Pi\left[ h^{(Q)}\left|\sum_{j\in\mc{J}_k}\mc{D}_k^{(j)}(Q)\right.\right]\right\rangle_Q\\
        =& \left\langle f, \Pi\left[ h^{(Q)}\left|\sum_{j\in\mc{J}_k}\mc{D}_k^{(j)}(Q)\right.\right]\right\rangle_Q\\
        =& 0
    \end{align*}
    where the last equality follows because $f\in \mc{T}(Q, \mc{Q})^{\perp}$ and \cref{lemma:li-luedtke-tangent} implies that for $h^{(Q)} \in \mc{T}(Q, \mc{Q})$, $\Pi\left[ h^{(Q)}\left|\sum_{j\in\mc{J}_k}\mc{D}_k^{(j)}(Q)\right.\right] \in \mc{T}(Q, \mc{Q})$. Hence, $\widetilde{\psi}^1_Q$ is an influence function of $\psi$ at $Q$ in model $\mc{Q}$. 

    Next we show that $\widetilde{\psi}^1_Q = \sum_{j\in[J]}\sum_{k\in[K^{(j)}]}m_k^{(j)}$ $\varphi^1_P$ is an influence function of $\varphi$ at $P$ in model $\mc{P}$ by part \ref{item:influence functions1} of \cref{Theorem:influence functions} because 
    \begin{align*}
        \widetilde{\psi}^1_Q(z) =& \sum_{k\in[{K^{[J]}}]}I\left(\overline{z}_{k-1} \in \mc{I}_k\right) \left\{E_Q\left[\psi^1_Q(W)|\overline{z}_{k}\right] - E_Q\left[\psi^1_Q(W)|\overline{z}_{k-1}\right]\right\}\\
        =& \sum_{k\in[K^{[J]}]}\sum_{j\in\mc{J}_k}m_k^{(j)}(\overline{z}_{k})\\
        =&\sum_{j\in[J]}\sum_{k\in[K^{(j)}]}m_k^{(j)}(\overline{z}^{(j)}_{k}).
    \end{align*}
    where the second equality follows from \eqref{eq:m_k-decomp-li-luedtke}. 
\end{proof}

\begin{proof}[Proof of \cref{cor:li-eif}]

    Recall from \cref{theorem:eif} that $\varphi^1_P$ will be the efficient influence function of $\varphi$ at $P$ in model $\mc{P}$ if it satisfies 
    \begin{align}
        \label{eq:li-eif-1}
        \varphi^1_P(o) = \sum_{j\in[J]}I(s=j) \sum_{k\in[K^{(j)}]}\Pi[h^{(Q)}|\mc{D}_k^{(j)}(Q)](\overline{z}_k^{(j)})
    \end{align}
    for some $h^{(Q)} \in \mc{T}(Q, \mc{Q})$ satisfying 
    \begin{align}
        \label{eq:li-eif-2}
        \psi^1_{Q, eff} = \sum_{j\in[J]}\sum_{k\in[K^{(j)}]}\Pi\left\{\frac{dP(\cdot|S=j)}{dQ}(\overline{Z}_{k-1}^{(j)})P(S=j)\Pi[h^{(Q)}|\mc{D}_k^{(j)}(Q)](\overline{Z}_k^{(j)})\Big\vert\mc{T}(Q;\mc{Q})\right\}.
    \end{align}

    Let 
    \begin{align}
        \label{eq:li-eif-3}
        h^{(Q)} \coloneqq \sum_{k\in[K^{[J]}]}&\frac{1}{P(S \in \mc{S}_k(\overline{z}_{k-1}))}\frac{dQ}{dP(\cdot|S \in \mc{S}_k(\overline{z}_{k-1}))}(\overline{z}_{k-1})\\
        &\times I\left(\overline{z}_{k-1}\in \bigcup_{j\in \mc{J}_k}\overline{\mc{Z}}_{k-1}^{(j)}\right)\{E_Q\left[\psi^1_{Q, eff}(W)|\overline{z}_{k}\right] - E_Q\left[\psi^1_{Q, eff}(W)|\overline{z}_{k-1}\right]\}.
    \end{align}
    
    We first show that with $h^{(Q)}$ so defined, \eqref{eq:li-eif-1} holds. To do so, note that
    \begin{align*}
        \Pi[h^{(Q)}|\mc{D}_k^{(j)}(Q)](\overline{z}_k) =& \frac{1}{P(S \in \mc{S}_k(\overline{z}_{k-1}))}\frac{dQ}{dP(\cdot|S \in \mc{S}_k(\overline{z}_{k-1}))}(\overline{z}_{k-1})\\
        &\times I\left(\overline{z}_{k-1}\in \overline{\mc{Z}}_{k-1}^{(j)}\right)\{E_Q\left[\psi^1_{Q, eff}(W)|\overline{z}_{k}\right] - E_Q\left[\psi^1_{Q, eff}(W)|\overline{z}_{k-1}\right]\}.
    \end{align*}
    Therefore \eqref{eq:li-eif-1} holds because  
    \begin{align*}
        \sum_{j\in[J]}I(s=j)& \sum_{k\in[K^{(j)}]}\Pi[h^{(Q)}|\mc{D}_k^{(j)}(Q)](\overline{z}_k)\\
        =\sum_{j\in[J]} \sum_{k\in[K^{(j)}]}&\Bigg[I(s=j)\frac{1}{P(S \in \mc{S}_k(\overline{z}_{k-1}))}\frac{dQ}{dP(\cdot|S \in \mc{S}_k(\overline{z}_{k-1}))}(\overline{z}_{k-1})\\
        &\times I\left(\overline{z}_{k-1}\in \overline{\mc{Z}}_{k-1}^{(j)}\right)\{E_Q\left[\psi^1_{Q, eff}(W)|\overline{z}_{k}\right] - E_Q\left[\psi^1_{Q, eff}(W)|\overline{z}_{k-1}\right]\}\Bigg]\\
        = \sum_{k\in[K^{[J]}]}\sum_{j\in\mc{J}_k}&\Bigg[I(s=j)\frac{1}{P(S \in \mc{S}_k(\overline{z}_{k-1}))}\frac{dQ}{dP(\cdot|S \in \mc{S}_k(\overline{z}_{k-1}))}(\overline{z}_{k-1})\\
        &\times I\left(\overline{z}_{k-1}\in \overline{\mc{Z}}_{k-1}^{(j)}\right)\{E_Q\left[\psi^1_{Q, eff}(W)|\overline{z}_{k}\right] - E_Q\left[\psi^1_{Q, eff}(W)|\overline{z}_{k-1}\right]\}\Bigg]\\
        = \sum_{k\in[K^{[J]}]}&\Bigg[\frac{I(s\in\mc{S}_k(\overline{z}_{k-1}))}{P(S \in \mc{S}_k(\overline{z}_{k-1}))}\frac{dQ}{dP(\cdot|S \in \mc{S}_k(\overline{z}_{k-1}))}(\overline{z}_{k-1})\\
        &\times \{E_Q\left[\psi^1_{Q, eff}(W)|\overline{z}_{k}\right] - E_Q\left[\psi^1_{Q, eff}(W)|\overline{z}_{k-1}\right]\}\Bigg]\\
        =\; \varphi^1_P(o)&.
    \end{align*}
    
    Next, we show that for $h^{(Q)}$ defined as in \eqref{eq:li-eif-3}, the equality \eqref{eq:li-eif-2} holds. To see this, first note that 
    \begin{align*}
        &\sum_{j\in\mc{J}_k}\frac{dP(\cdot|S=j)}{dQ}(\overline{z}_{k-1}^{(j)})P(S=j)\Pi[h^{(Q)}|\mc{D}_k^{(j)}(Q)](\overline{z}_k)\\
        =& \sum_{j\in\mc{J}_k}\frac{P(S=j)}{P(S \in \mc{S}_k(\overline{z}_{k-1}))}\frac{dP(\cdot|S=j)}{dP(\cdot|S \in \mc{S}_k(\overline{z}_{k-1}))}(\overline{z}_{k-1})\\
        &\times I\left(\overline{z}_{k-1}\in \overline{\mc{Z}}_{k-1}^{(j)}\right)\{E_Q\left[\psi^1_{Q, eff}(W)|\overline{z}_{k}\right] - E_Q\left[\psi^1_{Q, eff}(W)|\overline{z}_{k-1}\right]\}\\
        =& I\left(\overline{z}_{k-1}\in \mc{I}_{k-1}\right)\{E_Q\left[\psi^1_{Q, eff}(W)|\overline{z}_{k}\right] - E_Q\left[\psi^1_{Q, eff}(W)|\overline{z}_{k-1}\right]\}\\
        =& \Pi\left[\psi^1_{Q, eff}\Big|\sum_{j\in\mc{J}_k}\mc{D}_k^{(j)}(Q)\right](\overline{z}_k).
    \end{align*}
    We then have 
    \begin{align*}
        \sum_{j\in[J]}\sum_{k\in K^{(j)}}\frac{dP(\cdot|S=j)}{dQ}(\overline{z}_{k-1}^{(j)})P(S=j)\Pi[h^{(Q)}|\mc{D}_k^{(j)}(Q)](\overline{z}_k) =& \sum_{k\in[K^{[J]}]} \Pi\left[\psi^1_{Q, eff}\Big|\sum_{j\in\mc{J}_k}\mc{D}_k^{(j)}(Q)\right]\\
        =&\Pi\left[\psi^1_{Q, eff}\Big|\bigoplus_{k\in[K^{[J]}]}\sum_{j\in\mc{J}_k}\mc{D}_k^{(j)}(Q)\right]\\
        =& \psi^1_{Q, eff}
    \end{align*}
    where the second equality follows because $\sum_{j\in\mc{J}_k}\mc{D}_k^{(j)}(Q) \perp \sum_{j\in\mc{J}_{k'}}\mc{D}_{k'}^{(j)}(Q)$ for $k \not=k'$ and the third from the proof of part \ref{item:lemma:li-ifs1} of \cref{lemma:li-ifs} where we showed $\psi^1_{Q, eff} \in \bigoplus_{k\in[K^{[J]}]}\sum_{j\in\mc{J}_k}\mc{D}_k^{(j)}(Q)$. Because $\psi^1_{Q, eff} \in \mc{T}(Q, \mc{Q})$, the above display implies \eqref{eq:li-eif-2}

    Hence, $\varphi^1_P$ will be the efficient influence function for $\varphi$ at $P$ in model $\mc{P}$ if $h^{(Q)} \in \mc{T}(Q;\mc{Q})$. This follows immediately from the assumption that \eqref{eq:cor:li-eif-corrigendum} holds and the fact that $\mc{T}(Q;\mc{Q})$ is a linear space. 
\end{proof}

    


\section{Additional results and derivations for the examples in the main text}
\label{app:examples-extras}

\subsection{\texorpdfstring{\cref{example:disease-prevalence}}{Example 1} and shadow variables}
\label{app-subsec:shadow-variables}
In this subsection, we elaborate on the connection between \cref{example:disease-prevalence} and the missing data setting with shadow variables. We highlight subtleties of identification that arise when $V$ and $Y$ are polytomous or continuous, and discuss the resulting implications for inference. For concreteness, we compare the results of \cite{li_non-parametric_2023} with those obtained from the general theory developed in the main text, as the structure of that work aligns most closely with our framework.

Shadow variables were introduced by \cite{dhaultfoeuille_new_2010} and have since been studied extensively \citep{wang_instrumental_2014, miao_varieties_2016, li_non-parametric_2023, wang_identification_2024, park_single_2024}. This setting represents a class of missing-not-at-random problems in which an auxiliary variable provides information about the unobserved outcome. Under appropriate conditions to be described below, such a shadow variable can suffice to identify certain functionals of the full data distribution despite nonignorable missingness. In what follows, we first argue that the statistical model of shadow variables is a fused-data model. In lemma \cref{lemma:shadow-variables:tangent-space} we derive the observed data tangent space of this fused-data model corresponding to any ideal data model $\mc{Q}$. We also provide a simple condition on the ideal data law $Q$ under which the observed data extended model $\mc{P}^{ext}$ is strictly semiparametric, locally at any $P\in \mc{P}^{ext}$ such that $P\alignswith Q$. We then comment on the distinct targets of inference between the shadow-variables analysis and the corresponding fused-data analysis, and argue that characterizing inference about the parameters of interest in shadow-variables analysis essentially reduces to characterizing inference of a corresponding parameter under the associated fused-data framework. We conclude by discussing some subtleties of identification that arise in the resulting fused-data framework, and discuss the implications for inference. 

\subsubsection{The shadow variables model as a fused-data model}

The shadow-variables framework conceptualizes the random vector $(X, Y, V, S)$ where $Y$ is an outcome that is not observed in all sampled units, $X$ is a vector of fully observed covariates, and $S$ is the missingness indicator ($S=1$ if $Y$ is missing and $S=2$ if $Y$ is observed). The variable $V$ is fully observed and is assumed to satisfy $V \not\perp Y | X$ and
\begin{align*}
    V \perp S | (Y, X).
\end{align*}
The latter display can be equivalently written as 
\begin{align}
    \label{eq:shadow-variables-P*-alignment}
    P^*(V \leq v|Y, X, S=1) = P^*(V \leq v|Y, X, S=2) \text{ for all $v \in \textsf{Supp}[V;P^*]$ a.e. $P^*$}
\end{align}
where $P^*$ denotes the joint law of the random vector $(X, Y, V, S)$ had $Y$ been always observed. The corresponding statistical model for observed data distribution $P$ is exactly the shadow-variables model considered in \cite{li_non-parametric_2023}. 
 
To connect this framework with fused data, notice that for any joint law $P^*$ satisfying \eqref{eq:shadow-variables-P*-alignment} there exists laws $Q, U^{(2)}$ for $(X, Y, V)$ and $\lambda$ for $S$ such that $P^*$ can be written as
\begin{align}
    \label{eq:shadow-variables:likelihood}
    p^*(x, y, v, s) =& p^*(s)\times \{p^*(v|y, x, S=1)p^*(y|x, S=1)p^*(x|S=1)\}^{I(s=1)}\\
    \nonumber &\times\{p^*(v|y, x, S=2)p^*(y|x, S=2)p^*(x|S=2)\}^{I(s=2)}\\
    \nonumber =& \lambda(s) \times \{q(v|y, x) \times q(y|x) \times q(x)\}^{I(s=1)} \\
   \nonumber  &\times \{q(v|y, x)u^{(2)}(y|x)u^{(2)}(x)\}^{I(s=2)}
\end{align}
where we have assumed the existence of a dominating measure. The shadow-variable condition $V \perp S \mid (Y, X)$ is encoded through the equality
\begin{align*}
    p^*(v| y, x, S=1) =& p^*(v| y, x, S=2)\\
    \coloneqq&  q(v| y, x)
\end{align*}
Conversely, given $(Q, U^{(2)}, \lambda)$, we can construct a law $P^*$ satisfying \eqref{eq:shadow-variables-P*-alignment} whose density is equal to the right hand side of \eqref{eq:shadow-variables:likelihood}. The conditional independence $V \not\perp Y|X$  holds under $P^*$ whenever it holds under $Q$. 

The likelihood for the observed data law $P$ is then given by 
\begin{align}
    \label{eq:shadow-variables:observed-data-liklihood}
    p(o) =  & \lambda(s) \times \{q(v, x) \}^{I(s=1)}\times \{q(v|y, x)u^{(2)}(y|x)u^{(2)}(x)\}^{I(s=2)}.
\end{align}
From this likelihood, we see that the shadow-variables model assumed by \cite{li_non-parametric_2023} is a fused-data model $(\mc{Q}_{sv}, \mc{P}_{sv}, \mc{C})$ where $W = (Y, V, X)$ with $\mc{Q}_{sv}$ the collection of all laws for $W$ satisfying that $Y \not \perp V|X$ under $Q\in\mc{Q}_{sv}$, $W^{(1)} = (V, X)$, and $W^{(2)} = (Y, V, X)$, with alignments
\begin{align}
\label{eq:shadow-variables:alignments}
P(V \le v, X \le x | S=1) &= Q(V \le v, X \le x)\\
\nonumber P(V \le v | Y, X, S=2) &= Q(V \le v | Y, X) \quad \text{a.e. }Q,
\end{align}
for all $(x, v)$ in the support of $(X, V)$ under $Q$, and where the supports of $(X, Y)$ coincide under $Q$ and $P(\cdot | S=2)$. This model coincides with that of \cref{example:disease-prevalence} in the special case where $V$ and $Y$ are both binary. 

In what follows we call any fused-data model $(\mc{Q}, \mc{P}, \mc{C})$ where $\mc{C}$ encodes the alignments \eqref{eq:shadow-variables:alignments} and where $\mc{Q}$ is an ideal data model not necessarily equal to $\mc{Q}_{sv}$, a shadow variables model. 

\subsubsection{The tangent space of the shadow variables model}

Because the shadow variables model is a special case of a fused-data model, we are now in a position to exploit the results in the main text to derive the tangent space of the model. Specifically, the next result, which follows directly from \cref{lemma:score-operator}, characterizes the tangent space of $\mathcal{P}$ under the shadow variables model $(\mc{Q}, \mc{P}, \mc{C})$ for an arbitrary ideal data model $\mathcal{Q}$.
\begin{slemma}
    \label{lemma:shadow-variables:tangent-space}
    Let $(\mc{Q}, \mc{P}, \mc{C})$ be a shadow-variables fused-data model. Let $(Q, P) \in \mc{Q} \times \mc{P}$ where $(Q, P)$ is strongly aligned. Let $\mc{T}(Q, \mc{Q})$ be the ideal data tangent space at $Q$. Then, the observed data tangent space is given by 
    \begin{align*}
        \mc{T}(P, \mc{P}) \coloneqq& cl[{\{I(S=1)E_Q[h^{(Q)}(X, Y, V)|X, V]} \\
        & +I(S=2)\{h^{(Q)}(X, Y, V) - E_Q[h^{(Q)}(X, Y, V)|X, Y]\}:h^{(Q)} \in \mc{T}(Q, \mc{Q})\}]\\
        &\oplus\{I(S=2)g(X, Y): g \in L^2_0((X, Y);P(\cdot|S=2))\}\\
        &\oplus L^2_0(S;P)
    \end{align*}
    where $cl[\cdot]$ denotes closure with respect to $L^2(P)$.
\end{slemma}
It follows from \cref{lemma:shadow-variables:tangent-space} and \cref{lemma:extended-model-nonparametric} in the main text that for the tangent space $\mc{T}(P, \mc{P})$ to be a strict subset $L^2_0(P)$, at least one of the following two conditions must hold:
\begin{enumerate}
    \item $\mc{T}(Q, \mc{Q})$ is strictly included in $L^2_0(Q)$\\
    \item $L^2_0((X, V);Q) \cap \{g \in L^2_0((X, Y, V);Q): E_Q[g(X, Y, V)|X, Y] = 0\} \not= \{0\}$.
\end{enumerate}
Condition (b) is equivalent to the statement that the mere fact that the alignments \eqref{eq:shadow-variables:alignments} hold restricts the observed data tangent space, i.e. even when $\mc{T}(Q, \mc{Q}) = L^2_0(Q)$. Condition (b) is equivalent to the following bounded linear operator having a nontrivial null space
\begin{align}
    \label{eq:shadow-variables:operator}
    T_Q^*: L^2((X, V);Q) \to L^2((X, Y);Q), \qquad (T_Q^* f)(X, Y) \coloneqq  E_Q[f(X, V)| X, Y]
\end{align}
As we shall soon see, this bounded linear operator is a key element in establishing the identification of functionals of the ideal data distribution $Q$ that are of the form $\psi(Q) \coloneqq E_Q[\tau(X, Y)]$ for some known function $\tau \in L^2((X, Y);Q)$ from the observed data distribution. 

\subsubsection{Contrasting the parameters of interest in shadow variables and fused-data analysis}

While the shadow-variables model is a fused data model, the target of inference in the setting of shadow variables does not correspond to any functional $\psi(Q)$ of the ideal data distribution $Q$. This is because under shadow variables, the distribution $P^*$ is the target distribution, not $Q$. However, as mentioned at the beginning of this subsection, the problem of inference for parameters in shadow‐variable models effectively reduces to the problem of inference for the associated parameter under the associated fused‐data framework.

Specifically, in shadow-variables analysis we are interested in estimating $\phi(P^*) \coloneqq E_{P^*}[\tau(X, Y)]$ for some given $\tau \in L^2((X, Y);P^*)$. Note that
\begin{align*}
    \phi(P^*) = &P^*(S=1)E_{P^*}[\tau(X, Y)|S=1] + P^*(S=2)E_{P^*}[\tau(X, Y)|S=2]\\
    \coloneqq &\lambda(S=1)E_{Q}\left[\tau(X, Y)\right] + \lambda(S=2)E_{U^{(2)}}[\tau(X, Y)].
\end{align*}
Taking $\psi(Q) \coloneqq E_Q[\tau(X, Y)]$ as our target functional under the fused data model $(\mc{Q}, \mc{P}, \mc{C})$, we see that $\phi(P^*)$ is a weighted average of the fused-data target of interest and the $U^{(2)}$-mean of $\tau$. Both $\lambda(S=s) \coloneqq P^*(S=s)$ and $E_{U^{(2)}}[\tau(X, Y)] \coloneqq E_{P^*}[\tau(X, Y)|S=2]$ can be directly estimated from the observed data using sample means. Their observed data influence functions are trivial. 

Furthermore, because $\lambda$, $Q$, and $U^{(2)}$ are separate factors in the observed data likelihood \eqref{eq:shadow-variables:observed-data-liklihood}, the efficient influence function of functionals that depend only on $Q$, $U^{(2)}$, and $\lambda$ respectively will be $L^2(P)$ orthogonal. As such, we can construct an efficient RAL estimator of $\phi(P^*)$ by separately constructing efficient RAL estimators of $\psi(Q)$, $\lambda(S=s)$, and $E_{U^{(2)}}[\tau(X, Y)]$, and combining them with the above relation. The efficient influence function of $\phi(P^*)$ can then easily be computed via the delta method. In fact, the influence functions of all RAL estimators of $\phi(P^*)$ can be constructed via the delta method from the set of all observed data influence functions of $E_Q[\tau(X, Y)]$ which can be computed using \cref{Theorem:influence functions} in the main text when $E_Q[\tau(X, Y)]$ is identified. Then, the analysis of the shadow variables target of inference $\phi(P^*)$ essentially reduces to the analysis of the fused-data target of inference $\psi(Q)$. 

To apply the results of \cref{subsec:characterize-ifs} and \cref{subsec:compute-EIF} in the main text, we require that $\psi(Q) = E_Q[\tau(X, Y)]$ be identified under the shadow variables model $(\mc{Q}, \mc{P}, \mc{C})$, i.e. that \cref{as:identification} in the main text holds over the model $\mc{Q}$. \cite{li_non-parametric_2023} derive a necessary and sufficient condition on $Q$ such that $\psi$ agrees and therefore is identified over all $\widetilde{Q}$'s such that $Q \overset{\mc{C}}{\sim}\widetilde{\mc{Q}}$. Consequently, the ideal data model comprised by laws $Q$ that satisfy the necessary and sufficient condition is the maximal ideal data model over which $\psi(Q)$ is identified and where the results of our paper can be applied. 

In the next Theorem, we restate the identification result of \cite{li_non-parametric_2023} using our notation. 
\begin{stheorem}[\cite{li_non-parametric_2023}]
    \label{thrm:shadow-variables:identification}
    Consider the fused-data model $(\mc{Q}, \mc{P}, \mc{C})$ of shadow variables where $\mc{Q}$ is unrestricted except for the condition that for any $Q, \widetilde{Q} \in \mc{Q}$, $L^2((X, Y, V);Q) = L^2((X, Y, V);\widetilde{Q})$. Let $\psi(Q) \coloneqq E_Q[\tau(X, Y)]$. Fix $Q \in \mc{Q}$. Then, $\psi(Q) = \psi(\widetilde{Q})$ for all $\widetilde{Q} \in \mc{Q}$ such that $Q \overset{\mc{C}}{\sim} \widetilde{Q}$ if and only if 
    \begin{align}
        \label{eq:shadow-variables:identifying-condition}
        \tau \in \overline{\text{Range}(T^*_Q)}.
    \end{align}
    where $T_Q^*$ is defined in \eqref{eq:shadow-variables:operator}. 
\end{stheorem}
\cite{li_non-parametric_2023} state \cref{thrm:shadow-variables:identification} in terms of the observed data law $P$ and the $L^2(P)$-closure of the range of the observed data conditional expectation operator $T^*_P:L^2((X, V);P) \rightarrow L^2((X, Y);P)$, $(T_{P}^*f)(X, Y) = E_P[f(X, V)|X, Y, S=2]$. The operators $T^*_P$ and $T^*_Q$ are equivalent by the strong alignment of $(Q, P)$. By framing this condition in terms of laws $Q$ we see how to define the maximal ideal data model $\mc{Q}$ directly as the set of laws satisfying \eqref{eq:shadow-variables:identifying-condition} and the conditions of the lemma.

When, like in \cref{example:disease-prevalence}, $Y$ and $V$ binary, the condition \eqref{eq:shadow-variables:identifying-condition} is satisfied for any $\tau \in L^2((X, Y);Q)$ at any law $Q$ such that $Y$ and $V$ are dependent given $X$. More generally it is also true at every $Q$ such that $T_Q$ is injective, where $T_Q$ is the $L^2(Q)$ adjoint of $T_Q^*$. However whenever $Y$ is non-binary, a non-zero conditional correlation of $Y$ and $V$ given $X$ is not sufficient for $T_Q$ to be injective. In particular, this reveals that for $\tau$ a non-trivial function of $Y$, $\psi(Q)$ is not identified from observed data over the model $\mc{Q}_{sv}$ when $Y$ is non-binary, where we recall that $\mc{Q}_{sv}$ is the set of mutually absolutely continuous laws $Q$ for $(X, Y, V)$ restricted only by the condition that $Y$ and $V$ are conditionally correlated given $X$. To ensure identification, we must replace the restriction on our ideal data model that $Y$ and $V$ are conditionally correlated given $X$ with the restriction that all laws $Q$ in the ideal data model satisfy \eqref{eq:shadow-variables:identifying-condition} for the given $\tau$ of interest. We denote this model $\mc{Q}_{\tau}$ and refer to $(\mc{Q}_{\tau}, \mc{P}_{\tau}, \mc{C})$ as the identified fused-data model. 

\subsubsection{Inference under the identified fused-data model}

Because $\psi$ satisfies \cref{as:identification} in the main text under the model $\mc{Q}_{\tau}$, we then know that by \cref{thrm:identifiability} there exists a well defined observed data functional $\varphi:\mc{P} \rightarrow \mathbb{R}$ such that $\varphi(P) = \psi(Q)$ for all $(Q, P) \in \mc{Q} \times \mc{P}$ such that $P \alignswith Q$. 
Inference about $\varphi$ in the observed data model $\mc{P}_{\tau}$ is subtle for the following reason. As illustrated in \cref{example:shadow-variables:local-identification} below, it may happen that one can find laws in $\mc{Q}_{np}\setminus \mc{Q}_{\tau}$ arbitrarily close to certain $Q^* \in \mc{Q}_{\tau}$, thus locally restricting model $\mc{Q}_{\tau}$. Here $\mc{Q}_{np}$ is the nonparametric model comprised by all laws for $(X, Y, V)$. When $Y$ is discrete and finitely valued, it is possible to characterize $\mc{T}(Q, \mc{Q}_{\tau})$ for any $Q \in \mc{Q}$ and any $\tau \in L^2((X, Y);Q)$. However, when $Y$ is continuous, to the best of our knowledge, characterizing $\mc{T}(Q, \mc{Q}_{\tau})$ is an open problem. In spite of this, we can generically characterize the set of all influence functions of RAL estimators, equivalently gradients of $\varphi$, by invoking \cref{prop:two-source-solution}. Specifically, $\varphi^1_P$ is a gradient of $\varphi$ at $P $ in $\mc{P}_{\tau}$ if and only if it is of the form
\begin{align}
    \label{eq:shadow-variables:all-ifs}
    \varphi^1_P(O) \coloneqq \frac{I(S=1)}{P(S=1)}\{m_Q(X, V) - \varphi(P)\} + \frac{I(S=2)}{P(S=2)}\frac{q(X, Y)}{p(X, Y|S=2)}\{\psi^1_Q(X, Y) - m_Q(X, Y)\}
\end{align}
where $Q$ is any $Q \in \mc{Q}_{\tau}$ such that $(Q, P)$ is strongly aligned, $\psi^1_Q$ is a gradient of $\psi$ in model $\mc{Q}_{\tau}$ such that the equation
\begin{align}
    \label{eq:shadow-variables:T_Q^*-eq}
    T_Q^*m = \psi^1_Q
\end{align}
has a solution $m \in L^2((X, V);Q)$, and $m_Q$ is any such solution. 

While \eqref{eq:shadow-variables:all-ifs} gives a generic characterization of the set of all observed data gradients, we are unable to go beyond this characterization because as indicated above calculation of the tangent space $\mc{T}(Q, \mc{Q}_{\tau})$, and therefore of all ideal data gradients $\psi^1_Q$, is an open problem. For the same reason, we cannot characterize the canonical gradient beyond the generic form of \cref{theorem:eif}. However, $\psi^1_Q \coloneqq \tau(X, Y) - \psi(Q)$ is still an ideal data gradient, and therefore if \eqref{eq:shadow-variables:T_Q^*-eq} has a solution $m_Q$ the set of all observed data gradients $\varphi^1_P$ at $P$ in model $\mc{P}_{\tau}$ contains the set 
\begin{align}
    \label{eq:shadow-variables:Q-set-of-ifs}
    \Bigg\{&\frac{I(S=1)}{P(S=1)}\{m_Q(X, V) + f(X, Y) - \varphi(P)\} \\
    \nonumber &+ \frac{I(S=2)}{P(S=2)}\frac{q(X, Y)}{p(X, Y|S=2)}\{\tau(X, Y) - [m_Q(X, Y) + f(X, Y)]\}: f \in \text{Null}(T_Q^*)\Bigg\}
\end{align}
for any $Q \in \mc{Q}_{\tau}$ such that $(Q, P)$ is strongly aligned. 

Because when $T_Q$ is not injective, $Q$ need not be fully identified by $P$ even if $\psi(Q)$ is identified by $P$, the set in display \eqref{eq:shadow-variables:Q-set-of-ifs} gives a set of observed data gradients for any choice $Q \in \mc{Q}_{\tau}$ such that $(Q, P)$ is strongly aligned. Thus, taking the union of all sets \eqref{eq:shadow-variables:Q-set-of-ifs} over all $Q\in \mc{Q}_{\tau}$ such that $(Q, P)$ is strongly aligned gives a potentially larger set of observed data gradients. In the special case in which all $Q \in \mc{Q}_{\tau}$ have equivalent norms, the sets $\{m_Q(X, Y) + f(X, Y): f \in \text{Null}(T_Q^*), T_Q^*m_Q = \tau(X, Y) - \psi(Q)\}$ are equal for all $Q$'s that strongly align with $P$. Then, taking $Q \in \mc{Q}_{\tau}$ as fixed with $(Q, P)$ strongly aligned, the aforementioned union can be explicitly characterized as 
\begin{align*}
    \Bigg\{&\frac{I(S=1)}{P(S=1)}\{m_Q(X, V) + f(X, Y) - \varphi(P)\} \\
    \nonumber &+ \frac{I(S=2)}{P(S=2)}\frac{\widetilde{q}(X, Y)}{p(X, Y|S=2)}\{\tau(X, Y) - [m_Q(X, Y) + f(X, Y)]\}: f \in \text{Null}(T_Q^*), \widetilde{Q} \in \mc{Q}_{\tau}, P\alignswith\widetilde{Q}\Bigg\}
\end{align*}

\cite{li_non-parametric_2023} provided an estimator of $\phi(P^*) = E_{P^*}[\tau(X, Y)]$ that is asymptotically linear at any $P^*$ with observed data law $P \in \mc{P}_{\tau}$. They did not comment on regularity of this estimator except in the special case where $T_Q$ and $T_Q^*$ are both injective. As discussed earlier, observed data gradients for the shadow-variables target parameter determine observed data gradients of the fused-data target parameter and vise-versa. In fact, one can show, using the characterization in the last display, that the influence function of their proposed asymptotically linear estimator is indeed a gradient, and so their proposed estimator is regular with respect to the model $\mc{P}_{\tau}$. Recall from \cref{app:semiparametric-theory} that asymptotically linear estimators are regular if their influence functions are gradients of the target parameter. 

\subsubsection{Failure of local identification under \texorpdfstring{$\mc{Q}_{sv}$}{Q}}

We now give the aforementioned example that illustrates that $\psi(Q)$ may not be locally identified in the model $\mc{Q}_{sv}$, in the sense that given a law $Q \in \mc{Q}_{sv}$ where \eqref{eq:shadow-variables:identifying-condition} holds, there may be laws in $\mc{Q}_{sv}$ arbitrarily close to $Q$ in Hellinger distance where \eqref{eq:shadow-variables:identifying-condition} does not hold. 

\begin{sexample}
    \label{example:shadow-variables:local-identification}
    Suppose $Y$ and $V$ both take values in $\{1, 2, 3\}$, and that $X$ is degenerate at $0$. Let $\mc{Q}_{sv}$ be the collection of all laws on $(X, Y, V)$ such that $V$ and $Y$ are correlated. Suppose $Q$ is such that 
    \begin{align*}
        Q_{V|Y} \coloneqq \begin{pmatrix}
            q(V=1|Y=1) & q(V=2|Y=1) & q(V=3|Y=1)\\
            q(V=1|Y=2) & q(V=2|Y=2) & q(V=3|Y=2)\\
            q(V=1|Y=3) & q(V=2|Y=3) & q(V=3|Y=3)
        \end{pmatrix} \coloneqq 
        \begin{pmatrix}
            0.5 & 0.2 & 0.3\\
            0.5 & 0.3 & 0.2\\
            0.5 & 0.3 & 0.2\\
        \end{pmatrix}
    \end{align*}
    so that
    \begin{align*}
        \begin{pmatrix}
            (T^*_Qf)(1)\\
            (T^*_Qf)(2)\\
            (T^*_Qf)(3)
        \end{pmatrix} = Q_{V, Y}\begin{pmatrix}
            f(1), f(2), f(3)
        \end{pmatrix}.
    \end{align*}
    Hence, $\text{Range}(T^*_Q) = \text{Range}(Q_{V|Y})$ when we view functions of $V$ and $Y$ as three dimensional vectors. Because $Y$ and $V$ take values in a finite dimensional space, these ranges are closed. It follows from \cref{thrm:shadow-variables:identification} above that $\psi(Q)$ is identified if and only if $\tau \in \text{Range}(Q_{V|Y})$. Suppose $\tau(Y) \coloneqq (0.7, 0.8, 0.8)^T$. Then, $\tau(Y) = T^*_{Q} f$ where $f(V) \coloneqq (1, 1, 0)^T$, and so $\psi(Q) \coloneqq E_Q[\tau(Y)]$ is identified from any aligned observed data law. To complete defining the ideal and observed data laws, we let $Q(Y) \coloneqq (.5, .3, .2)$, $U^{(2)}(Y) \coloneqq (.5, .2, .3)$, and $\lambda(S=1) = 0.2$. Then, 
    \begin{align*}
        P_{Q, U, \lambda}(V|S=1)^T = Q(Y)^T \times Q_{V|Y}.
    \end{align*}
    
    We now construct a DQM submodel $\{Q_t: \|t\| \in (-\delta, \delta), Q_0 = Q\}$ of $\mc{Q}_{sv}$ where $\psi$ is not identified from the observed data $P_t$ no matter how small $\|t\|$ is. Let $t$ be in the $\delta$ ball around 0 in $\mathbb{R}^2$ and define
    \begin{align*}
        Q_{V|Y;t} \coloneqq \begin{pmatrix}
            q_t(V=1|Y=1) & q_t(V=2|Y=1) & q_t(V=3|Y=1)\\
            q_t(V=1|Y=2) & q_t(V=2|Y=2) & q_t(V=3|Y=2)\\
            q_t(V=1|Y=3) & q_t(V=2|Y=3) & q_t(V=3|Y=3)
        \end{pmatrix} \coloneqq 
        B(t_1 + t_2) \times Q_{V|Y}.
    \end{align*}
    where 
    \begin{align*}
        B(\epsilon)  = \begin{pmatrix}
            1 & \epsilon & -\epsilon\\
            -\epsilon & 1 & \epsilon\\
            \epsilon & -\epsilon & 1
        \end{pmatrix} 
    \end{align*}
    Let $Q_{t}(Y) \coloneqq  Q(Y)+ \{B(t_2)^{-1}\}^T(0, -t_2, t_2)^T$. Then, 
    \begin{align*}
        P_{Q_t, U, \lambda}(V|S=1)^T = (Q(Y)+ (0, -t_2, t_2)\times B(t_2)^{-1}) \times B(t_1 + t_2)\times Q_{V|Y}.
    \end{align*}
    Then, $P_{Q_{t = (\epsilon, 0)}, U, \lambda}(V|S=1)^T = P_{Q_{t = (0, \epsilon)}, U, \lambda}(V|S=1)^T$ because $(0, -t, t) \times Q_{V|Y} = 0$. As such, $Q_{t = (t_1, 0)} \overset{\mc{C}}\sim Q_{t = (0, t_2)}$ for all $t_1=t_2 = \epsilon$ small enough that the above are valid probability distributions. But, 
    \begin{align*}
        E_{Q_{t = (0, \epsilon)}}[\tau(Y)] = (Q(Y)^T+ (0, -\epsilon, \epsilon)\times B(\epsilon)^{-1} )\times \tau(Y) \not = Q(Y)^T \times \tau(Y) = E_{Q_{t = (\epsilon, 0)}}[\tau(Y)].
    \end{align*}
    This submodel for $Q$ is DQM and so by choosing $\delta$ sufficiently small all laws in the submodel can be made arbitrary close to $Q$ in Hellinger norm, but \cref{as:identification} does not hold for any $\delta > 0$.
\end{sexample}

\subsection{Extending the the examples to broader fused-data frameworks}
\label{app-subsec:gen-examples}
In this subsection we discuss classes of fused-data frameworks that include the frameworks of the examples as special cases. For each of these classes, we provide characterizations of the set of all observed data influence functions. We then provide a lemma that establishes, for each of these classes, whether or not the observed data tangent space $\mc{T}(P, \mc{P})$ and the extended tangent space is equal to $L^2_0(P)$. We conclude this subsection by introducing \cref{example:transporting} scenario (iv) that is like scenario (iii) except that we make weaker alignment assumptions that still allow identification of the average treatment effect in the target population of interest under the assumptions discussed in \cref{app-subsec:causal-identification}. We then discuss in more depth two classes of frameworks, one that includes \cref{example:transporting} scenarios (ii) and (iv), and another that includes \cref{example:transporting} scenario (iii). For the latter class we provide the efficient observed data influence function. As in the main text, we assume the existence of dominating product measures. 

\begin{sproposition}
\label{prop:example-ifs}
Let $\left( \mathcal{Q},\mathcal{P},%
\mathcal{C}\text{\thinspace },\psi ,\varphi \right) $ be a fused-data
framework with respect to $\left( Q_{0},P_{0}\right) $. Let $P\in \mathcal{P}
$. Suppose there exists $Q$ in $\mc{Q}$ such that $\left(
Q,P\right) $ is strongly aligned with respect to $\mathcal{C}$ and $\psi $
is pathwise differentiable at $Q$ in model $\mathcal{Q}$. Suppose there exists product measure $\mu$ that dominates $Q$ and let $q = \frac{dQ}{d\mu}$ and suppose there exist product measures $\mu^{(j)}$ that dominate $P(\cdot|S=j)$ and let $p(\cdot|S=j) = \frac{dP(\cdot|S=j)}{d\mu^{(j)}}$ for $j\in[J]$. 

\begin{enumerate}
    \item\label{item:prop:example-ifs1} Suppose $(\mc{Q}, \mc{P}, \mc{C})$ is the fused-data model of \cref{example:disease-prevalence}. Let 
    \begin{align*}
        m_Q(x, v) \coloneqq \frac{v - E_Q[V|Y=0, x]}{E_Q[V|Y=1, x] - E_Q[V|Y=0, x]}
    \end{align*}
    and suppose $E_Q[V|Y=1, x] - E_Q[V|Y=0, x] \geq \widetilde{\delta}$ for some $\widetilde{\delta} > 0$. Write the influence function $\psi^1_Q$ for $\psi$ at $Q$ in model $\mc{Q}$ as 
    \begin{align*}
        \psi^1_Q(x, v, y) = \psi^1_{Q; X}(x) + \psi^1_{Q;X, V}(x)v + \psi^1_{Q; X, Y}(x)y + \psi^1_{Q;X, Y, V}(x)vy
    \end{align*}
    for some $\psi^1_{Q; X}, \psi^1_{Q; X, V}, \psi^1_{Q; X, Y}, \psi^1_{Q; X, Y, V} \in L^2(X;Q)$.
    Then $\varphi$ is pathwise differentiable at $P$ in model $\mc{P}$ and 
    \begin{align*}
        \varphi^1_P(o) = \frac{I(s=1)}{P(S=1)}&\{\psi^1_{Q;X}(x) + \psi^1_{Q;X, V}(x)v  \\
        &+\psi^1_{Q; X, Y}(x)m_Q(x, v) + \psi^1_{Q;X, V, Y}(x)m_Q(x, v)E_Q[V|Y=1, x]\}\\
        + \frac{I(s=2)}{P(S=2)}&\frac{q(x, y)}{p(x, y|S=2)}\big\{\psi^1_{Q;X, Y}(x)\{y - m_Q(x, v)\}\\
        &+ \psi^1_{Q;X, V, Y}(x)\{vy - m_Q(x, v)E_Q[V|Y=1, x]\}\big\}
    \end{align*}
    is its unique influence function. In particular, the influence function for $\varphi$ in \cref{example:disease-prevalence} follows by applying the preceding formula to $\psi^1_Q = y - \psi(Q)$ since $\psi^1_{Q;X} = -\psi(Q)$, $\psi^1_{Q;X, V} = \psi^1_{Q;X,Y,V} = 0$ and $\psi^1_{Q;X,Y} = 1$. 
    \item \label{item:prop:example-ifs2} Suppose $\mc{C}$ is as in \cref{example:tsiv-lsm}. An influence function for $\psi$ at $Q$ in $\mc{Q}$ corresponds to an influence function for $\varphi$ at $P$ in model $\mc{P}$ if and only if $\psi^1_Q(l, x, y) = \psi^1_{Q;L, Y}(l, y) + \psi^1_{Q;L, X}(l, x)$ for some $\psi^1_{Q;L, Y} \in L^2(L, Y;Q)$, $\psi^1_{Q;L, X} \in L^2(L, X;Q)$, and 
    \begin{align*}
        E_Q[\psi^1_{Q;L, X}(L, X)|L] + E_Q[\psi^1_{Q;L, Y}(L, Y)|L] = 0.
    \end{align*}
    In such case, 
    \begin{align*}
        \varphi^1_P(o) =& \frac{I(s=1)}{P(S=1)}\frac{q(l)}{p(l|S=1)}\left\{\psi^1_{Q;L, Y}(l, y) - E_Q[\psi^1_{Q;L, Y}(L, Y)|l]\right\}\\
        &+\frac{I(s=2)}{P(S=2)}\frac{q(l)}{p(l|S=2)}\left\{\psi^1_{Q;L, Y}(l, x) + E_Q[\psi^1_{Q;L, Y}(L, Y)|l]\right\}
    \end{align*}
    is the unique influence function of $\varphi$ that corresponds to $\psi^1_Q$. 
    \item \label{item:prop:example-ifs3} Suppose $W = (U, B)$ with $U$ a finitely valued discrete random variable and $B$ a euclidean random vector of dimension $r$. Suppose $\mc{Q}$ is a nonparametric model. Suppose $J=2$ and the alignments in $\mc{C}$ are 
    \begin{align*}
        Q(U\leq u|B) &= P(U \leq u|B, S=1) \text{ a.e.-$Q$}\\
        Q(B\leq b|U=u_0) &= P(B \leq b|U=u_0, S=2)
    \end{align*}
    for all $u \in \textsf{Supp}[U;Q]$, $b \in \mathbb{R}^r$, and for some $u_0 \in \textsf{Supp}[U;Q]$. Additionally, suppose there exists $\widetilde{\delta} < \infty$ such that $q(u_0|B)^{-1}\leq \widetilde{\delta}$ a.e.-$Q$. Let $\psi^1_Q$ be the influence function of $\psi$ at $Q$ in model $\mc{Q}$. Then, $\varphi^1_P(o)$ is pathwise differentiable at $P$ in model $\mc{P}$ and 
    \begin{align*}
        \varphi^1_P(o) =& \frac{I(s=1)}{P(S=1)}\frac{q(b)}{p(b|S=1)}\left\{\psi^1_Q(u, b) - \frac{I(u=u_0)}{q(u_0|b)}E_Q\left[\psi^1_Q(U, B)|b\right]\right\}\\
        &+\frac{I(s=2)}{P(S=2)}\frac{q(u_0)}{p(u_0|S=2)}\frac{I(u=u_0)}{q(u_0|b)}E_Q\left[\psi^1_Q(U, B)|b\right]
    \end{align*}
    is its unique influence function.
    \item \label{item:prop:example-ifs4} Suppose $W = (U, B)$ with $U$ and $B$ random euclidean vectors of dimension $p$ and $r$ respectively. Suppose there exists no measurable maps $g_1$ such that $U=g_1\left(
B\right) $ a.e. - $Q$ and likewise there exists no measurable maps $g_2$ such that $B=g_2\left(
U\right) $ a.e. - $Q$. Suppose $\textsf{Supp}[(U, B);Q] = \textsf{Supp}[U;Q]\times \textsf{Supp}[B;Q]$, $\mc{Q}$ is a nonparametric model, and suppose there exists $\widetilde{\delta} < \infty$ such that $\frac{q(U)q(B)}{q(U, B)} \leq \widetilde{\delta}$ a.e.-$Q$. Suppose $J=2$ and the alignments in $\mc{C}$ are
    \begin{align*}
        Q(U \leq u|B) =& P(U\leq u|B, S=1) \text{ a.e.-$Q$}\\
        Q(B\leq b|U) =& P(B \leq b|U, S=2) \text{ a.e.-$Q$}
    \end{align*}
    for all $u \in \mathbb{R}^p$ and $b \in \mathbb{R}^r$. Let $\psi^1_Q$ be the influence function of $\psi$ at $Q$ in model $\mc{Q}$. Then $\varphi$ is pathwise differentiable at $P$ in model $\mc{P}$ and the set of its influence functions is given by
    \begin{align*}
        \left\{\varphi^1_P\right\} + \left\{\left[\frac{I(s=1)}{P(S=1)}\frac{q(b)}{p(b|S=1)} - \frac{I(s=2)}{P(S=2)}\frac{q(u)}{p(u|S=2)}\right]f(u, b):f \in \mc{F}\right\}
    \end{align*}
    where 
    \begin{align*}
        \varphi^1_P =& \frac{I(s=1)}{P(S=1)}\frac{q(b)}{p(b|S=1)}\left\{\psi^1_Q(u, b) - \frac{q(u)q(b)}{q(u, b)}E_Q[\psi^1_Q(U, B)|b]\right\} \\
        &+ \frac{I(s=2)}{P(S=2)}\frac{q(u)}{p(u|S=2)}\left\{\frac{q(u)q(b)}{q(u, b)}E_Q[\psi^1_Q(U, B)|b]\right\}
    \end{align*}
    and 
    \begin{align*}
        \mc{F} = \left\{f \in L^2_0(Q):E_Q[f(U, B)|U] = E_Q[f(U, B)|B] = 0 \text{ a.e. $Q$}\right\}.
    \end{align*}
\end{enumerate}
\end{sproposition}

\begin{sremark}
    \cite{tchetgen_tchetgen_doubly_2010} showed that the set $\mc{F}$ in part \ref{item:prop:example-ifs4} of the above lemma coincides with the set  
    \begin{align*}
        \left\{\frac{q(u)q(b)}{q(u, b)}\left\{t(u, b) - E_{Q^*}[t(U, B)|b] - E_{Q^*}[t(U, B)|u] + E_{Q^*}[t(U, B)]\right\}:t \in L_0^2(Q)\right\}
    \end{align*}
    where $E_{Q^*}$ denotes expectation under the law $Q^*$ with density $q^*(u, b) = q(u)q(b)$.
\end{sremark}

The fused-data frameworks in \cref{example:disease-prevalence} and \cref{example:tsiv-lsm} are special cases of the frameworks in parts \ref{item:prop:example-ifs1} and \ref{item:prop:example-ifs2} of the above proposition respectively, which generalize these examples by allowing $\psi$ to be any pathwise differentiable parameter satisfying \cref{as:identification}. \cref{example:tsiv-lsm} is further generalized by allowing $\mc{Q}$ to be any model. The fused-data frameworks in part \ref{item:prop:example-ifs3} generalize the frameworks of scenarios (ii). In that scenario, $U = Y$, $u_0 = y_0$, and $B = (L, A)$. The fused-data frameworks of part \ref{item:prop:example-ifs4} generalize the framework of scenario (iii) in \cref{example:transporting}. In this scenario, $U = Y$ and $B = (A, L)$. 

The next lemma states the tangent spaces of the observed data models for the fused-data models considered in Examples \ref{example:disease-prevalence}, \ref{example:tsiv-lsm}, and \ref{example:transporting}. 
\begin{slemma}
    \label{lemma:example-tangent-spaces}
    Let $\left( \mathcal{Q},\mathcal{P},%
    \mathcal{C}\text{\thinspace }\right) $ be a fused-data
    model with respect to $\left( Q_{0},P_{0}\right)$. Let $P\in \mathcal{P}$. Suppose there exists $Q$ in $\Phi \left( P;\mc{C}\right) $ such that $\left(
    Q,P\right) $ is strongly aligned with respect to $\mathcal{C}$. Suppose there exists product measure $\mu$ that dominates $Q$ and let $q = \frac{dQ}{d\mu}$ and suppose there exist product measures $\mu^{(j)}$ that dominate $P(\cdot|S=j)$ and let $p(\cdot|S=j) = \frac{dP(\cdot|S=j)}{d\mu^{(j)}}$. 
    \begin{enumerate}
        \item\label{item:lemma:example-tangent-spaces1} Suppose $(\mc{Q}, \mc{P}, \mc{C})$ is the fused-data model of \cref{example:disease-prevalence}. Then $\mc{T}(P, \mc{P}) = \mc{T}(P, \mc{P}^{ext}) = L_0^2(P)$.
        \item\label{item:lemma:example-tangent-spaces2} Suppose $(\mc{Q}, \mc{P}, \mc{C})$ is a fused-data model where $\mc{C}$ is as in \cref{example:tsiv-lsm}. Suppose that $\mc{Q}$ is such that $\{\mc{D}_2^{(1)}(Q)\cap \mc{T}(Q, \mc{Q})^\perp\}\cup\{\mc{D}_2^{(2)}(Q)\cap \mc{T}(Q, \mc{Q})^\perp\} \not=\{0\}$. Then $\mc{T}(P, \mc{P}) \subsetneq \mc{T}(P, \mc{P}^{ext}) = L_0^2(P)$.
        \item\label{item:lemma:example-tangent-spaces3} Suppose $(\mc{Q}, \mc{P}, \mc{C})$ is the fused-data model of \cref{example:transporting} scenario (i). Then $\mc{T}(P, \mc{P}) = \mc{T}(P, \mc{P}^{ext}) = L_0^2(P)$.
        \item\label{item:lemma:example-tangent-spaces4} Suppose $(\mc{Q}, \mc{P}, \mc{C})$ is the fused-data model of part \ref{item:prop:example-ifs3} of \cref{prop:example-ifs}, which includes \cref{example:transporting} scenario (ii) and scenario (iii) as special cases. Then $\mc{T}(P, \mc{P}) = \mc{T}(P, \mc{P}^{ext}) = L_0^2(P)$.
        \item\label{item:lemma:example-tangent-spaces5} Suppose $(\mc{Q}, \mc{P}, \mc{C})$ is the fused-data model of part \ref{item:prop:example-ifs4} of \cref{prop:example-ifs}. Then $\mc{T}(P, \mc{P}) = \mc{T}(P, \mc{P}^{ext}) \subsetneq L^2_0(P)$.
    \end{enumerate}
\end{slemma}

From the above lemma, the unique observed data influence function is efficient in parts \ref{item:prop:example-ifs1} and \ref{item:prop:example-ifs3} of \cref{prop:example-ifs} because the observed data model is non-parametric. However, in part \ref{item:prop:example-ifs2}, when $\mc{Q}$ is semiparametric and places restrictions on the aligned components, there will be infinitely many observed data influence functions, and the efficient one may be found using the technique of \cref{prop:only-Q-restricted-eif} in \cref{app:model-types-eif}. We demonstrate the application of this lemma for this example in \cref{app-subsec:additional-derivations}. In part \ref{item:prop:example-ifs4} there are also infinitely many observed data influence functions because the model $\mc{P}$ is strictly semiparametric. However, in contrast to \ref{item:prop:example-ifs3}, $\mc{P}$ is strictly semiparametric because the assumed alignments in $\mc{C}$ place equality constraints on the laws $P$ in the model $\mc{P}$. We provide a lemma giving the efficient influence function in \cref{app-subsubsec:generalizing-ex-3} below. 

\subsubsection{Additional discussion of  \texorpdfstring{\cref{example:transporting}}{Example 3}}
\label{app-subsubsec:generalizing-ex-3}

We first describe the aforementioned scenario (iv) for \cref{example:transporting}. 

\textbf{Scenario (iv):} Suppose we again perform a case-control study and measure that outcomes, covariates, and treatments. Additionally suppose the covariates $L$ are discrete. Under the assumptions in \cref{app-subsec:causal-identification}, the following alignments hold for $Q=Q_0$, $P=P_0$
\begin{gather}
    \label{eq:scenario-4-align-1} Q\left( Y=1|L=l_{0},A=0\right) =P\left( Y=1|L=l_{0},A=0,S=1\right)\\
    \label{eq:scenario-4-align-2} Q\left( L\leq l,A=a|Y\right) =P\left( L\leq l,A=a|Y,S=2\right) \text{ for all $l \in \mathbb{R}^p, a \in \{0, 1\}$ a.e.-$Q$}.
\end{gather}

This first alignment implies that only the prevalence among the unexposed with
covariate level $l_{0}$ agree between the source and target population. The second alignment is the same as in scenario (iii). Then, assuming that $0<Q\left( Y=1|L=l_{0},A=0\right) <1,$ the identity 
\begin{equation*}
\left. Q\left( Y=1\right) =\left\{ 1+\left. \frac{Q\left(
L=l_{0},A=0|Y=1\right) }{Q\left( Y=1|L=l_{0},A=0\right) }\right/ \frac{%
Q\left( L=l_{0},A=0|Y=0\right) }{1-Q\left( Y=1|L=l_{0},A=0\right) }\right\}
^{-1}\right.
\end{equation*}%
shows that $Q$ is identified by the coarsened data law $P$ due to
the alignments \eqref{eq:scenario-4-align-1} and \eqref{eq:scenario-4-align-2}. In particular, $\psi(Q)$ is equal to 
\begin{align}
    \label{eq:phi-scenario-4}\varphi \left( P\right) \coloneqq\widetilde{\alpha} \left( P\right)&\\
    \nonumber\times\sum_{a=0}^{1}\sum_{l}\left( -1\right) ^{1-a}&\widetilde{\beta} \left(
a,l;P\right) \left[ P\left( L=l|Y=1,S=2\right) \widetilde{\alpha} \left(
P\right) +P\left( L=l|Y=0,S=2\right) \left\{ 1-\widetilde{\alpha} \left(
P\right) \right\} \right]
\end{align}
where%
\begin{equation*}
\widetilde{\beta} \left( a,l;P\right) \coloneqq\frac{P\left(
A=a,L=l_0|Y=1,S=2\right) }{P\left( A=a,L=l_0|Y=1,S=2\right) \widetilde{\alpha} \left(
P\right) +P\left( A=a,L=l_0|Y=0,S=2\right) \left\{ 1-\widetilde{\alpha} \left(
P\right) \right\} }
\end{equation*}%
and 
\begin{equation*}
\widetilde{\alpha} \left( P\right) \coloneqq\left\{ 1+\left. \frac{P\left(
L=l_0,A=0|Y=1,S=2\right) }{P\left( Y=1|L=l_0,A=0,S=1\right) }\right/ 
\frac{P\left( L=l_0,A=0|Y=0,S=2\right) }{1-P\left(
Y=1|L=l_0,A=0,S=1\right) }\right\} ^{-1}.
\end{equation*}

As in scenario (ii), the mere existence of $Q_{0}$ satisfying \eqref{eq:scenario-4-align-1} and \eqref{eq:scenario-4-align-2} does not impose restrictions
on the coarsened data distribution $P_{0}$ other than support constraints.
Consequently, this new scenario fits the fused-data framework $\left( 
\mathcal{Q},\mathcal{P},\mathcal{C},\psi ,\varphi \right) $
with $\mathcal{Q}$ and $\mathcal{P}$ non-parametric, $\psi \left(
Q_{0}\right) $ defined in \eqref{eq:psi-example3}, $\varphi \left(
P_{0}\right) $ as defined in the preceding display, and $\mathcal{C}$ the
collection of alignments \eqref{eq:scenario-4-align-1} and \eqref{eq:scenario-4-align-2}.
\medskip

In \cref{example:transporting} scenarios (ii) and (iv) the components of $Q_{0}$ that align with the sources can be expressed as the conditional distributions of $B$ given $%
U=u_{0}$ and the conditional distribution of $U$ given $B$ for certain
variables $U$ and $B,$ with $U$ discrete. Specifically, $U=Y,B=\left( L,A\right) $ and $u_{0}=1$ for scenario (ii) and $U=(L, A),B=Y$ and $u_{0}=(l_0, 0)$ for scenario (iv). \cite{arnold_specification_1996} showed that given a pair of probability laws $%
P_{1}$ and $P_{2}$ on $\left( U,B\right) $ there exists a law $Q$ on $\left(
U,B\right) $ such that the conditional distribution of $U$ given $B$ under $%
P_{1}$ and $Q$ agree, and the conditional distribution of $B$ given $U=u_{0}$
under $P_{2}$ and $Q$ agree, provided certain conditions on the supports of $%
P_{1}$ and $P_{2}$ hold. Identifying $P_{1}$ with $P\left( \cdot|S=1\right) $ and 
$P_{2}$ with $P\left( \cdot|S=2\right) ,$ this result establishes that the model $%
\mathcal{P}$ for the law $P$ of coarsened data $O=\left( c\left( W,S\right)
,S\right) $ with $S\in \left\{ 1,2\right\} ,c\left( W,1\right) =c\left(
W,2\right) =W=\left( U,B\right) $ and alignments $Q\left( U\leq u|B\right)
=P\left( U\leq u|B,S=1\right) $ and $Q\left( B\leq b|U=u_{0}\right) =P\left(
B\leq b|U=u_{0},S=2\right) $ does not impose equality constraints on $P.$ This discussion is summarized in the following lemma.

\begin{slemma}
\label{lemma:case-control-ji} 
Let $\mathcal{U}\times \mathcal{B}\subseteq \mathbb{R}^{r}\times \mathbb{R}^{m}$. Let $Q$, $P_1$, $P_2$ be probability measures over $\mathcal{U}\times \mathcal{B}$ each dominated by some product measure $\lambda\times \mu$. Let $q, p_1, p_2:\mathcal{U}\times \mathcal{B}\rightarrow \mathbb{R}$ be versions of the densities of $Q, P_1, P_2$ with respect to $\lambda\times \mu$, respectively. Let $\mathcal{B}^{\ast }\coloneqq\left\{ b\in \mathcal{B}:q\left( b\right)
>0\right\}$ where $q\left( b\right) \coloneqq\int q\left( u,b\right)
d\lambda \left( u\right)$.
\begin{enumerate}
    \item \label{item:lemma:case-control-ji1} Suppose that there exists $u_{0}\in \mathcal{U}$
satisfying 
\begin{equation*}
q\left( u_{0}|b\right) >0\text{ for all }b\in \mathcal{B}^{\ast }
\end{equation*}%
Then, for any $u\in \mathcal{U}$ and $b\in \mathcal{B}^{\ast }$ it holds
that 
\begin{equation}
q\left( u,b\right) =q\left( u|b\right) \times \left. \frac{q\left(
b|u_{0}\right) }{q\left( u_{0}|b\right) }\right/ \int_{b:q\left( b\right) >0}%
\frac{q\left( b|u_{0}\right) }{q\left( u_{0}|b\right) }d\mu \left( b\right).
\label{q-formula}
\end{equation}%
\item \label{item:lemma:case-control-ji2}Let $\mc{B}^*_1 \coloneqq \{b \in \mc{B}:p_1(b) > 0\}$. Suppose that for any $b\in 
\mathcal{B}^{\ast }$ it holds that $b\in \mathcal{B}_{1}^{\ast }$, $p_{1}\left( u_{0}|b\right) >0$, and $p_{2}\left( b|u_{0}\right)
=q\left( b|u_{0}\right) $ and $p_{1}\left( u_{0}|b\right) =q\left(
u_{0}|b\right)$. Then, for all $b\in \mathcal{B}^{\ast }$%
\begin{equation*}
q\left( u,b\right) =p_{1}\left( u|b\right) \times \left. \frac{p_{2}\left(
b|u_{0}\right) }{p_{1}\left( u_{0}|b\right) }\right/ \int_{\mathcal{B}}\frac{%
p_{2}\left( b|u_{0}\right) }{p_{1}\left( u_{0}|b\right) }d\mu \left(
b\right). 
\end{equation*}%
\item \label{item:lemma:case-control-ji3} Suppose there exists $u_{0}\in \mathcal{U}$
satisfying $p_{2}\left( u_{0}\right) >0$ and $p_{1}\left( u_{0}|b\right) >0$
for all $b$ such that $p_{1}\left( b\right) >0.$ Then, there exists a
measure $\widetilde{Q}$ over $\mathcal{U}\times \mathcal{B}$ dominated by $\lambda
\times \mu $ with a density $\widetilde{q}:\mathcal{U}\times \mathcal{B}$ $\rightarrow 
\mathbb{R}$ satisfying $p_{2}\left( b|u_{0}\right) =\widetilde{q}\left( b|u_{0}\right) $
and $p_{1}\left( u_{0}|b\right) =\widetilde{q}\left( u_{0}|b\right) $ for all $b\in 
\mathcal{B}$ such that $p_{1}\left( b\right) >0$.
\end{enumerate}
\end{slemma}

We now turn to scenario (iii). As previously discussed, the alignments in this scenario lead to a strictly semiparametric observed data model $\mc{P}$ even though $\mc{Q}$ is nonparametric. The following proposition gives a closed form expression for the efficient influence function in the fused-data frameworks of part \ref{item:prop:example-ifs4} in \cref{prop:example-ifs} above where we additionally assume $U$ is discrete, which includes the fused-data framework of scenario (iii) as a special case.
\begin{sproposition}
    \label{prop:example-eifs}
    Let $\left( \mathcal{Q},\mathcal{P},%
    \mathcal{C}\text{\thinspace },\psi ,\varphi \right) $ be a fused-data
    framework. Let $P\in \mathcal{P}
    $ and suppose $\varphi$ is pathwise differentiable at $P$ in model $\mc{P}$. Suppose there exists $Q$ in $\mc{Q}$ such that $\left(
    Q,P\right) $ is strongly aligned and $\psi $
    is pathwise differentiable at $Q$ in model $\mathcal{Q}$. Suppose $(\mc{Q}, \mc{P}, \mc{C}, \psi, \varphi)$ is as in part \ref{item:prop:example-ifs4} of \cref{prop:example-ifs} except where $U$ is finite discrete whose support takes $T\geq 2$ values. Suppose there exists product measure $\mu$ that dominates $Q$ and let $q = \frac{dQ}{d\mu}$ and suppose there exist product measures $\mu^{(j)}$ that dominate $P(\cdot|S=j)$ and let $p(\cdot|S=j) = \frac{dP(\cdot|S=j)}{d\mu^{(j)}}$ for $j\in\{1, 2\}$. Suppose there exists $u_0$ such that $\frac{1}{q(u_0|B)}\leq \widetilde{\delta}$ a.e.-$Q$ for some $\widetilde{\delta}<\infty$. Let $\psi^1_Q$ be the influence function of $\psi$ at $Q$ in model $\mc{Q}$. Define $r_{Q, P}(u, b) \coloneqq \frac{q(u, b)}{p(u, b)}\psi^1_Q(u, b)$. Let 
        \begin{align*}
            \underline{r}_{Q, P}(b) \coloneqq& \left(r_{Q, U}(u_1, b), \dots, r_{Q, U}(u_T, b)\right)^T\\
            \underline{\pi}(b) \coloneqq& \left(p(S=2|b, u_1), \dots, p(S=2|b, u_T)\right)^T\\
            \underline{\beta}(b) \coloneqq& \left(p(U=u_1|b, S=1), \dots, p(U=u_T|b, S=1)\right)^T,
        \end{align*}
        $diag\{\underline{1}- \underline{\pi}(b)\}$ be the $T\times T$ diagonal matrix with $t^{th}$ diagonal element equal to $1 - \underline{\pi}(b)_t$, $R(U)$ be the $T\times T$ diagonal matrix with $t^{th}$ diagonal element equal to $I(U=u_t)P(U=u_t|S=2)^{-1}$, $Id$ be the $T\times T$ identity matrix, and for any matrix $D$, $D^{-}$ is a generalized inverse of $D$. Then, the efficient influence function $\varphi^1_{P, eff}$ of $\varphi$ at $P$ in model $\mc{P}$ is
        \begin{align}
            \label{eq:example-3-eif-general-form}\varphi^1_{P, eff}(o) =& I(s=1)\{h_{P, eff}^{(Q)}(u, b) -E_P[h_{P, eff}^{(Q)}(U, B)|b, S=1]\}\\
            \nonumber&+ I(s=2)\{h_{P, eff}^{(Q)}(u, b) -E_{P}[h_{P, eff}^{(Q)}(U, B)|u, S=2]\}
        \end{align}
        with $\left(h_{P, eff}^{(Q)}(u_1, b), \dots, h_{P, eff}^{(Q)}(u_T, b)\right)^T \eqqcolon \underline{h}_{P, eff}^{(Q)}(b)$ where 
        \begin{align*}
            \underline{h}_{P,eff}^{\left( Q\right) }\left( b\right) =&\left\{ Id-%
            \underline{\pi }\left( b\right)\underline{\beta} \left( b\right) ^{\prime
            }\right\} ^{-1}\times \Big\{\underline{r}_{Q,P}\left( b\right) + \\
            &diag\left\{ \underline{1}-\underline{\pi }\left( b\right) \right\} \times
            \left\{ Id-E_{P}\left[ R\left( U\right) \left\{ Id-\underline{\pi }\left(
            B\right) \underline{\beta} \left( B\right) ^{\prime }\right\}
            ^{-1}diag\left\{ \underline{1}-\underline{\pi }\left( B\right) \right\} |S=2%
            \right] \right\} ^{-} \\
            &\times E_{P}\left[ R\left( U\right) \left\{ Id-\underline{\pi }%
            \left( B\right) \underline{\beta} \left( B\right) ^{\prime }\right\} ^{-1}%
            \underline{r}_{Q,P}\left( B\right) \right] \Big\}.
        \end{align*}
\end{sproposition}

\subsection{Causal Assumptions in \texorpdfstring{\cref{example:transporting}}{Example 3}}
\label{app-subsec:causal-identification}

In this subsection, we discuss the causal assumptions made for the various scenarios of \cref{example:transporting} that justify the assumed alignments made in that example. To give a precise definition of the average causal effect of interest, we
assume that, for each unit in the union of both populations, there exists a
full-data vector $\left( Y\left( 0\right) ,Y\left( 1\right) ,A,L,T\right) $
where $Y\left( a\right) $ is the counterfactual outcome if, possibly
contrary to fact, treatment were set to $a$, $a=0,1;$ $T=1$ if the unit is
from the source population and $T=2$ if the unit is from the target
population. We let $H_{0}$ denote the distribution of a random draw of the
full-data vector from the combined population, which we assume has a density 
$h_{0}$ with respect to some dominating product measure. The average treatment
effect in the $t$-th population, $t=1,2,$ a.k.a. the causal risk difference,
is defined as 
\begin{equation*}
ATE\left( t\right) =H_{0}\left[ Y\left( 1\right) =1|T=t\right] -H_{0}\left[
Y\left( 0\right) =1|T=t\right].
\end{equation*}%
Hereafter, we will make the consistency assumption

\begin{sassumption}[Consistency]
    \label{as:consistency}
    $Y=AY\left( 1\right) +\left(
1-A\right) Y\left( 0\right) $
\end{sassumption}

Assuming that participants in $T=1$ study are randomly selected from the source population, we have that
for all $l\in \mathbb{R}^{p},a,y\in \left\{ 0,1\right\} $ : 
\begin{equation}
P_{0}\left( L\leq l,A=a,Y=y|S=1\right) =H_{0}\left( L\leq
l,A=a,Y=y|T=1\right).  \label{SwithTzero}
\end{equation}

Table 1 gives the definition of $c\left( W,2\right) \coloneqq W^{\left(
2\right) }$ in each study design scenario from the target population, along
with the alignments with respect to $H_0$ ensured by the sampling designs. 
\begin{table}[ht]
\caption{Coarsening and alignment structure for study two ($S=2$)}
\label{tab:example-4-setup}\centering
\resizebox{\textwidth}{!}{\begin{tabular}{lll}
\toprule \textbf{Scenario} & $c\left( W,2\right) \coloneqq W^{\left(
2\right) }$ & \textbf{Alignments} (for all $l\in\mathbb{R}^p, a, y\in\{0, 1\}$)\\ 
\midrule (i) Random sample of $L$ & $L$ & $P_{0}\left( L\leq l|S=2\right)
=H_{0}\left( L\leq l|T=2\right)$ \\ 
(ii) Random sample of cases & $\left( L,A,Y=1\right) $ & $P_{0}\left(
L\leq l,A=a|Y=1,S=2\right) =H_{0}\left( L\leq l,A=a|Y=1,T=2\right) $ \\ 
(iii) Case-control study & $\left( L,A,Y\right) $ & $P_{0}\left(
L\leq l,A=a|Y=y,S=2\right) =H_{0}\left(L\leq l, A=a|Y=y,T=2\right)$ \\ 
\bottomrule &  & 
\end{tabular}}
\end{table}

We now discuss assumptions that suffice to identify $ATE(2)$, the average
treatment effect in the target population, from the observed data. The
first scenario has been extensively studied (\cite%
{pearl_transportability_2011, rudolph_robust_2017, dahabreh_extending_2020,
dahabreh_towards_2020, shi_data_2023, li_efficient_2023}). As we shall see, in all three scenarios we can recast $%
ATE\left( 2\right) $ as the evaluation at a law $Q_{0}$ for $W$ of the
functional $\psi :\mathcal{Q}\rightarrow \mathbb{R}$ defined as 
\begin{equation}
\psi \left( Q\right) \coloneqq E_{Q}\left[ Q\left( Y=1|L,A=1\right) -Q\left(
Y=1|L,A=0\right) \right]  \label{eq:psi-example3-appendix}
\end{equation}%
where $\mathcal{Q}$ is a collection of probability laws on $W$ such that $%
\mathcal{T}\left( \mathcal{Q};Q_{0}\right) =L_{0}^{2}\left( Q_{0}\right) $
and the definition of $Q_{0}$ varies depending on the scenario.

\subsubsection{Scenario (i)}

\cite{dahabreh_towards_2020} showed that under \cref{as:consistency} and Assumptions \ref{as:no-unmeasured-confounding-in-source} - \ref{as:additive-effect-exchangeability} below, the average treatment effect
in the target population satisfies 
\begin{equation}
ATE\left( 2\right) =E_{H_{0}}\left[ H_{0}\left( Y=1|L,A=1,T=1\right)
-H_{0}\left( Y=1|L,A=0,T=1\right) |T=2\right].  \label{a4}
\end{equation}

\begin{sassumption}[No unmeasured confounding in source]
\label{as:no-unmeasured-confounding-in-source}
$H_{0}\left( Y\left( j\right) =1|A=a,L,T=1\right) =H_{0}\left( Y\left(
j\right) =1|L,T=1\right) $ a.e.- $H_{0}\left( \cdot |T=1\right) ,a,j\in
\left\{ 0,1\right\} $.
\end{sassumption}

\begin{sassumption}[Treatment positivity in source]
    \label{as:treatment-positivity-in-source}
    $%
    0<H_{0}\left( A=1|L,T=1\right) <1,$ a.e.\medskip - $H_{0}\left( \cdot
    |T=1\right) $.
\end{sassumption}

\begin{sassumption}[Absolute continuity]
    \label{as:absolute-continuity}
    $H_{0}\left( L\in
B|T=1\right) =0\Rightarrow H_{0}\left( L\in B|T=2\right) =0$ for any Borel
set $B$ of $\mathbb{R}^{p}$.
\end{sassumption}

\begin{sassumption}[Additive effect exchangeability]
    \label{as:additive-effect-exchangeability}
    \begin{align*}
        E_{H_{0}}%
\left[ Y\left( 1\right) -Y\left( 0\right) |L,T=1\right] =E_{H_{0}}\left[
Y\left( 1\right) -Y\left( 0\right) |L,T=2\right]
    \end{align*}
    a.e.-$H_{0}\left( \cdot |T=2\right)$.
\end{sassumption}

Let $Q_{0}$ be any law of $W=\left( L,A,Y\right) $ such that $%
\mathsf{Supp}\left[ \left( L,A\right) ;Q_{0}\right] =\mathsf{Supp}\left[
\left( L,A\right) ;H_{0}\left( \cdot |T=1\right) \right] $ and such that $%
Q_{0}\left( Y=1|L,A\right) =H_{0}\left( Y=1|L,A,T=1\right) $ a.e.- $Q_{0}$
and $Q_{0}\left( L\leq l\right) =H_{0}\left( L\leq l|T=2\right) $ for all $%
l\in \mathbb{R}^{p}$. We can express the right hand side of $\left( \ref{a4}%
\right) $ as $\psi \left( Q_{0}\right) $ for the functional $\psi $ defined
in \eqref{eq:psi-example3-appendix}. In addition, because the identity $\left( \ref%
{SwithTzero}\right) $ and the alignment in the first row of \cref%
{tab:example-4-setup} hold, then for $Q=Q_{0}$ and $P=P_{0}$ it holds that 
\begin{align*}
Q\left( Y=1|L,A\right) =&P\left( Y=1|L,A,S=1\right) \text{ a.e.- }Q\\
Q\left( L\leq l\right) =&P\left( L\leq l|S=2\right) \text{ for all }l\in 
\mathbb{R}^{p}
\end{align*}%
with the support of $\left( L,A\right) $ in source 1 equal to the support of 
$\left( L\,A\right) $ under $Q_{0}$ as stated in \cref{subsec:overview}. These are exactly the alignment assumptions \eqref{eq:scenario-1-align-1} and \eqref{eq:scenario-1-align-2} respectively in \cref{sec:inferential-problem} in the main text. In this example, alignment \eqref{eq:scenario-1-align-1} is justified on the basis of the structural assumptions \ref{as:consistency}-\ref{as:additive-effect-exchangeability} and alignment \eqref{eq:scenario-1-align-2} holds because of simple random
sampling from source 2.

\subsubsection{Scenario (ii)}

In this setting $h_{0}\left( L|T=2\right) $ is not
identified because we lack access to a random sample of $L$ from the target
population $T=2,$ so Assumptions \ref{as:consistency} - \ref{as:additive-effect-exchangeability} do not suffice to identify ATE(2)
from the observed data. However, suppose that instead of Assumptions \ref{as:no-unmeasured-confounding-in-source} -
\ref{as:additive-effect-exchangeability} we now assume:

\begin{sassumption}[No unmeasured confounding in target]
\label{as:no-unmeasured-confounding-in-target}
$\left(
Y\left( 0\right) ,Y\left( 1\right) \right) \bot A|L,T=2$
\end{sassumption}
\begin{sassumption}[Treatment positivity in target]
    \label{as:treatment-positivity-in-target}
    $0<H_{0}\left( A=1|L,T=2\right) <1,$ a.e.-$H_{0}\left( \cdot |T=2\right)$.
\end{sassumption}

\begin{sassumption}[Strong absolute continuity]
    \label{as:strong-absolute-continuity}
    $h_{0}\left(L,A|T=2\right) <<h_{0}\left( L,A|T=1\right) $.
\end{sassumption}
\begin{sassumption}[Equal conditional prevalence in source and
target]
    \label{as:equal-conditional-prevalence-in-source-and-target}
    \begin{equation*}
H_{0}\left( Y=1|L,A,T=1\right) =H_{0}\left( Y=1|L,A,T=2\right) \,\ \text{%
a.e.- }H_{0}\left( \cdot |T=2\right).
\end{equation*}
\end{sassumption}
Defining $Q_{0}$ to be equal to the law of $W$ under $H_{0}\left( \cdot
|T=2\right)$, under Assumptions \ref{as:consistency}, \ref{as:no-unmeasured-confounding-in-target} and \ref{as:treatment-positivity-in-target} ATE(2) is
equal to $\psi \left( Q_{0}\right) $ defined as in \eqref{eq:psi-example3-appendix}. Furthermore, the substantive assumptions \ref{as:strong-absolute-continuity} and \ref{as:equal-conditional-prevalence-in-source-and-target} and the
alignment assumption $\left(\ref{SwithTzero}\right)$  implies that the alignment assumption \eqref{eq:scenario-2-align-1} in \cref{sec:inferential-problem} holds with $Q=Q_0$ and $P=P_0$, i.e. that
\begin{align*}
    Q(Y=1|L, A) = P(Y=1|L, A, S=1) \text{ a.e.-$Q$}.
\end{align*}

Additionally, under these assumptions the support of $\left( L,A\right) $ in
source 1 includes the support of $\left( L,A\right) $ under $Q_{0}$.
Finally, under the alignment in the second row of Table 1 (justified
by the random sampling of cases), the alignment assumption \eqref{eq:scenario-2-align-1} in \cref{sec:inferential-problem} holds, i.e. that
\begin{equation*}
Q\left( L\leq l,A=a|Y=1\right) =P\left( L\leq l,A=a|Y=1,S=2\right) \text{for all $l\in \mathbb{R}%
^{p}$ and $a\in \left\{ 0,1\right\}$ a.e.-$Q$}. 
\end{equation*}

\subsubsection{Scenario (iii)}

Due to the biased sampling of the outcome $Y$, the
available data constitutes a sample of $W^{\left( 2\right) }=\left(
L,A,Y\right) $, drawn from a law $P_{0}\left( \cdot |S=2\right) $ satisfying
that for all $l\in \mathbb{R}^{p}$ and $a\in \left\{ 0,1\right\} ,y\in
\left\{ 0,1\right\} $ 
\begin{equation}
P_{0}\left( L\leq l,A=a|Y=y,S=2\right) =H_{0}\left( L\leq
l,A=a|Y=y,T=2\right)  \label{cc1}
\end{equation}%
but such that $P_{0}\left( Y=1|S=2\right) $ is not equal to $H_{0}\left(
Y=1|T=2\right)$. Suppose that, as in scenario (ii), we make \cref{as:consistency} and assumptions \ref{as:no-unmeasured-confounding-in-target}-\ref{as:equal-conditional-prevalence-in-source-and-target}. Then, just as in scenario (ii), the causal risk difference is
equal to the right hand side of $\left( \ref{idenphi}\right) $ because $%
\left( \ref{cc1}\right) $ holds for $Y=1$. There is, however, an important
distinction with scenario (ii) in that now, Assumptions \ref{as:strong-absolute-continuity} and \ref{as:equal-conditional-prevalence-in-source-and-target}, and the
alignments $\left( \ref{SwithTzero}\right) $ and $\left( \ref{cc1}\right) $
impose restrictions on the law $P_{0}$ of the observed data. Specifically,
these assumptions imply that $\left( \ref{cc1}\right) $ holds, and
additionally, the following identity holds for $y\in \left\{ 0,1\right\} $ 
\begin{equation}
P_{0}\left( Y=y|L,A,S=1\right) =H_{0}\left( Y=y|L,A,T=2\right) \text{ a.e. }%
H_{0}\left( \cdot |T=2\right). \label{cc2}
\end{equation}

Assume additionally that $0<H_{0}\left( Y=1|L,A,T=2\right) <1$ a.e.- $%
H_{0}\left( |T=2\right)$. Then the model $\mathcal{P}$ is
semiparametric in the sense that $\mathcal{T}\left( \mathcal{P},P_{0}\right)
\varsubsetneqq L_{0}^{2}\left( P_{0}\right)$ because the mere fact that a single law exists that aligns the conditionals $p_{0}\left( y|l,a,S=1\right) $ and $p_{0}\left(
l,a|y,S=2\right) $ implies equality constraints for certain components of
the observed data law $P_{0}$. To see this, observe that
under $\left( \ref{cc1}\right)$ and $\left( \ref{cc2}\right)$, 
\begin{eqnarray*}
h_{0}\left( l,a|T=2\right) &=&\left. \left\{ \frac{p_{0}\left(
l,a|Y=0,S=2\right) }{P_{0}\left( Y=0|l,a,S=1\right) }\right\} \right/
\left\{ \sum_{a^{\prime }=0}^{1}\int \frac{p_{0}\left( l,a^{\prime
}|Y=0,S=2\right) }{P_{0}\left( Y=0|l,a^{\prime },S=1\right) }dl\right\}
 \\
&=&\left. \left\{ \frac{p_{0}\left( l,a|Y=1,S=2\right) }{P_{0}\left(
Y=1|l,a,S=1\right) }\right\} \right/ \left\{ \sum_{a^{\prime }=0}^{1}\int 
\frac{p_{0}\left( l,a^{\prime }|Y=1,S=2\right) }{P_{0}\left( Y=1|l,a^{\prime
},S=1\right) }dl\right\}.  \notag
\end{eqnarray*}%
The second equality is an equality restriction on $P_{0}.$ As in scenario
(ii), letting $Q_{0}\left( \cdot \right) $ to be equal to $H_{0}\left( \cdot
|T=2\right) ,$ we can recast the target of inference as $\psi \left(
Q_{0}\right) $ defined in \eqref{eq:psi-example3-appendix}. Then, under these assumptions with $Q = Q_0$ and $P=P_0$, aligments \eqref{eq:scenario-3a-align-1} and \eqref{eq:scenario-3a-align-2} hold, i.e. 
\begin{gather*}
    Q(Y=1|L, A) = Q(Y=1|L, A) \text{ a.e.-$Q$}\\
    Q\left( L\leq l,A=a|Y\right) =P\left( L\leq l,A=a|Y,S=2\right)
\end{gather*}
for all $l \in \mathbb{R}^p$, $a\in \{0, 1\}$ a.e.-$Q$. Notice that the second
alignment is strictly stronger than alignment \eqref{eq:scenario-2-align-2} in scenario (ii). 

\subsubsection{Scenario (iv)}

We now consider scenario (iv) as in \cref{app-subsubsec:generalizing-ex-3}. As with scenario (iii), the available data constitutes a sample of $W^{\left( 2\right) }=\left(L,A,Y\right) $, drawn from a law $P_{0}\left( \cdot |S=2\right) $ satisfying
that for all $l\in \mathbb{R}^{p}$ and $a\in \left\{ 0,1\right\} ,y\in
\left\{ 0,1\right\} $ 
\begin{equation}
P_{0}\left( L\leq l,A=a|Y=y,S=2\right) =H_{0}\left( L\leq
l,A=a|Y=y,T=2\right)  \label{cc-iv}
\end{equation}%
but such that $P_{0}\left( Y=1|S=2\right) $ is not equal to $H_{0}\left(
Y=1|T=2\right)$. Suppose now that $L$ is discrete. Suppose we make
Assumptions \ref{as:consistency}, \ref{as:no-unmeasured-confounding-in-target}, and \ref{as:treatment-positivity-in-target} but we replace Assumptions \ref{as:strong-absolute-continuity} and \ref{as:equal-conditional-prevalence-in-source-and-target} with the significantly weaker substantive assumption
\begin{sassumption}
    \label{as:scenario-iiib-weakening}
    $H_{0}\left( Y=1|L=l_{0},A=0,T=1\right) =H_{0}\left(
Y=1|L=l_{0},A=0,T=2\right)$.
\end{sassumption}
This assumption states that only the prevalence among the unexposed with
covariate level $l_{0}$ agree between the source and target population. Once
again, letting $Q_{0}\left( \cdot \right) $ to be equal to $H_{0}\left(
\cdot |T=2\right) ,$ we have that under \ref{as:consistency}, \ref{as:no-unmeasured-confounding-in-target}, and \ref{as:treatment-positivity-in-target}, ATE(2) is equal
to $\psi (Q_{0})$ with $\psi (Q)$ defined as in \eqref{eq:psi-example3-appendix}. Furthermore, alignment $\left( \ref{SwithTzero}\right) $ and
Assumption \ref{as:scenario-iiib-weakening} imply, for $Q=Q_{0}$ and $P=P_{0}$, alignments \eqref{eq:scenario-4-align-1}, i.e. that
\begin{equation*}
Q\left( Y=1|L=l_{0},A=0\right) =P\left( Y=1|L=l_{0},A=0,S=1\right).
\end{equation*}%
Additionally, by definition of $Q_{0},$ $\left( \ref{cc-iv}\right) $ for $Q=Q_{0}$ and $P=P_{0}$ is the
same as alignment \eqref{eq:scenario-4-align-2}, i.e. that 
\begin{equation*}
Q\left( L\leq l,A=a|Y\right) =P\left( L\leq l,A=a|Y,S=2\right) \text{ a.e.-$Q$}
\end{equation*}
for all $l\in \mathbb{R}^{p},a\in \left\{0,1\right\}$. 

\subsection{Defining the fused-data frameworks}
\label{app-subsec:define-fdf-examples}
Here we provide a rigorous definition of each fused-data framework in Examples \ref{example:disease-prevalence}-\ref{example:transporting}. 

\subsubsection{Example 1 (Continuation)}
Model $\mathcal{Q}$ is the collection of all laws on $W=\left( X,V,Y\right) $
that are mutually absolutely continuous and dominated by some measure $\mu $
where $Y$ and $V$ are binary and where $Y, V$ are correlated given $X$ a.e.-$Q$ for all $Q\in \mc{Q}$.
The collection $\mathcal{C}$ is comprised of 
\begin{equation*}
\left\{ \left( W^{\left( 1\right) },\left\{Z_1^{(1)}\right\},%
\{\overline{\mathcal{Z}}_{0}^{\left( 1\right) }\}\right) ,\left( W^{\left(
2\right) },\left\{ Z%
_{1}^{\left( 2\right) },Z_{2}^{\left( 2\right)
}\right\},\left\{ \overline{\mathcal{Z}}%
_{0}^{\left( 2\right) },\overline{\mathcal{Z}}_{1}^{\left( 2\right)
}\right\} \right) \right\}
\end{equation*}%
where $W^{\left( 1\right) }=Z_{1}^{\left( 1\right)
}=\left( X,V\right)$ and $\overline{\mathcal{Z}}_{0}^{\left(
1\right) }=\left\{ \ast \right\};$ $W^{\left( 2\right) }=\left(
X,V,Y\right)$, $Z_1^{(2)} = (X, Y)$, $Z_2^{(2)} = (V)$, $\overline{%
\mathcal{Z}}_{0}^{\left( 2\right) }=\emptyset$, and $\overline{\mathcal{Z}}%
_{1}^{\left( 2\right) }=$ $\mathsf{Supp}\left[ \left( X,Y\right) ;\mathcal{Q}%
\right]$. Thus, a law $P$ is in $\mathcal{P}$ if and only if $P$ is
mutually absolutely continuous with $P_{0}$, there exists $Q\in \mathcal{Q}$
satisfying the alignments $\left( \ref{eq:prevalence-source1}\right) $ and $%
\left( \ref{eq:prevalence-source2}\right) $ for $v\in \left\{ 0,1\right\}
,y\in \left\{ 0,1\right\} $ and $x\in \mathbb{R}^{p}$ and such that \textsf{%
Supp}$\left[ \left( Y,X\right) ;\mathcal{Q}\right] =$\textsf{Supp}$\left[
\left( Y,X\right) ;P\left( \cdot |S=2\right) \right] $. The functionals $%
\psi $ and $\varphi $ are $\psi(Q) \coloneqq E_{Q}\left[ Y\right] $ and 
\begin{align*}
    \varphi(P) \coloneqq E_{P}\left[\frac{V-E_{P}\left[ V|Y=0,X, S=2\right] }{E_{P}\left[
V|Y=1,X, S=2\right] -E_{P}\left[ V|Y=0,X, S=2\right] }\Big|S=1\right].
\end{align*}

\subsubsection{Example 2 (Continuation)}
Model $\mathcal{Q}$ is the collection of all laws on $W=\left( L,X,Y\right) $
that are mutually absolutely continuous and dominated by some product
measure such that there exist unique scalars $\alpha \left( Q\right) $ and $%
\psi \left( Q\right) $ solving the moment equation $\left( \ref{a1}\right) $
a.e.-$Q$. The collection $\mathcal{C}$ is comprised of 
\begin{equation*}
\left\{ \left( W^{\left( 1\right) },\left\{ Z_{1}^{\left( 1\right) },Z_{2}^{\left( 1\right) }\right\},\left\{ \overline{\mathcal{Z}}_{0}^{\left( 1\right) },\overline{\mathcal{Z}%
}_{1}^{\left( 1\right) }\right\} \right) ,\left( W^{\left( 2\right) },\left\{ Z_{1}^{\left(
2\right) },Z_{2}^{\left( 2\right) }\right\},\left\{ \overline{\mathcal{Z}}_{0}^{\left(
2\right) },\overline{\mathcal{Z}}_{1}^{\left( 2\right) }\right\} \right)
\right\}
\end{equation*}%
where $W^{\left( 1\right) }=\left( L,Y\right)$, $Z_{1}^{\left( 1\right)
}=L,Z_{2}^{\left( 1\right) }=Y$, $\overline{%
\mathcal{Z}}_{0}^{\left( 1\right) }=\emptyset $, and $\overline{\mathcal{Z}}%
_{1}^{\left( 1\right) }=$ \textsf{Supp} $\left[ L;\mathcal{Q}\right]$; $%
W^{\left( 2\right) }=\left( L,X\right)$, $Z_{1}^{\left( 2\right)
}=L,Z_{2}^{\left( 2\right) }=X$, $\overline{%
\mathcal{Z}}_{0}^{\left( 2\right) }=\emptyset $, and $\overline{\mathcal{Z}}%
_{1}^{\left( 2\right) }=$ \textsf{Supp} $\left[ L;\mathcal{Q}\right] .$
Thus, a law $P$ is in $\mathcal{P}$ if and only if $P$ is mutually
absolutely continuous with $P_{0},$ and there exists $Q\in \mathcal{Q}$
satisfying the alignments $\left( \ref{a2}\right) $ for all $y\in \mathbb{R}$
and $x\in \mathbb{R}$ and such that \textsf{Supp} $\left[ L;\mathcal{Q}%
\right] \subseteq $\textsf{Supp} $\left[ L;P\left( \cdot |S=1\right) \right] 
$ and \textsf{Supp} $\left[ L;\mathcal{Q}\right] \subseteq $\textsf{Supp} $%
\left[ L;P\left( \cdot |S=2\right) \right] .$ The functional $\psi \left(
Q\right) $ is the unique solution to $\left( \ref{a1}\right) $ and $\varphi
\left( P\right) $ is the unique solution to 
\begin{equation*}
\text{ }E_{P}\left[ Y|L,S=1\right] -\tau -\phi E_{P}\left[ X|L,S=2\right]
=0\text{ a.e.- }Q.
\end{equation*}%
In this example, $\xi \left( Q;\mathcal{C}\right) $ contains more than one
element for each $Q\in \mathcal{Q}$, because the marginal distribution of $L$
under $Q$ remains unrestricted.

\subsubsection{Example 3 (Continuation)}\

\textbf{Scenario (i).} The model $\mathcal{Q}$ is the collection of all laws
on $W=\left( L,A,Y\right) $ that are mutually absolutely continuous and
dominated by some product measure $\mu $. The collection $\mathcal{C}$ is comprised
of 
\begin{equation*}
\left\{ \left( W^{\left( 1\right) },\left\{ Z_{1}^{\left( 1\right) },Z_{2}^{\left( 1\right) }\right\},\left\{ \overline{\mathcal{Z}}_{0}^{\left( 1\right) },\overline{\mathcal{Z}%
}_{1}^{\left( 1\right) }\right\} \right) ,\left( W^{\left( 2\right) },\left\{Z_{1}^{\left( 2\right)
}\right\},\left\{\overline{\mathcal{Z}}_{0}^{\left( 2\right)
}\right\}\right) \right\}
\end{equation*}%
where $W^{\left( 1\right) }=\left( L,A,Y\right)$, $Z_{1}^{\left( 1\right) }=\left( L,A\right) ,Z_{2}^{\left( 1\right)
}=Y$, $\overline{\mathcal{Z}}_{0}^{\left( 1\right) }=\emptyset$, and $\overline{%
\mathcal{Z}}_{1}^{\left( 1\right) }= \textsf{Supp}\left[ \left(
L,A\right) ;\mathcal{Q}\right]$; $W^{\left( 2\right) }=L$, $
Z_{1}^{\left( 2\right) }=L$, and $\overline{\mathcal{Z}}_{0}^{\left( 2\right)
}=\left\{ \ast \right\}.$ Thus, a law $P$ is in $\mathcal{P}$ if and only
if $P$ is mutually absolutely continuous with $P_{0},$ and there exists $%
Q\in \mathcal{Q}$ satisfying alignments \eqref{eq:scenario-1-align-1} and \eqref{eq:scenario-1-align-2} and such that \textsf{Supp} $\left[ \left(
L,A\right) ;\mathcal{Q}\right] =\textsf{Supp}\left[ \left( A,L\right)
;P\left( \cdot |S=1\right) \right]$. The functionals $\psi $ and $\varphi$
are defined as in displays \eqref{eq:psi-example3} and \eqref{phi-def-example3}. We draw special attention to the fact that
under solely these alignments the conditional distribution of $A$ given $L$
is not identified by the observed data law, so $\xi \left( Q;\mathcal{C}%
\right) $ contains more than one element for each $Q\in \mathcal{Q}$.

\textbf{Scenario (ii)} The model $\mathcal{Q}$ is the collection of all laws
on $W=\left( L,A,Y\right) $ that are mutually absolutely continuous and
dominated by some product measure $\mu$. The collection $\mathcal{C}$ is comprised
of 
\begin{equation*}
\left\{ \left( W^{\left( 1\right) },\left\{ Z_{1}^{\left( 1\right) },Z_{2}^{\left( 1\right) }\right\},\left\{ \overline{\mathcal{Z}}_{0}^{\left( 1\right) },\overline{\mathcal{Z}%
}_{1}^{\left( 1\right) }\right\} \right) ,\left( W^{\left( 2\right) },\left\{ Z_{1}^{\left(
2\right) },Z_{2}^{\left( 2\right) }\right\},\left\{ \overline{\mathcal{Z}}_{0}^{\left(
2\right) },\overline{\mathcal{Z}}_{1}^{\left( 2\right) }\right\} \right)
\right\}
\end{equation*}%
where $W^{\left( 1\right) }=\left( L,A,Y\right)$, $\overline{\mathcal{Z}}_{0}^{\left( 1\right) }=\emptyset $, $Z_{1}^{\left( 1\right) }=\left( L,A\right)$, $Z_{2}^{\left( 1\right)}=Y$, and $\overline{\mathcal{Z}}_{1}^{\left( 1\right) }=\mathsf{Supp}\left[ \left( L,A\right);\mathcal{Q}\right]$; $W^{\left( 2\right) }=\left( L,A,Y\right)$,$Z_{1}^{\left( 2\right) }=Y$, $Z_{2}^{\left( 2\right)
}=\left( L,A\right)$, $\overline{\mathcal{Z}}%
_{0}^{\left( 2\right) }=\emptyset $, and $\overline{\mathcal{Z}}_{1}^{\left(
2\right) }=\left\{ 1\right\}$. Thus, a law $P$ is in $\mathcal{P}$ if and
only if is mutually absolutely continuous with $P_{0},$ and there exists $%
Q\in \mathcal{Q}$ satisfying alignments \eqref{eq:scenario-2-align-1} and \eqref{eq:scenario-2-align-2} and such that \textsf{Supp} $\left[ \left(
L,A\right) ;\mathcal{Q}\right] \subseteq $\textsf{Supp} $\left[ \left(
A,L\right) ;P\left( \cdot |S=1\right) \right]$. Note that in this example, 
\textsf{Supp} $\left[ Y;P\left( \cdot |S=2\right) \right] =\left\{ 1\right\}$. The functional $\psi $ is as defined in \eqref{eq:psi-example3} and $\varphi(P)$ for $P\in \mc{P}$ is defined as in 
\eqref{idenphi} but with $P$ replacing $P_0$.

\textbf{Scenario (iii)} The model $\mathcal{Q}$ is the collection of all
laws on $W=\left( L,A,Y\right) $ that are mutually absolutely continuous and
dominated by some product measure $\mu $. The collection $\mathcal{C}$ is comprised
of 
\begin{equation*}
\left\{ \left( W^{\left( 1\right) },\left\{ Z_{1}^{\left( 1\right) },Z_{2}^{\left( 1\right) }\right\},\left\{ \overline{\mathcal{Z}}_{0}^{\left( 1\right) },\overline{\mathcal{Z}%
}_{1}^{\left( 1\right) }\right\} \right) ,\left( W^{\left( 2\right) },\left\{ Z_{1}^{\left(
2\right) },Z_{2}^{\left( 2\right) }\right\},\left\{ \overline{\mathcal{Z}}_{0}^{\left(
2\right) },\overline{\mathcal{Z}}_{1}^{\left( 2\right) }\right\} \right)
\right\}
\end{equation*}%
where $W^{\left( 1\right) }=\left( L,A,Y\right)$, $Z_{1}^{\left( 1\right) }=\left( L,A\right)$, $Z_{2}^{\left( 1\right)}=Y$, $\overline{\mathcal{Z}}_{0}^{\left( 1\right) }=\emptyset $, and $\overline{%
\mathcal{Z}}_{1}^{\left( 1\right) }=\mathsf{Supp}\left[ \left( L,A\right) ;%
\mathcal{Q}\right]$; $W^{\left( 2\right) }=\left( L,A,Y\right)$, $Z_{1}^{\left( 2\right) }=Y$, $Z_{2}^{\left( 2\right)
}=\left( L,A\right)$, $\overline{\mathcal{Z}}%
_{0}^{\left( 2\right) }=\emptyset $, and $\overline{\mathcal{Z}}_{1}^{\left(
2\right) }=\left\{ 0,1\right\}$. Thus, a law $P$ is in $\mathcal{P}$ if and
only if is mutually absolutely continuous with $P_{0}$ and there exists $Q\in \mathcal{Q}$ satisfying the alignments \eqref{eq:scenario-3a-align-1} and \eqref{eq:scenario-3a-align-2} and such that \textsf{Supp} $\left[ \left(
L,A\right) ;\mathcal{Q}\right] \subseteq $\textsf{Supp} $\left[ \left(
A,L\right) ;P\left( \cdot |S=1\right) \right]$. As argued earlier, $\mathcal{P}$ is a semiparametric model. The functional $\psi $ is as defined
in \eqref{eq:psi-example3} and $\varphi \left( P\right) $ for any $%
P\in \mathcal{P}$ is defined as in \eqref{idenphi} but with $%
P $ replacing $P_{0}$. Note that for any $P\in \mathcal{P}$ , $\varphi
\left( P\right) $ is also equal to 
\begin{equation*}
\int \left\{ P\left( Y=1|l,A=1,S=1\right) -P\left( Y=1|l,A=0,S=1\right)
\right\} \sum_{a=0}^1\beta\left(a, l, 0;P\right) dl.
\end{equation*}

\textbf{Scenario (iv)}

The model $\mathcal{Q}$ is the collection of all laws on $W=\left(
L,A,Y\right) $ that are mutually absolutely continuous and dominated by some
product measure $\mu $. The collection $\mathcal{C}$ is comprised of 
\begin{equation*}
\left\{ \left( W^{\left( 1\right) },\left\{ Z_{1}^{\left( 1\right) },Z_{2}^{\left( 1\right) }\right\},\left\{ \overline{\mathcal{Z}}_{0}^{\left( 1\right) },\overline{\mathcal{Z}%
}_{1}^{\left( 1\right) }\right\} \right) ,\left( W^{\left( 2\right) },\left\{ Z_{1}^{\left(
2\right) },Z_{2}^{\left( 2\right) }\right\},\left\{ \overline{\mathcal{Z}}_{0}^{\left(
2\right) },\overline{\mathcal{Z}}_{1}^{\left( 2\right) }\right\} \right)
\right\}
\end{equation*}%
where $W^{\left( 1\right) }=\left( L,A,Y\right)$, $Z_{1}^{\left( 1\right) }=\left( L,A\right)$, $Z_{2}^{\left( 1\right)}=Y$, $\overline{\mathcal{Z}}_{0}^{\left( 1\right) }=\emptyset $, and $\overline{%
\mathcal{Z}}_{1}^{\left( 1\right) }=\{l_0\}\times \{0\}$; $W^{\left( 2\right) }=\left( L,A,Y\right)$, $Z_{1}^{\left( 2\right) }=Y$, $Z_{2}^{\left( 2\right)
}=\left( L,A\right)$, $\overline{\mathcal{Z}}%
_{0}^{\left( 2\right) }=\emptyset $, and $\overline{\mathcal{Z}}_{1}^{\left(
2\right) }=\left\{ 0,1\right\}$. Thus, a law $P$ is in $%
\mathcal{P}$ if and only if is mutually absolutely continuous with $P_{0},$
and there exists $Q\in \mathcal{Q}$ satisfying the alignments \eqref{eq:scenario-4-align-1} and \eqref{eq:scenario-4-align-2}. These alignments do not
impose equality constraints on the coarsened data law. The functional $\psi $
is as defined in \eqref{eq:psi-example3} and $\varphi$ is defined as in \eqref{eq:phi-scenario-4}.

\subsection{Additional Derivations}
\label{app-subsec:additional-derivations}

\subsubsection{Example 1}
\label{app-subsubsec:additional-derivations-example1}
Let 
\begin{align*}
    m_Q(x, v) \coloneqq \frac{v-E_{Q}\left[ V|Y=0,X=x\right] }{E_{Q}\left[
V|Y=1,X=x\right] -E_{Q}\left[ V|Y=0,X=x\right] }
\end{align*}
as in \eqref{eq:mq}. We first demonstrate $E_Q[m_Q(X, V)|x, y] = y$. 
\begin{align*}
    E_Q[m_Q(X, V)|x, y] =& E_Q\left[\left. \frac{v-E_{Q}\left[ V|Y=0,X\right] }{E_{Q}\left[
V|Y=1,X\right] -E_{Q}\left[ V|Y=0,X\right]}\right|x, y\right]\\
=& \frac{E_{Q}\left[ V|y,x\right]-E_{Q}\left[ V|Y=0,x\right] }{E_{Q}\left[
V|Y=1,x\right] -E_{Q}\left[ V|Y=0,x\right]}\\
=&y.
\end{align*}
Then, by the tower law, $\psi(Q) = E_Q[Y] = E_Q[E_Q[m_Q(X, V)|X, Y]] = E_Q[m_Q(X, V)]$.

We now demonstrate $Q$ is identified by $P$ under alignments \eqref{eq:prevalence-source1} and \eqref{eq:prevalence-source2}. Let 
\begin{align*}
    m_{P(\cdot|S=2)}(x, v) \coloneqq \frac{v-E_{P}\left[ V|Y=0,X=x, S=2\right] }{E_{P}\left[
V|Y=1,X=x, S=2\right] -E_{P}\left[ V|Y=0,X=x, S=2\right]}.
\end{align*}
The alignments imply that $m_{P(\cdot|S=2)} = m_{Q}$ a.e.-$Q$. We additionally know that $q(y=1|X) = E_Q[m_Q(X, V)|X]=E_P[m_{P(\cdot|S=2)}(X, V)|X, S=1]$ a.e.-$Q$. But then, because $q(x) = p(x|S=1)$ and $q(v|y, x) = p(v|y, x, S=2)$, $Q$ is identified by $P$. 

\subsubsection{Example 2}

We use \cref{prop:two-source-solution} to derive the observed data
influence functions for this example. Applying that proposition,
we aim to solve the integral equation 
\begin{align*}
& m^{(2)}(l,x)-E_{Q}[m^{(2)}(L,X)|Y=y,L=l]+E_{Q}[m^{(2)}(L,X)|L=l] \\
=& \gamma _{Q}^{1}(l,x,y)-E_{Q}[\gamma _{Q}^{1}(L,X,Y)|Y=y,L=l]+E_{Q}[\gamma
_{Q}^{1}(L,X,Y)|L=l] \\
=& B_{Q}(g)^{-1}g(l)\left\{ -\psi \left( Q\right) x+\psi \left( Q\right)
E_{Q}[X|Y=y,L=l]-E_{Q}[Y-\alpha \left( Q\right) -\psi \left( Q\right)
X|L=l]\right\} \\
=& -B_{Q}(g)^{-1}g(l)\left\{ \psi \left( Q\right) x-\psi \left( Q\right)
E_{Q}[X|Y=y,L=l]\right\} .
\end{align*}%
Taking expectations given $L$ on both sides gives that 
\begin{equation*}
E_{Q}[m^{(2)}(L,X)|Y=y,L=l]=B_{Q}(g)^{-1}g(l)\left\{ \psi \left( Q\right)
E_{Q}[X|L=l]-\psi \left( Q\right) E_{Q}[X|Y=y,L=l]\right\} .
\end{equation*}%
Combined with the fact that $E_{Q}[m^{(2)}(L,X)|L=l]=0$, we arrive at 
\begin{align*}
m^{(2)}(l,x)=& B_{Q}(g)^{-1}g(l)\left\{ \psi \left( Q\right)
E_{Q}[X|L=l]-\psi (Q)x\right\} \\
=& B_{Q}(g)^{-1}g(l)\left\{ E_{Q}[Y|L=l]-\alpha (Q)-\psi (Q)x\right\} .
\end{align*}%
Now, set $m^{(1)}(l,y)\coloneqq\gamma _{Q}^{1}-m^{(2)}\left( l,x\right)
=B_{Q}(g)^{-1}g(l)\left\{ y-E_{Q}[Y|L=l]\right\} $. It follows from parts \ref{item:two-source-solution1}
and \ref{item:two-source-solution2} of \cref{prop:two-source-solution} that $\nu $ is pathwise
differentiable at $P$ and 
\begin{align*}
\nu _{P}^{1}\left( o\right) =& B_{Q}\left( g\right) ^{-1}g(l)\frac{q\left(
l\right) }{p(l|S=2)}\Bigg[\frac{I(s=1)}{P(S=1)}\frac{p(l|S=2)}{p(l|S=1)}%
\{y-E_{Q}\left( Y|L=l\right) \} \\
& +\frac{I(s=2)}{P(S=2)}\{E_{Q}\left( Y|L=l\right) -\alpha \left( Q\right)
-\psi \left( Q\right) x\}\Bigg] \\
=& B_{P\left( \cdot |S=2\right) }\left( t_{g,q}\right) t_{g,q}\left(
l\right) \varepsilon _{P}(o)
\end{align*}%
is one of its influence functions, where $t_{g,q}\left( l\right) \coloneqq %
g(l)\frac{q\left( l\right) }{p(l|S=2)}$ and $\varepsilon _{P}(o)\coloneqq%
\frac{I(s=1)}{P(S=1)}\frac{p(l|S=2)}{p(l|S=1)}\{y-E_{P\left( \cdot|S=1\right)
}\left( Y|L=l\right) \}+\frac{I(s=2)}{P(S=2)}\{E_{P\left( \cdot|S=1\right)
}\left( Y|L=l\right) -\tau \left( P\right) -\varphi \left( P\right) x\}$.
The second equality follows from 
\begin{eqnarray*}
B_{Q}\left( g\right) &=&E_{Q}\left[ g(L)%
\begin{pmatrix}
1, & X%
\end{pmatrix}%
\right] \\
&=&E_{P}\left[ \left. g(L)\frac{q\left( L,X\right) }{p(L,X|S=2)}%
\begin{pmatrix}
1, & X%
\end{pmatrix}%
\right\vert S=2\right] \\
&=&E_{P}\left[ \left. g(L)\frac{q\left( L\right) }{p(L|S=2)}%
\begin{pmatrix}
1, & X%
\end{pmatrix}%
\right\vert S=2\right] \\
&=&B_{P\left( \cdot |S=2\right) }\left( t_{g,q}\right)
\end{eqnarray*}%
since $q\left( x|l\right) =p(x|l,S=2).$ Now, the sets $\left\{ g:B_{Q}\left(
g\right) \text{ is non-singular}\right\} $ and \\$\{ t:B_{P\left(\cdot
|S=2\right) }\left( t\right) \text{ is non-singular}\} $ are equal
because $g$ and $t_{g,q}$ are in one to one correspondence for any given $q.$
Furthermore, since the alignment assumptions alone place no restrictions on $%
P$ other than inequality constraints, $\mathcal{P}^{ext}$ is nonparametric
and so by \cref{Theorem:influence functions} part \ref{item:influence functions2}, each ideal data
influence function $\gamma _{Q}^{1}$ corresponds to a single observed data
influence function $\nu _{P}^{1}$. We then conclude that the set 
\begin{equation}
\left\{ \nu _{P}^{1}\left( o\right) =B_{P\left( \cdot|S=2\right) }\left( t\right)
^{-1}t\left( l\right) \varepsilon _{P}(o):t\text{ such that }B_{P\left(
\cdot|S=2\right) }\left( t\right) \text{ is non-singular}\right\}  \label{IFset}
\end{equation}%
comprises the set of all observed data influence functions of $\nu .$ Note
that although the specific observed data influence function $B_{P\left(
\cdot |S=2\right) }\left( t_{g,q}\right) t_{g,q}\left( l\right) \varepsilon
_{P}(o)$ corresponding to a particular ideal data influence function depends
on the marginal distribution of $L$ under $Q$, the set of all observed data
influence functions does not.

Following \cref{prop:only-Q-restricted-eif} in \cref{app:model-types-eif}, the efficient influence
function $\nu _{P,eff}^{1}$ is the element of the set $\left( \ref{IFset}%
\right) $ with the smallest variance. Letting $U\coloneqq$ $\frac{I(S=2)}{P(S=2)}(1, X)^\prime$ we can write 
\begin{equation*}
B_{P\left( \cdot|S=2\right) }\left( t\right) =E_{P}\left[ t(L)U^{\prime }\right].
\end{equation*}%
Thus, letting $\sigma ^{2}\left( L\right) \coloneqq var_{P}\left(
\varepsilon _{P}|L\right) $ we have that%
\begin{eqnarray*}
&&\left. var_{P}\left[ \nu _{P}^{1}\left( O\right) \right] =var_{P}\left[
B_{P\left( \cdot|S=2\right) }\left( t\right) t\left( L\right) \varepsilon _{P}%
\right] \right. \\
&=&E_{P}\left[ \left\{ t(L)\sigma \left( L\right) \right\} \left\{ \sigma
^{-1}\left( L\right) U\right\} ^{\prime }\right] ^{-1}E_{P}\left[ \sigma
^{2}\left( L\right) t(L)t\left( L\right) ^{\prime }\right] E_{P}\left[
\left\{ \sigma ^{-1}\left( L\right) U\right\} \left\{ t(L)\sigma \left(
L\right) \right\} ^{\prime }\right] ^{-1} \\
&\geq &E_{P}\left[ \sigma ^{-1}\left( L\right) U^{\prime }\right] ^{-1}
\end{eqnarray*}%
by the Cauchy-Schwartz inequality. The lower bound is then achieved at $%
t_{P,eff}\left( L\right) \coloneqq\sigma ^{-2}\left( L\right) E_{P}\left(
U|L\right) $ rendering the efficient influence function of $\nu $ at $P:$ 
\begin{equation*}
\nu _{P,eff}^{1}\left( o\right) =B_{P\left( \cdot|S=2\right) }\left(
t_{P,eff}\right) ^{-1}t_{P,eff}\left( l\right) \varepsilon _{P}(o).
\end{equation*}%
The second component of each $\nu _{P}^{1}$ and $\nu _{P,eff}^{1}$
corresponds to the observed data influence function and the efficient
influence function of $\varphi $. Our results agree with those of \cite%
{zhao_two-sample_2019} who derived a class of estimating equations whose
solutions are, up to asymptotic equivalence, all RAL estimators of $\varphi
. $ Thus, the collection of influence functions of all of their RAL
estimators of $\nu $ coincides with the set of influence functions derived
here, and the efficient one in their collection must therefore have
influence function equal to $\nu _{P,eff}^{1}$.

\subsubsection{Example 3}

\textbf{Scenario (i):} Here we compute the unique observed data influence function. Since the terms in the expression for $\psi _{Q}^{1}$ satisfy $E_{Q}%
[ E_{Q}\left( Y|L,A=1\right) -E_{Q}\left( Y|L,A=0\right) -\psi \left(
Q\right) ] =0$ and $E_{Q}\left[ \frac{2A-1}{q\left( A|L\right) }%
\left\{ Y-E_{Q}\left( Y|L,A\right) \right\} |A,L\right] =0\,$\ then, by part
\ref{item:pathwise-diff-3} of \cref{lemma:pathwise-differentiabilty} $\varphi $ is pathwise
differentiable at $P$ and, by part \ref{item:influence functions1} of \cref{Theorem:influence functions}
its unique influence function is%
\begin{eqnarray*}
\varphi _{P}^{1}(o) &=&\frac{I(s=1)}{P(S=1)}\frac{q(a,l)}{p(a,l|S=1)}\frac{%
2a-1}{q\left( a|l\right) }\left\{ y-E_{Q}\left( Y|a,l\right) \right\} \\
&&+\frac{I(s=2)}{P(S=2)}\left\{ E_{Q}\left( Y|A=1,l\right) -E_{Q}\left(
Y|A=0,l\right) -\psi \left( Q\right) \right\} \\
&=&\frac{I(s=1)}{P(S=1)}\frac{p(l|S=2)}{p(l|S=1)}\frac{2a-1}{p\left(
a|l,S=1\right) }\left\{ y-E_{P\left( \cdot |S=1\right) }\left( Y|a,l\right)
\right\} \\
&&+\frac{I(s=2)}{P(S=2)}\left\{ E_{P\left( \cdot |S=1\right) }\left(
Y|A=1,l\right) -E_{P\left( \cdot |S=1\right) }\left( Y|A=0,l\right) -\varphi
\left( P\right) \right\}
\end{eqnarray*}%
where the second equality holds because of the alignments in the assumed
fused-data framework.

\textbf{Scenario (ii):} We use \cref{prop:two-source-solution} to derive the unique observed data
influence function for this example. We aim to solve the integral equation 
\begin{align*}
& m^{(1)}(l,a,y)-y\left\{ m^{(1)}(l,a,y)+E_{Q}[m^{(1)}(L,A,Y)|Y=y]\right\} \\
=& \psi _{Q}^{1}(l,a,y)-y\left\{ \psi _{Q}^{1}(l,a,y)+E_{Q}[\psi
_{Q}^{1}(L,A,Y)|Y=y]\right\}
\end{align*}%
When $y=0$, the above reduces to $m^{(1)}(l,a,0)=\psi _{Q}^{1}(l,a,0)$. We
also know that $E_{Q}[m^{(1)}(L,A,Y)|L=l,A=a]=0$. This implies that 
\begin{equation*}
m^{(1)}(l,a,1)Q(Y=1|l,a)=-m^{(1)}(l,a,0)Q(Y=0|l,a)
\end{equation*}%
which in turn yields 
\begin{equation*}
m^{(1)}(l,a,1)=-\psi _{Q}^{1}(l,a,0)\frac{Q(Y=0|l,a)}{Q(Y=1|l,a)}
\end{equation*}%
We thus arrive at
\begin{equation*}
m^{(1)}(l,a,y)=\psi _{Q}^{1}(l,a,y)-\frac{y}{Q(Y=1|l,a)}E_{Q}[\psi
_{Q}^{1}(L,A,Y)|l,a].
\end{equation*}
One can readily show $m^{(1)} \in L^2(Q)$ because $\frac{1}{Q(Y=1|L, A)} \leq \widetilde{\delta}$.
Now, set $m^{(2)}(l,a,y)\coloneqq\psi _{Q}^{1}(l,a,y)-m^{(1)}(l,a,y)=\frac{y%
}{Q(Y=1|l,a)}E_{Q}[\psi _{Q}^{1}(L,A,Y)|l,a]$. It follows from parts \ref{item:two-source-solution1} and \ref{item:two-source-solution2}
of \cref{prop:two-source-solution} that $\varphi $ is pathwise
differentiable at $P$ in model $\mathcal{P}$ and 
\begin{align*}
\varphi _{P}^{1}\left( o\right) & =\frac{I\left( s=1\right) }{P(S=1)}\frac{%
q\left( l,a\right) }{p\left( l,a|S=1\right) }\left\{ \psi _{Q}^{1}(l,a,y)-%
\frac{y}{Q(Y=1|l,a)}E_{Q}\left[ \psi _{Q}^{1}\left( L,A,Y\right) |l,a\right]
\right\} \\
& +\frac{I\left( s=2\right) }{P(S=2)}\frac{Q(Y=1)}{P(Y=1|S=2)}\frac{y}{%
Q(Y=1|l,a)}E_{Q}\left[ \psi _{Q}^{1}\left( L,A,Y\right) |l,a\right]
\end{align*}%
is its unique influence function at $P$ because model $\mathcal{P}$ is
non-parametric. Noticing that $E_{Q}\left[ \psi _{Q}^{1}\left( L,A,Y\right)
|l,a\right] =E_{Q}\left( Y|l,A=1\right) -E_{Q}\left( Y|l,A=0\right) -\psi
\left( Q\right) ,$ after some algebra we arrive at the expression%
\begin{align}
\varphi _{P}^{1}\left( o\right) =& \frac{I\left( s=1\right) }{P(S=1)}\frac{%
q\left( l,a\right) \left\{ y-E_{Q}\left(
Y|a,l\right) \right\} }{p\left( l,a|S=1\right) E_{Q}\left( Y|a,l\right)}\left( 2a-1\right) \left[
\left\{ \frac{1}{q\left( a|l\right) }-1\right\} E_{Q}\left( Y|l,a\right)
+E_{Q}\left( Y|l,1-a\right) \right]  \notag \\
& +\frac{I\left( s=2\right) }{P(S=2)}\frac{Q(Y=1)}{P(Y=1|S=2)}\frac{y}{%
E_{Q}\left( Y|a,l\right) }\left\{ E_{Q}\left( Y|l,A=1\right) -E_{Q}\left(
Y|l,A=0\right) -\psi \left( Q\right) \right\} .  \label{IF3ii}
\end{align}

Replacing $q(l,a,y)$ in the right hand side of $\left( \ref{IF3ii}\right) $
with the right-hand side 
\begin{equation*}
q(l,a,y)=p(y|l,a,S=1)\frac{p(l,a|Y=1,S=2)}{p(Y=1|l,a,S=1)}\left\{
\sum_{a^{\prime }=0}^{1}\int \frac{p(l^{\prime },a^{\prime }|Y=1,S=2)}{%
p(Y=1|l^{\prime },a^{\prime },S=1)}dl^{\prime }\right\} ^{-1}
\end{equation*}
yields the expression of $\varphi _{P}^{1}$ as a function of the observed data law $P$.

\textbf{Scenario (iii):} We use \cref{prop:two-source-solution} to derive the observed data influence functions in this example. We aim to solve the
integral equation 
\begin{equation*}
E_{Q}[m^{(1)}(L,A,Y)|Y=y]=E_{Q}[\psi _{Q}^{1}(L,A,Y)|Y=y]
\end{equation*}%
where $m^{(1)}\in L^{2}(L,A,Y;Q)$ such that $E_{Q}[m^{(1)}(L,A,Y)|L,A]=0$
a.e.- $Q$. To solve this equation, we set 
\begin{align*}
    m^{(1)}(l, a, y) \coloneqq \psi _{Q}^{1}(l,a,y) - \frac{q(l, a)q(y)}{q(l, a, y)}E_Q[\psi _{Q}^{1}(L,A,Y)|l, a].
\end{align*}
Notice that $\frac{q(l, a)q(y)}{q(l, a, y)}$ is the Radon-Nikodym derivative of $Q^*$ with respect to $Q$ where $Q^*$ is the law with density $q^*(l, a, y) = q(l, a)q(y)$. Hence, under $Q^*$, $(L, A)$ and $Y$ are independent. Then, 
\begin{align*}
    E_{Q*}[E_Q[\psi^1_Q(L, A, Y)|L, A]|y] =& E_{Q^*}[\psi^1_Q(L, A, Y)|L, A]]\\
    =& E_{Q}[E_Q[\psi^1_Q(L, A, Y)|L, A]]\\
    =& E_{Q}[\psi^1_Q(L, A, Y)]\\
    =& 0
\end{align*}
where the first equality follows because $(L, A)$ and $Y$ are independent under $Q^*$ and the second because distribution of $L, A$ under $Q$ and $Q^*$ are equal. But then 
\begin{align*}
    E_Q[m^{(1)}(L, A, Y)|y] =& E_Q[\psi^1_Q(L, A, Y)|y] -E_{Q*}[E_Q[\psi^1_Q(L, A, Y)|L, A]|y] \\
    =& E_Q[\psi^1_Q(L, A, Y)|y]
\end{align*}
as desired. Additionally, 
\begin{align*}
    E_Q[m^{(1)}(L, A, Y)|l, a] =& E_Q[\psi^1_Q(L, A, Y)|l, a] -E_{Q*}[E_Q[\psi^1_Q(L, A, Y)|L, A]|l, a]\\
    =&E_Q[\psi^1_Q(L, A, Y)|l, a] -E_Q[\psi^1_Q(L, A, Y)|l,a]\\
    =& 0
\end{align*}
and so $m^{(1)}\in L^{2}(L,A,Y;Q)$ such that $E_{Q}[m^{(1)}(L,A,Y)|L,A]=0$. It follows from parts \ref{item:two-source-solution1} and \ref{item:two-source-solution2} of \cref%
{prop:two-source-solution} that $\varphi $ is pathwise differentiable at $P$
in model $\mathcal{P}$ and 
\begin{align*}
\varphi _{P}^{1}\left( o\right) =& \frac{I\left( s=1\right) }{P(S=1)}\frac{%
q\left( l,a\right) }{p\left( l,a|S=1\right) }\left\{ \psi _{Q}^{1}(l,a,y)-%
\frac{q\left( l,a\right) q\left( y\right) }{q\left( l,a,y\right) }E_{Q}\left[
\psi _{Q}^{1}\left( L,A,Y\right) |l,a\right] \right\}
\\
& +\frac{I\left( s=2\right) }{P(S=2)}\frac{q(y)}{p(y|S=2)}\frac{q\left(
l,a\right) q\left( y\right) }{q\left( l,a,y\right) }E_{Q}\left[ \psi
_{Q}^{1}\left( L,A,Y\right) |l,a\right]  \notag
\end{align*}%
is an influence function of $\varphi $ at $P$ in model $\mathcal{P}$.

Recall that the given alignments imply a strictly semiparametric model for $%
\mathcal{P}$, even though $\mathcal{Q}$ is non-parametric, so there exist
infinitely many observed data influence functions. Part \ref{item:two-source-solution3} of \cref%
{prop:two-source-solution} tells us we can derive the set of all
influence functions by adding to $\varphi _{P}^{1}\left( o\right) $ any
element of the set described in part \ref{item:two-source-solution3} of that proposition. For this
example, that set is equal to
\begin{equation*}
\left\{ r\left( l,a,y,s\right) =\left[ \frac{I\left( s=1\right) }{P(S=1)}%
\frac{q\left( l,a\right) }{p\left( l,a|S=1\right) }-\frac{I\left( s=2\right) 
}{P(S=2)}\frac{q(y)}{p(y|S=2)}\right] f\left( l,a,y\right) :f\in \mathcal{F}%
\right\}
\end{equation*}%
where%
\begin{equation*}
\mathcal{F}=\left\{ f\in L_{0}^{2}(Q):E_{Q}[f(L,A,Y)|L,A]=E_{Q}[f(L,A,Y)|Y]=0%
\text{ a.e. - }Q\right\} .
\end{equation*}
When $\mathsf{Supp}\left[ \left( L,A,Y\right) ;Q\right] =\mathsf{Supp}\left[
\left( L,A\right) ;Q\right] \times \mathsf{Supp}\left[ Y;Q\right] ,$
\cite{tchetgen_tchetgen_doubly_2010} showed that $f\in \mathcal{F}$ if and only if
there exists $t\in L_{0}^{2}(Q)$ such that 
\begin{align*}
f\left( L,A,Y\right) =\frac{q\left( L,A\right) q\left( Y\right) }{q\left(
L,A,Y\right) }\{& t\left( L,A,Y\right) -E_{Q^{\ast }}\left[ \left.
t\left( L,A,Y\right) \right\vert L,A\right] \\
&-E_{Q^{\ast }}\left[ \left.
t\left( L,A,Y\right) \right\vert Y\right] +E_{Q^{\ast }}\left[ t\left(
L,A,Y\right) \right] \} \text{ a.e. - }Q
\end{align*}%
where $E_{Q^{\ast }}$ denotes expectation under the law $Q^{\ast }$ with
density $q^{\ast }\left( l,a,y\right) \coloneqq q\left( l,a\right) q\left(
y\right) .$ Thus, the collection of all functions of the form 
\begin{align}
\left( l,a,y\right) & \mapsto \frac{I\left( s=1\right) }{P(S=1)}\frac{%
q\left( l,a\right) }{p\left( l,a|S=1\right) }\left\{ \psi _{Q}^{1}(l,a,y)-%
\frac{q\left( l,a\right) q\left( y\right) }{q\left( l,a,y\right) }E_{Q}\left[
\psi _{Q}^{1}\left( L,A,Y\right) |l,a\right] \right\}
\label{eq:example3iiia-if-f} \\
& +\frac{I\left( s=2\right) }{P(S=2)}\frac{q(y)}{p(y|S=2)}\left\{ \frac{%
q\left( l,a\right) q\left( y\right) }{q\left( l,a,y\right) }E_{Q}\left[ \psi
_{Q}^{1}\left( L,A,Y\right) |l,a\right] \right\}  \notag \\
& +\left[ \frac{I\left( s=1\right) }{P(S=1)}\frac{q\left( l,a\right) }{%
p\left( l,a|S=1\right) }-\frac{I\left( s=2\right) }{P(S=2)}\frac{q(y)}{%
p(y|S=2)}\right]  \notag \\
& \times \frac{q\left( l,a\right) q\left( y\right) }{q\left( l,a,y\right) }%
\left\{ t\left( l,a,y\right) -E_{Q^{\ast }}\left[ \left. t\left(
L,A,Y\right) \right\vert l,a\right] -E_{Q^{\ast }}\left[ \left. t\left(
L,A,Y\right) \right\vert y\right] +E_{Q^{\ast }}\left[ t\left( L,A,Y\right) %
\right] \right\}  \notag
\end{align}%
for any $t\in L_{0}^{2}(Q)$ comprises the set of all observed data influence
functions. Similarly to scenario (ii), in this example, $Q$ is determined by
the aligned conditionals as is seen by replacing the right hand side of $%
\left( \ref{q-formula}\right) $ with the aligned conditionals: 
\begin{equation}
q(l,a,y)=p(l,a|y,S=2)\frac{p(y|l_{0},a=0,S=1)}{p(l_{0},a=0|y,S=2)}\left\{
\sum_{y^{\prime }=0}^{1}\frac{p(y^{\prime }|l_{0},a=0,S=1)}{%
p(l_{0},a=0|y^{\prime },S=2)}\right\} ^{-1}  \label{idenq}
\end{equation}%
In fact, $\left( l_{0},a=0\right) $ in the right hand side of $\eqref{idenq}$
can be replaced by $\left( l^{\ast },a^{\ast }\right) $ for any $l^{\ast }$
in $\mathsf{Supp}\left( L;Q\right) $ and $a^{\ast }\in \left\{ 0,1\right\} .$
Replacing $q(l,a,y)$ with the right hand side of the last equality in the
right hand side of \eqref{eq:example3iiia-if-f} yields the expression of any
observed data influence function in terms of the observed data law $P$.

Turn now to the computation of the efficient influence function $\varphi
_{P,eff}^{1}.$ It follows from \cref{prop:example-eifs} that the observed data efficient influence function is 
\begin{align*}
\varphi _{P,eff}^{1}\left( o\right) =& I(s=1)\left\{ h^{\left( Q\right)
}(l,a,y)-E_{Q}[h^{\left( Q\right) }(l,a,Y)|l,a]\right\}
\\
& -I(s=2)\left\{ h^{\left( Q\right) }(l,a,y)-E_{Q}[h^{\left( Q\right)
}(L,A,y)|y]\right\}  \notag
\end{align*}%
where $h^{\left( Q\right) }$ solves the integral equation in $h\in
L^{2}\left( Q\right) $%
\begin{align*}
\frac{q(l,a,y)}{p(l,a,y)}\psi _{Q}^{1}(l,a,y)=&
h(l,a,y)-p(S=1|l,a,y)E_{Q}[h(l,a,Y)|l,a] \\
& -p(S=2|l,a,y)E_{Q}[h(L,A,y)|y]
\end{align*}
Furthermore, since $Y$ is binary, $h^{\left( Q\right) }$ admits a closed
form expression given in \cref{prop:example-eifs}.

\textbf{Scenario (iv):} We do not give the derivation of the unique observed data influence function in this scenario because it is isomorphic to scenario (ii) as discussed in \cref{app-subsubsec:generalizing-ex-3}. 

\subsection{Efficiency Gains in \texorpdfstring{\cref{example:transporting}}{Example 3}}

When the alignments in $\mathcal{C}$ impose equality
constraints on the laws $P$ in $\mathcal{P}$, it may be possible to relax
some alignment assumptions while maintaining parameter identification. In
this subsection, we study this concept using \cref{example:transporting} scenarios (ii), (iii) and (iv) (see \cref{app-subsubsec:generalizing-ex-3} for a description of scenario (iv)). In scenario (iii) of that example, we could consider relaxing the alignments to match the
alignments in either scenarios (ii) or (iv), both of which suffice to
identify the ideal data distribution $Q$ and consequently the target
parameter $\psi \left( Q\right)$. In contrast to the alignments in 
scenario (iii), the alignments in scenarios (ii) or (iv) do not place equality constraints on the law $P$ in $\mathcal{P}$. This lack of equality constraints induces a decrease in the efficiency with which $\psi \left( Q\right)$ can be estimated. To demonstrate this phenomenon, we
computed the asymptotic variance of semiparametric efficient estimators of
the average treatment effect $\psi(Q) = \varphi(P)$ in the ideal population $Q$ at a particular law $P$ in model $\mathcal{P}$ that aligns with $Q$ 
under the fused-data framework of all of scenarios (ii), (iii), or (iv) simultaneously. We used a data generating process in which treatment and outcome
are both binary and the covariate $L$ was a vector $\left(
L_{1},L_{2}\right) $ with $L_{1}$ and $L_{2}$ discrete with two
and three levels respectively. \cref{app-subsubsec:dgp} below describes the data-generating
process in detail. \cref{fig:asym-var} below presents the results from
this analysis.

\begin{figure}[ht]
\caption{Asymptotic relative efficiency of efficient estimators of the ATE under the scenarios (ii), (iii), and (iv) of \cref{example:transporting}}
\label{fig:asym-var}\centering
\includegraphics[width=1\textwidth]{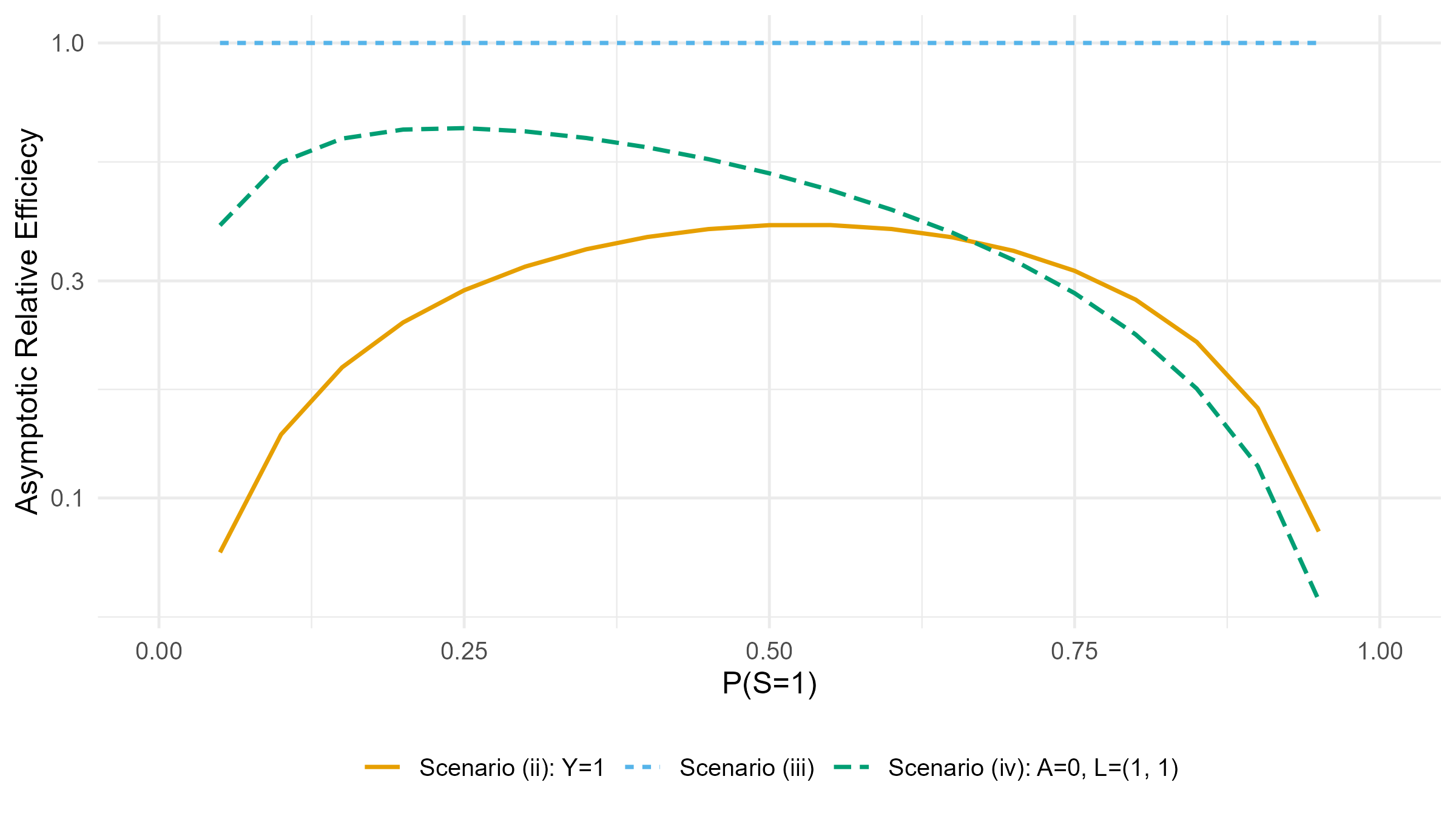}
\end{figure}

The plot in \cref{fig:asym-var} depicts the asymptotic relative efficiencies of
efficient estimators of $\varphi \left( P\right) $ under scenarios (ii) and
(iv) with respect to an efficient estimator of $\varphi\left( P\right) $
under scenario (iii) as a function of $P\left( S=1\right) ,$ the
probability of observing data from the prospective cohort study. The degree
of variance reduction under scenario (iii) illustrates that the alignment
assumptions for this scenario impose strong restrictions on the observed
data model. Recall that the observed data models in scenarios (ii) and
(iv) do not impose equality constraints. As usual, relaxing assumptions
broadens the set of data-generating processes under which efficient
estimators of $\varphi \left( P\right)$ are asymptotically unbiased. 

\subsubsection{Data-Generating process for \texorpdfstring{\cref{fig:asym-var}}{Figure 2}}
\label{app-subsubsec:dgp}
To produce \cref{fig:asym-var}, we chose a observed data law $P$ that belonged to the observed data model $\mc{P}$ for the three fused-data models of \cref{example:transporting} scenarios (ii), (iii), and (iv) simultaneously. We now describe the chosen law $P$ for $(L, A, Y, S)$. $A, Y$ are binary, and $L = (L_1, L_2)$ is a two dimensional random vector where $L_1$ and $L_2$ are discrete taking values in $\{1, 2\}$ and $\{1,2,3\}$ respectively. 

First, we define the ideal data law $Q$ that $P$ aligns with. $L_1$ and $L_2$ are independent under $Q$ and are both discrete uniformly distributed. We then define
\begin{align*}
    \text{logit}[Q(Y=1|a,l_1, l_2)] &\coloneqq 0.5 + 0.5a + 0.25l_1 - 0.25l_2\\
    \text{logit}[Q(A=1|l_1, l_2)] &\coloneqq -0.2 -0.15l_1 + 0.25l_2
\end{align*}
for $a \in \{0, 1\}$, $l_1 \in \{1, 2\}$ and $l_2 \in \{1, 2, 3\}$ which fully specifies the ideal data law $Q$. 

We then define the observed data law $P$. For $a \in \{0, 1\}$, $l_1 \in \{1, 2\}$ and $l_2 \in \{1, 2, 3\}$, 
\begin{align*}
    P(Y=1|A, L_1, L_2, S=1) &\coloneqq  Q(Y|A, L_1, L_2)\\
    P(A=a, L_1=l_1, L_2=l_2|Y, S=2) &\coloneqq  Q(A=a, L_1=l_1, L_2=l_2|Y)\\
    P(A=1|L_1, L_2, S=1) &\coloneqq 0.1 -0.2L_1 + 0.2L_2 \text{ a.e. $Q$}\\
    P(L_1|S=1) &\coloneqq 0.4I(L_1=1) + 0.6I(L_1=2)  \\
    P(L_2|S=1) &\coloneqq 0.3I(L_2=1) + 0.33I(L_2=2) + 0.37I(L_2=3)  \\
    P(Y=1|S=2) &\coloneqq 0.4
\end{align*}
and $L_1\perp L_2|S=1$ under $P$ where all equalities are a.e.-$Q$. We varied $P(S=1)$ over the set $\{0.05, 0.1, \dots, 0.9, 0.95\}$. This fully specifies the observed data law $P$. $P$ satisfies the alignment assumptions of scenarios (ii), (iii), and (iv) simultaneously because 
\begin{align*}
    P(Y=1|A, L_1, L_2, S=1) =  Q(Y=1|A, L_1, L_2) \text{ a.e.-$Q$}\\
    P(A, L_1, L_2|Y, S=2) =  Q(A, L_1, L_2|Y)  \text{ a.e.-$Q$}
\end{align*}
and because $L, A, Y$ were each discrete. Additionally, because the random vector $(O)$ is finite discrete, we can easily compute the asymptotic variances $E_P[\varphi^1_{P, eff}(O)^2]$ of the observed data efficient influence functions under the assumptions of the fused-data frameworks of each of the three scenarios.

\subsection{Proofs for \texorpdfstring{\cref{app:examples-extras}}{Supplement D}}
\label{subsec:example-proofs}
\begin{proof}[Proof of \cref{prop:example-ifs}]
$ $\newline
    \textbf{Proof of part \ref{item:prop:example-ifs1}:}
    We first note that because $Y, V$ are both binary, any $f \in L^2(Q)$ may be written as 
    \begin{align*}
        f(x, v, y) = f_{X}(x) + f_{X, V}(x)v + f_{ X, Y}(x)y + f_{X, Y, V}(x)vy
    \end{align*}
    for some $f_{X}, f_{X, V}, f_{X, Y}, f_{X, Y, V} \in L^2(X;Q)$. Define 
    \begin{align*}
        m_1^{(1)}(x, v) \coloneqq& \psi^1_{Q;X}(x) + \psi^1_{Q;X, V}(x)v + \psi^1_{Q; X, Y}(x)m_Q(x, v) + \psi^1_{Q;X, V, Y}(x)m_Q(x, v)E_Q[V|Y=1, x]\\
        m_2^{(2)}(x, y, v) \coloneqq& \psi^1_{Q;X, Y}(x)\{y - m_Q(x, v)\} + \psi^1_{Q;X, V, Y}(x)\{vy - m_Q(x, v)E_Q[V|Y=1, x]\}
    \end{align*}
    $\varphi$ will be pathwise differentiable at $P$ in model $\mc{P}$ with unique influence function $\varphi^1_P$ if 
    $m_1^{(1)} \in \mc{D}_1^{(1)}(Q)$, $m_2^{(2)} \in \mc{D}_2^{(2)}(Q)$
    and $m_1^{(1)} + m_2^{(2)} = \psi^1_Q$ by \cref{lemma:pathwise-differentiabilty} and part \ref{item:influence functions1} of \cref{Theorem:influence functions}. 
    
    Clearly $m_1^{(1)} + m_2^{(2)} = \psi^1_Q$ by construction. We now show that $m_1^{(1)} \in \mc{D}_1^{(1)}(Q)$. 
    \begin{align*}
        &E_Q[\psi^1_{Q;X}(X) + \psi^1_{Q;X, V}(X)V + \psi^1_{Q; X, Y}(X)m_Q(X, V) + \psi^1_{Q;X, V, Y}(X)m_Q(X, V)E_Q[V|Y=1, X]]\\
        =&E_Q[E_Q[\psi^1_{Q;X}(X) + \psi^1_{Q;X, V}(X)V \\
        &+ \psi^1_{Q; X, Y}(X)m_Q(X, V) + \psi^1_{Q;X, V, Y}(X)m_Q(X, V)E_Q[V|Y=1, X]|X, Y]]\\
        =&E_Q[\psi^1_{Q;X}(X) + \psi^1_{Q;X, V}(X)V +  \psi^1_{Q; X, Y}(X)Y + E_Q[V|X, Y]\psi^1_{Q;X, V, Y}(X)Y]\\
        =&E_Q[\psi^1_{Q;X}(X) + \psi^1_{Q;X, V}(X)V +  \psi^1_{Q; X, Y}(X)Y + \psi^1_{Q;X, V, Y}(X)YV]\\
        =&E_Q[\psi^1_Q(X, Y, V)]\\
        =&0
    \end{align*}
    through repeated use of the tower law. The second equality follows because $E_Q[m_Q(X, V)|x, y] = y$. Additionally, $m^{(1)}_1\in L^2((X, V);Q)$ because $\psi^1_{Q;X}, \psi^1_{Q;X, V},\psi^1_{Q; X, Y}, \psi^1_{Q;X, V, Y} \in L^2(Q)$, $E_Q[V|Y=1, X] \leq 1$ and $m_Q(X, V) \leq \widetilde{\delta}^{-1}$ a.e.-$Q$.
    
    Next we show that $m_2^{(2)} \in \mc{D}_2^{(2)}(Q)$: 
    \begin{align*}
        &E_Q\left[\psi^1_{Q;X, Y}(X)\{Y - m_Q(X, V)\} + \psi^1_{Q;X, V, Y}(X)\{VY- m_Q(X, V)E_Q[V|Y=1, X]\}|x, y\right]\\
        =&\psi^1_{Q;X, Y}(x)\{y - E_Q[m_Q(X, V)|x, y]\} +\psi^1_{Q;X, V, Y}(x)\{yE_Q[V|x, y] - E_Q[m_Q(X, V)|x, y]E_Q[V|Y=1, x]\}\\
        =&0
    \end{align*}
    where in the second equality we used $E_Q[m_Q(X, V)|x, y] = y$. Clearly, $m_2^{(2)} \in L^2(Q)$. Setting $m_1^{(2)} \coloneqq 0$ and applying \cref{lemma:pathwise-differentiabilty} and part \ref{item:influence functions1} of \cref{Theorem:influence functions} gives the desired result. $\varphi^1_P$ is unique by part \ref{item:lemma:example-tangent-spaces1} of \cref{lemma:example-tangent-spaces}. 

    \textbf{Proof of part \ref{item:prop:example-ifs2} ($\Rightarrow$):}
    Suppose that $\psi^1_Q$ corresponds to an influence function $\varphi^1_P$ at $P$ in model $\mc{P}$. Then, by part \ref{item:influence functions1} of \cref{Theorem:influence functions} $\psi^1_Q = m_2^{(1)} + m_2^{(2)}$ where $m_2^{(j)} \in \mc{D}_2^{(j)}(Q)$, $j\in\{1, 2\}$ and we have used that the alignments in $\mc{C}$ imply that $m_1^{(j)} = 0$ for $j\in\{1, 2\}$. Let $\psi^1_{Q; L, Y} \coloneqq m_2^{(1)}\in \mc{D}_2^{(1)}(Q)\subseteq L^2((L, Y);Q)$ and $\psi^1_{Q; L, X} \coloneqq m_2^{(2)}\in \mc{D}_2^{(2)}(Q)\subseteq L^2((L, X);Q)$. Then, $E_Q[\psi^1_{Q; L, Y}(L, Y)|L] = E_Q[\psi^1_{Q; L, X}(L, X)|L] = 0$ and so $E_Q[\psi^1_{Q; L, Y}(L, Y)|L] + E_Q[\psi^1_{Q; L, X}(L, X)|L] =0$. 
    
    Additionally, 
    \begin{align*}
        \varphi^1_P(o) \coloneqq& \frac{I(s=1)}{P(S=1)}\frac{q(l)}{p(l|S=1)}\left\{\psi^1_{Q;L, Y}(l, y) - E_Q[\psi^1_{Q;L, Y}(L, Y)|l]\right\}\\
        &+\frac{I(s=2)}{P(S=2)}\frac{q(l)}{p(l|S=1)}\left\{\psi^1_{Q;L, Y}(l, x) + E_Q[\psi^1_{Q;L, Y}(L, Y)|l]\right\}
    \end{align*}
    is an influence function for $\varphi$ that corresponds to $\psi^1_Q$ by part \ref{item:influence functions1} of \cref{Theorem:influence functions}. By part \ref{item:influence functions2} of that Theorem and part \ref{item:lemma:example-tangent-spaces2} of \cref{lemma:example-tangent-spaces},  $\varphi^1_P$ is the unique observed data influence function that corresponds to $\psi^1_Q$ because $\mc{T}(P, \mc{P}^{ext})$ is nonparametric. 
    
    \textbf{Proof of part \ref{item:prop:example-ifs2} ($\Leftarrow$):}
    Suppose that $\psi^1_Q(l, x, y) = \psi^1_{Q;L, Y}(l, y) + \psi^1_{Q;L, X}(l, x)$ for some $\psi^1_{Q;L, Y} \in L^2(L, Y;Q)$, $\psi^1_{Q;L, X} \in L^2(L, X;Q)$, and 
    \begin{align*}
        E_Q[\psi^1_{Q;L, X}(L, X)|L] + E_Q[\psi^1_{Q;L, Y}(L, Y)|L] = 0.
    \end{align*}
    Let $m_2^{(1)}(l, y) \coloneqq \psi^1_{Q;L, Y}(l, y) - E_Q[\psi^1_{Q;L, Y}(L, Y)|l]$,  $m_2^{(2)}(l, x) \coloneqq \psi^1_{Q;L, X}(l, x) + E_Q[\psi^1_{Q;L, Y}(L, Y)|l]$, and $m_1^{(j)} \coloneqq 0$ for $j\in\{1, 2\}$. Then, $\sum_{j\in[J]}\sum_{k\in[K^{(j)}]} m_k^{(j)} = \psi^1_Q$. Additionally, $E_Q[m_2^{(1)}(L, Y)|l] = E_Q[m_2^{(1)}(L, X)|l] = 0$ and furthermore $m_2^{(1)}(L, Y) \in L^2((L, Y);Q)$ and $m_2^{(2)}(L, X) \in L^2((L, X);Q)$. Hence, $m_k^{(j)} \in \mc{D}_k^{(j)}(Q)$ for $k\in[K^{(j)}]$, $j\in[J]$ and so $\varphi$ is pathwise differentiable at $P$ in model $\mc{P}$.

    We also have that
    \begin{align*}
        \varphi^1_P(o) \coloneqq& \frac{I(s=1)}{P(S=1)}\frac{q(l)}{p(l|S=1)}\left\{\psi^1_{Q;L, Y}(l, y) - E_Q[\psi^1_{Q;L, Y}(L, Y)|l]\right\}\\
        &+\frac{I(s=2)}{P(S=2)}\frac{q(l)}{p(l|S=1)}\left\{\psi^1_{Q;L, Y}(l, x) + E_Q[\psi^1_{Q;L, Y}(L, Y)|l]\right\}
    \end{align*}
    is an influence function for $\varphi$ that corresponds to $\psi^1_Q$ by part \ref{item:influence functions1} of \cref{Theorem:influence functions}. By part \ref{item:influence functions2} of that Theorem and part \ref{item:lemma:example-tangent-spaces2} of \cref{lemma:example-tangent-spaces},  $\varphi^1_P$ is the unique observed data influence function that corresponds to $\psi^1_Q$ because $\mc{T}(P, \mc{P}^{ext})$ is nonparametric. 

    \textbf{Proof of part \ref{item:prop:example-ifs3}:}
    Let $m_1^{(j)} \coloneqq 0$ for $j\in\{1, 2\}$. Let $m_2^{(1)}(u, b) \coloneqq \psi^1_Q(u, b) - \frac{I(u=u_0)}{q(u_0|b)}E_Q\left[\psi^1_Q(U, B)|b\right]$ and $m_2^{(2)}(u, b) \coloneqq \frac{I(u=u_0)}{q(u_0|b)}E_Q\left[\psi^1_Q(U, B)|b\right]$. Clearly $\sum_{j\in[J]}\sum_{k\in[K^{(j)}]}m_k^{(j)} = \psi^1_Q$. By \cref{lemma:pathwise-differentiabilty} and part \ref{item:influence functions1} of \cref{Theorem:influence functions} the proof of this part will be completed if we show that $m_2^{(j)} \in \mc{D}_2^{(j)}(Q)$ for $j\in\{1, 2\}$.

    We first show that $m_2^{(1)} \in \mc{D}_2^{(1)}(Q)$. 
    \begin{align*}
        &E_Q[m_2^{(1)}(U, B)|b]\\
        =&E_Q\left[\left.\psi^1_Q(U, B) - \frac{I(U=u_0)}{q(u_0|B)}E_Q\left[\psi^1_Q(U, B)|B\right]\right|b\right]\\
        =&E_Q\left[\psi^1_Q(U, B)|b\right] - E_Q\left[\psi^1_Q(U, B)|b\right]\frac{E_Q\left[I(U=u_0)|b\right]}{q(u_0|b)}\\
        =&0.
    \end{align*}
    Additionally, $\frac{I(U=u_0)}{q(u_0|B)}E_Q\left[\psi^1_Q(U, B)|b\right] \in L^2(Q)$ because $\frac{1}{Q(u_0|B)} \leq \widetilde{\delta}$ a.e.-$Q$. But then $m_2^{(1)} \in L^2(Q)$ because $E_Q[\psi^1_Q(U, B)|b] \in L^2(Q)$. Hence $m_2^{(1)} \in \mc{D}_2^{(1)}(Q)$. 

    Now we show $m_2^{(2)} \in \mc{D}_2^{(2)}(Q)$. First, $m_2^{(2)}(u, b) = 0$ for $u \not = u_0$. 
    \begin{align*}
        &E_Q[m_2^{(2)}(U, B)|u]\\
        =&E_Q\left[\left.\frac{I(U=u_0)}{q(u_0|B)}E_Q\left[\psi^1_Q(U, B)|B\right]\right|u\right]\\
        =& \int\frac{I(U=u_0)}{q(u|b)}E_Q\left[\psi^1_Q(U, B)|b\right]q(b|u)db\\
        =& \frac{I(U=u_0)}{q(u_0)}\int E_Q\left[\psi^1_Q(U, B)|b\right]q(b)db\\
        =& \frac{I(U=u_0)}{q(u_0)}E_Q[\psi^1_Q(U, B)]\\
        =&0.
    \end{align*}
    Again $\frac{I(U=u_0)}{q(u_0|b)}E_Q\left[\psi^1_Q(U, B)|b\right] \in L^2(Q)$ because $\frac{1}{Q(u_0|B)} \leq \widetilde{\delta}$ a.e.-$Q$. Hence, $m_2^{(2)} \in \mc{D}_2^{(2)}(Q)$, proving the desired result. Part \ref{item:lemma:example-tangent-spaces4} of \cref{lemma:example-tangent-spaces} shows that $\varphi^1_P$ is unique. 

    \textbf{Proof of part \ref{item:prop:example-ifs4}:}
    We first show $\varphi^1_P$ is an influence function of $\varphi$ at $P$ in model $\mc{P}$. Let $m_1^{(j)} \coloneqq 0$ for $j\in\{1, 2\}$. Let $m_2^{(1)} \coloneqq \psi^1_Q(u, b) - \frac{q(u)q(b)}{q(u,b)}E_Q\left[\psi^1_Q(U, B)|b\right]$ and $m_2^{(2)} \coloneqq \frac{q(u)q(b)}{q(u,b)}E_Q\left[\psi^1_Q(U, B)|b\right]$. Clearly $\sum_{j\in[J]}\sum_{k\in[K^{(j)}]}m_k^{(j)} = \psi^1_Q$. Then, by \cref{lemma:pathwise-differentiabilty} and part \ref{item:influence functions1} of \cref{Theorem:influence functions} $\varphi$ is pathwise differentiable with influence function $\varphi^1_P$ if $m_2^{(j)} \in \mc{D}_2^{(j)}(Q)$ for $j\in\{1, 2\}$. 

    To see that $m_2^{(1)} \in \mc{D}_2^{(1)}(Q)$,
    \begin{align*}
        &E_Q[m_2^{(1)}(U, B)|b]\\
        =&E_Q\left[\left.\psi^1_Q(U, B) - \frac{q(U)q(B)}{q(U,B)}E_Q\left[\psi^1_Q(U, B)|B\right]\right|b\right]\\
        =&E_Q\left[\psi^1_Q(U, B)|b\right] - E_Q\left[\psi^1_Q(U, B)|b\right]E_Q\left[\left.\frac{q(U)q(B)}{q(U,B)}\right|b\right]\\
        =& E_Q\left[\psi^1_Q(U, B)|b\right] - E_Q\left[\psi^1_Q(U, B)|b\right]\int \frac{q(u)q(b)}{q(u,b)}q(u|b)db\\
        =& E_Q\left[\psi^1_Q(U, B)|b\right] - E_Q\left[\psi^1_Q(U, B)|b\right]\int q(b)db\\
        =& 0.
    \end{align*}
    We have that $\frac{q(u)q(b)}{q(u,b)}E_Q\left[\psi^1_Q(U, B)|b\right] \in L^2(Q)$ because $\frac{q(U)q(B)}{q(U, B)}\leq \widetilde{\delta}$ a.e.-$Q$. Additionally, $m_2^{(1)} \in L^2(Q)$ because $E_Q[\psi^1_Q|b] \in L^2(Q)$. Hence $m_2^{(1)} \in \mc{D}_2^{(1)}(Q)$. 

    Next
    \begin{align*}
        &E_Q[m_2^{(2)}(U, B)|u]\\
        =&E_Q\left[\left.\frac{q(U)q(B)}{q(U,B)}E_Q\left[\psi^1_Q(U, B)|B\right]\right|u\right]\\
        =&\int \frac{q(u)q(b)}{q(u,b)}E_Q\left[\psi^1_Q(U, B)|b\right]q(b|u)db\\
        =&\int E_Q\left[\psi^1_Q(U, B)|b\right]q(b)db\\
        =& E_Q[\psi^1_Q(U, B)]\\
        =&0.
    \end{align*}
    Again $\frac{q(U)q(B)}{q(U,B)}E_Q\left[\psi^1_Q(U, B)|B\right] \in L^2(Q)$ because $\frac{q(U)q(B)}{q(U, B)}\leq \widetilde{\delta}$ a.e.-$Q$. Hence $m_2^{(2)} \in \mc{D}_2^{(2)}(Q)$. Then $\varphi^1_P$ is an influence function of $\varphi$ at $P$ in model $\mc{P}$. 

    That all influence functions take the form 
    \begin{align*}
        \varphi^1_P =& \frac{I(s=1)}{P(S=1)}\frac{q(b)}{p(b|S=1)}\left\{\psi^1_Q(u, b) - \frac{q(u)q(b)}{q(u, b)}E_Q[\psi^1_Q(U, B)|b]\right\} \\
        &+ \frac{I(s=2)}{P(S=2)}\frac{q(u)}{p(u|S=2)}\left\{\frac{q(u)q(b)}{q(u, b)}E_Q[\psi^1_Q(U, B)|b]\right\}
    \end{align*}
    for 
    \begin{align*}
        \mc{F} = \left\{f \in L^2_0(Q):E_Q[f(U, B)|U] = E_Q[f(U, B)|B] = 0 \text{ a.e. $Q$}\right\}
    \end{align*}
    is a direct corollary of part \ref{item:compute-ifs-algorithm3} of \cref{prop:two-source-solution} with the insight that $\mc{F} = \mc{D}_2^{(1)}(Q) \cap \mc{D}_2^{(2)}(Q)$ and $\mc{D}_1^{(j)}(Q) = \{0\}$ for $j\in\{1, 2\}$. 
\end{proof}

\begin{proof}[Proof of \cref{prop:example-eifs}]
    Recall that by \cref{theorem:eif}, because $\mc{Q}$ is nonparametric, the efficient influence function is 
    \begin{align*}
        \varphi^1_{P, eff}(o) \coloneqq \sum_{j\in[J]}I(s=j)\sum_{k\in[K^{(j)}]}\lim_{n\rightarrow \infty}^{L^2(Q)} I(\overline{z}_{k-1}^{(j)} \in \overline{\mc{Z}}_{k-1}^{(j)})\{E_Q[h_n(W)|\overline{z}_k^{(j)}] - E_Q[h_n(W)|\overline{z}_{k-1}^{(j)}]\}
    \end{align*}
    where the sequence ($h_n$) with $h_n\in \mc{T}(Q, \mc{Q}) = L^2_0(Q)$ solves 
    \begin{align*}
        \psi^1_{Q}(w) = \lim^{L^2(Q)}_{n\rightarrow\infty}\sum_{j\in[J]}P(S=j)\sum_{k\in[K^{(j)}]}\frac{p(\overline{z}_{k-1}^{(j)}|S=j)}{q(\overline{z}_{k-1}^{(j)})}I(\overline{z}_{k-1}^{(j)} \in \overline{\mc{Z}}_{k-1}^{(j)})\{E_Q[h_n(W)|\overline{z}_k^{(j)}] - E_Q[h_n(W)|\overline{z}_{k-1}^{(j)}]\}
    \end{align*} 
    for $\psi^1_Q$ the unique influence function of $\psi$ at $Q$ in model $\mc{Q}$ and $\lim^{L^2(Q)}_{n\rightarrow \infty}$ denotes a limit with respect to the $L^2(Q)$ norm. 
    
    For the fused-data framework of this part, the above equations can be rewritten as 
    \begin{align*}
        \varphi^1_{P, eff}(o) =& I(s=1)\lim^{L^2(Q)}_{n\rightarrow \infty}\{h_n(u, b) - E_Q[h_n(U, B)|B=b]\}\\
        &+I(s=2)\lim^{L^2(Q)}_{n\rightarrow \infty}\{h_n(u, b) - E_Q[h_n(U, B)|U=u]\}
    \end{align*}
    where the sequence $h_n$ with $h_n\in \mc{T}(Q, \mc{Q}) = L^2_0(Q)$ for $n\in\{1,2, \dots\}$ solves 
    \begin{align*}
        \psi^1_{Q}(u, b) = \lim^{L^2(Q)}_{n\rightarrow\infty}&P(S=1)\frac{p(b|S=1)}{q(b)}\{h_n(u, b) - E_Q[h_n(U, B)|B=b]\}\\
        &+P(S=2)\frac{p(u|S=2)}{q(u)}\{h_n(u, b) - E_Q[h_n(U, B)|U=u]\}.
    \end{align*}

    In the proof of \cref{prop:example-ifs} we showed that any element $f\in L^2_0(Q) = \mc{T}(Q, \mc{Q})$ can be written as $f = m_2^{(1)} + m_2^{(2)} + m_1^{(1)} + m_2^{(1)}$ with $m_2^{(1)} \in \mc{D}_2^{(1)}(Q)$ given by 
    \begin{align*}
        m_2^{(1)}(u, b) \coloneqq f(u, b) - \frac{I(u=u_0)}{q(u_0|b)}E_Q[f(U, B)|b],
    \end{align*}
    $m_2^{(2)} \in \mc{D}_2^{(2)}(Q)$ given by
    \begin{align*}
        m_2^{(2)}(u, b) \coloneqq \frac{I(u=u_0)}{q(u_0|b)}E_Q[f(U, B)|b],
    \end{align*}
    and $m_1^{(1)}\coloneqq m_1^{(2)} \coloneqq 0 \in \mc{D}_1^{(1)}(Q) = \mc{D}_1^{(2)}(Q) = \{0\}$. Hence, $L^2_0(Q) = \sum_{j=1}^2\bigoplus_{k=1}^2\mc{D}_k^{(j)}(Q)$. 
    
    Let $U^{(j)} = P(\cdot|S=j)$ for each $j\in[J]$. Then, part \ref{item:lemma:operator-ranges3} of \cref{lemma:operator-ranges} in \cref{sec:appendix-score-operator} tells us that $A^*_{Q}$ has closed range equal to $L^2_0(Q)$. But then, parts the Closed Range Theorem (\cite{rudin_functional_1991} Theorem 4.14)  and \ref{item:lemma:operator-ranges5} of that lemma reveals that we may replace the sequence $h_n \in L^2(Q)$ in the above expressions with a single $h \in L^2(Q)$ for all $n$, leading to 
     \begin{align*}
        \varphi^1_{P, eff}(o) =& I(s=1)\{h(u, b) - E_Q[h(U, B)|b]\}\\
        &+I(s=2)\{h(u, b) - E_Q[h(U, B)|u]\}
    \end{align*}
    where $h\in L^2(Q)$ solves
    \begin{align}
        \label{eq:eif-fredholm-eq-precursor}
        \psi^1_{Q}(u, b) = &P(S=1)\frac{p(b|S=1)}{q(b)}\{h(u, b) - E_Q[h(U, B)|b]\}\\
        &+P(S=2)\frac{p(u|S=2)}{q(u)}\{h(u, b) - E_Q[h(U, B)|u]\}.\nonumber
    \end{align}
    
    Now,
    \begin{align*}
        \frac{q(u,b)}{p(u,b)}P(S=1)\frac{p(b|S=1)}{q(b)} = p(S=1|u,b).
    \end{align*}
    and
    \begin{align*}
        \frac{q(u,b)}{p(u,b)}P(S=2)\frac{p(u|S=2)}{q(u)} = p(S=2|u,b).
    \end{align*}
    because $P\alignswith Q$. Hence, multiplying both sides of \eqref{eq:eif-fredholm-eq-precursor} by $\frac{q(u,b)}{p(u,b)}$ we arrive at 
    \begin{align*}
        \frac{q(u, b)}{p(u, b)}\psi^1_{Q}(u, b) = &p(S=1|u, b)\{h(u, b) - E_Q[h(U, B)|B=b]\}\\
        &+p(S=2|u, b)\{h(u, b) - E_Q[h(U, B)|U=u]\}\\
        =&h(u, b) - p(S=1|u, b)E_Q[h(U, B)|B=b] - p(S=2|u, b)E_Q[h(U, B)|U=u].
    \end{align*}
    Define 
    \begin{align*}
        r_{Q, P}(u,b) \coloneqq \frac{q(u,b)}{p(u,b)}\psi^1_Q(u,b).
    \end{align*}
    We aim to find a solution to 
    \begin{align*}
        r_{Q,P}(u, b) = h(u, b) - p(S=1|u, b)E_Q[h(U, B)|B=b] - p(S=2|u, b)E_Q[h(U, B)|U=u]
    \end{align*}
    for $h\in L^2_0(Q)$. In particular, because $U$ is finite discrete with $T$ levels, the above holds for each $u \in \{u_1, \dots, u_T\}$. Note that $E_Q[h(U, B)|u] = E_P\left[\frac{I(U=u)}{P(U=u|S=2)}h(U, B)\right]$ for each $u\in \{u_1, \dots, u_T\}$. Hence we may rewrite the above display as the matrix equation 
    \begin{align}
        \label{eq:fredholm-intermediate1}
        \underline{r}_{Q,P}(b) =& \underline{h}(b) - \underline{\pi}(b)\underline{\beta}(b)'\underline{h}(b) - \text{diag}(1-\underline{\pi}(b))E_P[R(U)\underline{h}(B)|S=2]\\
        =&(Id - \underline{\pi}(b)\underline{\beta}(b)')\underline{h}(b) - \text{diag}(1-\underline{\pi}(b))E_P[R(U)\underline{h}(B)|S=2]\nonumber
    \end{align}
     $(Id - \underline{\pi}(b)\underline{\beta}(b)')$ is invertible with inverse 
    \begin{align*}
        (Id - \underline{\pi}(b)\underline{\beta}(b)')^{-1} = Id + \frac{1}{1 - \underline{\pi}(b)'\underline{\beta}(b)}\underline{\pi}(b)\underline{\beta}(b)'
    \end{align*}
    whenever $\underline{\pi}(b)'\underline{\beta}(b) \not 
= 1$. But 
    \begin{align*}
        \underline{\pi}(b)'\underline{\beta}(b) =& \sum_{t\in[T]}p(S=1|u_t, b)p(U=u_t|B=b, S=1)\\
        <&\sum_{t\in[T]}p(U=u_t|B=b, S=1)\\
        =&1
    \end{align*}
    because $p(S=1|U, B) < 1$ a.e.-$Q$. Hence \eqref{eq:fredholm-intermediate1} is equivalent to
    \begin{align}
        \underline{h}(b) = (Id - \underline{\pi}(b)\underline{\beta}(b)')^{-1}\left\{\underline{r}_{Q, P}(b) + \text{diag}(1-\underline{\pi}(b))E_P[R(U)\underline{h}(B)|S=2]\right\}.\label{eq:h-precursor}
    \end{align}
   Multiplying by $R(u)$ and taking expectations under $P(\cdot|S=2)$ on both sides of the above display we arrive at
    \begin{align*}
        E_P[R(U)\underline{h}(B)|S=2] =&E_P\left[R(U)(Id - \underline{\pi}(b)\underline{\beta}(b)')^{-1}\underline{r}_{Q, P}(B)|S=2\right] \\
        &+ E_P\left[R(U)(Id - \underline{\pi}(b)\underline{\beta}(b)')^{-1}\text{diag}(1-\underline{\pi}(B))|S=2\right]E_P[R(U)\underline{h}(B)|S=2].
    \end{align*}
    We rewrite the above expression as 
    \begin{align*}
        &\{Id - E_P\left[R(U)(Id - \underline{\pi}(b)\underline{\beta}(b)')^{-1}\text{diag}(1-\underline{\pi}(B))|S=2\right]\}E_P[R(U)\underline{h}(B)|S=2]\\
        =&E_P\left[R(U)(Id - \underline{\pi}(b)\underline{\beta}(b)')^{-1}\underline{r}_{Q, P}(B)|S=2\right].
    \end{align*}
    Now $\varphi$ has an efficient influence function because $\varphi$ is pathwise differentiable by part \ref{item:prop:example-ifs3} of \cref{prop:example-ifs}. Hence, the above display must have solution and as such \\$E_P\left[R(U)(Id - \underline{\pi}(b)\underline{\beta}(b)')^{-1}\underline{r}_{Q, P}(B)|S=2\right]$ is necessarily in the range of 
    \begin{align*}
        \{Id - E_P\left[R(U)(Id - \underline{\pi}(b)\underline{\beta}(b)')^{-1}\text{diag}(1-\underline{\pi}(B))|S=2\right]\}.
    \end{align*}
    Thus we may write 
    \begin{align*}
        &E_P[R(U)\underline{h}(B)|S=2] \\
        =& \{Id - E_P\left[R(U)(Id - \underline{\pi}(b)\underline{\beta}(b)')^{-1}\text{diag}(1-\underline{\pi}(B))|S=2\right]\}^-\\
        &\times E_P\left[R(U)(Id - \underline{\pi}(b)\underline{\beta}(b)')^{-1}\underline{r}_{Q, P}(B)|S=2\right].
    \end{align*}
    Combining the above display with \eqref{eq:h-precursor} gives that 
    \begin{align*}
        \underline{h}(b) = (Id - \underline{\pi}(b)\underline{\beta}(b)')^{-1}\big\{&\underline{r}_{Q, P}(b) \\
        +&\text{diag}(1 - \underline{\pi}(b))\{Id - E_P\left[R(U)(Id - \underline{\pi}(b)\underline{\beta}(b)')^{-1}\text{diag}(1-\underline{\pi}(B))|S=2\right]\}^-\nonumber\\
        \times& E_P\left[R(U)(Id - \underline{\pi}(b)\underline{\beta}(b)')^{-1}\underline{r}_{Q, P}(B)|S=2\right]\big\}\nonumber.
    \end{align*}
    Then, $\left(h(u_1, b), \dots, h(u_T, b)\right)^{T} = \underline{h}(b)$ is the solution to \eqref{eq:eif-fredholm-eq-precursor}. Hence, 
    \begin{align*}
        \varphi^1_{P, eff}(o) \coloneqq& I(s=1) \{h^{(Q)}_{P, eff}(u, b) - E_P\left[h^{(Q)}_{P, eff}(U, B)|B=b, S=1\right]\}\\
        &+ I(s=2) \{h^{(Q)}_{P, eff}(u, b) - E_P\left[h^{(Q)}_{P, eff}(U, B)|U=u, S=2\right]\}
    \end{align*}
    where 
    \begin{align*}
        \underline{h}_{P, eff}^{(Q)}(b) \coloneqq &(Id - \underline{\pi}(b)\underline{\beta}(b)')^{-1}\big\{\underline{r}_{Q, P}(b) \\
        &+\text{diag}(1 - \underline{\pi}(b))\{Id - E_P\left[R(U)(Id - \underline{\pi}(b)\underline{\beta}(b)')^{-1}\text{diag}(1-\underline{\pi}(B))|S=2\right]\}^-\nonumber\\
        &\times E_P\left[R(U)(Id - \underline{\pi}(b)\underline{\beta}(b)')^{-1}\underline{r}_{Q, P}(B)|S=2\right]\big\}\nonumber.
    \end{align*}
    is the efficient influence function of $\varphi$ at $P$ in model $\mc{P}$. 
\end{proof}

\begin{proof}[Proof of \cref{lemma:example-tangent-spaces}]
    Let $U^{(j)} \coloneqq P(\cdot|S=j)$ for each $j\in[J]$. Then $(Q, U, P)$ is strongly aligned. Throughout this proof, we will use that $\mc{T}(P, \mc{P}) = \overline{A_{Q}\mc{T}(Q, \mc{Q})}\oplus \bigoplus_{j\in[J]} A_{U^{(j)}}L^2_0(U^{(j)}) \oplus L^2_0(\lambda) = \text{Null}(A^*_{Q})^\perp \oplus A_{U^{(j)}}L^2_0(U^{(j)}) \oplus L^2_0(\lambda)$ and $\mc{T}(P, \mc{P}^{ext}) = \overline{A_{Q}L^2_0(Q)}\oplus \bigoplus_{j\in[J]} A_{U^{(j)}}L^2_0(U^{(j)}) \oplus L^2_0(\lambda) = \text{Null}(A^{ext,*}_{Q})^\perp \oplus A_{U^{(j)}}L^2_0(U^{(j)}) \oplus L^2_0(\lambda)$ where $A^{ext}_{Q}$ and $A^{ext,*}_{Q}$ are defined as in the proof of part \ref{item:influence functions2} of \cref{Theorem:influence functions}. We will also use that $g \in \text{Null}(A^*_{Q})$ if and only if there exists $m_k^{(j)} \in \mc{D}_k^{(j)}(Q)$, $k\in[K^{(j)}], j\in[J]$ such that
    \begin{align*}
        g = \sum_{j\in[J]}I(S=j)\sum_{k\in[K^{(j)}]}m_k^{(j)}
    \end{align*} 
    and 
    \begin{align*}
        0=\sum_{j\in[J]}\sum_{k\in[K^{(j)}]}\frac{dP(\cdot|S=j)}{dQ}(\overline{z}_{k-1}^{(j)})P(S=j)m_k^{(j)}(\overline{z}_k^{(j)})
    \end{align*}
    which follows from the fact that $A^*_{Q, U, \lambda}g = (A^*_Qg, A^*_{U^{(1)}}g,\dots, A^*_{U^{(J)}}g, A^*_\lambda g)$ and the expressions for $A^*_Q$, $A^*_{U^{(j)}}, j\in [J]$, and $A^*_\lambda$.
    
    \textbf{Proof of part \ref{item:lemma:example-tangent-spaces1}:} Suppose $g \in \text{Null}(A^*_{Q})$. Write $g = \sum_{j=1}^2I(S=j)\sum_{k\in[K^{(j)}]}m_k^{(j)}$ for some $m_k^{(j)} \in \mc{D}_k^{(j)}(Q)$, $k\in[K^{(j)}]$, $j\in[2]$. Note that $\mc{D}_1^{(2)}(Q) = \{0\}$ and so $m_1^{(2)}=0$ by the alignments in $\mc{C}$. Let $f_k^{(j)} \coloneqq \frac{dP(\cdot|S=j)}{dQ}(\overline{z}_{k-1}^{(j)})P(S=j)m_k^{(j)}(\overline{z}_k^{(j)})$ for $k\in[K^{(j)}]$, $j\in\{1, 2\}$. Then 
    \begin{align*}
        f_1^{(1)}  = -f_2^{(2)}
    \end{align*}
    by the expression of $A^*_Q$ in \cref{lemma:score-operator} because $g \in \text{Null}(A^*_{Q})$. Let $f \coloneqq f_1^{(1)} = -f_2^{(2)}$. 
    We have that $f\in \mc{D}_1^{(1)}(Q) \cap \mc{D}_2^{(2)}(Q)$ because $f_1^{(1)} \in \mc{D}_1^{(1)}(Q)$ and $-f_2^{(2)} \in \mc{D}_2^{(2)}(Q)$ by \cref{lemma:f-in-Q}. But then $f$ is a function of $X$ and $V$ alone because $f\in \mc{D}_1^{(1)}(Q)$. On the other hand 
    \begin{align*}
        E_Q[f(X, V)|X, Y] = 0
    \end{align*}
    because $f \in \mc{D}_2^{(2)}(Q)$. Hence 
    \begin{align*}
        f(X, 0)q(0|X, Y) = -f(X, 1)q(1|X, Y)
    \end{align*}
    because $V$ is binary which means that 
    \begin{align*}
        \frac{f(X, 0)}{f(X, 1)} = -\frac{q(1|X, Y)}{q(0|X, Y)}.
    \end{align*}
    if $f(X, 1) \not =0$ where all statements are a.e.-$Q$. The right-hand side of the above display must be a function of $X$ alone because the left-hand side is. But this could only occur if $V\perp Y|X$. However, by assumption $V\not\perp Y|X$. Hence, $f(X, 1) = 0$ which in turn implies that $f(X, 0) = 0$ and so $f(X, V) = 0$ a.e.-$Q$. We then have that $g=0$ and so $\text{Null}(A^*_{Q, U, \lambda}) = \{0\}$. But this means that $\mc{T}(P, \mc{P}) = L^2_0(P)$, concluding the proof of this part.

    \textbf{Proof of part \ref{item:lemma:example-tangent-spaces2}:} Suppose that $g \in \text{Null}(A^{ext,*}_{Q})$. Write $g = \sum_{j=1}^2I(S=j)\sum_{k=1}^2m_k^{(j)}$ for some $m_k^{(j)} \in \mc{D}_k^{(j)}(Q)$, $k\in\{1, 2\}$, $j\in\{1, 2\}$. Note that $\mc{D}_1^{(j)}(Q) = \{0\}$ and so $m_1^{(j)}=0$ for $j\in\{1, 2\}$ by the alignments in $\mc{C}$. Let $f_k^{(j)} \coloneqq \frac{dP(\cdot|S=j)}{dQ}(\overline{z}_{k-1}^{(j)})P(S=j)m_k^{(j)}(\overline{z}_k^{(j)})$ for $k\in\{1, 2\}$, $j\in\{1, 2\}$. Then 
    \begin{align*}
        f_2^{(1)}  = -f_2^{(2)}
    \end{align*}
    by the expression of $A^{ext,*}_{Q}$ in \cref{lemma:score-operator} because $g \in \text{Null}(A^{ext,*}_{Q})$. Let $f \coloneqq f_2^{(1)}  = -f_2^{(2)}$. Then $f \in \mc{D}_2^{(1)}(Q) \cap \mc{D}_2^{(2)}(Q)$ because $f_2^{(j)} \in \mc{D}_2^{(j)}(Q)$ by \cref{lemma:f-in-Q}. As $f \in \mc{D}_2^{(1)}(Q)$, it is a function of $Y$ and $L$ alone. On the other hand, because $f \in \mc{D}_2^{(2)}(Q)$, it is a function of $X$ and $L$ alone. Hence $f$ is a function of $L$ alone because we have assumed throughout that neither $Y$ nor $X$ is a deterministic function of the other. But $E_Q[f(L)|L] = f(L) = 0$ a.e.-$Q$ because $f \in \mc{D}_2^{(1)}(Q)$. This implies $g=0$ proving that $\text{Null}(A^{ext,*}_{Q, U, \lambda}) = \{0\}$. Hence, $\mc{T}(P, \mc{P}^{ext}) = L^2_0(P)$. 

    We now show that $\mc{T}(P, \mc{P})\subsetneq L^2_0(P)$. By the assumptions of the lemma, there exists $f_2^{(j)} \in \mc{D}_2^{(j)}(Q)\cap \mc{T}(Q, \mc{Q})^\perp$ such that at least one of $f_2^{(1)}$, $f_2^{(2)}$ is non-zero. Let 
    \begin{align*}
        g(o) \coloneqq \sum_{j=1}^2I(S=j)\sum_{k=1}^2\frac{dQ}{dP(\cdot|S=j)}(\overline{z}_{k-1}^{(j)})P(S=j)^{-1}f_k^{(j)}(\overline{z}_{k-1}^{(j)}).
    \end{align*}
    Then, $g\not=0$ and $g \in \text{Null}(A^{*}_Q)$ because $A^*_Qg = \Pi[A^{ext,*}_Qg|\mc{T}(Q, \mc{Q})] = \Pi[f_1^{(1)} + f_2^{(2)}|\mc{T}(Q, \mc{Q})]$. This shows $\{0\} \subsetneq \text{Null}(A^{*}_Q)$. But then $\mc{T}(P, \mc{P})\subsetneq L^2_0(P)$ which proves this part of the lemma.  

    \textbf{Proof of part \ref{item:lemma:example-tangent-spaces3}:} Suppose that $g \in \text{Null}(A^{*}_{Q, U, \lambda})$. Write $g = \sum_{j=1}^2I(S=j)\sum_{k\in[K^{(j)}]}m_k^{(j)}$ for some $m_k^{(j)} \in \mc{D}_k^{(j)}(Q)$, $k\in[K^{(j)}]$, $j\in\{1,2\}$. Note that $\mc{D}_1^{(1)}(Q) = \{0\}$ and so $m_1^{(1)}=0$ by the alignments in $\mc{C}$. Let $f_k^{(j)} \coloneqq \frac{dP(\cdot|S=j)}{dQ}(\overline{z}_{k-1}^{(j)})P(S=j)m_k^{(j)}(\overline{z}_k^{(j)})$ for $k\in[K^{(j)}]$, $j\in\{1,2\}$. Then 
    \begin{align*}
        f_2^{(1)}  = -f_1^{(2)}.
    \end{align*}
    by the expression of $A^*_Q$ in \cref{lemma:score-operator} because $g \in \text{Null}(A^*_{Q})$. Let $f \coloneqq f_2^{(1)} = -f_1^{(2)}$. Then $f \in \mc{D}_2^{(1)}(Q) \cap \mc{D}_1^{(2)}(Q)$ because $f_2^{(1)} \in \mc{D}_2^{(1)}(Q)$ and $f_1^{(2)} \in \mc{D}_1^{(2)}(Q)$.
    But $\mc{D}_2^{(1)}(Q) \perp \mc{D}_1^{(2)}(Q)$, so $f = 0$. This means $g=0$ and so $\text{Null}(A^*_{Q, U, \lambda}) = \{0\}$. Hence, $\mc{T}(P, \mc{P}) = L^2_0(P)$. 

    \textbf{Proof of part \ref{item:lemma:example-tangent-spaces4}:} Suppose that $g \in \text{Null}(A^{*}_{Q})$. Write $g = \sum_{j=1}^2I(S=j)\sum_{k=1}^2m_k^{(j)}$ for some $m_k^{(j)} \in \mc{D}_k^{(j)}(Q)$, $k\in\{1, 2\}$, $j\in\{1, 2\}$. Note that $\mc{D}_1^{(j)}(Q) = \{0\}$ and so $m_1^{(j)}=0$ for $j\in\{1, 2\}$ by the alignments in $\mc{C}$. Let $f_k^{(j)} \coloneqq \frac{dP(\cdot|S=j)}{dQ}(\overline{z}_{k-1}^{(j)})P(S=j)m_k^{(j)}(\overline{z}_k^{(j)})$ for $k\in\{1, 2\}$, $j\in\{1, 2\}$. Then 
    \begin{align*}
        f_2^{(1)}  = -f_2^{(2)}.
    \end{align*}
    by the expression of $A^*_Q$ in \cref{lemma:score-operator} because $g \in \text{Null}(A^*_{Q})$. Let $f \coloneqq f_2^{(1)} = -f_2^{(2)}$. Then $f \in \mc{D}_2^{(1)}(Q) \cap \mc{D}_2^{(2)}(Q)$ because $f_2^{(j)} \in \mc{D}_2^{(j)}(Q)$ for $j\in\{1, 2\}$ by \cref{lemma:f-in-Q}. We have that $f(u, b) = 0$ for $u \not = u_0$ because $f \in \mc{D}_2^{(2)}(Q)$. On the other hand $E_Q[f(U, B)|B] = 0$ a.e. $Q$ because $f \in \mc{D}_2^{(1)}(Q)$. But these two facts combined give that 
    \begin{align*}
        E_Q[f(U, B)|B] = f(u_0, B) = 0 
    \end{align*}
    which means that $f(u_0, B) = 0$ a.e.-$Q$ and so $f = 0$. This implies $g=0$ and so $\text{Null}(A^*_{Q, U, \lambda}) = \{0\}$. Hence, $\mc{T}(P, \mc{P}) = L^2_0(P)$. 

    \textbf{Proof of part \ref{item:lemma:example-tangent-spaces5}:} Suppose that $g \in \text{Null}(A^{*}_{Q})$. Write $g = \sum_{j=1}^2I(S=j)\sum_{k=1}^2m_k^{(j)}$ for some $m_k^{(j)} \in \mc{D}_k^{(j)}(Q)$, $k\in\{1, 2\}$, $j\in\{1, 2\}$. Note that $\mc{D}_1^{(j)}(Q) = \{0\}$ and so $m_1^{(j)}=0$ for $j\in\{1, 2\}$ by the alignments in $\mc{C}$. Let $f_k^{(j)} \coloneqq \frac{dP(\cdot|S=j)}{dQ}(\overline{z}_{k-1}^{(j)})P(S=j)m_k^{(j)}(\overline{z}_k^{(j)})$ for $k\in\{1, 2\}$, $j\in\{1, 2\}$. Then 
    \begin{align*}
        f_2^{(1)}  = -f_2^{(2)}.
    \end{align*}
    by the expression of $A^*_Q$ in \cref{lemma:score-operator} because $g \in \text{Null}(A^*_{Q})$. Let $f \coloneqq f_2^{(1)} = -f_2^{(2)}$. Then $f \in \mc{D}_2^{(1)}(Q) \cap \mc{D}_2^{(2)}(Q)$ because $f_2^{(j)} \in \mc{D}_2^{(j)}(Q)$ for $j\in\{1, 2\}$ by \cref{lemma:f-in-Q}. It follows from  \cite{tchetgen_tchetgen_doubly_2010} that if neither $U$ nor $B$ is a measurable map of the other, then $\{0\}\subsetneq \mc{D}_2^{(1)}(Q) \cap \mc{D}_2^{(2)}(Q)$. This follows from those authors characterization of the intersection as  
    \begin{align*}
        &\mc{D}_2^{(1)}(Q) \cap \mc{D}_2^{(2)}(Q)  \\
        &=\left\{\frac{q(u)q(b)}{q(u, b)}\left\{t(u, b) - E_{Q^*}[t(U, B)|b] - E_{Q^*}[t(U, B)|u] - E_{Q^*}[t(U, B)]\right\}:t \in L_0^2(Q)\right\}.
    \end{align*}
    Then there exists $g \not=0$ such that $g \in \text{Null}(A^*_{Q})$. Hence, $\mc{T}(P, \mc{P}) \subsetneq L^2_0(P)$.
\end{proof}

\begin{proof}[Proof of \cref{lemma:case-control-ji}]
$ $\newline
    \textbf{Proof of part \ref{item:lemma:case-control-ji1}:} 
    
    Let $b\in \mathcal{B}^{\ast }.$ Then from 
    \[
    q\left( u_{0}|b\right) q\left( b\right) =q\left( b|u_{0}\right) q\left(
    u_{0}\right) 
    \]%
    and the assumption $q\left( u_{0}|b\right) >0$ for any $b\in \mathcal{B}%
    ^{\ast }$, we have 
    \begin{equation}
    q\left( b\right) =\frac{q\left( b|u_{0}\right) }{q\left( u_{0}|b\right) }%
    q\left( u_{0}\right)   \label{ii1}
    \end{equation}%
    Consequently, 
    \begin{align*}
    1=&\int_{\mathcal{B}}q\left( b\right) d\mu \left( b\right) \\
    =&\int_{\mathcal{B}%
    ^{\ast }}q\left( b\right) d\mu \left( b\right)  \\
    =&\int_{\mathcal{B}^{\ast }}\frac{q\left( b|u_{0}\right) }{q\left(
    u_{0}|b\right) }q\left( u_{0}\right) d\mu \left( b\right)  \\
    =&q\left( u_{0}\right) \int_{\mathcal{B}^{\ast }}\frac{q\left(
    b|u_{0}\right) }{q\left( u_{0}|b\right) }d\mu \left( b\right) 
    \end{align*}%
    Therefore,%
    \[
    q\left( u_{0}\right) =\left[ \int_{\mathcal{B}^{\ast }}\frac{q\left(
    b|u_{0}\right) }{q\left( u_{0}|b\right) }d\mu \left( b\right) \right] ^{-1}
    \]%
    Replacing in $\left( \ref{ii1}\right) $ we arrive at 
    \[
    q\left( b\right) =\frac{q\left( b|u_{0}\right) }{q\left( u_{0}|b\right) }%
    \left[ \int_{\mathcal{B}^{\ast }}\frac{q\left( b|u_{0}\right) }{q\left(
    u_{0}|b\right) }d\mu \left( b\right) \right] ^{-1}
    \]%
    and therefore 
    \begin{eqnarray*}
    q\left( u,b\right)  &=&q\left( u|b\right) q\left( b\right)  \\
    &=&q\left( u|b\right) \frac{q\left( b|u_{0}\right) }{q\left( u_{0}|b\right) }%
    \left[ \int_{\mathcal{B}^{\ast }}\frac{q\left( b|u_{0}\right) }{q\left(
    u_{0}|b\right) }d\mu \left( b\right) \right] ^{-1}
    \end{eqnarray*}
    
    \textbf{Proof of part \ref{item:lemma:case-control-ji2}:}
    
    By $\left\{ b\in \mathcal{B}:q\left( b|u_{0}\right) >0\right\} \subseteq 
    \mathcal{B}^{\ast }$ we have that 
    \[
    1\geq \int_{\mathcal{B}^{\ast }}q\left( b|u_{0}\right) d\mu \left( b\right)
    =\int_{\left\{ b\in \mathcal{B}:q\left( b|u_{0}\right) >0\right\} }q\left(
    b|u_{0}\right) d\mu \left( b\right) =1
    \]%
    Therefore, $\int_{\mathcal{B}^{\ast }}q\left( b|u_{0}\right) d\mu \left(
    b\right) =1.$ Then,  
    \[
    \int_{\mathcal{B}^{\ast }}p_{2}\left( b|u_{0}\right) d\mu \left( b\right) =1
    \]%
    because $p_2(b|u_0) = q(b|u_0)$ for $b\in \mc{B}^*$. We then conclude that for $\mu -$almost all $b$ in $\mathcal{B}%
    \backslash \mathcal{B}^{\ast },$ it holds that $p_{2}\left( b|u_{0}\right) =0
    $. Then,   
    \[
    \int_{\mathcal{B}^{\ast }}\frac{p_{2}\left( b|u_{0}\right) }{p_{1}\left(
    u_{0}|b\right) }d\mu \left( b\right) =\int_{\mathcal{B}}\frac{p_{2}\left(
    b|u_{0}\right) }{p_{1}\left( u_{0}|b\right) }d\mu \left( b\right) 
    \]%
    Furthermore, by assumption $p_{1}\left( u_{0}|b\right) $ is defined for any $b\in 
    \mathcal{B}^{\ast }$ and $p_{1}\left( u_{0}|b\right) =q\left(
    u_{0}|b\right) $ for any $b\in \mathcal{B}^{\ast }.$\ Then, replacing for
    any $b\in \mathcal{B}^{\ast },q\left( b|u_{0}\right) $ with $p_{2}\left(
    b|u_{0}\right) $ and $q\left( u_{0}|b\right) $ with $p_{1}\left(
    u_{0}|b\right) $ we arrive at     
    \begin{eqnarray*}
    q\left( u,b\right)  &=&q\left( u|b\right) \frac{q\left( b|u_{0}\right) }{%
    q\left( u_{0}|b\right) }\left[ \int_{\mathcal{B}^{\ast }}\frac{q\left(
    b|u_{0}\right) }{q\left( u_{0}|b\right) }d\mu \left( b\right) \right] ^{-1}
    \\
    &=&p_{1}\left( u|b\right) \frac{p_{2}\left( b|u_{0}\right) }{p_{1}\left(
    u_{0}|b\right) }\left[ \int_{\mathcal{B}^{\ast }}\frac{p_{2}\left(
    b|u_{0}\right) }{p_{1}\left( u_{0}|b\right) }d\mu \left( b\right) \right]
    ^{-1} \\
    &=&p_{1}\left( u|b\right) \frac{p_{2}\left( b|u_{0}\right) }{p_{1}\left(
    u_{0}|b\right) }\left[ \int_{\mathcal{B}}\frac{p_{2}\left( b|u_{0}\right) }{%
    p_{1}\left( u_{0}|b\right) }d\mu \left( b\right) \right] ^{-1}
    \end{eqnarray*}
    
    \textbf{Proof of part \ref{item:lemma:case-control-ji3}:} 
    
    Define $\mathcal{B}^{\ast }=\left\{ b\in 
    \mathcal{B}:p_{1}\left( b\right) >0\right\} $ and let $u_{0}$ be such that $%
    p_{1}\left( u_{0}|b\right) >0$ for all $b\in \mathcal{B}^{\ast }.$ For any $%
    b\in \mathcal{B}^{\ast }$ let  
    \[
    q\left( u,b\right) \coloneqq p_{1}\left( u|b\right) \frac{p_{2}\left(
    b|u_{0}\right) }{p_{1}\left( u_{0}|b\right) }\left[ \int_{\mathcal{B}^{\ast
    }}\frac{p_{2}\left( b|u_{0}\right) }{p_{1}\left( u_{0}|b\right) }d\mu \left(
    b\right) \right] ^{-1}
    \]%
    and let $q\left( u,b\right) =0$ otherwise, 
    
    We will first show that $q\left( u,b\right) $ is a density with respect to
    some prob. measure $Q$ on $\mathcal{U}\times \mathcal{B}$. To do so, it
    suffices to show that 
    \[
    \int_{\mathcal{B}}\int_{\mathcal{U}}q\left( u,b\right) d\lambda \left(
    u\right) d\mu \left( b\right) =1
    \]%
    Now,   
    \[
    \int_{\mathcal{B}}\int_{\mathcal{U}}q\left( u,b\right) d\lambda \left(
    u\right) d\mu \left( b\right) =\int_{\mathcal{B}^{\ast }}\int_{\mathcal{U}%
    }q\left( u,b\right) d\lambda \left( u\right) d\mu \left( b\right) 
    \]%
    because by definition $q\left( u,b\right) =0$ when $b\notin \mathcal{B}%
    ^{\ast }$. Furthermore, for any $b$ in $\mathcal{B}^{\ast }$,%
    \begin{eqnarray*}
    \int_{\mathcal{U}}q\left( u,b\right) d\lambda \left( u\right)  &=&\int_{%
    \mathcal{U}}p_{1}\left( u|b\right) \frac{p_{2}\left( b|u_{0}\right) }{%
    p_{1}\left( u_{0}|b\right) }\left[ \int_{\mathcal{B}^{\ast }}\frac{%
    p_{2}\left( b|u_{0}\right) }{p_{1}\left( u_{0}|b\right) }d\mu \left(
    b\right) \right] ^{-1}d\lambda \left( u\right)  \\
    &=&\left\{ \int_{\mathcal{U}}p_{1}\left( u|b\right) d\lambda \left( u\right)
    \right\} \frac{p_{2}\left( b|u_{0}\right) }{p_{1}\left( u_{0}|b\right) }%
    \left[ \int_{\mathcal{B}^{\ast }}\frac{p_{2}\left( b|u_{0}\right) }{%
    p_{1}\left( u_{0}|b\right) }d\mu \left( b\right) \right] ^{-1} \\
    &=&\frac{p_{2}\left( b|u_{0}\right) }{p_{1}\left( u_{0}|b\right) }\left[
    \int_{\mathcal{B}^{\ast }}\frac{p_{2}\left( b|u_{0}\right) }{p_{1}\left(
    u_{0}|b\right) }d\mu \left( b\right) \right] ^{-1}
    \end{eqnarray*}%
    Consequently 
    \[
    \int_{\mathcal{B}^{\ast }}\int_{\mathcal{U}}q\left( u,b\right) d\lambda
    \left( u\right) d\mu \left( b\right) =\int_{\mathcal{B}^{\ast }}\frac{%
    p_{2}\left( b|u_{0}\right) }{p_{1}\left( u_{0}|b\right) }\left[ \int_{%
    \mathcal{B}^{\ast }}\frac{p_{2}\left( b|u_{0}\right) }{p_{1}\left(
    u_{0}|b\right) }d\mu \left( b\right) \right] ^{-1}d\mu \left( b\right) =1
    \]%
    Next, we show that $p_{2}\left( b|u_{0}\right) =q\left( b|u_{0}\right) $ and 
    $p_{1}\left( u_{0}|b\right) =q\left( u_{0}|b\right) $ for all $b\in \mathcal{%
    B}^{\ast }$. Specifically, for any $b\in \mathcal{B}^{\ast }$  
    \begin{eqnarray*}
    &&q\left( u_{0}|b\right) \coloneqq\frac{q\left( u,b\right) }{\int q\left(
    u,b\right) d\lambda \left( u\right) } \\
    &&\coloneqq\frac{p_{1}\left( u|b\right) \frac{p_{2}\left( b|u_{0}\right) }{%
    p_{1}\left( u_{0}|b\right) }\left[ \int_{\mathcal{B}^{\ast }}\frac{%
    p_{2}\left( b|u_{0}\right) }{p_{1}\left( u_{0}|b\right) }d\mu \left(
    b\right) \right] ^{-1}}{\int_{\mathcal{U}}p_{1}\left( u|b\right) \frac{%
    p_{2}\left( b|u_{0}\right) }{p_{1}\left( u_{0}|b\right) }\left[ \int_{%
    \mathcal{B}^{\ast }}\frac{p_{2}\left( b|u_{0}\right) }{p_{1}\left(
    u_{0}|b\right) }d\mu \left( b\right) \right] ^{-1}d\lambda \left( u\right) }
    \\
    &=&\frac{p_{1}\left( u|b\right) \frac{p_{2}\left( b|u_{0}\right) }{%
    p_{1}\left( u_{0}|b\right) }\left[ \int_{\mathcal{B}^{\ast }}\frac{%
    p_{2}\left( b|u_{0}\right) }{p_{1}\left( u_{0}|b\right) }d\mu \left(
    b\right) \right] ^{-1}}{\frac{p_{2}\left( b|u_{0}\right) }{p_{1}\left(
    u_{0}|b\right) }\left[ \int_{\mathcal{B}^{\ast }}\frac{p_{2}\left(
    b|u_{0}\right) }{p_{1}\left( u_{0}|b\right) }d\mu \left( b\right) \right]
    ^{-1}} \\
    &=&p_{1}\left( u|b\right) 
    \end{eqnarray*}%
    On the other hand, 
    \begin{eqnarray*}
    &&q\left( b|u_{0}\right) \coloneqq\frac{q\left( u,b\right) }{\int q\left(
    u,b\right) d\mu \left( b\right) } \\
    &&\coloneqq\frac{p_{1}\left( u_{0}|b\right) \frac{p_{2}\left( b|u_{0}\right) 
    }{p_{1}\left( u_{0}|b\right) }\left[ \int_{\mathcal{B}^{\ast }}\frac{%
    p_{2}\left( b|u_{0}\right) }{p_{1}\left( u_{0}|b\right) }d\mu \left(
    b\right) \right] ^{-1}}{\int_{\mathcal{B}^{\ast }}p_{1}\left( u_{0}|b\right) 
    \frac{p_{2}\left( b|u_{0}\right) }{p_{1}\left( u_{0}|b\right) }\left[ \int_{%
    \mathcal{B}^{\ast }}\frac{p_{2}\left( b|u_{0}\right) }{p_{1}\left(
    u_{0}|b\right) }d\mu \left( b\right) \right] ^{-1}d\mu \left( b\right) } \\
    &=&\frac{p_{2}\left( b|u_{0}\right) \left[ \int_{\mathcal{B}^{\ast }}\frac{%
    p_{2}\left( b|u_{0}\right) }{p_{1}\left( u_{0}|b\right) }d\mu \left(
    b\right) \right] ^{-1}}{\left[ \int_{\mathcal{B}^{\ast }}\frac{p_{2}\left(
    b|u_{0}\right) }{p_{1}\left( u_{0}|b\right) }d\mu \left( b\right) \right]
    ^{-1}\int_{\mathcal{B}^{\ast }}p_{2}\left( b|u_{0}\right) d\mu \left(
    b\right) } \\
    &=&p_{2}\left( b|u_{0}\right) 
    \end{eqnarray*}%
    This concludes the proof. 
\end{proof}

\section{A more detailed review of semiparametric theory}

\label{app:semiparametric-theory}

Given a fused-data framework $\left( \mathcal{Q},\mathcal{P},\mathcal{C}%
,\psi ,\varphi \right) $ with respect to $\left( Q_{0},P_{0}\right) $, our
task is to provide a general template for conducting semiparametric
inference about $\varphi \left( P_{0}\right) .$ We will now review key
elements of semiparametric theory which highlight the pivotal role that
influence functions play in constructing semiparametric estimators; more
precisely in constructing regular, asymptotically linear, debiased machine
learning estimators. The results that we will establish will inform how to
compute influence functions, and in particular the efficient influence
function, of $\varphi :\mathcal{P}\mapsto \mathbb{R}$ from influence
functions of $\psi :\mathcal{Q}\mapsto \mathbb{R}.$

A parametric model $\{P_{\theta}:\theta \in \Theta\}$ of mutually absolutely continuous probability laws on $X$ with $\Theta\subseteq \mathbb{R}^d$ is differentiable in quadratic mean (DQM) at $\theta_0$ for $\theta_0$ in the interior of $\Theta$ with score $h \in \prod_{i=1}^d L^2_0(P)$ if and only if the mapping $\theta \mapsto \sqrt{\frac{dP_{\theta}}{dP_{\theta_0}}}$ is Frechet differentiable with derivative $\frac{1}{2}h$ at $\theta_0$ when viewed as a mapping from $\mathbb{R}^d$ to $L^2(P_{\theta_0})$. A parametric model is regular at ${\theta_0}$ and if and only if it is DQM and $E_{P_{\theta_0}}[h(X)h(X)']$ is non-singular. A parametric model is regular if and only if it is regular at every $\theta \in \Theta$. For a model $\mathcal{M}$ on the distribution of a random vector $X$, i.e. a
collection of probability laws $P$ on $X$, and a collection $\mathcal{A}%
\coloneqq\mathcal{A}\left( P\right) $ of regular parametric submodels of $%
\mathcal{M}$ indexed by a scalar parameter $t$, each containing $P$ and such
that $P_{t=0}=P,$ the tangent set $\mathcal{T}^{\circ }(P,\mathcal{A})$ is
the collection of scores at $t=0$ of all submodels in $\mathcal{A}.$ The
closed linear span of $\mathcal{T}^{\circ }(P,\mathcal{A})$ in $L^{2}(P),$
denoted as $\mathcal{T}(P,\mathcal{A}),$ is called the tangent space. When $%
\mathcal{A}$ includes all regular parametric submodels, $\mathcal{T}^{\circ
}(P,\mathcal{A})$ and $\mathcal{T}(P,\mathcal{A})$ are called the maximal
tangent set and space and are denoted as $\mathcal{T}^{\circ }(P,\mathcal{M}%
) $ and $\mathcal{T}(P,\mathcal{M}).$ When $\mathcal{M}$ is unrestricted or
constrained only by complexity or smoothness conditions on certain
infinite-dimensional functionals of $P$, such as conditional expectations or
densities, then $\mathcal{T}\left( P,\mathcal{M}\right) $ coincides with $%
L_{0}^{2}\left( P\right) ,$ when $\mathcal{M}$ imposes equality constraints
on $P$ then $\mathcal{T}\left( P,\mathcal{M}\right) $ is a strict subspace
of $L_{0}^{2}\left( P\right) $. Models $\mathcal{M}$ with $\mc{T}(P, \mc{P}) = L^2_0(P)$ are called (locally at $P$) nonparametric. When $\mc{T}(P, \mc{P})$ is finite dimensional, $\mc{M}$ is called parametric. Otherwise $\mc{M}$ is called semiparametric. A
mapping $\gamma :\mathcal{M}\rightarrow \mathbb{R}$ is said to be pathwise
differentiable at $P$ in $\mathcal{M}$ with respect to a a class $\mathcal{A}
$ of submodels through $P,$ if there exists $\gamma _{P}^{1}\in
L_{0}^{2}(P), $ called a gradient of $\gamma ,$ such that for any regular
parametric submodel of $\mathcal{A}$ indexed by $t$ with score at $t=0$
denoted by $h$ and with $P_{t=0}=P,$ it holds that $\left. \frac{d}{dt}%
\gamma (P_{t})\right\vert _{t=0}=\left\langle \gamma _{P}^{1},h\right\rangle
_{L^{2}\left( P\right) }$. The canonical gradient, a.k.a. efficient
influence function, of $\gamma $ at $P$ denoted as $\gamma _{P,eff,\mathcal{A%
}}^{1}$ with respect to $\mathcal{A}$, is the unique gradient of $\gamma $
that belongs to $\mathcal{T}(P,\mathcal{A}).$ When $\mathcal{T}(P,\mathcal{A}%
)$ is maximal we write $\gamma _{P,eff}^{1}.$

Consider an estimator $\gamma _{n}$ of $\gamma \left( P\right) \in \mathbb{R}
$ based on i.i.d. random draws $X_{i},i=1,...,n,$ from a probability law in
model $\mathcal{M}$. The estimator $\gamma _{n}$ is asymptotically linear at 
$P$ if there exists $\Gamma _{P}\in L_{0}^{2}(P),$ referred to as the
influence function of $\gamma _{n},$ such that $n^{1/2}\left\{ \gamma
_{n}-\gamma \left( P\right) \right\} =n^{-1/2}\sum_{i=1}^{n}\Gamma
_{P}\left( X_{i}\right) +o_{p}(1)$ where $o_{p}(1)$ is a sequence that
converges to 0 under $P$. Asymptotically linear estimators have a limiting
normal distribution with mean zero and variance $var_{P}\left( \Gamma
_{P}\right) .$ In particular, consistent estimation of the asymptotic
variance is readily available from the empirical variance of the estimated
influence function. The estimator $\gamma _{n}$\ is regular with respect to
submodel class $\mathcal{A}$ at $P,$ if its convergence to its limiting
distribution is locally uniform over laws contiguous to $P$. See \cite%
{bickel_efficient_1998} for a precise definition. Regularity is a desirable
property for estimators because Wald confidence intervals centered around
irregular estimators exhibit suboptimal performance due to high local bias.
Specifically, when these intervals are computed using the pointwise limiting
distribution of irregular estimators, their coverage probability does not
uniformly converge across all laws within the model to the nominal level.
Consequently, regardless of the sample size, there will always be some laws
within the model where the actual coverage probability will significantly
deviate from the nominal level.

The convolution theorem (Theorem 25.20 of \cite%
{van_der_vaart_asymptotic_2000}) states that if $\mathcal{T}(P,\mathcal{M})$
is convex and $\gamma _{n}$ is regular at $P$, then $n^{1/2}\left\{ \gamma
_{n}-\gamma \left( P\right) \right\} \,\ $converges in law under $P$ to $%
Z+U, $ where $Z\sim N\left( 0,var_{P}\left( \gamma _{P,eff}^{1}\right)
\right) $ and $U$ is independent of $Z.$

An asymptotically linear estimator $\gamma _{n}$ of $\gamma \left( P\right) $
at $P$ with influence function $\Gamma _{P}$ is regular with respect to $%
\mathcal{A}$ at $P,$ if and only if $\gamma :$ $\mathcal{M}$ $\mapsto 
\mathbb{R}$ is pathwise differentiable at $P$ with respect to $\mathcal{A}$
and $\Gamma _{P}$ is a gradient of $\gamma $ (\cite%
{van_der_vaart_asymptotic_2000}). This result explains why the terms
"gradients" and "influence functions" are often used interchangeably, a
practice we will continue throughout this paper. Importantly, locally
efficient estimators are both regular and asymptotically linear, as noted in
van der Vaart's Theorem 25.23 (\cite{van_der_vaart_asymptotic_2000}).

For a multivariate $\gamma :\mathcal{M\mapsto }\mathbb{R}^{b},$ $b>1,$ the
definitions and results are to be understood component-wise.

The procedure known as one-step estimation (\cite{bickel_adaptive_1982}) is
a strategy for computing a regular asymptotically linear (RAL) estimator
with a given influence function $\gamma _{P}^{1}.$ Specifically, given $\hat{%
P}_{n}$ and $P$ in $\mathcal{M}$ and a gradient $\gamma _{P}^{1}$, the
one-step estimator associated with $\gamma _{P}^{1}$ is defined as 
\begin{equation*}
\widehat{\gamma }_{n}\coloneqq\gamma \left( \hat{P}_{n}\right)
+n^{-1}\sum_{i=1}^{n}\gamma _{\hat{P}_{n}}^{1}\left( X_{i}\right) .
\end{equation*}

We then have 
\begin{equation*}
\widehat{\gamma }_{n}-\gamma (P)=(\mathbb{P}_{n}-P)\gamma _{P}^{1}+(\mathbb{P%
}_{n}-P)[\gamma _{\widehat{P}_{n}}^{1}-\gamma _{P}^{1}]+R(\widehat{P}_{n},P)
\end{equation*}%
where $R(\widehat{P}_{n},P)=\gamma \left( \hat{P}_{n}\right) -\gamma \left(
P\right) +P\gamma _{\widehat{P}_{n}}^{1}.$ Thus, by the Central Limit
Theorem and Slutsky's Lemma, the estimator $\widehat{\gamma }_{n}$ is
asymptotically linear provided $R(\widehat{P}_{n},P)=o_{p}\left(
n^{-1/2}\right) $ and $(P_{n}-P)[\gamma _{\widehat{P}_{n}}^{1}-\gamma
_{P}^{1}]=o_{p}\left( n^{-1/2}\right) .$ When $\widehat{P}_{n}$ and $P$ are
in $\mathcal{M}$ and $\gamma _{P}^{1}$ is a gradient of $\gamma $ at $P$ in
model $\mathcal{M}$, $R(\widehat{P}_{n},P)$ is a second order term that it
often takes the form of a single integral or a sum of integrals. These are
integrals of either squared differences of infinite-dimensional nuisance
parameters evaluated at $\widehat{P}_{n}$ and $P,$ or of products of two
such differences. Then, provided the nuisance parameters at $\widehat{P}_{n}$
converge to their counterparts evaluated at $P$ at rate $o_{p}\left(
n^{-1/4}\right) ,$ the term $R(\widehat{P}_{n},P)$ will be $o_{p}\left(
n^{-1/2}\right) $ . On the other hand, $(P_{n}-P)[\gamma _{\widehat{P}%
_{n}}^{1}-\gamma _{P}^{1}]$ will be of order $o_{p}\left( n^{-1/2}\right) $
when $P[\gamma _{\widehat{P}_{n}}^{1}-\gamma _{P}^{1}]^{2}=o_{P}(1)$ and $%
\gamma _{\widehat{P}_{n}}^{1}-\gamma _{P}^{1}$ falls in a Donsker class with
probability tending to one. Alternatively, the Donsker class requirement can
be avoided if cross-fitting is employed (\cite{klaassen_consistent_1987,
schick_asymptotically_1986}). The one-step estimator is often referred to as
a debiased machine learning estimator. This is because when flexible machine
learning estimation strategies are used to estimate the infinite-dimensional
nuisance parameters on which $\gamma \left( P\right) $ depends, the plug-in
estimator $\gamma \left( \hat{P}_{n}\right) $ typically converges to $\gamma
\left( P\right) $ at rates slower than order $o_{P}\left( n^{-1/2}\right) .$
This is because except for specially tailored estimators of the infinite
dimensional nuisance parameters, $\gamma \left( \hat{P}_{n}\right) $
inherits the bias of their estimation, which converges to zero at rates
slower than $o\left( n^{-1/2}\right) .\,\ $The term $n^{-1}\sum_{i=1}^{n}%
\gamma _{\hat{P}_{n}}^{1}\left( X_{i}\right) $ in the one-step estimator
acts as a bias correction term essentially because it is an estimator of the
first term of a functional Taylor's expansion of the map $\gamma :\mathcal{M}%
\mapsto \mathbb{R}$ around $\hat{P}_{n}$. Alternative approaches for
constructing regular, asymptotically linear, debiased machine learning
estimators with a given influence function include the targeted maximum
likelihood estimation of \cite{van_der_laan_targeted_2006}, and the
estimating equations approach (\cite{van_der_laan_unified_2003,
tsiatis_semiparametric_2006}). Critically, these approaches also require
knowing the expression for the dependence on $P$ and $X$ of a gradient of $%
\gamma $.

\section{Efficiency estimation in the four types of fused-data models}
\label{app:model-types-eif}

As discussed in \cref{subsec:compute-EIF}, fused-data models can be classified into four types, one where the observed data model is unrestricted, one where it is restricted only by the restrictions on the ideal data model, one where it is restricted only by the assumed alignments, and one where it is restricted by both the ideal data model and the alignments. This section discusses the challenges in conducting semiparametric efficient inference under each model type. 

The first type is a model in which $\mathcal{Q}$ is non-parametric and the alignments in $\mathcal{C}$ do not impose equality constraints on the laws $%
P $ in $\mathcal{P}.$ In this case $\mathcal{P}$ is non-parametric.
Therefore for any observed data pathwise differentiable parameter $\varphi $
there exists exactly one observed data influence function which then agrees
with the efficient influence function $\varphi _{P,eff}^{1}.$ This can be
computed by applying \cref{algorithm:decompose} \texttt{DECOMPOSE} with input the unique
ideal data influence function $\psi_{Q}^{1}$. \cref%
{example:disease-prevalence} and \cref{example:transporting} scenarios
(i) and (ii) illustrate fused-data frameworks of this type.

The second type is a model in which $\mathcal{Q}$ is non-parametric but the
alignments in $\mathcal{C}$ impose equality constraints on the laws $P$ in $%
\mathcal{P}.$ In this case the tangent space $\mathcal{T}(P,\mathcal{P)}$ is strictly included in $L_{0}^{2}(P)$ and therefore for any
given pathwise differentiable observed data parameter $\varphi $ there exist
infinitely many observed data influence functions. However, there exists exactly one ideal data
influence function because $\mathcal{Q}$ is non-parametric. Therefore, the conditions \ref{item:eif2} and \ref%
{item:eif3} \cref{theorem:eif} are exactly the same conditions. %
\cref{example:transporting} scenario (iii) illustrates a fused-data framework of this second type. In this example
there exist closed-form expressions for $\varphi_{P,eff}^{1}$. The
expression for $\varphi_{P,eff}^{1}$ follows as a special
case of the efficient influence function derived in \cref{prop:example-eifs} in \cref{app:examples-extras}.

The third type is a model in which $\mathcal{Q}$ is strictly semiparametric and restricts at least one of the aligned conditional distributions 
but the alignments in $\mathcal{C}$ do not impose equality constraints on
the laws $P$ in $\mathcal{P}$. In this case $\mathcal{P}^{ext},$ but not $%
\mathcal{P}$, is non-parametric. \cref{example:tsiv-lsm} illustrates a
fused-data framework of this type. For frameworks of this third type, by
part \ref{item:influence functions2} of \cref{Theorem:influence functions}
we know that for every ideal data influence function there corresponds at
most one observed data influence function. In this case, to compute the
observed data influence function one can attempt a strategy that avoids
directly applying parts \ref{item:eif2} or \ref{item:eif3} of \cref%
{theorem:eif}. The rationale for this strategy is as follows. By part \ref{item:influence functions1} of %
\cref{Theorem:influence functions}, the class of all observed data
gradients $\varphi _{P}^{1}$ is obtained by sweeping over the class of all
ideal data gradients $\psi _{Q}^{1}$ that decompose as %
\eqref{eq:if-decomposition}. This then implies that each decomposing $\psi
_{Q}^{1}$ gives rise to a unique observed data gradient and it implies that
among the set of single observed data gradients corresponding to each ideal
data gradient, the one with minimum variance is the efficient observed data
influence function. In \cref{example:tsiv-lsm} we found $\varphi
_{P,eff}^{1}$ by first characterizing the aforementioned set and then
finding the minimizer of the variance of the elements invoking the
Cauchy-Schwartz inequality. We are optimistic that in most data-frameworks
of this third type, one will be able to find shortcuts for computing $%
\varphi _{P,eff}^{1}$ that avoid solving integral equations. Intuitively, $\mathcal{P}^{ext}$ non-parametric implies that the alignments
alone do not entangle the observed data distributions from the different
sources, so roughly speaking, finding $\varphi _{P,eff}^{1}$ should entail
the same level of difficulty as finding $\psi _{Q,eff}^{1}.$ The following proposition summarises this discussion. 

\begin{sproposition}
    \label{prop:only-Q-restricted-eif}
    Let $\left( \mathcal{Q},\mathcal{P},%
    \mathcal{C}\text{\thinspace },\psi ,\varphi \right) $ be a fused-data
    framework. Let $P\in \mathcal{P}
    $ and suppose $\varphi$ is pathwise differentiable at $P$ in model $\mc{P}$. Suppose there exists $Q$ in $\mc{Q}$ such that $\left(
    Q,P\right) $ is strongly aligned and $\psi $
    is pathwise differentiable at $Q$ in model $\mathcal{Q}$. Suppose that $\mathcal{T}(Q,\mathcal{Q})\varsubsetneq L_{0}^{2}(Q)$ and $\mathcal{T}(P,\mathcal{P}^{ext})=L_{0}^{2}(P)$. Let $\mathcal{D}$ be an index set such that the
    collection of all ideal data influence functions of $\psi $ at $Q$ in model $%
    \mathcal{Q}$ is given by 
    \begin{equation*}
    \left\{ \psi _{Q,d}^{1}:d\in \mathcal{D}\right\}.
    \end{equation*}
    and for all $d, d^{\prime }\in \mathcal{D}$, $\psi _{Q,d}^{1} = \psi
    _{Q,d^{\prime }}^{1}$ a.e. $Q$ if and only if $d = d^{\prime }$. Let $%
    \widetilde{\mathcal{D}}\subseteq \mathcal{D}$ be such that for all $d \in 
    \widetilde{\mathcal{D}}$, $\psi^1_{Q, d}  = \sum_{j\in[J]}\sum_{k\in [K^{(j)}]} m_k^{(j)}$ for some collection $\{m_k^{(j)} \in \mc{D}_k^{(j)}(Q):k\in[K^{(j)}], j\in[J]\}$. Let $\varphi^1_{P, d}$ be the unique influence function for $\varphi$ at $P$ in model $\mc{P}$ that corresponds to $\psi^1_{Q, d}$ for $d \in \widetilde{\mc{D}}$. Then, the efficient observed data influence function $\varphi^1_{P, eff}$ is equal to $\varphi^1_{P, d^*}$ where 
    \begin{align*}
        d^* = \argmin_{d\in \widetilde{\mc{D}}}var_{P}\left(\varphi^1_{P, d}\right).
    \end{align*}
\end{sproposition}

The fourth and last type is a model in which $\mathcal{Q}$ is strictly
semiparametric and the alignments in $\mathcal{C}$ impose equality
constraints on the laws $P$ in $\mathcal{P}$. In general, we
expect the computational challenge for deriving $\varphi_{P,eff}^{1}$ for
fused-data frameworks of this fourth type will be greater than for
frameworks of the other three types. Part \ref{item:eif2} of \cref{theorem:eif} will not be helpful in general to derive $\varphi _{P,eff}^{1}$ because the
specific $\psi _{Q}^{1}$ in that part is unknown. Notice that such $\psi
_{Q}^{1}$ will, in general, not be equal to $\psi _{Q.eff}^{1}$.

\subsection{Proofs for \texorpdfstring{\cref{app:model-types-eif}}{Supplement F}}

\begin{proof}[Proof for \cref{prop:only-Q-restricted-eif}]
    Each influence function $\psi^1_{Q, d}$ for $\psi$ at $Q$ in $\mc{Q}$ such that $d\in \widetilde{\mc{D}}$ corresponds to exactly one influence function $\varphi^1_{P, d}$ for $\varphi$ at $P$ in $\mc{P}$ by \cref{Theorem:influence functions} part \ref{item:influence functions2}
    
    Then, we know that the set 
    \begin{align*}
        \left\{\varphi^1_{P, d}:d\in\widetilde{\mc{D}}\right\}
    \end{align*}
    is the set of all influence functions for $\varphi$ at $P$ in model $\mc{P}$ by \cref{Theorem:influence functions} parts \ref{item:influence functions1} and \ref{item:influence functions2}. Recall that the efficient influence function is the unique influence function with $L^2(P)$ norm, equivalently minimum variance. Hence, the efficient influence function is given by $\varphi^1_{P, eff} \coloneqq \varphi^1_{d*}$ where
    \begin{align*}
        d^* = \argmin_{d\in \widetilde{\mc{D}}}Var_P[\varphi^1_{P, d}]
    \end{align*}
    The above $\argmin$ is well defined because we have that $\varphi^1_{Q, d} = \varphi^1_{Q, d'}$ if and only if $d = d'$ since, for all $d, d'\in \widetilde{\mc{D}}$, $\psi^1_{Q, d} = \psi^1_{Q, d'}$ if and only if $d = d'$ and by part \ref{item:influence functions2} of \cref{Theorem:influence functions} only one observed data influence function corresponds to each ideal data influence function. 
\end{proof}
\section{Caveat about constructing one-step estimators}

\label{app:efficient-one-step-estimators}

In this section we discuss a potential challenge in constructing debiased machine learning estimators when the alignment assumptions themselves restrict the observed data model, and a general strategy to overcome that challenge. Suppose that $\mathcal{T}\left( P,\mathcal{P}\right) $ is strictly included
in $L_{0}^{2}\left( P\right) $ and let $\varphi _{P}^{1}$ be a given
gradient of $\varphi :\mathcal{P\rightarrow }\mathbb{R}$ at $P.$ Recall from %
\cref{app:semiparametric-theory} that a key condition for convergence of
the one-step estimator 
\begin{equation}
\widetilde{\varphi }_{n}\coloneqq\varphi \left( \widetilde{P}_{n}\right)
+n^{-1}\sum_{i=1}^{n}\varphi _{\widetilde{P}_{n}}^{1}\left( O_{i}\right)
\label{one-step}
\end{equation}%
to a normal distribution at rate$\sqrt{n}$ was that the term $R\left( 
\widetilde{P}_{n},P\right) =o_{p}\left( n^{-1/2}\right) $ where for any $%
\widetilde{P}\in \mathcal{P}$, the term $R\left( \widetilde{P},P\right) %
\coloneqq\varphi \left( \widetilde{P}\right) +\int \varphi _{\widetilde{P}%
}^{1}\left( o\right) dP\left( o\right) -\varphi \left( P\right) .$ The term $%
R\left( \widetilde{P},P\right) $ acts as the second-order remainder in the
first-order functional Taylor expansion of $\varphi $ around $\widetilde{P}$%
. However, for $\widetilde{P}\not\in \mathcal{P}$, $R\left( \widetilde{P}%
,P\right) $ is not well-defined because $\varphi $ and its influence
function are only defined for $P$ in $\mathcal{P}$. Even if one were to
extend $\varphi \left( P\right) $ and $\varphi _{P}^{1}$ to any $P$ in the
non-parametric model $\mathcal{P}^{np}$ for the observed data law, there
would be no reason to expect that the term $R\left( \widetilde{P},P\right) $
would be of second order, unless the extension of $\varphi $ to $\mathcal{P}%
^{np}$ were pathwise differentiable at $\widetilde{P}$ in model $\mathcal{P}%
^{np}$ and its unique gradient coincided with $\varphi _{P}^{1}$ for $P$ in $%
\mathcal{P}$. As such, there is no reason to expect the one-step estimator $%
\widetilde{\varphi }_{n}$ to be RAL with influence function $\varphi
_{P}^{1} $ at $P$ unless either (1) $\widetilde{P}_{n}\in \mathcal{P}$, or
(2) $\varphi \left( \widetilde{P}_{n}\right) $ $\coloneqq\varphi _{np}\left( 
\widetilde{P}_{n}\right) $ for $\varphi _{np}:\mathcal{P}^{np}$ $\rightarrow 
\mathbb{R}$ a pathwise differentiable extension of $\varphi $ over $\mathcal{%
P}^{np}$ whose unique gradient coincides with $\varphi _{P}^{1}$ at any $%
P\in \mathcal{P}$.

One instance in which it will be difficult to find a non-parametric
extension $\varphi _{np}$ whose unique gradient in model $\mathcal{P}^{np}$
is equal to the efficient influence function is when model $\mathcal{P}%
^{ext} $ is strictly semiparametric. As such, even when $\varphi
_{P,eff}^{1} $ exists in closed form as in \cref{example:transporting}
(iii), in constructing efficient one-step estimators it will often be
crucial that one evaluates all the components of $P$ on which $\varphi
\left( P\right) $ and $\varphi _{P}^{1}$ depend on at a model-obedient
estimator of $P.$ However, off-the-shelf machine learning estimators of such
components will typically not result in a model obedient estimator of $P$.
This has the negative consequence that naive one-step estimators $\widetilde{%
\varphi }_{n}$ in $\left( \ref{one-step}\right) $ constructed with $\varphi
\left( \widetilde{P}_{n}\right) $ equal to an arbitrary pathwise
differentiable extension of $\varphi $ to $\mathcal{P}^{np}$ evaluated at a
some off-the-shelf non-parametric estimators of $\widetilde{P}_{n}$ and with 
$\varphi _{\widetilde{P}_{n},eff}^{1}$ instead of $\varphi _{\widetilde{P}%
_{n}}^{1}$ will fail to be asymptotically efficient. Even worse, when $%
\varphi \left( P\right) $ depends on infinite dimensional components of $P,$
such naive construction will yield estimators that generally not converge at
rate $O_{p}\left( n^{-1/2}\right) $ because they will be based on the
incorrect influence function in the bias correction term.

To illustrate that naive one-step-efficient-like constructions based on an
estimator $\widetilde{P}_{n}$ that does not obey the model, need not yield
an asymptotically efficient estimator, consider the fused-data frameworks
described in \cref{app-subsubsec:generalizing-ex-3}, which recall, includes the
framework in \cref{example:transporting} scenario (iii). Suppose that $B$
and $U$ take values in finite sets $\mathcal{U}$ and $\mathcal{B}$ such that 
$q\left( u,b\right) \coloneqq Q\left( U=u,B=b\right) >0$ and $P\left(
S=j,U=u,B=b\right) >0$ for $\left( u,b\right) $ in $\mathcal{U}\times 
\mathcal{B}$ and $j=1,2.$ Since model $\mathcal{P}$ is a model for a
finitely valued vector $O$, it is finite dimensional. Letting $\widehat{p}%
_{n}^{ML}\left( o\right) \coloneqq\widehat{P}_{n}^{ML}\left( O=o\right) $
denote the maximum likelihood estimator (MLE) of $p\left( o\right) \coloneqq %
P\left( O=o\right) $ in model $\mathcal{P},$ we have that $\varphi \left( 
\widehat{P}_{n}^{ML}\right) $ is the MLE of $\varphi \left( P\right) $ and
consequently asymptotically efficient provided, as we assume throughout, $%
\varphi $ is pathwise differentiable at $P$.

We will now illustrate that a natural one-step-efficient-like construction
suggested by formula $\left( \ref{one-step}\right) $ with $\varphi _{%
\widetilde{P}_{n},eff}^{1}$ instead of $\varphi _{\widetilde{P}_{n}}^{1}$
for $\widetilde{P}_{n}$ the empirical law of $O,$ might yield an inefficient
estimator of $\varphi \left( P\right) $. Let
\begin{align*}
    \widetilde{p}_{n}\left(
j,u,b\right) \coloneqq\widetilde{P}_{n}\left( S=j,U=u,B=b\right) \coloneqq %
n^{-1}\sum\limits_{i=1}^{n}I\left( S_{i}=j,U_{i}=u,B_{i}=b\right) .
\end{align*} 
Note
that while consistent for $p\left( o\right) ,$ $\widetilde{p}_{n}\left(
o\right) $ ignores the constraints on $P$ imposed by model $\mathcal{P}$ and
it is not equal to the MLE $\widehat{p}_{n}^{ML}\left( o\right) .$
Furthermore, $\widetilde{P}_{n}$ is not in $\mathcal{P}$ with positive
probability. This is because there is a non-zero probability that the
following equality fails for some $u$ and $u^{\prime }$ in $\mathcal{U}$ and 
$b$ in $\mathcal{B}$,%
\begin{equation}
\left. \left\{ \frac{\widetilde{p}_{n}\left( b|u,S=2\right) }{\widetilde{p}%
_{n}\left( u|b,S=1\right) }\right\} \right/ \left\{ \sum_{b^{\prime }}\frac{%
\widetilde{p}_{n}\left( b^{\prime }|u,S=2\right) }{\widetilde{p}_{n}\left(
u|b^{\prime },S=1\right) }\right\} =\left. \left\{ \frac{\widetilde{p}%
_{n}\left( b|u^{\prime },S=2\right) }{\widetilde{p}_{n}\left( u^{\prime
}|b,S=1\right) }\right\} \right/ \left\{ \sum_{b^{\prime }}\frac{\widetilde{p%
}_{n}\left( b^{\prime }|u^{\prime },S=2\right) }{\widetilde{p}_{n}\left(
u^{\prime }|b^{\prime },S=1\right) }\right\} .  \label{Ptilde}
\end{equation}%
and this equality is a necessary condition for $\widetilde{P}_{n}$ to be in $%
\mathcal{P}$ since both the right and left hand sides of $\left( \ref{Ptilde}%
\right) $ agree with $\widetilde{q}_{n}\left( b\right) \coloneqq\widetilde{Q}%
_{n}\left( B=b\right) $ for some $\widetilde{Q}_{n}$ when $\widetilde{P}_{n}$
is in $\mathcal{P}$.

Since $\widetilde{P}_{n}$ is not necessarily in $\mathcal{P}$, to proceed
with a one-step construction, we must first define $\varphi \left( 
\widetilde{P}_{n}\right) .$ To do so, we must define a pathwise
differentiable extension of $\varphi $ to the non-parametric model $\mathcal{%
P}^{np}.$ We have many choices to do so because there exist many possible
such extensions. For instance, for any fixed $u_{0},$ the functional $%
\varphi \left( P\right) \coloneqq\psi \left( Q_{P}\right) ,$ defined on $%
\mathcal{P}^{np},$ where%
\begin{equation}
q_{P}\left( u,b\right) \coloneqq p\left( u|b,S=1\right) \left. \left\{ \frac{%
p\left( b|u_{0},S=2\right) }{p\left( u_{0}|b,S=1\right) }\right\} \right/
\left\{ \sum_{b^{\prime }}\frac{p\left( b^{\prime }|u_{0},S=2\right) }{%
p\left( u_{0}|b^{\prime },S=1\right) }\right\}  \label{qp-example}
\end{equation}%
is one such possible extension. Suppose that we define $\varphi \left( 
\widetilde{P}_{n}\right) \coloneqq\psi \left( Q_{\widetilde{P}_{n}}\right) $
with $q_{P}$ defined as in the last display for a given fixed $u_{0},$ and
we naively compute our one-step-efficient-like estimator as 
\begin{equation*}
\left. \widetilde{\varphi }_{n,naive}\coloneqq\right. \varphi \left( 
\widetilde{P}_{n}\right) +n^{-1}\sum_{i=1}^{n}\varphi _{\widetilde{P}%
_{n},eff}^{1}\left( O_{i}\right)
\end{equation*}%
where 
\begin{eqnarray*}
\varphi _{\widetilde{P}_{n},eff}^{1}\left( o\right) \coloneqq I(s
&=&1)\left\{ h_{\widetilde{P}_{n}}^{\left( Q_{\widetilde{P}_{n}}\right)
}\left( u,b\right) -E_{\widetilde{P}_{n}}[h_{\widetilde{P}_{n}}^{\left( Q_{%
\widetilde{P}_{n}}\right) }(b,U)|B=b,S=1]\right\} \\
+I(s &=&2)\left\{ h_{\widetilde{P}_{n}}^{\left( Q_{\widetilde{P}_{n}}\right)
}\left( u,b\right) -E_{\widetilde{P}_{n}}[h_{\widetilde{P}_{n}}^{\left( Q_{%
\widetilde{P}_{n}}\right) }(B,u)|U=u,S=2]\right\}
\end{eqnarray*}%
with $h_{\widetilde{P}_{n}}^{\left( Q_{\widetilde{P}_{n}}\right) }$ the
solution of equation $\left( \ref{new-big-eq}\right) $, which we know exists
by \cref{prop:example-eifs} in \cref{app:examples-extras}.

We have that $%
n^{-1}\sum_{i=1}^{n}\varphi _{\widetilde{P}_{n},eff}^{1}\left( O_{i}\right)
=0$ by virtue of $\widetilde{P}_{n}$ being the empirical law of $O$. Thus, $%
\widetilde{\varphi }_{n,naive}$ is equal to the plug-in estimator $\varphi
\left( \widetilde{P}_{n}\right) .$ While $\widetilde{\varphi }%
_{n,naive}=\varphi \left( \widetilde{P}_{n}\right) $ is a RAL estimator of $%
\varphi \left( P\right) $, its influence function is equal to the unique
gradient of the functional $\varphi :\mathcal{P}^{np}\mathcal{\rightarrow }%
\mathbb{R}$ defined as $\varphi \left( P\right) \coloneqq\psi \left(
Q_{P}\right) ,$ which need not equal $\varphi _{P,eff}^{1}.$

As an example, consider the estimation of $\psi \left( Q\right) =Q\left(
U=u^{\ast }|B=b^{\ast }\right) .$ For this functional, the estimator $%
\varphi \left( \widetilde{P}_{n}\right) $ is equal to $\widetilde{P}%
_{n}\left( U=u^{\ast }|B=b^{\ast },S=1\right) $ which has influence function 
\begin{equation}
\varphi _{P}^{1}\left( o\right) =I\left( s=1\right) \frac{I\left( b=b^{\ast
}\right) }{P\left( S=1,B=b\right) }\left\{ I\left( u=u^{\ast }\right)
-P\left( U=u^{\ast }|B=b^{\ast },S=1\right) \right\}  \label{exphi1p}
\end{equation}%
By \cref{prop:example-eifs} in \cref{app:examples-extras} we
know that 
\begin{equation*}
\varphi _{P,eff}^{1}=I\left( s=1\right) \left\{ h\left( u,b\right) -E_{Q} 
\left[ h\left( U,b\right) |B=b\right] \right\} +I\left( s=2\right) \left\{
h\left( u,b\right) -E_{Q}\left[ h\left( u,B\right) |U=u\right] \right\}
\end{equation*}%
for some $h\left( u,b\right) .$ The right hand side of $\left( \ref{exphi1p}%
\right) $ is not equal to $\varphi _{P,eff}^{1}.$ If it were, then equating
the terms in $I\left( s=1\right) $ and $I\left( s=2\right) ,$ we conclude
that $h\left( u,b\right) $ would have to simultaneously satisfy%
\begin{equation}
h\left( u,b\right) -E_{Q}\left[ h\left( U,b\right) |B=b\right] =\frac{%
I\left( b=b^{\ast }\right) }{P\left( S=1,B=b\right) }\left\{ I\left(
u=u^{\ast }\right) -P\left( U=u^{\ast }|B=b^{\ast },S=1\right) \right\}
\label{s1}
\end{equation}%
and 
\begin{equation}
h\left( u,b\right) -E_{Q}\left[ h\left( u,B\right) |U=u\right] =0  \label{s2}
\end{equation}%
Suppose there exists no invertible measurable map $g$ such that $U=g\left(
B\right) $ a.e. - $Q.$ Then no $h\left( u,b\right) $ exists that satisfies $%
\left( \ref{s1}\right) $ and $\left( \ref{s2}\right) $ simultaneously
because the equality $\left( \ref{s2}\right) $ implies that $h\left(
u,b\right) $ is a function of $u$ only, but equality $\left( \ref{s1}\right) 
$ implies that $h\left( u,b\right) $ is equal to a non-zero constant times $%
I\left( b=b^{\ast }\right) I\left( u=u^{\ast }\right) $ plus a function of $%
b $ only. This establishes that $\varphi _{P}^{1}\not=\varphi _{P,eff}^{1}$
and consequently that $\widetilde{\varphi }_{n,naive}=\varphi \left( 
\widetilde{P}_{n}\right) $ is inefficient. Note that the fact that $\varphi
_{P}^{1}\not=\varphi _{P,eff}^{1}$ in this example shows that data from
source $2\,\ $carries information about the conditional law $Q\left(
U=u^{\ast }|B=b^{\ast }\right) $ even though in source 2 only the
conditional distribution of $B$ given $U$ aligns with the corresponding
conditional of $Q.$ In this example, though inefficient, $\varphi \left( 
\widetilde{P}_{n}\right) $ remains RAL because $U$ and $B$ are finitely
valued. If $B$ and/or $U$ had been continuous, and we had followed the same
construction but now with $\widetilde{p}_{n}\left( u|b,S=1\right) $ and $%
\widetilde{p}_{n}\left( b|u,S=2\right) $ being some off-the-shelf
non-parametric estimators of the conditional densities $p\left(
u|b,S=1\right) $ and $p\left( b|u,S=2\right) $, the estimator $\widetilde{%
\varphi }_{n}$ would had not even converged at rate $O_{p}\left(
n^{-1/2}\right) $ because, as noted earlier, it would be based on the
incorrect influence function in the bias correction term.

We will now outline a general strategy to construct a model obedient
estimator $\widehat{P}_{n}$ that should preserve the consistency property
and can therefore be used to construct one-step RAL estimators with
influence function $\varphi _{P}^{1}$. Our presentation will be informal
because a rigorous analysis of the properties of our proposal is beyond the
scope of this paper.

Given an arbitrary influence function $\varphi _{P}^{1},$ possibly but not
necessarily equal to $\varphi _{P,eff}^{1},$ suppose $\widetilde{P}_{n}$ is
an estimator that is consistent for $P$ in the sense that it satisfies 
\begin{equation}
\int \left\{ \varphi _{\widetilde{P}_{n}}^{1}\left( o\right) -\varphi
_{P}^{1}\left( o\right) \right\} ^{2}dP\left( o\right) \overset{P}{%
\rightarrow }_{n\rightarrow \infty }0  \label{consistency}
\end{equation}
The estimator $\widetilde{P}_{n}$ need not be model obedient.

Suppose first that $\mathcal{T}\left( \mathcal{Q},Q\right) =L_{0}^{2}\left(
Q\right) $, $\mathcal{T}\left( \mathcal{P},P\right) =\mathcal{T}\left( 
\mathcal{P}^{ext},P\right) \varsubsetneq L_{0}^{2}\left( P\right) $ and $\xi
\left( Q,\mathcal{C}\right) =\left\{ Q\right\} $ for all $Q;$ equivalently,
the aligned components of $Q$ determine it. In this scenario, there exist
several distinct maps $P\longmapsto Q_{P}^{\left( m\right) },m=1,...,M,$
from $\mathcal{P}^{np}$ to $Q$, such that $P \alignswith Q_{P}^{\left( m\right) }$ for $P\in \mathcal{P}$.
For instance, in the preceding example, every distinct choice of $u_{0}$ in $%
\left( \ref{qp-example}\right) $ yields one different such map. Given $%
\widetilde{P}_{n},$ define 
\begin{equation*}
\widehat{Q}_{\widetilde{P}_{n}}\coloneqq m^{-1}\sum_{m=1}^{M}Q_{\widetilde{P}%
_{n}}^{\left( m\right) },
\end{equation*}%
and for all $k\in \left[ K^{\left( j\right) }\right] ,j\in \left[ J\right] $
and $z_{k}^{(j)}\in \mathbb{R}^{\dim \left( Z_{k}^{\left( j\right) }\right)
},$ define \ 
\begin{equation*}
\widehat{P}_{n}\left( \left. Z_{k}^{\left( j\right) }\leq
z_{k}^{(j)}\right\vert \overline{Z}_{k-1}^{\left( j\right) },S=j\right) %
\coloneqq\widehat{Q}_{\widetilde{P}_{n}}\left( \left. Z_{k}^{\left( j\right)
}\leq z_{k}^{(j)}\right\vert \overline{Z}_{k-1}^{\left( j\right) }\right) 
\text{ a.e.- }Q\text{ on }\overline{\mathcal{Z}}_{k-1}^{\left( j\right) }
\end{equation*}%
\begin{align*}
\widehat{P}_{n}\left( \left. Z_{k}^{\left( j\right) }\leq
z_{k}^{(j)}\right\vert \overline{Z}_{k-1}^{\left( j\right) },S=j\right) %
\coloneqq&\widetilde{P}_{n}\left( \left. Z_{k}^{\left( j\right) }\leq
z_{k}^{(j)}\right\vert \overline{Z}_{k-1}^{\left( j\right) },S=j\right) \\
&\text{ a.e.- }P\left( \cdot|S=j\right) \text{ on }\textsf{Supp}\left[ \overline{Z}%
_{k-1}^{\left( j\right) };P\left( \cdot|S=j\right) \right] \backslash \overline{%
\mathcal{Z}}_{k-1}^{\left( j\right) }
\end{align*}%
Finally, define for all $j\in \left[ J\right] $ 
\begin{equation*}
\widehat{P}_{n}\left( S=j\right) \coloneqq\widetilde{P}_{n}\left( S=j\right)
.
\end{equation*}

By construction, $\widehat{P}_{n}$ is in model $\mathcal{P}$. Furthermore,
if the maps $P\rightarrow Q_{P}^{\left( m\right) },m=1,...,M,$ are smooth in
some sense we expect $\widehat{P}_{n}$ to preserve the consistency property $%
\left( \ref{consistency}\right) $.

Suppose next that $\mathcal{T}\left( \mathcal{Q},Q\right) =L_{0}^{2}\left(
Q\right) $, $\mathcal{T}\left( \mathcal{P},P\right) =\mathcal{T}\left( 
\mathcal{P}^{ext},P\right) \varsubsetneq L_{0}^{2}\left( P\right) $ but now $%
\xi \left( Q,\mathcal{C}\right) $ strictly includes $\left\{ Q\right\} ,$
i.e. $Q$ is not entirely determined by $P.$ In this case, the preceding
construction still yields a model obedient estimator $\widehat{P}_{n}$ if
one defines $\widehat{Q}_{\widetilde{P}_{n}}$ as before, but replacing in
each $Q_{\widetilde{P}_{n}}^{\left( m\right) }$ the undetermined components
of $Q$ with arbitrary ones. The estimator $\widehat{P}_{n}$ so constructed
will not depend on the undetermined components of $Q$ arbitrarily imputed
and should preserve the consistency of $\widetilde{P}_{n}$.

Finally, suppose $\mathcal{T}\left( \mathcal{Q},Q\right) \varsubsetneq
L_{0}^{2}\left( Q\right) ,$ and let $\widetilde{Q}_{\widetilde{P}_{n}}$ be a
law in $\mathcal{Q}$ closest to $\widehat{Q}_{\widetilde{P}_{n}}$ according
to some distance or discrepancy measure $d,$ i.e. 
\begin{equation*}
\widetilde{Q}_{\widetilde{P}_{n}}\coloneqq\arg \min_{Q\in \mathcal{Q}%
}d\left( Q,\widehat{Q}_{\widetilde{P}_{n}}\right)
\end{equation*}%
Define $\widehat{P}_{n}$ as before but with $\widetilde{Q}_{\widetilde{P}%
_{n}}$ replacing $\widehat{Q}_{\widetilde{P}_{n}}.$ We expect that $\widehat{%
P}_{n}$ will preserve the consistency of $\widetilde{P}_{n}$, although this
might depend on the choice of $d$.

\section{Additional results for the score operator}

\label{sec:appendix-score-operator}

\label{subsec:add-score-op-results} 
\subsection{Counterexample: Failure of pathwise differentiability of an identified parameter}
We begin this section by showing that there exist fused data frameworks $(\mc{Q}, \mc{P}, \mc{C}, \psi, \varphi)$ where the ideal data parameter is identified and pathwise differentiable over the model $\mc{Q}$, but where $\varphi$ is not pathwise differentiable at all laws $P \in \mc{P}$. This result follows from the fact that the range of the adjoint of the score operator is not necessarily closed. Hence, in light of \cref{lemma:pathwise-differentiabilty} in the main text, it suffices to find an ideal data parameter that is identified, but whose ideal data influence function is on the boundary of the range of the adjoint of the score operator.
\begin{scounterexample}   
    Let $Q_0$ be the joint distribution on $(V, U)$ generate as follows. $W$ and $U$ are drawn independently from a $Unif(0, 1)$ distribution. Conditional on $U, W$, $\delta$ is Bernoulli with probability $\delta=1$ given $U, W$ equal to $U$. Finally $V = U + \delta W$.
    
    Consider the parameter 
    \begin{align*}
        \psi(Q) \coloneqq E_Q[U^{^{-1/2}} - V^{-1/2}].
    \end{align*}
    defined for $Q$ in 
    \begin{align*}
        \mathcal{Q} \coloneqq \{Q_t: \frac{dQ_t}{dQ_0}(v, u) \coloneqq (1 + th(v, u))I(u \in (0, 1), v \in (0, 2)):& t \in (-\epsilon, \epsilon), h \in L^2_0(V, W;Q_0), \|h\|_{\infty} < \epsilon/2, \\ 
        &Q_t, Q_0 \text{ mutually absolutely continuous }\}. 
    \end{align*} 
    Example A.4.1 of \cite{bickel_efficient_1998} shows that 
    \begin{align*}
        E_{Q_0}[\{U^{{-1/2}} - V^{-1/2}\}^2] < \infty.
    \end{align*} 
    $\psi(Q) < \infty$ for all $Q \in \mc{Q}$ by the bounded density ratio between $Q_0$ and $Q$ for any $Q \in \mc{Q}$. Hence, $\psi(Q)$ is pathwise differentiable in $\mathcal{Q}$ at $Q_0$ with unique influence function 
    \begin{align*}
        \psi^1_{Q_0}(U, V) \coloneqq U^{^{-1/2}} - V^{-1/2} - \psi(Q_0).
    \end{align*}
    
    Now, consider the alignments
    \begin{align*}
        Q(U \leq u) &= P(U \leq u|S=1)\\
        Q(V \leq v) &= P(V \leq v|S=2)
    \end{align*}
    for all $u \in \textsf{Supp}[U;Q]$ and $v \in \textsf{Supp}[V;Q]$. It follows that $\psi(Q_0)$ is identified under these alignments by 
    \begin{align*}
        \varphi(P) = E_{P}[U^{-1/2}|S=1] - E_{P}[V^{-1/2}|S=2],
    \end{align*}
    which completes the fused-data framework $(\mc{Q}, \mc{P}, \mc{C}, \psi, \varphi)$.
    
    By \cref{lemma:pathwise-differentiabilty} in the main text, the observed data parameter $\varphi$ is pathwise differentiable at $P_0 \in \mc{P}$ for $(P_0, Q_0)$ strongly aligned if and only if the ideal data influence function at $Q_0$ belongs to the space 
    \begin{align*}
        L^2_0(U;Q_0) + L^2_0(V;Q_0).
    \end{align*}
    Hence, we must decompose $\psi^1_{Q_0}$ into a $Q_0$ mean-zero function of $U$ alone and a $Q_0$ mean-zero function of $V$ alone. The mean zero function of $U$ alone is necessarily 
    \begin{align*}
        f_U(U) \coloneqq U^{-1/2} - 2
    \end{align*}
    leaving the mean zero function of $V$ alone to be
    \begin{align*}
        f_V(V) \coloneqq  - (V^{-1/2} - 2) - \psi(Q_0).
    \end{align*}
    But, \cite{bickel_efficient_1998} show in in Example A.4.1 that $E_{Q_0}[\{U^{-1/2} - 2\}^2] = \infty$. This reveals that there does not exist a decomposition of $\psi^1_{Q_0}$ into the sum space $L^2_0(U;Q_0) + L^2_0(V;Q_0)$. As such, $\varphi$ is not pathwise differentiable at $P_0$ in model $\mc{P}$. We note that $\psi^1_Q \in \overline{L^2_0(U;Q) + L^2_0(V;Q)}\setminus  L^2_0(U;Q_0) + L^2_0(V;Q_0)$, i.e. the boundary of the range of the adjoint score operator for $Q_0$ $A_{Q_0}^*$. 
\end{scounterexample}

\subsection{Counterexample: The information operator is not a contraction}

We now demonstrate that
there exist fused-data models where the information operator has a bounded
inverse when considered as a map from $\text{Null}(A^*_{Q, U, \lambda}A_{Q,
U, \lambda})^\perp$ to $\text{Null}(A^*_{Q, U, \lambda}A_{Q, U,
\lambda})^\perp$, but the identity minus the information operator is not a
contraction. As indicated in the main text, this is in contrast with coarsening at random models, where the
identity minus the information operator is a contraction under the assumption that the probability of observing the full data is bounded away from 0 (\cite{robins_estimation_1994, van_der_laan_unified_2003}).

\begin{scounterexample}
Suppose $W=(X, Y)$ with $Y$ binary, $X\in \mathbb{R}^p$. Suppose $\mathcal{Q}
$ is unrestricted beyond that $\mathsf{Supp}\left[(X,Y);\mathcal{Q}\right] = 
\mathbb{R}^p \times \{0, 1\}$, and the alignments in $\mathcal{C}$ are such
that 
\begin{align*}
P(X \leq x|Y=y, S=1) &= Q(X \leq x|Y=y) \text{ for $y\in\{0, 1\}$, $x\in 
\mathbb{R}^p$} \\
P(Y=y| S=2) &= Q(Y=y) \text{ for $y\in\{0, 1\}$}.
\end{align*}
Suppose there exists $\delta>0$ such that $\delta^{-1}\leq\frac{P(Y=y|S=1)}{%
P(Y=y|S=2)}\leq \delta$ for $y\in\{0, 1\}$. Let $U^{(j)} = P(\cdot|S=j)$ for
each $j\in[J]$. Then $(Q, U, P)$ are strongly aligned. For all $h^{(Q)}\in
L^2_0(Q)$, 
\begin{align*}
((I-A_{Q, U, \lambda}^*A_{Q, U, \lambda})h^{(Q)})(x, y) =& h^{(Q)}(x, y)
\\
&-\left\{\frac{P(Y=y|S=1)P(S=1)}{Q(Y=y)}\left\{h^{(Q)}(x, y) - E_Q[h^{(Q)}(X,
Y)|y]\right\}\right\} \\
&-P(S=2)E_Q[h(X, Y)|y].
\end{align*}
Fix $h^{(Q)}(x, y) = I(y=1)\{f(x, y) - E_Q[f(X, Y)|y]\}$ for some $f\in
L^2_0(Q)$. Then, 
\begin{align*}
((I-A_{Q, U, \lambda}^*A_{Q, U, \lambda})h^{(Q)})(x, y) =& h^{(Q)}(x, y) - 
\frac{P(Y=1|S=1)}{P(Y=2|S=1)}P(S=1)h^{(Q)}(x, y) \\
=& \left\{1 - \frac{P(Y=1|S=1)}{P(Y=2|S=1)}P(S=1)\right\}h^{(Q)}(x, y)
\end{align*}
and so 
\begin{align*}
\|(I-A_{Q, U, \lambda}^*A_{Q, U, \lambda})h^{(Q)})\|_{L^2(Q)} =&
\left\|\left\{1 - \frac{P(Y=1|S=1)}{P(Y=2|S=1)}P(S=1)\right\}h^{(Q)}\right%
\|_{L^2(Q)} \\
=&\left|\left\{1 - \frac{P(Y=1|S=1)}{P(Y=2|S=1)}P(S=1)\right\}\right|\left%
\|h^{(Q)}\right\|_{L^2(Q)}.
\end{align*}
The above display reveals that if $\left|1 - \frac{P(Y=1|S=1)}{P(Y=2|S=1)}%
P(S=1)\right| > 1$, $I-A_{Q, U, \lambda}^*A_{Q, U, \lambda}$ is not a
contraction. But $\frac{P(Y=1|S=1)}{P(Y=1|S=2)}P(S=1)$ is restricted only by 
$\delta^{-1}\leq\frac{P(Y=1|S=1)}{P(Y=1|S=2)}\leq \delta$ for some $\delta>0$%
, and so can be made arbitrarily large, proving that for this fused-data
model, $I-A_{Q, U, \lambda}^*A_{Q, U, \lambda}$ is not a contraction.

We now demonstrate that the information operator $A_{Q, U, \lambda}^*A_{Q,
U, \lambda}$ has a bounded inverse when considered
as a map from $\text{Null}(A^*_{Q, U, \lambda}A_{Q, U, \lambda})^\perp$ to $%
\text{Null}(A^*_{Q, U, \lambda}A_{Q, U, \lambda})^\perp$. We first provide
an expression of the space $\text{Null}(A^*_{Q, U, \lambda}A_{Q, U,
\lambda})^\perp$. \cref{lemma:info-operator} below establishes that $(h^{(Q)},
h^{(U)}, h^{(\lambda)}) \in \text{Null}(A^*_{Q, U, \lambda}A_{Q, U,
\lambda}) $ if and only if $h^{(Q)} \in \text{Null}(A^*_QA_Q)$, $%
h^{(U^{(j)})} \in \text{Null}(A^*_{U^{(j)}}A_{U^{(j)}})$ for $j\in\{1,2\}$, and $%
h^{(\lambda)} \in \text{Null}(A^*_\lambda A_\lambda)$. $h^{(\lambda)} \in 
\text{Null}(A^*_\lambda A_\lambda)$ if and only if $h^{(\lambda)} = 0$. $%
h^{(U^{(j)})}\in \text{Null}(A^*_{U^{(j)}}A_{U^{(j)}})$ if and only if $%
h^{(U^{(j)})} \in \left(\bigoplus_{k\in[K^{(j)}]}\mathcal{R}_k^{(j)}(P_{Q, U,
\lambda})\right)^\perp$. In addition, it follows from the expression of $A^*_QA_Q$ in \cref{lemma:info-operator} that $h^{(Q)} \in \text{%
Null}(A^*_QA_Q)$ if and only if $h^{(Q)} \in \left(\sum_{j=1}^2\bigoplus_{k\in[%
K^{(j)}]}\mathcal{D}_k^{(j)}(Q)\right)^\perp$ because $\left\{\mathcal{D}_k^{(j)}(Q):k\in[K^{(j)}], j\in\{1, 2\}\right\}$ are mutually orthogonal in this fused-data model. Then$\left(\sum_{j=1}^2\bigoplus_{k\in[K^{(j)}]}\mathcal{D}%
_k^{(j)}\right)^\perp = \{0\}$. Hence, 
\begin{align*}
\text{Null}(A_{Q, U, \lambda}^*A_{Q, U, \lambda})^\perp = L^2_0(Q)\times
\prod_{j\in[J]}\left\{\sum_{k\in[K^{(j)}]}\mathcal{R}_k^{(j)}(P_{Q, U,
\lambda})\right\}\times L^2_0(\lambda).
\end{align*}

Now, to show that $A_{Q, U, \lambda}^*A_{Q, U, \lambda}$ is invertible with bounded inverse, it suffices to show that $A^*_{Q, U, \lambda}A_{Q, U, \lambda}$ is a bijection by the Banach open mapping theorem (\cite{kress_linear_1999} Theorem 10.8). Clearly, the information operator is injective when the domain is 
$\text{Null}(A^*_{Q, U, \lambda}A_{Q, U, \lambda})^\perp$. It remains to
show $A^*_{Q, U, \lambda}A_{Q, U, \lambda}$ is surjective. Let $h =
(h^{(Q)}, h^{(U^{(1)})}, h^{(U^{(2)})}, h^{(\lambda)})\in \text{Null}%
(A^*_{Q, U, \lambda}A_{Q, U, \lambda})^\perp$. From the expression for $\text{Null}%
(A^*_{Q, U, \lambda}A_{Q, U, \lambda})^\perp$ in the preceding display and the
properties of $\mathcal{R}_k^{(j)}(P_{Q, U, \lambda})$ in this fused-data
model we know that $h^{(Q)} \in L^2_0(Q)$, $h^{(U^{(1)})} \in
L^2_0(Y;P(\cdot|S=1))$, $h^{(U^{(2)})}(x, y) = f(x, y) - E_P[f(X,Y)|y, S=2]$
for some $f\in L^2_0(P(\cdot|S=2))$, and $h^{(\lambda)} \in L^2_0(\lambda)$.

Let $\widetilde{h} \coloneqq (\widetilde{h}^{(Q)}, \widetilde{h}%
^{(U^{(1)})}, \widetilde{h}^{(U^{(2)})}, \widetilde{h}^{(\lambda)})$ be
given by 
\begin{align*}
\widetilde{h}^{(Q)} &\coloneqq \frac{Q(Y=y)}{P(Y=y|S=1)P(S=1)}\{h^{(Q)}(x,
y) - E_Q[h^{(Q)}(X, Y)|y]\} + \frac{1}{P(S=2)}E_Q[h^{(Q)}(X, Y)|y] \\
\widetilde{h}^{(U^{(1)})} &\coloneqq \frac{1}{P(S=1)}h^{(U)^{(1)}} \\
\widetilde{h}^{(U^{(2)})} &\coloneqq \frac{1}{P(S=2)}h^{(U)^{(2)}} \\
\widetilde{h}^{(\lambda)} &\coloneqq h^{(\lambda)}.
\end{align*}
It is easy to check that $\widetilde{h} \in \text{Null}(A^*_{Q, U,
\lambda}A_{Q, U, \lambda})^\perp$ and that by \cref{lemma:info-operator}
below $A^*_{Q, U, \lambda}A_{Q, U, \lambda}\widetilde{h} = h$. As $h$ was
arbitrary, $A^*_{Q, U, \lambda}A_{Q, U, \lambda}$ is surjective. This concludes the counterexample.
\end{scounterexample}

\subsection{Additional lemmas}

The following lemma provides the expression for the information operator.

\begin{slemma}
\label{lemma:info-operator} Let $\left( \mathcal{Q},\mathcal{P},\mathcal{C}%
\right) $ be a fused-data model with respect to $\left( Q_{0},P_{0}\right) .$
Let $\left( Q,U,P\right) $ be strongly aligned with respect to $\mathcal{C}.$
Let $\lambda(S=j) = P(S=j)$. Then, the information operator $%
A_{Q,U,\lambda }^*A_{Q,U,\lambda }:\mathcal{H}\rightarrow \mathcal{H}$
exists, is bounded and linear, and for any $h\coloneqq\left(
h^{(Q)},h^{(U^{(1)})},...,h^{(U^{(J)})},h^{(\lambda )}\right) \in \mathcal{H}
$, 
\begin{align*}
A_{Q,U,\lambda }^*A_{Q,U,\lambda }h =
\left(A_Q^*A_Qh^{(Q)},A_{U^{(1)}}^*A_{U^{(1)}}h^{(U^{(1)})}, \dots,
A_{U^{(J)}}^*A_{U^{(J)}}h^{(U^{(J)})}, A_\lambda^*A_\lambda
h^{(\lambda)}\right)
\end{align*}
where 
\begin{align*}
A_Q^*A_Qh^{(Q)}(w) =& \sum_{j\in[J]}\sum_{k\in[K^{(j)}]}\Pi\left[\left.\frac{%
dP(\cdot|S=j)}{dQ}(\overline{Z}_{k-1}^{(j)})\lambda(S=j)\Pi%
\left[\left.h^{(Q)}\right|\mathcal{D}_k^{(j)}(Q)\right](\overline{Z}_{k}^{(j)})\right|\mathcal{T}(Q, 
\mathcal{Q})\right](w),
\end{align*}
\begin{align*}
A_{U^{(j)}}^*A_{U^{(j)}}h^{(U^{(j)})}(z^{(j)}) =& \sum_{k\in[K^{(j)}]}\frac{%
dP(\cdot|S=j)}{dU^{(j)}}(\overline{z}_{k-1}^{(j)})%
\lambda(S=j)\Pi\left[\left.h^{(U^{(j)})}\right|\mathcal{R}_k^{(j)}(P)\right](\overline{z}_{k}^{(j)}),
\end{align*}
and $A_{\lambda}^*A_{\lambda}h^{(\lambda)}(s) = h^{(\lambda)}(s)$.
\end{slemma}

We conclude this section with a lemma summarizing several results on the range of the score operator $A_{Q, U, \lambda}$ and its adjoint $A^*_{Q, U, \lambda}$. In what follows, 
\begin{align*}
&\Pi\left[\left.\sum_{j\in[J]}\bigoplus_{k\in[K^{(j)}]}\mc{D}_k^{(j)}(Q)\right\vert \mc{T}(Q, \mc{Q})\right]\\
    &\coloneqq\left\{\Pi\left[\left.\sum_{j\in[J]}\sum_{k\in[K^{(j)}]}m_k^{(j)}\right\vert \mc{T}(Q, \mc{Q})\right]:m_k^{(j)} \in \mc{D}_k^{(j)}(Q), k\in[K^{(j)}], j\in[J]\right\}.
\end{align*}
\begin{slemma}
    \label{lemma:operator-ranges}
    Let $\left( \mathcal{Q},\mathcal{P},\mathcal{C}%
\right) $ be a fused-data model with respect to $\left( Q_{0},P_{0}\right) .$
Let $\left( Q,U,P\right) $ be strongly aligned with respect to $\mathcal{C}.$
Let $\lambda(S=j) = P(S=j)$. Then,
\begin{enumerate}
    \item \label{item:lemma:operator-ranges1} The range of $A_{Q, U, \lambda}^{*}$ is
    \begin{align*}
        A^*_{Q, U, \lambda}L^2_0(P) =& A^*_QL^2_0(P) \times \prod_{j\in[J]}A^*_{U^{(j)}}L^2_0(P) \times A^*_\lambda L^2_0(P) \\
        =&\left(\Pi\left[\left.\sum_{j\in[J]}\bigoplus_{k\in[K^{(j)}]}\mc{D}_k^{(j)}(Q)\right\vert \mc{T}(Q, \mc{Q})\right]\right) \times \prod_{j\in[J]} \left(\bigoplus_{k\in[K^{(j)}]}\mc{R}_k^{(j)}(P)\right) \times  L^2_0(\lambda).
    \end{align*}
    \item \label{item:lemma:operator-ranges2}$A_{Q, U, \lambda}^{*}$ will have closed range if and only if $\Pi\left[\left.\sum_{j\in[J]}\bigoplus_{k\in[K^{(j)}]}\mc{D}_k^{(j)}(Q)\right\vert \mc{T}(Q, \mc{Q})\right]$ is closed.
    \item \label{item:lemma:operator-ranges3}If  $\sum_{j\in[J]}\bigoplus_{k\in[K^{(j)}]}\mc{D}_k^{(j)}(Q)$ is closed then $A_{Q, U, \lambda}^{*}$ has a closed range.
    \item \label{item:lemma:operator-ranges5}Let $\left( \mathcal{Q},\mathcal{P},\mathcal{C}\text{%
\thinspace },\psi ,\varphi \right) $ be a fused-data framework. Suppose that $%
\varphi $ is pathwise differentiable at $P$ in $\mathcal{P}$. Suppose $\psi $ is pathwise
differentiable at $Q$ in $\mathcal{Q}$. If $A_{Q, U, \lambda}$ has closed range then
    $\varphi _{P,eff}^{1}(o)=\sum_{j\in[J]}I\left(
    s=j\right) \sum_{k\in \left[ K^{\left( j\right) }\right] }\Pi \left[ \left. h^{\left( Q\right) }\left( W\right)
    \right\vert \mathcal{D}_{k}^{(j)}\left( Q\right) \right](\overline{z}_k^{(j)}) $ is the efficient influence function for $\varphi$ at $P$ in model $\mc{P}$ where $h^{\left( Q\right) }$ $\in \mathcal{T}\left( Q;\mathcal{Q}\right)$ satisfies 
    \begin{equation}
    \label{eq:operator-ranges-eif}
    \psi _{Q,eff}^{1}=\sum_{j\in \lbrack
    J]}\sum_{k\in \lbrack K^{(j)}]}\Pi \left\{ \left. \frac{dP(\cdot|S=j)}{dQ}(\overline{Z}_{k-1}^{(j)})P(S=j)\Pi \left[
    h^{\left( Q\right) }|\mathcal{D}_{k}^{(j)}\left( Q\right) \right](\overline{Z}_{k}^{(j)})
    \right\vert \mathcal{T}\left( Q;\mathcal{Q}\right) \right\} .
    \end{equation}
    \item \label{item:lemma:operator-ranges6}If there
exists a regular parametric submodel $\{Q_t:t\in(-\varepsilon,\varepsilon)\}$ in $%
\Phi(P;\mc{C})$ with $Q_t\vert_{t=0}=Q$ then $%
A^*_{Q}$ is not surjective.
\end{enumerate}
\end{slemma}

\subsection{Proofs for \texorpdfstring{\cref{sec:appendix-score-operator}}{Supplement H}}
\begin{proof}[Proof of \cref{lemma:info-operator}]
    This lemma follows directly from \cref{lemma:score-operator} and the decomposition \eqref{decompose}. 
\end{proof}

\begin{proof}[Proof of \cref{lemma:operator-ranges}]
$ $\newline
    \textbf{Proof of part \ref{item:lemma:operator-ranges1}:}
    Let 
    \begin{align*}
        \mc{F} \coloneqq \left(\Pi\left[\left.\sum_{j\in[J]}\bigoplus_{k\in[K^{(j)}]}\mc{D}_k^{(j)}(Q)\right\vert \mc{T}(Q, \mc{Q})\right]\right) \times \prod_{j\in[J]} \left(\bigoplus_{k\in[K^{(j)}]}\mc{R}_k^{(j)}(P)\right) \times  L^2_0(\lambda).
    \end{align*}
    First, take $f \in \mc{F}$. Then, 
    \begin{align*}
        f = \left(\Pi\left[\left.\sum_{j\in[J]}\bigoplus_{k\in[K^{(j)}]}m_k^{(j)}\right\vert \mc{T}(Q, \mc{Q})\right], \sum_{k\in[K^{(1)}]}n_k^{(1)}, \dots, \sum_{k\in[K^{(J)}]}n_k^{(J)},  \gamma\right)
    \end{align*}
    for some $m_k^{(j)} \in \mc{D}_k^{(j)}(Q)$, $n_k^{(j)} \in \mc{R}_k^{(j)}(P)$ for $k\in[K^{(j)}]$, $j\in[J]$, and $\gamma \in L^2_0(\lambda)$. 
    Let $\widetilde{m}_k^{(j)}(\overline{z}_k^{(j)}) \coloneqq \frac{dQ}{dP(\cdot|S=j)}(\overline{z}_{k-1}^{(j)})P(S=j)m_k^{(j)}(\overline{z}_k^{(j)})$, $\widetilde{n}_k^{(j)}(\overline{z}_k^{(j)}) \coloneqq \frac{dU^{(j)}}{dP(\cdot|S=j)}(\overline{z}_{k-1}^{(j)})P(S=j)n_k^{(j)}(\overline{z}_k^{(j)})$, for $k\in[K^{(j)}]$, $j\in[J]$ and let $\widetilde{\gamma} \coloneqq \gamma$. By \cref{lemma:f-in-Q},  $\widetilde{m}_k^{(j)} \in \mc{D}_k^{(j)}(Q)$ and $\widetilde{n}_k^{(j)} \in \mc{R}_k^{(j)}(P)$ for $k\in[K^{(j)}]$, $j\in[J]$. Let 
    \begin{align*}
        g(o) \coloneqq \widetilde{\gamma}(s) + \sum_{j\in[J]}I(s=j)\sum_{k\in[K^{(j)}]}\{\widetilde{m}_k^{(j)}(\overline{z}_k^{(j)}) + \widetilde{n}_k^{(j)}(\overline{z}_k^{(j)})\}.
    \end{align*}
    Then $A^*_{Q, U, \lambda}g = f$, so $\mc{F} \subseteq A^*_{Q, U, \lambda}L_0^2(P)$. 
    
    Now, let $f \in A^*_{Q, U, \lambda}L^2_0(P)$. Let $g\in L^2_0(P)$ be such that $A^*_{Q, U, \lambda}g = f$. By \eqref{decompose} we may write 
    \begin{align*}
        g(o) = \gamma(s) + \sum_{j\in[J]}I(s=j)\sum_{k\in[K^{(j)}]}{m}_k^{(j)}(\overline{z}_k^{(j)}) + {n}_k^{(j)}(\overline{z}_k^{(j)}).
    \end{align*}
    for some $m_k^{(j)}\in \mc{D}_k^{(j)}(Q)$, $n_k^{(j)}\in \mc{R}_k^{(j)}(P)$, $\gamma \in L^2_0(\lambda)$. By \cref{lemma:f-in-Q}, $\frac{dP(\cdot|S=j)}{dQ}(\overline{z}_{k-1}^{(j)})P(S=j)m_k^{(j)}\overline{z}_{k}^{(j)} \in \mc{D}_k^{(j)}(Q)$ and $\frac{dP(\cdot|S=j)}{dU^{(j)}}(\overline{z}_{k-1}^{(j)})P(S=j)n_k^{(j)}\overline{z}_{k}^{(j)} \in \mc{R}_k^{(j)}(P)$. Then, $f \in \mc{F}$ by the expression of $A^*_{Q, U, \lambda}$. Hence $\mc{F} = A^*_{Q, U, \lambda}L^2_0(P)$. 

    \textbf{Proof of part \ref{item:lemma:operator-ranges2}:}
    The orthogonal sum $\bigoplus_{k\in[K^{(j)}]}\mc{R}_k^{(j)}(P_{Q,U,\lambda})$ is closed because $\mc{R}_k^{(j)}(P)$ are mutually orthogonal closed linear spaces for $k\in[K^{(j)}]$. $L^2_0(\lambda)$ is also closed. Hence, $A^*_{Q, U, \lambda}$ will have closed range if and only if $\Pi\left[\left.\sum_{j\in[J]}\bigoplus_{k\in[K^{(j)}]}\mc{D}_k^{(j)}(Q)\right\vert \mc{T}(Q, \mc{Q})\right]$ is closed by part \ref{item:lemma:operator-ranges1} of this lemma. 

    \textbf{Proof of part \ref{item:lemma:operator-ranges3}:} This result is a direct corollary of parts \ref{item:lemma:operator-ranges1} and \ref{item:lemma:operator-ranges2} of this Lemma and the fact that orthogonal projections of closed linear spaces are closed. 

    \textbf{Proof of part \ref{item:lemma:operator-ranges5}:} Recall from the discussion in \cref{subsec:score-operator-theory} that the efficient influence function for $\varphi$ at $P$ in model $\mc{P}$ is the unique element $\varphi^1_{P, eff}$ of $\mc{T}(P, \mc{P})$ that satisfies $A^*_{Q, U, \lambda}\varphi^1_{P, eff} = (\psi^1_{Q, eff}, \bs{0}_J, 0)$. We have that $\mc{T}(P, \mc{P}) = \overline{A_{Q, U, \lambda}\mc{H}} = A_{Q, U, \lambda}\mc{H}$ because $A_{Q, U, \lambda}$ has closed range. Hence $\varphi^1_{P, eff} = A_{Q, U, \lambda}h$ for some $h \in \mc{H}$. Then, $h = (h^{(Q)}, h^{(U^{(1)})}, \dots, h^{(U^{(J)})}, h^{(\lambda)})$ solves $A^*_{Q,U, \lambda}A_{Q, U, \lambda}h = (\psi^1_{Q, eff}, \bs{0}_J, 0)$. From \cref{lemma:info-operator}, this statement can alternatively be written as 
    \begin{align*}
        A^*_{Q}A_Q h^{(Q)} =& \psi^1_{Q, eff}\\
        A^*_{U^{(j)}}A_{U^{(j)}} h^{(U^{(j)})} =& 0 \text{ for $j\in[J]$}\\
        A^*_{\lambda}A_\lambda h^{(\lambda)} =& 0
    \end{align*}
    The first equality in the above expression is equivalent to \eqref{eq:operator-ranges-eif}. The second and third equalities will be satisfied if and only if $h^{(U^{(j)})}\in \text{Null}(A_{U^{(j)}})$ for $j\in[J]$ and $h^{(\lambda)} = 0$, which in turn implies that $\varphi _{P,eff}^{1}(o)=\sum_{j\in[J]}I\left(
    s=j\right) \sum_{k\in \left[ K^{\left( j\right) }\right] }\Pi \left[ \left. h^{\left( Q\right) }\left( W\right)
    \right\vert \mathcal{D}_{k}^{(j)}\left( Q\right) \right](\overline{z}_k^{(j)}) $ as desired. 
    
    \textbf{Proof of part \ref{item:lemma:operator-ranges6}}
    Let $\{Q_t:t\in(-\varepsilon,\varepsilon)\}$ be a regular parametric submodel in $\Phi(P;\mc{C})$ with $Q_t\vert_{t=0}=Q$. By the definition of a regular parametric submodel the score $h^{(Q)}$ of this submodel at $t=0$ is non-zero. Let $U$ be such that $(Q, U, P)$ is strongly aligned. Let $U_{t} \coloneqq U$ for $t\in(-\varepsilon,\varepsilon)$ and $\lambda_{t} \coloneqq \lambda$ for $t\in(-\varepsilon,\varepsilon)$. We have that $P_{Q_t, U_{t}, \lambda_{t}} = P_{Q, U, \lambda}$ because $Q_t\in \Phi(P;\mc{C})$ for all $t\in(-\varepsilon,\varepsilon)$. Hence, the score of $P_{Q_t, U_{t}, \lambda_{t}}$ is 0 at $t=0$. Thus there exists a score $(h^{(Q)}, \bs{0}_{J}, 0)\not=0$ that is in the null space of the operator $A_{Q, U, \lambda}$. We have that $h^{(Q)} \not\in \text{Range}(A_Q^*)$ because $\overline{\text{Range}(A_{Q, U, \lambda}^*)} = \text{Null}(A_{Q, U, \lambda})^\perp$, and $(h^{(Q)}, \bs{0}_{J}, 0) \in \text{Null}(A_{Q, U, \lambda})$. Hence, $A_Q^*$ is not surjective.
\end{proof}

\printbibliography[heading=subbibliography]
\end{refsection}

\end{document}